%% file: gkBALFarXiv.tex
\documentclass{gtart}
\usepackage{amsfonts}
\usepackage{amsmath}
\usepackage{amsthm}
\usepackage{graphicx}
\usepackage{color}
\usepackage{latexsym,pictex}

\newtheorem{theorem}{Theorem}
\newtheorem{lemma}[theorem]{Lemma}
\newtheorem{proposition}[theorem]{Proposition}
\newtheorem{corollary}[theorem]{Corollary}
\newtheorem{conjecture}[theorem]{Conjecture}

\theoremstyle{definition}
\newtheorem{definition}[theorem]{Definition}

\newtheorem{remark}[theorem]{Remark}

\newtheorem{question}[theorem]{Question}

\numberwithin{theorem}{section}

\numberwithin{equation}{section}

\font\thinlinefont=cmr5

\def\Z{\mathbb Z}

\def\R{\mathbb R}
\def\C{\mathbb C}

\def\S{\Sigma}

\newcommand{\tb}{\mathop{\rm tb}\nolimits}
\newcommand{\rot}{\mathop{\rm rot}\nolimits}

\newcommand{\pf}{\mathop{\rm pf}\nolimits}

\begin{document}

\title{Constructing Lefschetz-type fibrations on four-manifolds}
\shorttitle{}

\authors{David T. Gay
\footnote{Supported in part by NSF/DMS-0244558 and fellowships from
  CRM/ISM and CIRGET}
\\Robion Kirby
\footnote{Supported in part by NSF/DMS-0244558}
}
\address{Department of Mathematics and Applied Mathematics,
University of Cape Town,
Private Bag X3,
Rondebosch 7701,
South Africa}
\secondaddress{University of California, Berkeley, CA 94720, USA} 
\email{David.Gay@uct.ac.za}
\secondemail{kirby@math.berkeley.edu}
 
\begin{abstract}
We show how to construct broken, achiral Lefschetz fibrations on
arbitrary smooth, closed, oriented $4$-manifolds. These are
generalizations of Lefschetz fibrations over the $2$-sphere, where we
allow Lefschetz singularities with the non-standard orientation as well
as circles of singularities corresponding to round $1$-handles. We can
also arrange that a given surface of square $0$ is a fiber. The
construction is easier and more explicit in the case of doubles of
$4$-manifolds without $3$- and $4$-handles, such as the homotopy
$4$-spheres arising from nontrivial balanced presentations of the
trivial group.
\end{abstract}


\sloppy

\maketitlepage

\section{Introduction}
\begin{theorem} \label{T:Closed}
Let $X$ be an arbitrary closed $4$-manifold and let $F$ be a closed
surface in $X$ with $F \cdot F = 0$. Then there exists a broken,
achiral Lefschetz fibration (BALF) from $X$ to $S^2$ with $F$ as a
fiber. 
\end{theorem}

Recall that a (topological) Lefschetz fibration (LF) on a closed
$4$-manifold is a smooth map to a closed surface with all
singularities locally modelled by the complex map $(w,z) \mapsto w^2
+z^2$. (We call these ``Lefschetz singularities''.) An achiral LF
(ALF) is one in which we also allow singularities modelled by $(w,z)
\mapsto (\overline{w})^2 + z^2$, the same model as above but with the
opposite orientation on the domain. (We call these ``anti-Lefschetz
singularities''.) All Lefschetz and anti-Lefschetz singularities in
this paper will be allowable, see Definition~\ref{D:Allowable}.  A
broken LF (BLF) is one in which we also allow singularities modelled
by the map from $S^1 \times \R^3$ to $S^1 \times \R$ given by
$(\theta,x,y,z) \mapsto (\theta,-x^2+y^2+z^2)$. (We call these ``round
$1$-handle singularities''.) Such a fibration was called a ``singular
LF'' in~\cite{ADK}, and the singularities were called ``indefinite
quadratic singularities'' there. Finally, a broken achiral LF (BALF)
is one in which all three types of singularities are allowed.

This theorem can be compared to work of Auroux, Donaldson and
Katz\-arkov~\cite{ADK}, and of Etnyre and Fuller~\cite{EF}.  In the
first it is shown that if $X^4$ has a near-symplectic form (which it
does when $b^+_2 > 0$), then $X^4$ is a broken Lefschetz pencil (BLP).
This is a generalization of Donaldson's earlier results on Lefschetz
pencils and symplectic structures~\cite{Donaldson}.  In particular,
$X$ blown up some number of times is a Lefschetz fibration over each
hemisphere of $S^2$ with different genus fibers, and then over the
equator round $1$-handles are added (independently) to the side with
lower genus; also the Lefschetz singularities can all (topologically)
be placed over the high genus hemisphere. In our paper, round
$1$-handles can also be added independently; see the Addendum below.

Etnyre and Fuller show that $X^4$ connected sum with a $2$-sphere
bundle over $S^2$ is an achiral Lefschetz fibration (ALF); the
connected sum occurs as the result of surgery on a carefully chosen
circle in $X$. Baykur~\cite{Baykur} has results relating this
construction to folded symplectic structures.

\begin{conjecture}
Not all closed, smooth, oriented $4$-manifolds are BLFs. For example,
it is possible that $\C P^2$ is necessarily achiral as a fibration
(even though it does have a Lefschetz pencil structure).
\end{conjecture}

We also prove:

{\bf Addendum to Theorem~\ref{T:Closed}} {\it If we are given a
collection of embedded $2$-spheres $S_1, \ldots, S_n$, each
intersecting $F$ in a single positive intersection, then we can
construct the BALF so that each $S_i$ is a section. In particular, if
the initial ``fiber'' $F$ has positive self-intersection, we can blow
up its self-intersection points, make a BALF in which the exceptional
divisors are sections, and then blow down these sections, to get a
broken, achiral Lefschetz pencil (BALP) with $F$ as a fiber.

We can arrange that the round $1$-handle singularities all project to
the tropics of Cancer and Capricorn, with their high genus sides
towards the equator and with all Lefschetz and anti-Lefschetz
singularities over the equator. 
}


A significant section of this paper is devoted to proving a result
(Theorem~\ref{T:ConvexACS} and Corollary~\ref{C:ConvexOBD}) on the
existence of ``convex'' BLFs on $4$-manifolds built from $0$-, $1$-
and $2$-handles, with prescribed boundary conditions. This is
essential to the proof of Theorem~\ref{T:Closed}, but is also of
independent interest as a natural generalization of Loi and
Piergallini's result~\cite{LoiPiergallini} (see also~\cite{AkbOzb}) on
the existence of Lefschetz fibrations on Stein surfaces.

The virtues of Theorem~\ref{T:Closed} are:

\begin{enumerate}

\item It covers small $4$-manifolds such as homology $4$-spheres.  In
particular the Gluck construction on a knotted $2$-sphere $K$ in $S^4$
is a possibly exotic homotopy $4$-sphere which is a BALF with $K$ as a
fiber.  Also, the homotopy $4$-spheres arising from non-trivial
presentations of the trivial group (see Problems~5.1 and~5.2
of~\cite{kprob}) are seen by a simplified construction to be BALFs.
$\C P^2$ with either orientation can be seen as a simple example of a
BALF.

\item The proof is fairly constructive, with the least constructive part
coming from the use of Giroux's theorem that two open books on a
$3$-manifold are stably equivalent if their $2$-plane fields are
homotopic~\cite{Giroux} and Eliashberg's theorem that homotopic
overtwisted contact structures are isotopic~\cite{EliashOvertwisted}.

\item Conceivably these BALFs can be used as LFs are used in
Donaldson-Smith theory~\cite{DonSmith} (and BLFs in Perutz's
generalization~\cite{Perutz, Perutz1, Perutz2}) to find multisections
which are pseudoholomorphic curves, in the sense of Taubes'
program~\cite{Taubes1, Taubes2} on pseudoholomorphic curves in
near-symplectic $4$-manifolds.

\item In a philosophical sense, this paper complexifies Morse
functions as much as possible, in the sense that it produces maps from
arbitrary $4$-manifolds to $\C P^1$ which, locally, are as complex
analytic as possible. This continues the long line of results
(obtaining pencils) from Lefschetz ($X$ algebraic) to Donaldson ($X$
symplectic) to Auroux-Donaldson-Katzarkov ($X$ near-symplectic).

\end{enumerate}

This is an existence theorem, so of course there ought to be a
uniqueness theorem, which we hope will be the subject of a following
paper.

We would especially like to thank the African Institute of
Mathematical Sciences in Cape Town for their hospitality during the
final writing of this paper.

\subsection{Outline}
We begin in Section~\ref{S:BALPsandBALFs} by giving precise
definitions of the types of fibrations considered, including control
on behavior near singularities and along boundaries. While doing this,
we also show how to achieve the singularities and boundary behavior in
terms of handle additions, and we show how such handle additions
affect the monodromies of fibrations and open book decompositions
(OBDs) on the boundaries. The two important types of boundary behavior
we define are ``convexity'' and ``concavity'' along boundaries,
conditions which mean that the fibrations restrict to OBDs on the
boundary and that concave boundaries can be glued to convex boundaries
as long as the OBDs match. The proof of Theorem~\ref{T:Closed} then
boils down to constructing a concave piece and a convex piece and
arranging that the open books match.

In Section~\ref{S:ADKexample} we look in detail at an example
from~\cite{ADK} of a BLF on $S^4$, breaking it down into handles as in
Section~\ref{S:BALPsandBALFs}. The goal is to get the reader
accustomed to the tools and language we use in the rest of the paper,
and to see various ways to split the BLF into convex and concave
pieces. In particular we show (Lemma~\ref{L:ConcaveFXB2}) how to
construct a concave BLF on $F \times B^2$ for any closed surface $F$.

In Section~\ref{S:Doubles} we show how to construct a BALF on the
double of any $4$-dimensional $2$-handlebody. This construction is
more explicit than the general case because it does not depend on
Giroux's work on open books or Eliashberg's classification of
overtwisted contact structures. This section also includes a method
(Lemma~\ref{L:Concave1H}) for adding $1$-handles to a concave (BA)LF.
At the end of the section we discuss the relationship between doubles
and the Andrews-Curtis conjecture about balanced presentations of the
trivial group.

Then in Section~\ref{S:Convex} we show that a $4$-manifold $X$ built
from just $0$-, $1$- and $2$-handles is a convex BLF. Furthermore, if
we are given a homotopy class of plane fields on $\partial X$, we can
arrange that the induced OBD on $\partial X$ supports an overtwisted
contact structure in this homotopy class. (This is not true for ALFs.)
In order to achieve this, we need to be able to positively and
negatively stabilize the OBD on $\partial X$. (Stabilization means
plumbing on Hopf bands, positive being left-handed bands and negative
being right-handed bands.) Positive stabilization is easy to achieve;
negative stabilization is easy if we allow achirality, but to avoid
achirality as much as possible we show in Lemma~\ref{L:1right2lefts}
that we can negatively stabilize with round $1$-handles instead of
achiral vanishing cycles.  This section also includes a detailed
analysis of almost complex structures carried by BLFs.

Section~\ref{S:ClosedProof} finishes off the proof of
Theorem~\ref{T:Closed} and the addendum. We take the concave BLF on $F
\times B^2$ from Section~\ref{S:ADKexample} and add enough $1$-handles
(as in Section~\ref{S:Doubles}) so that the complement is built with
just $0$-, $1$- and $2$-handles. This induces a particular OBD on the
boundary of this concave piece.  We then construct a convex BLF on the
complement as in Section~\ref{S:Convex}, inducing an OBD on its
boundary which supports a contact structure homotopic to the contact
structure supported by the OBD coming from the concave piece. We
arrange that both contact structures are overtwisted, so by
Eliashberg's classification of overtwisted contact
structures~\cite{EliashOvertwisted} they are isotopic. By Giroux's
work on open books~\cite{Giroux} the two OBDs have a common positive
stabilization, which we already know we can achieve on the convex
piece without introducing achirality. (Note that at this point the two
pieces are BLFs, not BALFs.)  The only new tool developed in this
section is a trick for stabilizing OBDs on concave boundaries of
(BA)LFs; unfortunately, to achieve the positive stabilizations we are
forced to introduce anti-Lefschetz singularities (achirality).

Section~\ref{S:Questions} gives a list of questions.

\subsection{Notation and conventions}
Unless otherwise stated, all manifolds are smooth, compact, connected
and oriented (possibly with boundary), and all maps between manifolds
are smooth.  Whenever we specify a local model for the behavior of a
map, we imply that the local models respect all orientations involved.
All almost complex structures respect orientations and all contact
structures are positive and co-oriented.

For our purposes, an open book decomposition (OBD) on a closed
$3$-manifold $M$ is a smooth map $f\co M \rightarrow B^2$ such that
$f^{-1}(\partial B^2)$ is a compact $3$-dimensional submanifold on
which $f$ is a surface bundle over $S^1 = \partial B^2$ and such that
the closure of $f^{-1}(B^2 \setminus \partial B^2)$ is a disjoint
union of solid tori on each of which $f$ is the projection $S^1 \times
B^2 \rightarrow B^2$. The binding is $B = f^{-1}(0)$, and the page
over $z \in S^1$ is $\Sigma_z = f^{-1}\{ \lambda z | 0 \leq \lambda
\leq 1 \}$, with $B = \partial \Sigma_z$. The monodromy is the isotopy
class (rel. boundary) of the return map $h\co \Sigma_1 \rightarrow
\Sigma_1$ for any vector field transverse to the interiors of all the
pages and meridinal near the binding. We will usually blur the
distinction between the isotopy class and its
representatives. Positively (resp. negatively) stabilizing an OBD
$f\co M \rightarrow B^2$ means plumbing on a left-handed
(resp. right-handed) Hopf band. Thus if $f'\co M \rightarrow B^2$ is
the result of positively (resp. negatively) stabilizing $f\co M
\rightarrow B^2$, then $f'\co -M \rightarrow B^2$ is the result of
negatively (resp. positively) stabilizing $f\co -M \rightarrow B^2$.

When a knot $K$ lies in a page of an open book decomposition or a
fiber of a fibration over $S^1$, we call the framing induced by the
page the ``page framing'', and abbreviate it $\pf(K)$.

\section{Broken, achiral Lefschetz fibrations and pencils}
\label{S:BALPsandBALFs}

We will be constructing and working with smooth surjective maps from
compact $4$-manifolds to compact surfaces with controlled behavior at
singularities and along boundaries, this control to be discussed
below. When such a map $f\co X^4 \rightarrow \Sigma^2$ is defined on
all of $X$ we will call $f$ a ``fibration'', decorated with various
adjectives which characterize the allowed singularities and boundary
behavior. When $f$ is defined only on the complement of a discrete set
$B \subset X$, near each point of which $f$ is locally modelled by the
canonical map $\C^2 \setminus 0 \rightarrow \C P^1$, we will call $f$
a ``pencil'', decorated with the same adjectives; the points of $B$
are called ``base points''. Note that for a pencil the target surface
$\Sigma$ is necessarily $S^2$. Also note that blowing up each base
point turns a pencil into a fibration, with the exceptional divisors
becoming sections. Similarly, blowing down square $-1$ sections of a
fibration over $S^2$ yields a pencil.  If $f\co X \setminus B
\rightarrow \S^2$ is a pencil and $p \in S^2$, we abuse terminology
slightly to say that the ``fiber'' over $p$ is $f^{-1}(p) \cup B$, a
compact surface, so that any two fibers intersect transversely and
positively at each base point.

Now we describe the adjectives which characterize the singularities,
as well as interpretations of the singularities in terms of handlebody
decompositions and the effects of the various singularities on
monodromies of fibrations on boundaries. The relationships between
singularities, handles and monodromies are critical for all the
constructions in this paper.

Consider a general smooth map $f$ from a $4$-manifold
$X$ to a surface $\Sigma$.

\begin{definition} \label{D:LefAntiLef}
A critical point $p \in X$ of $f$ is a {\em Lefschetz singularity} if
$f$ is locally modelled near $p$ by the map $g\co (w,z) \mapsto w^2 +
z^2$ from $\C^2$ to $\C$. If instead $f$ is locally modelled
near $p$ by $g \circ \tau$, where $\tau(w,z) = (\bar{w},z)$ reverses
orientation, then $p$ is an {\em anti-Lefschetz singularity}.
\end{definition}
A Lefschetz singularity is the standard singularity in a Lefschetz
fibration, corresponding to the critical point of a vanishing
cycle. The following remark is a standard result and, if the reader
finds it confusing, a more detailed exposition can be found
in~\cite{GompfStipsicz}.
\begin{remark}[Vanishing cycles as $2$-handles.] \label{R:LefAntiLef2H}
If $[0,1] \times S^1$ is an annulus in $\Sigma$ with a single
Lefschetz singularity in $f^{-1}([0,1] \times S^1)$, then
$f^{-1}([0,1] \times S^1)$ is a cobordism from $M_0 = f^{-1}(0 \times
S^1)$ to $M_1 = f^{-1}(1 \times S^1)$ on which the projection to
$[0,1]$ is a Morse function with a single Morse critical point of
index $2$ (at the Lefschetz singularity). The corresponding $2$-handle
is attached along a knot $K$ in $M_0$ which in fact lies in a fiber of
the fibration of $M_0$ over $S^1$, and the framing is one less than
the framing induced by the fiber, i.e. $\pf(K)-1$. Conversely, suppose we
start with a fibration $f\co X^4 \rightarrow \Sigma^2$, where $\Sigma$
has nonempty boundary and $f$ has no singularities over $\partial
\Sigma$. Now attach a $2$-handle to $X$ along a knot $K$ in a fiber
of the fibration $f^{-1}(\partial \Sigma) \rightarrow \partial
\Sigma$, with framing $\pf(K) - 1$, to make a new $4$-manifold $X'
\supset X$. Then $f$ extends to a fibration of $X'$ over $\Sigma$ with
exactly one new singularity, a Lefschetz singularity, at the core of
the $2$-handle. Lastly, if the monodromy of the fibration on $\partial
X$ is $h$ and the monodromy of the fibration on $\partial X'$ is $h'$,
the relation is that $h' = \tau_K \circ h$, where $\tau_K$ is a
right-handed Dehn twist along $K$. 

If instead we started with an anti-Lefschetz singularity, the
$2$-handle would be attached with framing $\pf(K)+1$ and, conversely,
if we attach a $2$-handle as above but with framing $\pf(K)+1$ rather
than $\pf(K)-1$, we can extend the fibration creating a single new
anti-Lefschetz singularity, and the monodromy changes by a left-handed
Dehn twist (i.e. $h' = \tau_K^{-1} \circ h$).
\end{remark}

\begin{definition} \label{D:Allowable}
An (anti-)Lefschetz singularity is {\em allowable} if the attaching
circle of its vanishing is homologically nontrivial in the fiber.
\end{definition}

\begin{figure}
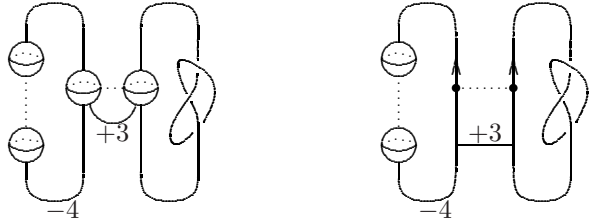

\begin{center}
\include{RoundExample}
\caption{Two drawings of a handlebody decomposition of $B^4$ involving
a $1$-handle, a round $1$-handle and a $2$-handle; on the left the
round $1$-handle is drawn as a $1$-handle and a $2$-handle.}
\label{F:RoundExample}
\end{center}
\end{figure}

As preamble to the next definition, recall that a ``round $k$-handle''
is $S^1$ times a $k$-handle. Thus a $4$-dimensional round $1$-handle
is $S^1 \times B^1 \times B^2$ attached along $S^1 \times S^0 \times
B^2$, i.e. attached along a pair of oriented framed knots. It is not
hard to see that the only important data is the relative orientation
of the pair (if we reverse one knot we should reverse the other) and
the relative framing (if we increase one framing by $k$ we should
decrease the other by $k$). A round $1$-handle can also be thought of
as a $1$-handle and a $2$-handle, with the attaching circle for the
$2$-handle running geometrically twice and algebraically zero times
over the $1$-handle. We will either draw round $1$-handles this way,
or shrink the balls of the $1$-handles down to small solid black
disks, so that we see two framed knots each decorated with a big black
dot, and a dashed line connecting the two dots. Drawn this latter way,
it is important to indicate the orientations with arrows.  Since only
the relative framing matters, we will only label one of the two knots
with a framing, implying that the other is $0$-framed. If a $2$-handle
runs over a round $1$-handle, we see its attaching circle as an arc or
sequence of arcs starting and ending on the attaching circles for the
round $1$-handle. Figure~\ref{F:RoundExample} gives two drawings of a
handlebody decomposition of $B^4$ involving a $1$-handle, a round
$1$-handle and a $2$-handle.

\begin{definition} \label{D:Round1H}
  An embedded circle $S \subset X$ of critical points of $f$ is a {\em
    round $1$-handle singularity} if $f$ is locally modelled near $S$ by
    the map $h\co(\theta,x,y,z) \mapsto (\theta,-x^2+y^2+z^2)$ from $S^1
    \times \R^3$ to $S^1 \times \R$.  Note that the genus of a fiber
    on one side of $f(S)$ is one higher than the genus on the other
    side; we will refer to these as the {\em high-genus side} and the
    {\em low-genus side}.
\end{definition}
This type of singularity is called an ``indefinite quadratic
singularity'' in~\cite{ADK}, which in principle also allows for a
local model which is a quotient of the above model by a $\Z/2$ action
so that the annulus $\{y=z=0\}$ becomes a M\"{o}bius band. In this
paper we do not need this nonorientable model. 

\begin{remark}[Attaching round $1$-handles.] \label{R:Round1H}
Let $[0,1] \times S^1$ be an annulus in $\Sigma$ with a single round
$1$-handle singularity $S$, and no other singularities, in
$f^{-1}([0,1] \times S^1)$, with $f(S) = 1/2 \times S^1$, and with the
low genus side over $0 \times S^1$ and the high genus side over $1
\times S^1$. Then $f^{-1}([0,1] \times S^1)$ is a cobordism from $M_0
= f^{-1}(0 \times S^1)$ to $M_1 = f^{-1}(1 \times S^1)$ which is the result
of attaching a round $1$-handle to $M_0$ along a framed, oriented pair
of knots $(K_1,K_2)$ each of which is a section of the fibration over
$S^1$, i.e. each one is transverse to all
the fibers and wraps once around the fibration in the positive
direction. Conversely, if we start with a fibration $f\co X
\rightarrow \Sigma$ with no singularities in $f^{-1}(\partial
\Sigma)$, if we choose any such pair $(K_1,K_2)$ in $f^{-1}(\partial
\Sigma)$, and if we attach a round $1$-handle along $(K_1,K_2)$ to
produce a new $4$-manifold $X' \supset X$, then $f$ extends to $f'\co
X' \rightarrow \Sigma$ with one new round $1$-handle singularity the
image of which is parallel to $\partial \Sigma$, and no other new
singularities. The fibers in $\partial X'$ are the result of
$0$-surgery on the fibers in $\partial X$ at the two points where
$K_1$ and $K_2$ intersect the fibers. To see how the monodromy
changes, consider a vector field transverse to the fibers in $\partial
X$ with $K_1$ and $K_2$ as closed orbits such that the return map $h$
on a fiber $F$ fixes a disk neighborhood $D_i$ of each $F \cap K_i$
and such that closed orbits close to $K_1$ and $K_2$ represent the
framings with which we are to attach the round $1$-handle. Let $F'$ be
the new fiber obtained by replacing $D_1 \cup D_2$ by $[0,1] \times
S^1$. Then the new monodromy is equal to $h$ on $F \setminus (D_1 \cup
D_2)$ and the identity on $[0,1] \times S^1$.

Since a round $1$-handle turned upside down is a round $2$-handle, we
could also understand constructions with round $1$-handle
singularities in terms of round $2$-handles. However, in our proofs we
do not seem to need this perspective.
\end{remark}

\begin{definition}
The adjective ``Lefschetz'' is used to mean that a given map
  (fibration or pencil) is allowed to have Lefschetz singularities. We
  add the adjective ``achiral'' to ``Lefschetz'' to indicate that we
  allow both Lefschetz and anti-Lefschetz singularities (recall that
  these are always allowable, as in Definition~\ref{D:Allowable}). The
  adjective ``broken'' means that round $1$-handle singularities are
  allowed.  (This term is due to Perutz~\cite{Perutz} and Smith and
  has been chosen to indicate that the non-singular fibers change
  genus when moving across the image in the base of a round $1$-handle
  singularity; since the singular circles disconnect the base, these
  singularities ``break'' the fibration in a certain sense.) If a type
  of singularity is not explicitly allowed then it is forbidden.

  To summarize and abbreviate, we have four kinds of ``fibrations'':
  Lefschetz fibrations (LFs), achiral Lefschetz fibrations (ALFs),
  broken Lefschetz fibrations (BLFs) and broken achiral Lefschetz
  fibrations (BALFs), with containment as follows: $LF \subset ALF$,
  $LF \subset BLF$, $ALF \subset BALF$ and $BLF \subset
  BALF$. Replacing ``fibration'' with ``pencil'' and ``F'' with ``P''
  in the preceding sentence also works.
\end{definition}

Now we describe the kind of boundary behavior we will allow for
fibrations and pencils on  $4$-manifolds with nonempty
boundary. Again consider a general smooth map $f\co X^4 \rightarrow
\Sigma^2$, and now let $M^3$ be a component of $\partial X$.
\begin{definition} \label{D:FlatConvexConcave}
We say that $f$ is ``flat'' along $M$ if $f(M)$ is a component of
  $\partial \Sigma$ and if $f|_M$ is an honest fibration over this
  component. We say that $f$ is ``convex'' along $M$ if $f(M) = \Sigma
  = B^2$ and if $f|_M\co M \to B^2$ is an open book decomposition of
  $M$. We say that $f$ is ``concave'' along $M$ if $f(M)$ is a disk
  $B^2$ in the {\em interior} of $\Sigma$ and if $f|_M$ is an open book
  decomposition of $M$. If $f$ is flat (resp. convex or concave) along
  each component of $\partial X$, we simply say that $f$ is flat
  (resp. convex or concave).
\end{definition}
Note that, for a convex fibration, the fibers are surfaces with
boundary. We use the term ``convex'' because a convex Lefschetz
fibration with ``allowable'' vanishing cycles (homologically
nontrivial in the fiber) naturally carries a symplectic structure (in
fact, a Stein structure) which has convex boundary. Likewise, a
concave Lefschetz pencil carries a symplectic structure with concave
boundary; in this case some fibers are closed and some are compact
with boundary. The term ``flat'' is similarly motivated; here the
fibers are all closed.

The typical example of a convex (BA)LF is $F \times B^2$ where $F$ is
a surface with nonempty boundary, together with vanishing cycles
(maybe of both kinds) and round $1$-handles.

\begin{remark}[Convex $1$-handles and concave $3$-handles.] 
\label{R:Convex1HConcave3H}
Suppose that $f\co X \rightarrow B^2$ is a convex fibration and that
$X'$ is the result of attaching a $1$-handle to $X$ at two balls
$B_0$, $B_1$ which are ``strung on the binding'' of the induced OBD on
$\partial X$ in the sense that $f|_{B_i}$ is the standard projection
$B^3 \rightarrow B^2$. Then $f$ extends to a convex fibration $f'\co
X' \rightarrow B^2$ with no new singularities. Each fiber $F'$ of $f'$
is diffeomorphic to a fiber $F$ of $f$ with a $2$-dimensional
$1$-handle attached along the two intervals $\partial F \cap B_0$ and
$\partial F \cap B_1$, and the same relation holds between the pages
of the new OBD on $\partial X'$ and the pages of the old OBD on
$\partial X$. The new monodromy is the old monodromy extended by the
identity across the $1$-handle. 

Dually, if $f \co X \rightarrow \Sigma^2$ is a concave fibration and
$X'$ is the result of attaching a $3$-handle to $X$ along a $2$-sphere
$S$ such that $f|_S$ is the standard projection $S^2 \rightarrow B^2$,
then $f$ extends to a concave fibration $f': X' \rightarrow \Sigma$
with no new singularities. Each page $F'$ of the new OBD on $\partial
X'$ is diffeomorphic to a page $F$ of $\partial X$ cut open along the
arc $S \cap F$. Implicit here is that the old monodromy was trivial in
a neighborhood of this arc, and so the new monodromy is just the old
monodromy restricted to $F'$. The fibers of $f'$ are related to the
fibers of $f$ as follows: If $f(\partial X) = B^2 \subset \Sigma$,
then the fibers over $\Sigma \setminus B^2$ do not change, while the
fibers of $f'$ over points in $B^2$ are obtained from the fibers of
$f$ over the same points by attaching $2$-dimensional $1$-handles. The
subtle point here is that each fiber of the fibration {\em inside} the
$4$-manifold gains a $1$-handle while each fiber of the OBD on the
{\em boundary} loses a $1$-handle.
\end{remark}

\begin{remark}
Some other handle attachments that are not used in this paper but that
can help develop the reader's intuition are as follows: If one
attaches $2$-handles to a convex (BA)LF, with one $2$-handle attached
along each component of the binding of the induced open book, with
framings $0$ relative to the pages, one produces a flat (BA)LF. Using
$+1$ framings instead produces a concave (BA)LP~\cite{GayConfig}.
\end{remark}

\begin{remark}[From flat to concave.] \label{R:Section2H}

One way to construct a concave (BA)LF is to start with a flat (BA)LF
and attach one or more $2$-handles along sections of the surface
bundle induced on the boundary. More concretely, suppose that $f\co X
\rightarrow \Sigma$ is flat along a boundary component $M \subset
\partial X$ and that $K_1, \ldots, K_n$ are framed knots in $M$ which
are sections of the induced fibration $f\co M \rightarrow S^1 \subset
\partial \Sigma$. Let $X' \supset X$ be the result of attaching
$2$-handles along $K_1, \ldots, K_n$ to $X$, and let $M'$ be the new
boundary component coming from surgery on $M$. Then $f$ extends to
$f'\co X' \rightarrow \Sigma'$, where $\Sigma'$ is the result of
attaching a disk $D$ to the relevant component of $\partial \Sigma$,
so that $f'$ is concave along $M'$. The cores of the $2$-handles
become sections of $f'$ over $D$, which extend as sections over all of
$X'$ as long as the knots $K_i$ extend as sections of $f$ over all of
$X$.  A concave (BA)LF which is used later in this paper is obtained
simply from $F \times B^2$, $F$ a closed surface, together with a
$2$-handle added to $point \times S^1$ with framing $0$.

In this process we transform a surface bundle over $S^1$ on $\partial
X$ into an OBD on $\partial X'$. Each page of the new OBD is
diffeomorphic to a fiber of the fibration on $\partial X$ with a disk
removed at each point of intersection with the sections $K_1, \ldots,
K_n$. If we choose a vector field $V$ transverse to the fibers in
$\partial X$ such that each $K_i$ is a closed orbit with a
neighborhood $\nu_i$ of closed orbits realizing the given framing of
$K_i$, and if $h$ is the return map on a fiber $F$ for flow along $V$,
then the monodromy of the new OBD on $\partial X'$ is precisely $h$
restricted to the new page $F \setminus (D_1 \cup \ldots \cup D_n)$,
where $D_i = \nu_i \cap F$.
\end{remark}

\begin{remark} [Glueing fibrations and pencils along boundaries.]
The point of spelling out the above boundary conditions is that it
should now be clear that fibrations and pencils can be glued along
common boundaries as long as we either
\begin{enumerate}
\item glue flat boundaries to flat boundaries via
  orientation-reversing diffeomorphisms respecting the induced
  fibrations over $S^1$ or
\item glue convex boundaries to concave boundaries via
  orientation-reversing diffeomorphisms respecting the induced open
  book decompositions.
\end{enumerate}
\end{remark}

\section{The Auroux-Donaldson-Katzarkov $4$-sphere example}

\label{S:ADKexample}

In section~8 of their paper~\cite{ADK} on singular (or broken) Lefschetz
fibrations, Auroux, Donaldson and Katzarkov construct a
BLF $f \co S^4 \to S^2$.  The fiber over the north pole is $S^2$, and over
the south pole is $T^2$. Over the polar caps are $S^2 \times B^2$ and
$T^2 \times B^2$. A round $1$-handle is attached to $S^2 \times B^2$,
giving a new boundary equal to $T^2 \times S^1 \to S^1$. Now this is
glued to $T^2 \times B^2 \to B^2$ by a diffeomorphism of $T^2 \times
S^1$ which rotates $T^2$ along a meridian as $S^1$ is traversed,
i.e. by a matrix of the form

$$
\left( \begin{matrix}
1&0&1 \\
0&1&0 \\
0&0&1
\end{matrix} \right)
\left(
\begin{matrix}
s \\
t \\
\theta
\end{matrix} \right)
= 
\left(\begin{matrix}
s+\theta \\
t   \\
\theta
\end{matrix} \right)
,\qquad  \theta \in S^1.$$   

The complement of the preimage $S^2 \times B^2$ of the arctic cap is
an interesting BLF for $B^3 \times S^1 \to B^2$ restricting to $S^2
\times S^1 \to S^1$ on the boundary; it is made from $T^2 \times B^2$
by adding a round $2$-handle in the right way.

However, it is more useful to describe the BLF in a somewhat different
way. If we pick the $0$-handle and one of the $1$-handles in $T^2$,
then its thickening gives $[0,1] \times S^1 \times B^2 \rightarrow
B^2$, a convex fibration with fiber an annulus. The base $B^2$ will
become the southern hemisphere $D_S$ of $S^2$. The complement in $S^4$
must be $S^2 \times B^2$, with a smaller $S^2 \times B^2$ in its
interior mapped by projection $S^2 \times B^2 \rightarrow B^2$ into
the northern hemisphere $D_N$. The fibration on this smaller $S^2
\times B^2$ is then flat along its boundary, inducing the fibration
$S^2 \times S^1 \rightarrow S^1$. The cobordism in between, $S^2
\times S^1 \times I$, will be mapped into $S^2$ in a way described
below, with one concave boundary component and one flat boundary
component which match the convex and flat boundaries of the two pieces
constructed above.

The cobordism $S^2 \times S^1 \times I$ can be written as a cancelling
$1$-$2$-handle pair and a cancelling $2$-$3$-handle pair, attached to
$S^2 \times B^2$ and not changing its diffeomorphism type. The
$1$-handle from the first pair and the $2$-handle from the second pair
will form a round $1$-handle, attached trivially along a pair of
circles $\{p_1,p_2\} \times \partial B^2 \subset \partial (S^2 \times
B^2)$, and mapping down to $D_N$. (The fibration extends over this
$1$-handle as in Remark~\ref{R:Round1H}).

\begin{figure}
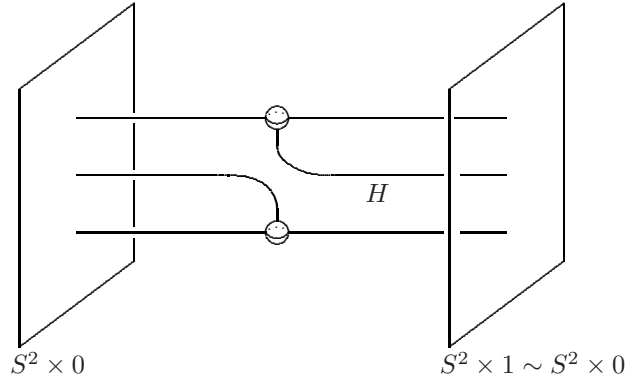

\begin{center}
\include{FindingAnH}
\caption{Finding an $H$.}
\label{F:FindingAnH}
\end{center}
\end{figure}

\begin{figure}
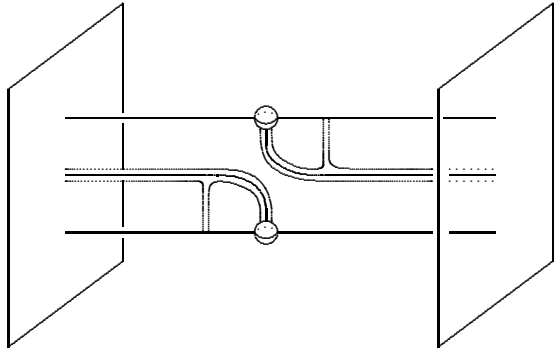

\begin{center}
\include{FindingAnHSlide}
\caption{Sliding twice over $H$.}
\label{F:FindingAnHSlide}
\end{center}
\end{figure}

The remaining $2$-handle and $3$-handle in fact make up a round
$2$-handle, or dually a round $1$-handle attached to the thickened
annulus $[0,1] \times S^1 \times B^2$, since the complement of an
annulus in $T^2$ is an annulus, and adding an annulus is the same as
adding a round $1$-handle. However, we do not use it as a round
$2$-handle here, but rather we map the $2$-handle and $3$-handle down
to $D_S$ as follows:

A handlebody picture of the process is given in
Figure~\ref{F:FindingAnH}. The $2$-handle labelled $H$ is the
$2$-handle of the round $2$-handle in the preceding paragraph, and in
the figure we see that its attaching map is a section of the fibration
over $S^1$, so that fibration will extend over $H$ exactly as in
Remark~\ref{R:Section2H}. Here, the framings of the $2$-handles are
chosen so that when the $2$-handle in the round $1$-handle is slid
twice over $H$ (see Figure~\ref{F:FindingAnHSlide}), then it becomes
an unknot, separated from the other components, with framing $0$, so
that it defines a $2$-sphere to which the $3$-handle (in the round
$2$-handle) is attached. $H$ then cancels the remaining $1$-handle.

There are several features about this construction that should be
noted. First, the concave piece has been constructed by adding the
$2$-handle $H = B^2 \times B^2$ along a section of $T^2 \times S^1
\rightarrow S^1 = \partial D_N$ which does not (in this case) extend
over $D_N$, and which maps to $D_S$ by projection on the first
factor. The fact that the section does not extend to $D_N$ is
necessary for otherwise $S^4$ would contain a hyperbolic pair, the
fiber and the global section. This is a key to finding BALFs for all
homology $4$-spheres. In particular, Theorem~\ref{T:Closed} shows that
any knotted $2$-sphere $K$ in $S^4$ can be made the fiber of a BALF on
$S^4$. Then, after performing the Gluck construction on $K$, the
resulting homotopy $4$-sphere is seen to be a BALF with fiber still
equal to $K$.

Second, the $3$-handle of the round $2$-handle (the $2$-handle being
$H$) is in a sense attached upside down to the concave side, as in
Remark~\ref{R:Convex1HConcave3H}; the attaching $2$-sphere consists of
a pair of disks parallel to $H$ and a cylinder $S^1 \times I$ which is
attached to a circle family of arcs in the fibers of $T^2 \times S^1
\rightarrow \partial D_N$.

Third, it is not necessary to begin building the concave piece with a
$2$-sphere fiber. Instead begin with $F^2 \times B^2 \rightarrow B^2 =
D_N$, where $F$ is a closed surface of genus $g$. Pick a pair of
points $p_1,p_2 \in F$ and attach a round $1$-handle along the
sections $\{p_1,p_2\} \times \partial B^2$ over $\partial D_N$. Now
add $H$ and the $3$-handle as before, and all the handles cancel
topologically. (Figures~\ref{F:FindingAnH} and~\ref{F:FindingAnHSlide}
are the same except that the squares at either end represent disks in
$F$.) We have thus proved:
\begin{lemma} \label{L:ConcaveFXB2}
Given any closed surface $F$ there exists a concave BLF  $f \co F \times
B^2 \to S^2$.
\end{lemma}
Note, however, that this statement of the result is deliberately vague
about the resulting OBD on $F \times S^1 = \partial (F \times B^2)$;
this is because we will not need to know anything about the OBD when
we use it later. However, in this $S^4$ example, it is important to
see the OBD, and it is instructive to think about what happens with
higher genus fibers.

\begin{figure}
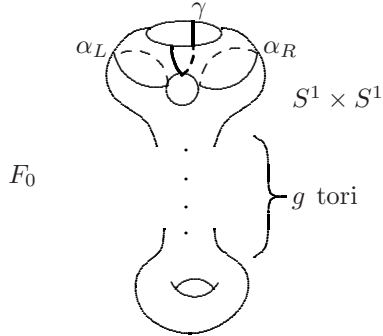

\begin{center}
\include{OBDAfterH}
\caption{Monodromy after attaching $H$.}
\label{F:OBDAfterH}
\end{center}
\end{figure}

Before adding the $3$-handle, the boundary is an open book with a
once-punctured fiber of genus $g+1$, called $F_0$. This open book does
not have trivial monodromy when $g+1 \geq 2$, a fact that needs
explaining.  It is easiest to understand the monodromy after attaching
$H$ if $H$ is added to a circle which corresponds to a fixed point of
the monodromy before attaching $H$; see Remark~\ref{R:Section2H}. In
this case, the initial monodromy is trivial, but $H$ is added to a
curve representing the sum of the class of $\{p\} \times S^1$ in $(F
\sharp (S^1 \times S^1)) \times S^1$ and the class of a curve running
over the first factor in $S^1 \times S^1$, which we call $\alpha$.  To
adjust for this fact, monodromy is introduced along two curves
$\alpha_L$ and $\alpha_R$ parallel to $\alpha$ which have the point
$p$ between them, with a left twist $\tau_{\alpha_L}^{-1}$ on one and
a right twist $\tau_{\alpha_R}$ on the other. Then the open book can
be represented, as in Figure~\ref{F:OBDAfterH}, by a fixed surface
$F_0$ (obtained by removing a disk neighborhood of $p$ from $F \sharp
(S^1 \times S^1)$) with twists along the curves $\alpha_L$ and
$\alpha_R$ drawn. When $g=0$ as in the case of $S^4$ above, then
$\alpha_L$ and $\alpha_R$ are isotopic in $F_0$ so that the two twists
cancel and the monodromy is still trivial after attaching $H$. But
when $g > 0$, $\alpha_L$ and $\alpha_R$ are not isotopic in $F_0$, so
this construction gives a concave BLF whose boundary is an open book
with non-trivial monodromy. The $3$-handle is then attached along the
$2$-sphere which intersects each page in the arc $\gamma$, so that
$\alpha_L$ and $\alpha_R$ become boundary parallel Dehn twists. It
follows that the convex piece, in order to fit with the concave piece,
cannot be just a $(g+1)$-genus surface minus an annulus, crossed with
$B^2$, for that has trivial monodromy on its boundary. However, if two
vanishing cycles were added to the convex side along $\alpha_L$ and
$\alpha_R$ (one framed $\pf+1$ and one framed $\pf-1$), this would
produce a convex piece which would ``dock'' into the concave piece.

Having given the construction for $S^4$, it is now easy to describe a
BALF on $\C P^2$: Simply take the above BLF for $S^4$ and add a
$+1$-framed $2$-handle to the $T^2 \times B^2$ along a nontrivial
circle in the fiber on the boundary. This produces a single
anti-Lefschetz singularity. The same construction with $-1$ gives us a
BLF on $\overline{\C P^2}$. This is interesting because $\C P^2$ is
symplectic (and is therefore a Lefschetz pencil) but seems to require
achirality when described as a fibration, while $\overline{\C P^2}$ is
far from symplectic but can be described as a fibration without using
anti-Lefschetz singularities.

\section{Doubled $4$-manifolds as BALFs}
\label{S:Doubles}

In this section we will prove a simpler version of
Theorem~\ref{T:Closed}, namely that the double $DX$ of any
$4$-dimensional $2$-handlebody $X$ is a BALF over $S^2$. Along the way
we prove some important lemmas needed for the full proof of
Theorem~\ref{T:Closed}, but this simpler result has the nice feature
of being more explicit than the full result in the sense that it does
not rely on Giroux's work on open books or Eliashberg's classification
of overtwisted contact structures.

The first tool we need is standard (see~\cite{AkbOzb}, for example).
\begin{lemma} \label{L:SimpleConvexStabilization}
Suppose that $f\co X \rightarrow B^2$ is a convex fibration and that $A$
is a properly embedded arc in a page of the induced OBD on $\partial
X$. First attach a $1$-handle to $X$ at the two endpoints of $A$ and
extend $f$ across the $1$-handle as in Remark~\ref{R:Convex1HConcave3H}. Let
$K$ be the knot lying in a page obtained by connecting the endpoints
of $A$ by going over the $1$-handle, and now attach a $2$-handle along
$K$ with framing $\pf(K)-1$ (resp. $\pf(K)+1$) and extend $f$ across
the $2$-handle as in Remark~\ref{R:LefAntiLef2H}. Since the $2$-handle
cancels the $1$-handle we get a new BALF on $X$ with one more
Lefschetz (resp. anti-Lefschetz) singularity (and different
fibers). Then the new OBD on $\partial X$ is the original OBD with a
left-handed (resp. right-handed) Hopf band plumbed on along $A$.
\end{lemma}

For clarification, recall that a Lefschetz singularity corresponds to
a right-handed Dehn twist, which in the lemma above corresponds to a
left-handed Hopf band (positive stabilization). Similarly, an
anti-Lefschetz singularity corresponds to a left-handed Dehn twist,
which in the lemma above corresponds to a right-handed Hopf band
(negative stabilization).

\begin{definition}
Given a handlebody decomposition of a manifold $X$, let $X_{(k)}$
denote the union of handles of index less than or equal to $k$. We
call $X$ a $k$-handlebody if $X = X_{(k)}$.
\end{definition}

We will make essential use of the following result:
\begin{proposition}[Harer~\cite{HarerLF}, Akbulut-Ozbagci~\cite{AkbOzb}]
  \label{P:1handlebodyLF} 
  Given a $4$-dimensional $2$-handlebody $X$, let $L$ be the attaching
  link for the $2$-handles in $\partial X_{(1)}$. Then there exists a
  convex LF $f\co X_{(1)} \rightarrow B^2$ such that $L$ lies in the
  interior of a single page $F$ of the induced open book
  decomposition of $\partial X_{(1)}$.  Furthermore, it can be
  arranged that each component $K$ of $L$ can be connected to
  $\partial F$ by an arc $A \subset F$ avoiding $L$ (i.e the
  interior of $A$ is disjoint from $L$).
\end{proposition}

\begin{proof} We do not need the full strength of the result
  in~\cite{AkbOzb}, so here we provide a streamlined proof of the
  result as we need it. The key fact we need is that if the page of an
  OBD of $S^3$ is obtained by plumbing left-handed Hopf bands onto a
  disk~\cite{Harer}, then this OBD is induced by a Lefschetz fibration
  on $B^4$. (Start with the fibration $B^4 = B^2 \times B^2
  \rightarrow B^2$ and plumb on the Hopf bands using
  Lemma~\ref{L:SimpleConvexStabilization}.)
  Figure~\ref{F:PALFconstruction} is an example illustrating the
  following construction:

\begin{figure}
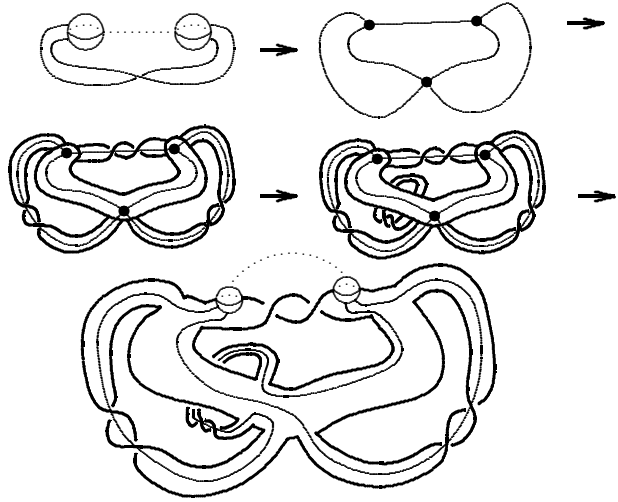

\begin{center}
\include{PALFconstruction}
\caption{Constructing a Lefschetz fibration as in Proposition~\ref{P:1handlebodyLF}.}
\label{F:PALFconstruction}
\end{center}
\end{figure}

Consider a standard balls-and-link diagram in $\R^3 = S^3 \setminus
\{\infty\}$ (balls for the $1$-handles, a link for the $2$-handles)
for the given handlebody decomposition of $X$.  Let $\Gamma$ be the
graph in $\R^2 \subset \R^3$ which is the projection of the diagram,
with crossings for $L$ and balls for $1$-handles made into vertices,
and with dotted lines for $1$-handles made into edges. By an isotopy
of $L$ we can always assume $\Gamma$ is connected.  Thus we have two
types of vertices: $4$-valent vertices for crossings and two
$(n+1)$-valent vertices for each $1$-handle which has $n$ strands of
the link running over it.

By plumbing left-handed Hopf bands onto a disk, one can easily
construct a surface $S$ which is made up of one disk neighborhood in
$\R^2$ of each vertex of $\Gamma$ and one (sometimes twisted) band
neighborhood in $\R^3$ of each edge of $\Gamma$. (Start with a disk
neighborhood of a spanning tree and then plumb on one Hopf band for
each remaining edge.)  At each $4$-valent vertex corresponding to a
crossing, plumb on an extra left-handed Hopf band along an arc at
right angles to one of the over-passing incident edges, underneath the
surface. Now $S$ is the page of an open book decomposition of $S^3$
induced by a Lefschetz fibration on $B^4$.

At this point, if there were no $1$-handles, we would be done, since
we could resolve the crossings of $\Gamma$ to reconstruct the link
simply by letting the under-crossing strand at each crossing go over
the extra Hopf band at that crossing. To deal with the $1$-handles, at
each $1$-handle vertex, string the foot of the $4$-dimensional
$1$-handle on the binding near that vertex (as in
Remark~\ref{R:Convex1HConcave3H}) and now pass all the strands
entering that vertex over the $1$-handle, remaining in the page the
whole time.

\end{proof}

To construct our BALF on the double $DX$ of a $2$-handlebody $X$, we
will use Proposition~\ref{P:1handlebodyLF} on $X_{(1)}$ so that the
$2$-handles lie in a page with some framing. The $2$-handles can now
be slid over their duals (small linking circles with framing $0$) so
that their framings are $\pf -1$ or $\pf -2$ (sliding over the dual
changes the framing by $\pm 2$ and does not change $L$ otherwise). If
$\pf -2$, then plumb on one more left-handed Hopf band along a short
boundary-parallel arc in a page and run the attaching circle over the
band so that the framing becomes $\pf-1$.

Note that we have now expressed $DX$ as equal to $DX'$ where $X'$ has
the same $0$- and $1$-handles as $X$ and has $2$-handles attached
along the same link but with different framings than $X$. We now
forget about the original $X$ and work with $X'$, which we simply call
$X$. In addition, the LF on $X_{(1)}$, with fiber $F$, in fact gives a
more complicated handlebody decomposition of $X_{(1)}$, where the
$1$-handles are those needed to build $F \times B^2$, and the
$2$-handles are the vanishing cycles needed to turn $F \times B^2$
into $X_{(1)}$. We now use this, together with the rest of the
$2$-handles needed to make $X$, as our handlebody decomposition of
$X$, and forget the previous handlebody decomposition. Thus $DX$ is
expressed as $F \times B^2$ together with $n$ $2$-handles attached
along knots in a page with framing $\pf -1$ and $n$ more dual
$2$-handles attached along small linking circles with framing $0$.

Now if we slide each dual over the $2$-handle it comes from, it
becomes a parallel $2$-handle, lying in a page with framing
$\pf+1$. Thus $(DX)_{(2)}$ is expressed as a convex ALF over $B^2$
with $n$ Lefschetz singularities and $n$ anti-Lefschetz singularities,
inducing an OBD on $\partial (DX)_{(2)}$ with trivial
monodromy, since each right-handed Dehn twist has a corresponding
parallel left-handed Dehn twist. (Note that at this stage we have not
used any round $1$-handles.)

To finish the construction, we will construct a concave BLF on the
union of the $3$- and $4$-handles of $DX$ inducing the same open book
as above. The concave structure we need, after
turning things upside down, is given by the following two results:

\begin{lemma} \label{L:Concave0H}
There exists a concave BLF $f\co B^4 \rightarrow S^2$ which restricts
to $S^3 = \partial B^4$ to give the standard OBD with disk pages.
\end{lemma}

\begin{proof}
Take the ADK $4$-sphere, discussed in Section~\ref{S:ADKexample}
above, and remove from $S^4$ a $4$-ball consisting of a neighborhood
of a section over $D_S$; that is, remove the $0$-handle of each torus
fiber over $D_S$. The result is the desired concave BLF.
\end{proof}

(We could equally well remove the $0$-handle of each sphere fiber over
$D_S$. However, the final BALF constructed on $DX$ will have the
undesirable feature that, as we move the torus fiber over the south
pole to the north pole, the genus of the fibers decreases from $1$ to
$0$, then increases. If we use the construction given in the proof above,
however, the genus will strictly increase as we move from one pole to
the other. When we finally prove Theorem~\ref{T:Closed}, the genus
will strictly increase as we move from each pole to the equator, but
will not have more than one ``local maximum''.)

\begin{lemma}[Attaching a $1$-handle to a concave boundary.] 
\label{L:Concave1H}
Suppose that $f\co X \rightarrow \Sigma$ is a concave fibration and
that $X'$ is the result of attaching a $1$-handle to $X$. Then, after
changing the handlebody decomposition of the cobordism from $\partial
X$ to $\partial X'$, we can extend $f$ to a concave fibration $f'\co
X' \rightarrow \Sigma$ with the following properties:
\begin{enumerate}
\item Each page of the new OBD on $\partial X'$ is diffeomorphic to a
  page of the OBD on $\partial X$ with a $2$-dimensional $1$-handle
  attached along two intervals in the binding. (The locations of these
  intervals can be chosen in advance.)
\item The monodromy of the new OBD is the monodromy of the old OBD
  extended by the identity across the $2$-dimensional $1$-handle.
\item The only singularities in $f'\co X' \rightarrow \Sigma$ that are
  not in $f\co X \rightarrow \Sigma$ are a single round $1$-handle
  singularity and a single Lefschetz singularity.
\end{enumerate}
\end{lemma}

\begin{proof}

Let $I_0$ and $I_1$ be the two intervals in the binding along which
the $2$-dimensional $1$-handle is to be attached. Move one foot of the
$4$-dimensional $1$-handle into a ball neighborhood $B_0$ of $I_0$ and
the other into a ball neighborhood $B_1$ of $I_1$. Inside $B_0$
introduce a cancelling $2$-$3$-handle pair so that the $2$-handle is
attached along a $0$-framed unknot $K$ and the $3$-handle is attached
along a $2$-sphere $S$ made of the Seifert disk for $K$ and the core
disk of the $2$-handle. Now slide an arc of $K$ over the $1$-handle so
that we see one unknotted loop of $K$ sticking out of the $1$-handle
in the ball $B_0$ and another unknotted loop sticking out of the
$1$-handle in the ball $B_1$. Now push each loop across the binding,
and the $1$-handle together with the $2$-handle becomes a round
$1$-handle as in Remark~\ref{R:Round1H}, across which the fibration
$f$ extends. This much is illustrated in
Figure~\ref{F:Concave1H.1}. 
\begin{figure}
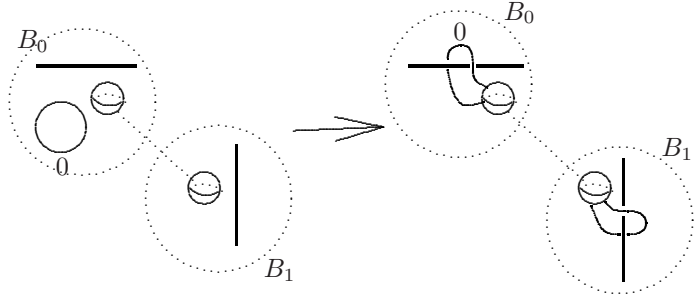

\begin{center}
\include{Concave1H.1}
\caption{Attaching a $1$-handle to a concave boundary.}
\label{F:Concave1H.1}
\end{center}
\end{figure}
Now observe that the page has changed by
removing a disk near $I_0$ and a disk near $I_1$ and replacing with
$[0,1] \times S^1$, with the monodromy extended by the identity across
$[0,1] \times S^1$, as illustrated in Figure~\ref{F:Concave1H.2}. 
\begin{figure}
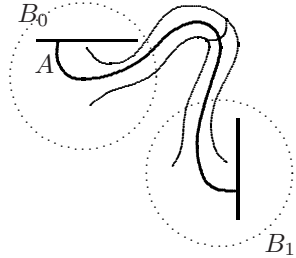

\begin{center}
\include{Concave1H.2}
\caption{How the page changes after attaching a $1$-handle to a
  concave boundary.}
\label{F:Concave1H.2}
\end{center}
\end{figure}
The
$3$-handle can then be seen to be attached along the $2$-sphere which
intersects each page in the arc $A$ drawn in
Figure~\ref{F:Concave1H.2}. Thus the fibration extends across the
$3$-handle as in Remark~\ref{R:Convex1HConcave3H}. The page has now
changed by cutting open along $A$, which amounts to attaching a
$2$-dimensional $1$-handle to the original page of $f$ at the two
intervals $I_0$ and $I_1$.
\end{proof}

Using these two lemmas, build a concave BLF on the union of the $3$-
and the $4$-handles which has an OBD on its boundary with trivial
monodromy and with pages diffeomorphic to the pages coming from the
convex BLF on $(DX)_{(2)}$. This can be done because the number of
$3$-handles in $DX$ equals the number of $1$-handles in $DX$ which
equals the number of $1$-handles in each page on $\partial
(DX)_{(2)}$; also Lemma~\ref{L:Concave1H} gives us the freedom to
attach the $2$-dimensional $1$-handles to the pages so as to get the
right number of boundary components.

Now glue the two pieces together using the diffeomorphism we get by
identifying their open books. This
gives $X$ for the following reason:  The $4$ and $3$-handles of $X$
form a boundary connected sum of $S^1 \times B^3$'s.  By a classical
theorem of Laudenbach and Poenaru~\cite{LP} it does not matter which
diffeomorphism of a connected sum of $S^1\times S^2$'s is used to glue
on the $4$- and $3$-handles; the resulting $4$-manifolds are
diffeomorphic.  Thus, in gluing trivial open book to trivial open book
above, we must obtain $X$.

Thus we have proved:

\begin{proposition}
If $X^4$ is a $2$-handlebody then its double $DX$ has a BALF $f\co DX
\rightarrow S^2$.
\end{proposition}

\subsection{The Andrews-Curtis conjecture}

One way of constructing smooth, homotopy $4$-spheres which may not be
diffeomorphic to $S^4$ (they are homeomorphic~\cite{Freedman}) is to
use balanced presentations of the trivial group which are not known to
satisfy the Andrews-Curtis Conjecture.

If a finite presentation of a group $G$ is described by attaching $1$-
and $2$-handles to an $n$-ball, then sliding handles over handles and
introducing or cancelling $1$-$2$-handle pairs correspond to what are
called Andrews-Curtis moves on the presentation. The Andrews-Curtis
conjecture is that any balanced presentation of the trivial group can
be reduced to the trivial presentation using only these moves. The
point is that you ``can't remember'', meaning that at any moment the
only relations available for use are those of current
$2$-handles. (When one $2$-handle slides over another the old relation
represented by the old $2$-handle is lost.)

A balanced presentation $P = \{x_1, \ldots , x_n \| r_1, \ldots ,
r_n\}$ of the trivial group determines uniquely a homotopy $4$-sphere
by attaching $n$ $1$-handles to the $5$-ball, and then $n$ $2$-handles
whose attaching maps read off the relations $\{ r_1, \ldots , r_n\}$.
If two attaching maps represent the same relation, then they are
homotopic, and homotopic circles in dimension $4$ are isotopic.  Hence
this $5$-manifold $V^5$ is unique up to diffeomorphism and is
contractible.  Its boundary $\partial V = S_P$, is the homotopy
$4$-sphere associated with the presentation $P$.

Given $P$, we can also build $4$-manifolds $X$ which are contractible by
adding $n$ $1$-handles to the $4$-ball and then $n$ $2$-handles
corresponding to the relations.  This involves choices because different
attaching maps which are homotopic are not necessarily isotopic, so there
are many possible choices of $X$ corresponding to $P$.  However in all
cases, $X \times I$ is diffeomorphic to $V^5$ because with the extra
dimension homotopic attaching maps are isotopic.  We have shown:

\begin{lemma}Our homotopy $4$-sphere $\partial V$ is diffeomorphic to $\partial
(X \times I)$ and hence diffeomorphic to the double $DX$ which is known to
be a BALF.
\end{lemma}

\begin{question}
Is the fact that $\partial V$ is known to be a BALF helpful in
showing that $\partial V$ is, or is not, diffeomorphic to $S^4$?
\end{question}

\begin{remark}
If a presentation $P$ can be reduced to the trivial presentation
by Andrews-Curtis moves, then these moves can be mirrored geometrically in
handle slides, and then $V^5 = B^5$ so $DX$ is $S^4$.

But it is possible that $DX$ is diffeomorphic to $S^4$ even though $P$
cannot be reduced to the trivial presentation by Andrews-Curtis moves.
This would have to be the case if the Andrews-Curtis Conjecture is
false (as is expected by many experts) and the smooth $4$-dimensional
Poincar\'{e} Conjecture is true.

The authors know of only one presentation $P$, namely $\{x,y \| xyx=yxy,
x^4 = y^5\}$, which is not known to satisfy the Andrews-Curtis Conjecture
but is known to give $S^4$.  The latter was shown in~\cite{ak1, ak2} with
a beautiful denouement by Gompf in~\cite{gompf-ak}.
\end{remark}

There are many tantalizing presentations to play with.  A full
discussion appears in~\cite{metzler}.

\section{The general construction of convex $2$-handlebodies}
\label{S:Convex}

To prove Theorem~\ref{T:Closed} we will need a general construction of
convex BLFs on $2$-handlebodies, with prescribed boundary conditions.
As a warm-up we prove a simple version without the boundary conditions.

\begin{proposition}[Quick and easy recipe for convex BLFs]
\label{P:QuickNEasy}
Every $4$-dimensional $2$-handlebody $X$ can be given the structure of
a convex BLF.
\end{proposition}

\begin{proof}
  Let $f \co 
  X_{(1)} \rightarrow D^2$ be the LF whose existence is
  asserted by Proposition~\ref{P:1handlebodyLF}. The idea now is to turn each
  $2$-handle (whose attaching circle lies in a page of the open book
  on $\partial X_{(1)}$) into a round $1$-handle whose attaching circles
  are transverse to the pages of the open book. For each such
  attaching circle $K$ of a $2$-handle $H$, consider a neighborhood
  $U$ of the arc $A$ mentioned in Proposition~\ref{P:1handlebodyLF}, in which
  we see only an arc of the binding $B$ and an arc of $K$ lying in a
  half-disk of the page $F$. The following construction is
  illustrated in Figure~\ref{F:2HandleToRound1Handle}.

First introduce a cancelling $1$-$2$-handle pair inside $U$ so that
the feet of the $1$-handle intersect $F$ in small disks, so that the
attaching circle of the cancelling $2$-handle runs from one foot
straight to the other staying in $F$ with framing $-1$ with respect to
this picture.  Next, slide a small loop of $K$ over the $1$-handle,
and now $H$ together with the $1$-handle form a round $1$-handle $H'$;
the two attaching circles of $H'$ are a small unknot $U$ near $B$ and
a copy $K'$ of the original knot $K$. Now push $U$ across $B$ to
become a small meridinal loop, hence transverse to the pages of the
open book. Likewise, push a small finger out from $K'$ and across $B$
and then tilt the rest of $K'$ out of the page $F$ so that $K'$ also
becomes transverse to the pages. Thus the two feet of this round
$1$-handle wrap once around the binding and the broken Lefschetz
fibration extends across the round $1$-handle.  Lastly note that the
cancelling $-1$-framed $2$-handle now lies in the extended page (after
attaching the round $1$-handle) and has framing $\pf-1$, so the
fibration also extends across these $2$-handles.

\begin{figure}
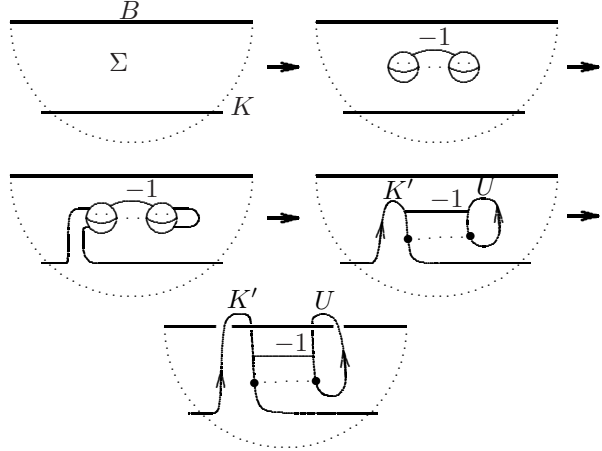

\begin{center}
\include{2HandleToRound1Handle}
\caption{Turning a $2$-handle into a round $1$-handle (proof of
  Proposition~\ref{P:QuickNEasy}).}
\label{F:2HandleToRound1Handle}
\end{center}
\end{figure}

\end{proof}

For the more general result we need to keep track of almost complex
structures and homotopy classes of plane fields associated to
fibrations and OBDs.

Given a B(A)LF $f\co X \rightarrow \Sigma$, we will use $R_f$ to
denote the union of the round $1$-handle singularities.  Given a
$4$-manifold $X$, let $\mathcal{J}(X)$ be the set of all almost
complex structures on $X$ modulo homotopy. Given a $3$-manifold $M$
let $\mathcal{Z}(M)$ be the set of all co-oriented plane fields on $M$
modulo homotopy. (This is of course equivalent to the set of all
nowhere-zero vector fields modulo homotopy, but we take the plane
field perspective because of the connections with contact topology.)
First note the following facts relating Lefschetz fibrations, almost
complex structures, open book decompositions and homotopy classes of
plane fields.
\begin{enumerate}
\item A BLF $f$ on a $4$-manifold $X$
  determines a homotopy class $j(f) \in \mathcal{J}(X \setminus
  R_f)$, characterized by having a representative $J \in j(f)$
  such that the fibers of $f$ are $J$-holomorphic curves. 
\item An OBD $f$ on a $3$-manifold $M$
  determines a homotopy class $z(f) \in \mathcal{Z}(M)$,
  characterized by having a representative which is positively
  transverse to a vector field $V$ which in turn is positively
  transverse to the pages of $f$ and positively tangent to the
  binding of $f$. This is the same as the homotopy class of the
  unique isotopy class of positive contact structures supported by
  $f$ in the sense of Giroux~\cite{Giroux}. 
\item If $X$ is a $4$-manifold and $M = \partial X$, then a homotopy
  class $j \in \mathcal{J}(X)$ determines a homotopy class $z(j) \in
  \mathcal{Z}(M)$ characterized by having a representative $\xi$ which
  is the field of $J$-complex tangencies to $M$ for some $J \in j$.
\item If $f$ is a convex BLF on a $4$-manifold $X$, inducing the
  OBD $f|_M$ on $M = \partial X$, then
  $z(j(f)) = z(f|_M)$.
\end{enumerate}

\begin{theorem} \label{T:ConvexACS} 
  Let $X$ be a $4$-dimensional $2$-handlebody, let $C$ be a
  (possibly empty) finite disjoint union of points and circles in the
  interior of $X$ and let $J$ be an almost complex structure on $X
  \setminus C$. Let $N$ be a given open neighborhood of $C$. Then
  there exists a convex BLF $f\co X \rightarrow B^2$ with the following
  properties:
\begin{itemize}
\item The union of the round $1$-handle singularities $R_f$ is
  contained in $N$.
\item For any almost complex structure $J' \in j(f)$, $J$ and $J'$
  will be homotopic on $X \setminus N$.
\item The positive contact structure supported by $f|_{\partial X}$
  is overtwisted.
\end{itemize}
\end{theorem}

At this point it is worth emphasizing that, to prove
Theorem~\ref{T:Closed}, we would be satisfied if
Theorem~\ref{T:ConvexACS} produced a BALF. However, we feel it is of
independent interest that we are able to avoid achirality on the
convex half of the construction.  Before we prove
Theorem~\ref{T:ConvexACS}, the corollary that we will actually use is:
\begin{corollary} \label{C:ConvexOBD}
Given any $4$-dimensional $2$-handlebody $X$ and any OBD $g\co
\partial X \rightarrow B^2$ which supports an overtwisted contact
structure, there exists a convex BLF $f\co X \rightarrow B^2$ such
that the open book $f|_{\partial X}$ is obtained from $g$ by a
sequence of positive stabilizations.
\end{corollary}

\begin{proof}
  The homotopy class of plane fields $z(f)$ on $\partial X$ determines
  a homotopy class of almost complex structures on a collar
  neighborhood of $\partial X$, which extends across all of $X$ except
  perhaps a finite disjoint union of points and circles. (This is
  because the space of almost complex structures on $\R^4$ respecting
  a given metric is $S^2$, so we only see obstructions to extending
  almost complex structures when we reach the $3$-skeleton.) Then
  Theorem~\ref{T:ConvexACS} produces the BLF $f'$, such that $\xi \in
  z(f'|_{\partial X})$. Eliashberg's classification of overtwisted
  contact structures~\cite{EliashOvertwisted} tells us that the
  contact structure supported by $f'|_{\partial X}$ is isotopic to
  $\xi$, and Giroux's results on contact structures and open books
  tell us that $f'|_{\partial X}$ and $f$ have a common positive
  stabilization (where stabilization is plumbing on left-handed Hopf
  bands). Lastly each stabilization of $f'|_{\partial X}$ can be
  implemented using Lemma~\ref{L:SimpleConvexStabilization}.
\end{proof}

\begin{figure}
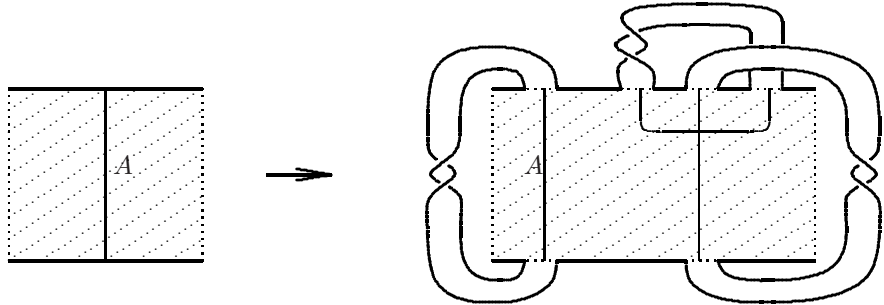

\begin{center}
\include{1right2lefts}
\caption{Plumbing on one right-handed Hopf band along an arc $A$,
  together with a left-handed Hopf band plumbed along a parallel copy
  of $A$ and a left-handed Hopf band plumbed along a short arc
  transverse to this parallel copy.}
\label{F:1right2lefts}
\end{center}
\end{figure}

To prove Theorem~\ref{T:ConvexACS} (producing a BLF rather than a
BALF) we need a way of negatively stabilizing OBDs on convex
boundaries without introducing anti-Lefschetz singularities.
Figure~\ref{F:1right2lefts} shows a modification of an OBD involving
plumbing one {\em right}-handed Hopf band along an arc $A$ in a page,
one left-handed Hopf band along a parallel copy of $A$, and one more
left-handed Hopf band along a short arc transverse to this parallel
copy. We will now show that this modification can be achieved using
round $1$-handles but no anti-Lefschetz singularities. One should
think of the following lemma as giving us the freedom to plumb on
right-handed Hopf bands wherever we want, avoiding achirality, at the
expense of introducing extraneous left-handed Hopf bands.

\begin{lemma} \label{L:1right2lefts}
  Given a convex (BA)LF $f \co X \rightarrow B^2$ and an arc $A$ in a page
  of the OBD on $\partial X$, there exists a B(A)LF $f' \co X
  \rightarrow B^2$ inducing
  the OBD indicated in Figure~\ref{F:1right2lefts}, which agrees
  with $f$ outside a neighborhood $U$ of $A$ and which has one
  Lefschetz singularity and one round 1-handle singularity inside $U$.
\end{lemma}

\begin{proof}
  Attach two cancelling $1$-$2$-handle pairs as on the left in
  Figure~\ref{F:1R2Lhandles} (so we have not changed the
  $4$-manifold).  Then observe, as on the right in
  Figure~\ref{F:1R2Lhandles}, that this configuration can also be seen
  as a $1$-handle with feet strung on the binding, a round $1$-handle
  with feet wrapping once around the binding and a $2$-handle whose
  foot is a knot in a page running over the round $1$-handle, with
  framing $\pf-1$. The monodromy of the new open book decomposition is
  indicated on the left in Figure~\ref{F:1R2Lmonodromy}. To see this,
  note that we would like to see both feet of the round $1$-handle as
  given by fixed points of the monodromy, but the left foot goes over
  the $1$-handle. However, if we introduce a left-handed Dehn twist
  and a right-handed Dehn twist along parallel curves that go along
  the arc $A$ and over the $1$-handle (the product of which is
  isotopic to the identity), the section determined by a fixed point
  in between the two twists is in fact the same as the left foot of
  the round $1$-handle. The extra right-handed Dehn twist in
  Figure~\ref{F:1R2Lmonodromy} comes from the $-1$ framed vanishing
  cycle $2$-handle. Figure~\ref{F:1R2Lmonodromy} then shows a two-step
  isotopy so that we see that the resulting monodromy agrees with the
  monodromy for Figure~\ref{F:1right2lefts}. (In these figures the
  indicated monodromy should be understood to be composed with any
  pre-existing monodromy coming from the initial open book
  decomposition.) Thus the new page is isotopic to that in
  Figure~\ref{F:1right2lefts}. (To go from a statement about the
  monodromy of an open book to a statement about the isotopy class of
  an open book is not safe in general. Here, however, we have the fact
  that the operation in question amounts to a Murasugi sum with an
  open book decomposition of $S^3$, and in $S^3$ open book
  decompositions are completely determined up to isotopy by their
  monodromy, since the mapping class group of $S^3$ is trivial.)

\begin{figure}
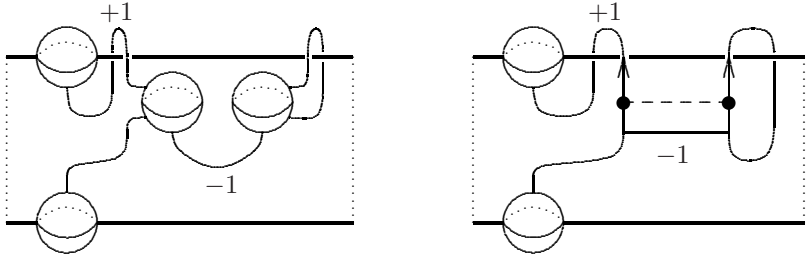

\begin{center}
\include{1R2Lhandles}
\caption{Two cancelling $1$-$2$-handle pairs becoming a $1$-handle, a
  round $1$-handle and a $2$-handle.}
\label{F:1R2Lhandles}
\end{center}
\end{figure}

\begin{figure}
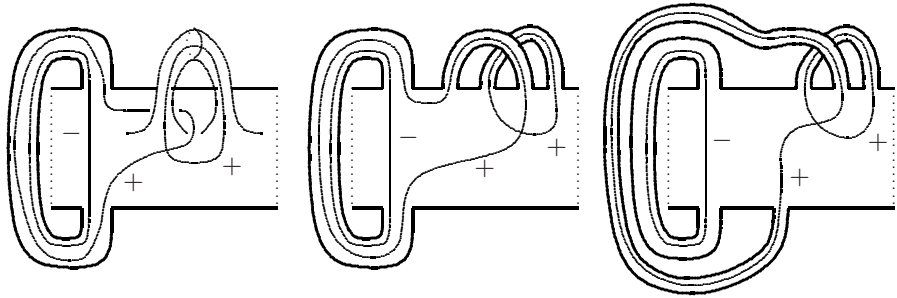

\begin{center}
\include{1R2Lmonodromy}
\caption{Three equivalent descriptions of the monodromy corresponding
  to Figure~\ref{F:1R2Lhandles}.}
\label{F:1R2Lmonodromy}
\end{center}
\end{figure}

\end{proof}

The last techniques we need to develop before the proof of
Theorem~\ref{T:ConvexACS} are techniques for computing Chern classes
of almost complex structures and invariants of co-oriented plane
fields in terms of BLFs and OBDs. We begin by collecting some relevant facts.

First, note that there is a well-defined connected sum operation $\#$
on $\{ (M^3,z) | z \in \mathcal{Z}(M)\}$, since any two plane fields
are locally homotopic. Likewise there is a well-defined boundary
connected sum operation $\natural$ on $\{ (X^4,j) | j \in
\mathcal{J}(X)\}$ which induces the connected sum on the boundary.
(These extend in the obvious way to self-connect sums and boundary
self-connect sums.) If one attaches a $1$-handle from one convex BLF
$(X_1,f_1)$ to another one $(X_2,f_2)$ such that the feet are strung
on the bindings as in Remark~\ref{R:Convex1HConcave3H}, giving a BLF
$f$ on $X_1 \natural X_2$, the resulting $j(f) \in \mathcal{J}((X_1
\natural X_2) \setminus R_f)$ is equal to the boundary connect sum
$j(f_1) \natural j(f_2)$, and $z(f) \in \mathcal{Z}(\partial X_1
\sharp \partial X_2)$ is equal to $z(f_1) \sharp z(f_2)$.

Next, we summarize some results from~\cite{Gompf}; a useful exposition
can also be found in~\cite{DingGeigesStipsicz}: 

There are two invariants $d_2$ and $d_3$ of $\mathcal{Z}(M)$, which as
a pair constitute a complete invariant. The ``$2$-dimensional
invariant'' $d_2$ of a given $z \in \mathcal{Z}(M)$ is simply the
spin$^\C$ structure determined by $z$; in the case where $H^2(M;\Z)$
has no $2$-torsion, this is completely characterized by $c_1(z) \in
H^2(M;\Z)$. In general, the set $\mathcal{S}(M)$ of spin$^\C$
structures on $M$ is an affine space for $H^2(M;\Z)$, and the action
of $H^2(M;\Z)$ on $\mathcal{S}(M)$ has the property that $c_1(a \cdot
s) = 2a + c_1(s)$ for $a \in H^2(M;\Z)$ and $s \in
\mathcal{S}(M)$. The ``$3$-dimensional invariant'' $d_3(z)$ lies in an
affine space for a cyclic group; the key properties of $d_3$ that we
need are summarized in the following two items.

Focusing on the case of $M=S^3$ (in which case there is only one
spin$^\C$ structure and so we need only pay attention to $d_3$)
suppose that $z_1, z_2 \in \mathcal{Z}(S^3)$ and that $S^3 = \partial
X_1 = \partial X_2$, with $j_i \in \mathcal{J}(X_i)$ such that
$j_i|_{S^3} = z_i$, for $i=1,2$. Then $c_1(j_1)^2 - 2\chi(X_1) -
3\sigma(X_1) = c_1(j_2)^2 - 2\chi(X_2) - 3\sigma(X_2)$ if and only if
$d_3(z_1) = d_3(z_2)$ (i.e. if and only if $z_1 = z_2$). Hence, in the
case of $S^3$, we identify $d_3(z)$ with $(c_1(j)^2 - 2\chi(X) -
3\sigma(X))/4 \in \Z - 1/2$, where $j$ is any extension of $z$ over a
$4$-manifold $X$. Now, for a general $3$-manifold $M$, if $z_1, z_2
\in \mathcal{Z}(M)$ and $d_2(z_1) = d_2(z_2)$, then there exists a $z
\in \mathcal{Z}(S^3)$ such that $d_3(z_2) = d_3(z_1 \# z)$. In
particular, $(M,z_2) = (M,z_1) \# (S^3,z)$. If $M = S^3$, then
$d_3(z_2) = d_3(z_1) + d_3(z) + 1/2$.

Now we summarize a discussion in~\cite{GayKirby} on constructing
almost complex structures in prescribed homotopy classes: 

Given an almost complex structure $J$ on a smooth manifold $X$, $J$
can always be trivialized over the $1$-skeleton of $X$. Then $c_1(J)$
is represented by the cocycle whose value on a $2$-cell $e$ is the
obstruction to extending this trivialization across $e$, as an element
of $\pi_1(GL_2 (\C)) = \Z$. Any two almost complex structures can be
made, via a homotopy, to agree on the $1$-skeleton. Given two almost
complex structures $J_1$ and $J_2$ over the $2$-skeleton which agree
on the $1$-skeleton, if their obstruction cocycles are equal for a
given trivialization over the $1$-skeleton then $J_1$ is homotopic to
$J_2$ on all of the $2$-skeleton. Thus, if we wish to construct a
given almost complex structure up to homotopy on a $2$-handlebody, we
must be able to construct an almost complex structure $J_1$ on the
$1$-skeleton with a trivialization and then, for any given cocycle
$c$, be able to extend $J_1$ to an almost complex structure $J$ on the
$2$-skeleton with $c$ as its obstruction cocycle. In the absence of
$2$-torsion in $H^2(X;\Z)$, this just amounts to getting $c_1(J)$
correct, but when there is $2$-torsion, there will be different
cocycles representing a fixed $c_1$ but corresponding to different
almost complex structures.

Next we combine some standard contact and symplectic topology and some
results from~\cite{DingGeigesStipsicz} to relate the above facts to
surgery and handle addition: 

Given a $3$-manifold $M$ and a homotopy class $z \in \mathcal{Z}(M)$,
suppose that $\xi \in z$ and that $K$ is a knot in $M$ tangent to
$\xi$. Then $K$ comes with a canonical framing $c$ given by $\xi$; let
$M'$ be the result of $c \pm 1$ surgery on $K$. Then there is a
well-defined $z' \in \mathcal{Z}(M')$ which can be characterized in
either of the following two equivalent ways:
\begin{enumerate}
\item Homotope $\xi$, remaining fixed along $K$, to be positive
  contact in a neighborhood of $K$. Then there is a unique contact
  $\pm 1$ surgery along $K$, producing $\xi'$ on $M'$, and we let $z'$
  be the homotopy class of $\xi'$.
\item Express a neighborhood $N$ of $K$ as $S^1 \times [-1,1] \times
  [-1,1]$, with $K = S^1 \times 0 \times 0$ and with $\xi$ tangent to
  $S^1 \times [-1,1] \times 0$ along $K$. Now homotope $\xi$,
  remaining fixed along $K$, to be tangent to the foliation $S^1
  \times [-1,1] \times t$ on all of $N$. As in~\cite{Lickorish}, $c
  \pm 1$ surgery along $K$ can be viewed as cutting $M$ open along
  $S^1 \times [-1,1] \times 0$ and reglueing via a left/right-handed
  Dehn twist along $K$. Thus the surgered neighborhood $N'$ in $M'$
  naturally inherits a foliation by annuli, and we let $\xi'$ be
  tangent to this foliation inside $N'$ and be equal to $\xi$ outside
  the surgery. Then we define $z'$ to be the homotopy class of $\xi'$.
\end{enumerate}

Now suppose that $X$ is a $4$-manifold with $\partial X = M$ and with
a given $j \in \mathcal{J}(M)$ restricting to $z \in
\mathcal{Z}(M)$. Let $\xi \in z$ with $K$ tangent to $\xi$ as above,
with canonical framing $c$.  Let $X'$ be the result of attaching a
$2$-handle $H$ along $K$ with framing $c \pm 1$, so that $\partial X'
= M'$ as above, and let $z' \in \mathcal{Z}(M')$ be as above.  Then,
in the case of $c - 1$ framing, there is a canonical extension $j'$ of
$j$ across $H$ so that $j'|_{M'} = z'$, and in the case of $c + 1$
framing, there is a canonical extension $j'$ of $j$ across $H
\setminus B$, where $B$ is a small ball in the interior of $H$, so
that $j'|_{M'} = z'$.  These extensions can be characterized as
follows:
\begin{enumerate}
\item In the case of $c - 1$ framing, identify $H = D^2 \times D^2$ as
  a subset of $\C^2$ via the orientation-preserving map $D^2 \times
  D^2 \ni ((x_1,x_2),(y_1,y_2)) \mapsto (x_1 + i y_1, x_2 - i y_2) =
  (z_1,z_2) \in \C^2$. Then $j'$ is represented by an almost complex
  structure $J' \in j'$ which equals the standard integrable complex
  structure on $H \subset \C^2$ and, when restricted to $X = X'
  \setminus H$, represents $j$. In particular, the fibers of the map
  $(z_1,z_2) \mapsto z_1^2 + z_2^2$ in $H \subset \C^2$ are
  $J'$-holomorphic.
\item In the case of $c + 1$ framing, identify $H = D^2 \times D^2$ as
  a subset of $\C^2$ via the orientation-{\em reversing} map $D^2
  \times D^2 \ni ((x_1,x_2),(y_1,y_2)) \mapsto (x_1 + i y_1, x_2 + i
  y_2) = (z_1,z_2) \in \C^2$. Then $j'$ is represented by an almost
  complex structure $J' \in j'$ defined everywhere except at $(0,0)
  \in H$ which, when restricted to $X = X' \setminus H$, represents
  $j$ and which is characterized on $H$ by the fact that the fibers of
  the map $(z_1,z_2) \mapsto z_1^2 + z_2^2$ are $J'$-holomorphic
  except at $(0,0)$. The ball $B$ is then a small ball around
  $(0,0)$. Although $j'$ does not extend across $B$, if we replace $B$
  with $\C P^2 \setminus B^4$ (i.e.  connect sum with $\C P^2$), then
  $j'$ does extend across $\C P^2 \setminus B^4$ so as to agree with
  the standard complex structure on $\C P^2$.
\end{enumerate}

Now suppose that, in the setting of the preceding paragraph, we are
also given a trivialization of $\xi$ in a neighborhood of $K$ (i.e.  a
non-vanishing section $v$ of $\xi$). This gives $K$ a rotation number
$\rot(K)$ (the winding number of $TK$ inside $\xi$ relative to the
trivialization). Suppose that $J \in j$ so that $\xi$ is the field of
$J$-complex tangencies to $M$; then we naturally get a trivialization
$(v,n)$ of $J$ in a neighborhood of $K$, where $n$ is the outward
normal to $M$. Let $J' \in j'$ agree with $J$ on $X$.  Then, in both
the case of $c - 1$ framing and $c + 1$ framing, the obstruction to
extending this trivialization of $J$ to a trivialization of $J'$, as
an element of $\pi_1(GL_2(\C)) = \Z$, is precisely
$\rot(K)$. (In~\cite{DingGeigesStipsicz} this is proved in the case
where $X = B^4$, $\xi$ is the standard contact structure on $S^3$, $v$
is defined on all of $S^3$, and $c \neq 0$. Note, however, that our
assertion is purely local to $K$ and $H$, and that, given any $\xi$ on
$S^1 \times B^2$ which is tangent to $S^1 \times \{0\}$, with any
trivialization $v$ of $\xi$, after a homotopy of $\xi$ fixed along $K$
there exists an embedding of $S^1 \times B^2$ into $S^3$ carrying
$\xi$ to the standard contact structure on $S^3$, taking $S^1 \times
\{0\}$ to a Legendrian knot with $\tb \neq 0$, and taking $v$ to a
trivialization which extends over all of $S^3$.)

Finally, if $X$ is equipped with a convex (BA)LF $f \co X \to B^2$ and
if $(X',f')$ is the (BA)LF resulting from attaching a $2$-handle along a
knot in a page of the induced OBD on $\partial X$ with framing $\pf
\pm 1$, then $j(f') = j(f)'$ in the sense that $j(f')$ is precisely
the canonical extension of $j(f)$ discussed above.

This gives us the following algorithm for computing the invariants of
a homotopy class $z \in \mathcal{Z}(M)$ associated to an open book
decomposition on a closed $3$-manifold $M$ in terms of a factorization
$h=\tau_1 \circ \ldots \circ \tau_n$ of the monodromy $h$ into Dehn
twists $\tau_i$ along curves $\gamma_i$ in the page $F$. (We hope some
readers may find this algorithm useful in other contexts; a similar
algorithm is spelled out in~\cite{EtnOzb}.)
\begin{enumerate}
\item Begin with a standard immersion of the page $F$ in $\R^2$ as a
disk with $2$-dimensional $1$-handles attached around the boundary.
\item This gives a trivialization of $TF$ coming from the
standard trivialization of $T\R^2$. Together with the standard
trivialization of $TB^2$, we get a trivialization of $T(F \times
B^2)$ which yields a trivialization of the standard almost complex
structure on $F \times B^2$.
\item Each Dehn twist curve $\gamma_i$ can be thought of as a curve in
$F \times p_i$, where $p_i \in S^1$; with respect to the above
trivialization, we get a rotation number $\rot(\gamma_i)$ which is
precisely the winding number of $\gamma_i$ as an immersed curve in
$\R^2$, seen via the immersion of $F$ in $\R^2$.
\item Now interpret the Dehn twist curves as attaching circles for
$2$-handles attached to $F \times D^2$, with framing $\pf-1$ for
each right-handed Dehn twist and framing $\pf+1$ for each left-handed
Dehn twist. This describes an ALF on a $4$-manifold $X$ with an almost complex
structure $J$ on the complement of $q$ points, where $q$ is the number
of left-handed Dehn twists, and $J|_{\partial X}$ induces the required
homotopy class $z$  of plane fields on $M = \partial X$. 
\item Then $J$ extends to an almost complex structure $J'$ on all of
$X' = X \sharp^q \C P^2$ which is standard on each $\C P^2$ summand,
and we still have $J'|_{\partial X'}=J|_{\partial X}$ inducing $z$ on
$M = \partial X'$. 
\item Now read off $c_1(J')$ as a cocycle from the rotation numbers of
each $\gamma_i$ and the fact that $c_1$ evaluates to $3$ on each generator of
$H_2(X')$ coming from a $\C P^2$ summand.
\item Now use the intersection form on $X'$ to identify the
Poincar\'{e} dual of $c_1(J')$ and hence compute $c_1(J')|_{\partial
X'}$ to get $d_2(z)$ and compute $\chi(X')$, $\sigma(X')$ and
$c_1(J')^2$ to get $d_3(z) = (c_1(J')^2 - 2\chi(X') - 3\sigma(X'))/4$.
\item The last two steps are equivalent to the following shortcut:
  Read off $c_1(J)$ from the rotation numbers of each $\gamma_i$. Use
  the intersection form on $X$ to identify $c_1(J)$ and
  $c_1(J)|_{\partial X}$ to get $d_2(z)$. Then compute $\chi(X)$,
  $\sigma(X)$ and $c_1(J)^2$ to get $d_3(z) = (c_1(J')^2 - 2\chi(X') -
  3\sigma(X'))/4+ q$.
\end{enumerate}

\begin{proof}[Proof of Theorem~\ref{T:ConvexACS}]

  First we will prove the theorem when $X=B^4$ and $C$ is
  a point. Then we will prove it when $X=S^1 \times B^3$ and
  $C=S^1 \times \{0\}$. Finally we will prove the general case. Before
  beginning, however, note that the ability to plumb on right-handed
  Hopf bands, as in Lemma~\ref{L:1right2lefts}, immediately gives us
  the last assertion of the Theorem, that we can arrange for our
  contact structures to be overtwisted.

  \textbf{Trivial case: $X=B^4$ and $C=\emptyset$.} Here there is only
  one almost complex structure, achieved by the fibration $B^2 \times
  B^2 \rightarrow B^2$.

  \textbf{Simplest nontrivial case: $X=B^4$ and $C$ is a single
  point.} In this case all we need to do is to construct a broken
  Lefschetz fibration on $B^4$ inducing a given homotopy class of
  plane fields on $S^3$. Recall that $\mathcal{Z}(S^3)$ is in
  one-to-one correspondence with $\Z - 1/2$ via the formula $d_3(z) =
  (c_1(j)^2 - 2\chi(X) - 3\sigma(X))/4$, where $j$ is an extension of
  $z$ over a $4$-manifold $X$. Thus, suppose we are given $n \in \Z$,
  and we wish to construct a convex BLF on $B^4$ inducing a given $z
  \in \mathcal{Z}(S^3)$ with $d_3(z) = n - 1/2$.

  It is well known that plumbing on a left-handed Hopf band will not
  change $d_3$ (in fact it does not change the isotopy class of the
  contact structure~\cite{Giroux, Torisu}), while plumbing on a
  right-handed Hopf band increases $d_3$ by one~\cite{Torisu}.
  Furthermore the trivial fibration $B^4 = B^2 \times B^2 \rightarrow
  B^2$ yields $d_3 = -1/2$. Thus, using Lemma~\ref{L:1right2lefts} we
  can achieve our goal for any $n \geq 0$. By the comments on
  connected sums above, we now need only perform the construction for
  some negative value of $n$ to complete the proof when $X=B^4$. 

  We give a construction explicitly in Figure~\ref{F:B4nm1} for
  $n=-1$, i.e. $d_3 = -3/2$; the figure should be interpreted as
  follows: The topmost diagram shows a page of an open book
  decomposition of $S^3$ involving $2$ left-handed Hopf bands and $2$
  right-handed Hopf bands plumbed in sequence onto a disk, so that the
  page is a $4$-punctured disk.  Each right-handed Hopf band should
  really have an extra pair of left-handed Hopf bands immediately
  adjacent, as in Figure~\ref{F:1right2lefts}, but we have suppressed
  this extra pair since they play no further role in the
  construction. This open book decomposition of $S^3$ (including the
  $4$ extra left-handed Hopf bands not drawn) is thus the boundary of
  a convex BLF, using Lemmas~\ref{L:1right2lefts}
  and~\ref{L:SimpleConvexStabilization}. To this we add a $1$-handle
  strung along the binding, a round $1$-handle which wraps around the
  binding once, and a $2$-handle on a page with framing $\pf-1$
  running over the round $1$-handle, as in the Figure. This gives a
  more complicated convex BLF on $B^4$.

  To see that we have achieved $n=-1$, we first analyze the monodromy
  of the new open book decomposition of $S^3$, exactly as in the proof
  of Lemma~\ref{L:1right2lefts}, Figure~\ref{F:1R2Lmonodromy}. We need
  to introduce pairs of right- and left-handed Dehn twists parallel to
  and on either side of the two feet of the round $1$-handle to
  compensate for the fact that the feet are not initially described as
  fixed points of the monodromy. This is indicated in the middle
  diagram in Figure~\ref{F:B4nm1}. The Dehn twist curves are labelled
  and oriented for use in the calculation to come. We now use this
  factorization of the monodromy into Dehn twists to compute $d_3$ as
  in the algorithm explained above. This describes a new $4$-manifold
  shown in the bottom diagram in the figure; each right- (resp. left-)
  handed Dehn twist has become a $2$-handle on a page with framing
  $-1$ (resp. $+1$), attached to an open book with page a
  $6$-punctured torus and monodromy equal to the identity. We note
  that $H_2$ is generated by $A = c - a_1 -b_1 - a_2 - b_2$ and $B =
  d+f-g-a_1-b_1 - a_2 - b_2$, with $A^2 = 1$ and $B^2 = -1$ and $A
  \cdot B = 0$. Reading off rotation numbers we see that $c_1(A) = -3$
  and $c_1(B) = -5$ so that $c_1$ is Poincar\'{e} dual to $-3 A + 5 B$
  and $c_1^2 = -16$. Also, $\sigma = 0$, $\chi = 3$ and the number of
  left-handed Dehn twists is $q=4$. Thus a final calculation gives
  $d_3 = (c_1^2 - 2\chi - 3\sigma)/4 + q = (-16-6)/4+4 = -3/2$.

\begin{figure}
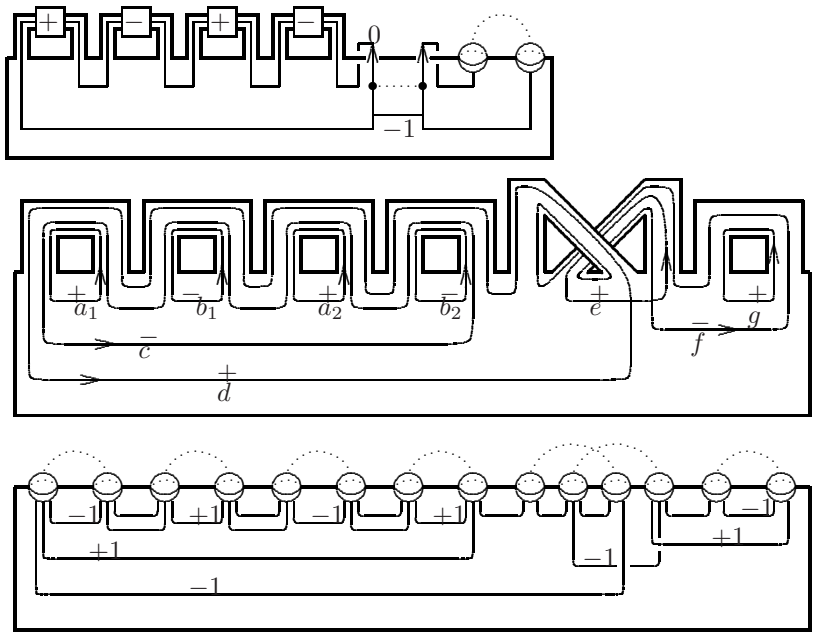

\begin{center}
\include{B4nm1}
\caption{A broken Lefschetz fibration on $B^4$ for the case $n=-1$,
  i.e. $d_3 = -3/2$.}
\label{F:B4nm1}
\end{center}
\end{figure}

  \textbf{Next simplest case: $X=S^1 \times B^3$ and $C=S^1 \times
  \{0\}$.} Now we need to construct a convex BLF on
  $X = S^1 \times B^3$ inducing a given homotopy class of plane fields
  on $S^1 \times S^2$. By the comments earlier on the $3$-dimensional
  invariant and connected sums of broken Lefschetz fibrations, if we
  get the $2$-dimensional invariant correct then we can use the case
  above for $B^4$ to get the $3$-dimensional invariant correct. Thus
  we need to construct a convex BLF $f$ on $X$
  such that $c_1(j(f)|_{\partial X}) = 2n$ for any given $n \in \Z =
  H^2(S^1 \times S^2)$. Note that we do not need to worry about the
  potential sign ambiguities associated with the identification of
  $\Z$ with $H^2(S^1 \times S^2)$ because there is an
  orientation-preserving automorphism of $S^1 \times B^3$ which
  induces multiplication by $-1$ on $H^2(S^1 \times S^2)$. So we can
  simplify the problem slightly to say that, given any non-negative
  integer $n$, we should construct $f$ so that
  $|c_1(z(f|_{\partial X}))| = 2n$. If $n=0$ the fibration is $S^1
  \times [0,1] \times D^2 \rightarrow D^2$.  
\begin{figure}
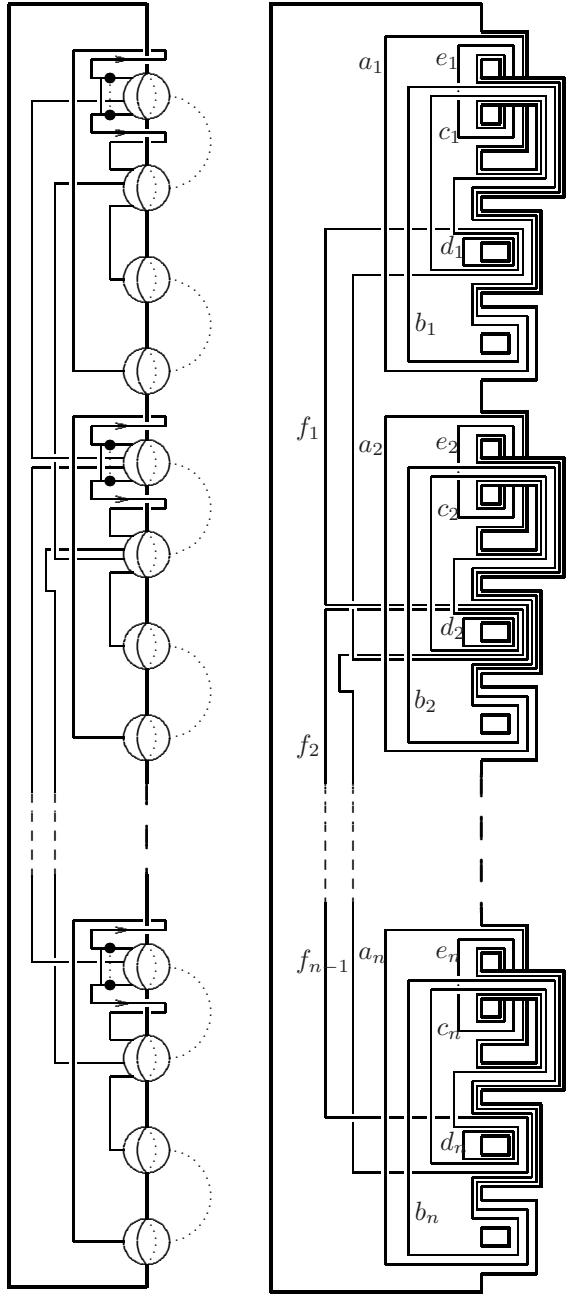

\begin{center}
\include{S1B3case}
\caption{A broken Lefschetz fibration on $S^1 \times B^3$.}
\label{F:S1B3case}
\end{center}
\end{figure}
\normalsize
  Figure~\ref{F:S1B3case} is an explicit example for $n > 0$, and
  should be interpreted as follows: The left diagram is a (round)
  handlebody decomposition of $S^1 \times B^3$, starting with $S^3 =
  \partial B^4$ with the binding of the standard open book indicated
  by the heavy lines, and involving $2n$ $1$-handles strung on this
  binding, $n$ round $1$-handles each wrapping once around the
  binding, $n$ $2$-handles with framing $\pf-1$ each running over one
  of the round $1$-handles, and $n-1$ $2$-handles with framing $\pf-1$
  each connecting one $1$-handle to another. (The framings are not
  indicated in the diagram.) This describes a convex BLF $f\co X
  \rightarrow D^2$. Again, using the techniques of
  Lemma~\ref{L:1right2lefts}, we compute the monodromy which is
  indicated in the right diagram. Here the curves indicate Dehn twists
  but their handedness is not indicated in the diagram. Their
  handedness is as follows: The curves labelled $a_i$, $c_i$, $e_i$
  and $f_i$ are right-handed Dehn twists and the curves labelled $b_i$
  and $e_i$ are left-handed Dehn twists.

  We will now compute $c = c_1(z(f|_{\partial X}))$.  Orient each
  curve so that its lowermost straight line segment is oriented
  left-to-right. With this orientation, we have that $\rot(a_i)=-1$
  and $\rot(b_i)=\rot(c_i)=\rot(d_i)=\rot(f_i) = +1$. Now we convert
  the monodromy diagram into a handlebody decomposition of a new
  $4$-manifold so that the $1$-handles of the page become
  $4$-dimensional $1$-handles and the Dehn twist curves become
  attaching circles for $2$-handles, with $a_i$, $c_i$, $e_i$ and
  $f_i$ framed $-1$ and $b_i$ and $d_i$ framed $+1$. Then we see that
  $H_2$ is generated by $A_1, \ldots, A_n$ and $F_1, \ldots, F_{n-1}$
  where $A_i = a_i-b_i+c_i-d_i$ and $F_i = f_i-d_i+d_{i+1}$. All
  intersections between generators are $0$ except for $F_i \cdot F_i =
  1$, $A_i \cdot F_i = 1$ and $F_{i+1} \cdot A_i = -1$. Thus $H_2$ of
  the boundary $S^1 \times S^2$ is generated by $A = A_1 + A_2 +
  \ldots + A_n$. Finally, we evaluate $c$ on $A$ using the rotation
  numbers above to get $c(A) = c(A_1)+c(A_2)+\ldots + c(A_n)$, with
  $c(A_i) = \rot(a_i)-\rot(b_i)+\rot(c_i)-\rot(d_i) = -2$, so that
  $c(A) = -2n$. Thus $|c| = 2n \in H^2(S^1 \times S^2) = \Z$.

  \textbf{General case:} Now we are given a general $2$-handlebody $X$
  and a collection $C \subset (X \setminus \partial X)$ of $m$ points
  and $n$ circles. Choose a handlebody decomposition of $X$ involving
  $m$ $0$-handles, $n$ copies of $S^1 \times B^3$ (i.e. $n$
  $0$-handles and $n$ $1$-handles), and then some more $1$- and
  $2$-handles, so that each of the $m$ points is contained in one of
  the $m$ $0$-handles and each of the $n$ circles is the core of one
  of the copies of $S^1 \times B^3$. We apply the above cases to each
  of the $m$ $0$-handles and each of the $n$ copies of $S^1 \times
  B^3$. There is only one way to extend an almost complex structure
  across a $1$-handle, so now we need to get the almost complex
  structure right as we extend across each $2$-handle. By comments
  above, this is simply a matter of getting the rotation number
  correct for each $2$-handle's attaching circle, relative to a given
  trivialization of a page. Having the freedom to plumb on both left-
  and right-handed Hopf bands is precisely what makes this
  possible. (See~\cite{EF}, for example.) The final subtlety is that,
  each time we plumb on a right-handed Hopf band, we introduce more
  round $1$-handle singularities. However, each one lies in a ball,
  and so at the very end we have some new balls across which the
  almost complex structure may not extend. This can easily be
  compensated for, however, by performing one more connected sum with
  a judiciously chosen convex BLF on $B^4$ as described at the
  beginning of this proof.
\end{proof}
\section{Proof of the main result}
\label{S:ClosedProof}

We need one more trick to complete the proof of our main result,
namely the ability to stabilize OBDs on concave boundaries.

\begin{lemma}  \label{L:ConcaveStabilization}
Given a concave fibration $f: X^4 \rightarrow \Sigma$ and an arc $A$
in a page of the induced OBD $f: \partial X \rightarrow B^2$. Let $f':
\partial X \rightarrow B^2$
be the result of positively (resp. negatively) stabilizing this OBD
along $A$. Then $f'$ extends to a concave fibration $f' : X
\rightarrow \Sigma$ which agrees with $f$ outside a ball neighborhood of
$A$, inside which the only singularities are a round $1$-handle
singularity and a Lefschetz (resp. anti-Lefschetz) singularity.
\end{lemma}

\begin{proof}  
Let  $\alpha_0$ and $\alpha_1$ be the endpoints of $A$.
Now add a cancelling pair of $1-$ and $2$-handles where the feet of
the one handle lie on $A$ near $\alpha_0$ and $\alpha_1$ and where
the attaching map of the $2$-handle goes over the $1$-handle once and
follows $A$, with framing $\pf -1$ (resp. $\pf +1$).

Now add a cancelling $2$-$3$-handle pair and proceed as in the proof
of Lemma~\ref{L:Concave1H} to turn the one handle into a round $1$-handle
(adding a $T^2$ to the fiber) whereupon the $3$-handle removes an
annulus from the page leaving the following: the page has had a
$1$-handle attached with a Dehn twist along $\alpha$, right (left)
handed if the framing was page framing -1 (+1).  This is exactly the
stabilization that was desired.
\end{proof}

\begin{proof}[Proof of Theorem~\ref{T:Closed}] 

Split $X$ as $A \cup B$ where $A$ is the result of attaching some
number of $1$-handles to a neighborhood $F \times B^2$ of $F$ and $B$
is a $2$-handlebody. Lemma~\ref{L:ConcaveFXB2} gives a concave BLF $f
\co F \times B^2 \to S^2$, which we extend across the rest of the
$1$-handles of $A$ using Lemma~\ref{L:Concave1H}. Use
lemma~\ref{L:ConcaveStabilization} to positively stabilize
the induced OBD $f|_{\partial A} \co \partial A \to B^2$. Now shift
attention to $B$, where $\partial B = -\partial A$, and consider the
problem of extending the given OBD on $\partial A$ across $B$. First
note that the (positive) contact structure supported by this OBD on
$\partial B$ is in fact overtwisted, precisely because it resulted
from a positive stabilization on $\partial A$, which is a negative
stabilization on $\partial B = -\partial A$. Thus
Corollary~\ref{C:ConvexOBD} gives us a convex BLF $g \co B \to B^2$
which induces on OBD on $\partial B$ which is the result of positively
stabilizing the given OBD coming from $-\partial A$. Note that at this
point we have a concave BLF on $A$ and a convex BLF on $B$, which do
not quite match because we need to positively stabilize (in the sense
of the orientation coming from $B$) the OBD coming from $A$, i.e. we
need to negatively stabilize the OBD on $\partial A$. Thus we
use Lemma~\ref{L:ConcaveStabilization} one more time, but finally we
are forced to introduce achirality, in order to achieve these negative
stabilizations, and then the BALF on $A$ can be glued to the BLF on $B$.
\end{proof}

Note that if we could find a trick for negatively stabilizing the
concave side without introducing anti-Lefschetz singularities, we
would be able to avoid achirality completely. The authors did find
some tricks analogous to Lemma~\ref{L:1right2lefts} that work on the
concave side, but these always involved extraneous extra positive
stabilizations which we could not control.

\begin{proof}[Proof of addendum to Theorem~\ref{T:Closed}]
Here we are given the extra data of some $2$-spheres $S_1, \ldots,
S_n$ which should become sections of the BALF. In this case split $X$
as $A \cup B$ where $A$ is $F \times B^2$ together with $n$
$2$-handles attached along sections $p_i \times S^1$ of $F \times
S^1$, and some extra $1$-handles so that the complement is a
$2$-handlebody, and so that the $2$-handles give the spheres
$S_i$. Start with the flat fibration $F \times B^2 \to B^2$, then
attach the section $2$-handles as in Remark~\ref{R:Section2H}, to get
concave boundary, and then attach the $1$-handles as in
Lemma~\ref{L:Concave1H}. This gives the concave piece, and proceed as
in the proof above to put a convex BLF on the complement, and to make
the open books match by appropriate stabilizations.

To arrange that the round $1$-handle singularities all lie over the
tropics of Cancer and Capricorn, notice that the only place in our
construction where the attaching circles for a round $1$-handle might
run over another round $1$-handle is in the negative stabilizations on
the convex side (Lemma~\ref{L:1right2lefts}). However, if we do not
try to keep the convex side chiral, we can achieve these
stabilizations with anti-Lefschetz singularities rather than round
$1$-handle singularities. Then the round $1$-handle singularities on
the convex side are independent and therefore can lie over the tropic
of Capricorn, and those on the concave side are also independent and
can lie over the tropic of Cancer. Finally, on each side, the
vanishing cycle $2$-handles can always be attached after attaching the
round $1$-handles, so we can arrange for them to project to the
equator.

\end{proof}

\section{Questions}

\label{S:Questions}

\begin{question}
The most basic question is, ``What is the uniqueness theorem?'' Many
choices are made in the construction of a BALF; if different choices
are made, what is the set of moves relating the two BALFs? These
should include, for example, the positive and negative stabilizations
used to match the convex and concave pieces, and adding cancelling
round $1$-$2$-handle pairs. A critical ingredient would be a
uniqueness statement for the sequences of stabilizations coming from
Giroux's results.
\end{question}

\begin{question}
Another question is whether achirality can be avoided. By the results
in~\cite{ADK}, if $b_2^+(X) > 0$ and we blow up enough, then this can
be done; but even in this case we do not have a constructive proof.

Achirality could be avoided in general if we could find a way to
positively and negatively stabilize the concave side using only
Lefschetz and round $1$-handle singularities. If this cannot be done,
there ought to be an obstruction which lies in the set of equivalence
classes of OBDs on connected sums of $S^1 \times S^2$'s, where the
equivalence relation is derived from the basic moves in the uniqueness
question above. Even better would be to push this obstruction to the
$S^3$ boundary of the $4$-handle.
\end{question}

\begin{question}
In~\cite{ADK} it is shown that a BLF supports a near-symplectic form
as long as there is a $2$-dimensional cohomology class evaluating
positively on the fibers. (This is a closed $2$-form vanishing
identically along the round $1$-handle singularities, symplectic in
their complement, and satisfying a certain transversality along the
circles.) Does this generalize meaningfully to the case of a BALF?
What control on the $2$-form can we expect near the anti-Lefschetz
singularities? Baykur~\cite{Baykur} has used ALFs to construct
folded symplectic structures.
\end{question}

\begin{question}
To what extent does achirality destroy Perutz's program~\cite{Perutz,
  Perutz1, Perutz2} to generalize the
  Donaldson-Smith-Usher~\cite{DonSmith, Usher} results relating smooth
  $4$-manifold invariants to counts of multisections of Lefschetz
  fibrations? Perutz proposes to count multisections of BLFs (some of
  which may limit on round $1$-handle singularities); see
  also~\cite{Baykur2}.
\end{question}

\begin{question}
A smooth, simply-connected $5$-dimensional $h$-cobordism is a product
off of a contractible manifold, which is an $h$-cobordism between two
contractible $4$-manifolds, $A_0$ and $A_1$~\cite{CFHS,
KirbyCorks, Matveyev}. These $4$-manifolds, called Akbulut's corks, are
constructed from a symmetric link of $0$-framed unknots by changing
half the unknots to $1$-handles (replacing the $0$ by a dot), or the
other half to $1$-handles. There is an involution $h \co \partial A_0
\to \partial A_1 = \partial A_0$ which does not extend to a
diffeomorphism from $A_0$ to $A_1$.

The question here is whether each of $A_0$ and $A_1$ are convex B(A)LFs such
that the involution $h$ preserves the induced OBD on the boundary, so
that the process of surgering out $A_0$ and replacing with $A_1$ can
be carried out on closed B(A)LFs without changing the fibration outside
$A_0$.
\end{question}

\bibliographystyle{gtart}
\bibliography{singPALF}
\end{document}

%% file: RoundExample.tex
\font\thinlinefont=cmr5
\mbox{\beginpicture
\small
\setcoordinatesystem units <0.8cm,0.8cm>
\unitlength=1.04987cm
\linethickness=1pt
\setplotsymbol ({\makebox(0,0)[l]{\tencirc\symbol{'160}}})
\setshadesymbol ({\thinlinefont .})
\setlinear
%
%
\linethickness= 0.500pt
\setplotsymbol ({\thinlinefont .})
{\color[rgb]{0,0,0}\circulararc 73.740 degrees from  1.143 22.384 center at  1.429 22.765
}%
%
%
\linethickness= 0.500pt
\setplotsymbol ({\thinlinefont .})
\setdots < 0.0953cm>
{\color[rgb]{0,0,0}\circulararc 73.740 degrees from  1.715 22.384 center at  1.429 22.003
}%
%
%
\linethickness= 0.500pt
\setplotsymbol ({\thinlinefont .})
\setsolid
{\color[rgb]{0,0,0}\ellipticalarc axes ratio  0.286:0.286  360 degrees 
	from  1.715 22.384 center at  1.429 22.384
}%
%
%
\linethickness= 0.500pt
\setplotsymbol ({\thinlinefont .})
{\color[rgb]{0,0,0}\circulararc 73.740 degrees from  1.143 20.955 center at  1.429 21.336
}%
%
%
\linethickness= 0.500pt
\setplotsymbol ({\thinlinefont .})
\setdots < 0.0953cm>
{\color[rgb]{0,0,0}\circulararc 73.740 degrees from  1.715 20.955 center at  1.429 20.574
}%
%
%
\linethickness= 0.500pt
\setplotsymbol ({\thinlinefont .})
\setsolid
{\color[rgb]{0,0,0}\ellipticalarc axes ratio  0.286:0.286  360 degrees 
	from  1.715 20.955 center at  1.429 20.955
}%
%
%
\linethickness= 0.500pt
\setplotsymbol ({\thinlinefont .})
{\color[rgb]{0,0,0}\circulararc 73.740 degrees from  2.095 21.907 center at  2.381 22.289
}%
%
%
\linethickness= 0.500pt
\setplotsymbol ({\thinlinefont .})
\setdots < 0.0953cm>
{\color[rgb]{0,0,0}\circulararc 73.740 degrees from  2.667 21.907 center at  2.381 21.526
}%
%
%
\linethickness= 0.500pt
\setplotsymbol ({\thinlinefont .})
\setsolid
{\color[rgb]{0,0,0}\ellipticalarc axes ratio  0.286:0.286  360 degrees 
	from  2.667 21.907 center at  2.381 21.907
}%
%
%
\linethickness= 0.500pt
\setplotsymbol ({\thinlinefont .})
{\color[rgb]{0,0,0}\circulararc 73.740 degrees from  3.048 21.907 center at  3.334 22.289
}%
%
%
\linethickness= 0.500pt
\setplotsymbol ({\thinlinefont .})
\setdots < 0.0953cm>
{\color[rgb]{0,0,0}\circulararc 73.740 degrees from  3.620 21.907 center at  3.334 21.526
}%
%
%
\linethickness= 0.500pt
\setplotsymbol ({\thinlinefont .})
\setsolid
{\color[rgb]{0,0,0}\ellipticalarc axes ratio  0.286:0.286  360 degrees 
	from  3.619 21.907 center at  3.334 21.907
}%
%
%
\linethickness= 0.500pt
\setplotsymbol ({\thinlinefont .})
{\color[rgb]{0,0,0}\circulararc 73.740 degrees from  7.334 22.384 center at  7.620 22.765
}%
%
%
\linethickness= 0.500pt
\setplotsymbol ({\thinlinefont .})
\setdots < 0.0953cm>
{\color[rgb]{0,0,0}\circulararc 73.740 degrees from  7.906 22.384 center at  7.620 22.003
}%
%
%
\linethickness= 0.500pt
\setplotsymbol ({\thinlinefont .})
\setsolid
{\color[rgb]{0,0,0}\ellipticalarc axes ratio  0.286:0.286  360 degrees 
	from  7.906 22.384 center at  7.620 22.384
}%
%
%
\linethickness= 0.500pt
\setplotsymbol ({\thinlinefont .})
{\color[rgb]{0,0,0}\circulararc 73.740 degrees from  7.334 20.955 center at  7.620 21.336
}%
%
%
\linethickness= 0.500pt
\setplotsymbol ({\thinlinefont .})
\setdots < 0.0953cm>
{\color[rgb]{0,0,0}\circulararc 73.740 degrees from  7.906 20.955 center at  7.620 20.574
}%
%
%
\linethickness= 0.500pt
\setplotsymbol ({\thinlinefont .})
\setsolid
{\color[rgb]{0,0,0}\ellipticalarc axes ratio  0.286:0.286  360 degrees 
	from  7.906 20.955 center at  7.620 20.955
}%
%
%
\linethickness= 0.500pt
\setplotsymbol ({\thinlinefont .})
\setdots < 0.0953cm>
{\color[rgb]{0,0,0}\plot  7.620 22.098  7.620 21.241 /
}%
%
%
\linethickness=3pt
\setplotsymbol ({\makebox(0,0)[l]{\tencirc\symbol{'162}}})
\setsolid
{\color[rgb]{0,0,0}\plot  8.572 21.907  8.572 21.907 /
%
%
\linethickness=3pt
\setplotsymbol ({\makebox(0,0)[l]{\tencirc\symbol{'162}}})
\color[rgb]{0,0,0}\plot  9.525 21.907  9.525 21.907 /
%
%
\linethickness= 0.500pt
\setplotsymbol ({\thinlinefont .})
\color[rgb]{0,0,0}\putrule from  8.572 20.955 to  9.525 20.955
}%
%
%
\linethickness= 0.500pt
\setplotsymbol ({\thinlinefont .})
\setdots < 0.0953cm>
{\color[rgb]{0,0,0}\plot  8.572 21.907  9.525 21.907 /
}%
%
%
\linethickness= 0.500pt
\setplotsymbol ({\thinlinefont .})
\setsolid
{\color[rgb]{0,0,0}\putrule from  8.572 21.907 to  8.572 22.384
%
%
\plot  8.636 22.130  8.572 22.384  8.509 22.130 /
}%
%
%
\linethickness= 0.500pt
\setplotsymbol ({\thinlinefont .})
{\color[rgb]{0,0,0}\putrule from  9.525 21.907 to  9.525 22.384
%
%
\plot  9.588 22.130  9.525 22.384  9.462 22.130 /
}%
%
%
\linethickness= 0.500pt
\setplotsymbol ({\thinlinefont .})
{\color[rgb]{0,0,0}\putrule from  7.620 22.669 to  7.620 22.672
\putrule from  7.620 22.672 to  7.620 22.680
\putrule from  7.620 22.680 to  7.620 22.699
\putrule from  7.620 22.699 to  7.620 22.729
\putrule from  7.620 22.729 to  7.620 22.761
\putrule from  7.620 22.761 to  7.620 22.792
\putrule from  7.620 22.792 to  7.620 22.822
\putrule from  7.620 22.822 to  7.620 22.849
\putrule from  7.620 22.849 to  7.620 22.879
\putrule from  7.620 22.879 to  7.620 22.907
\putrule from  7.620 22.907 to  7.620 22.930
\putrule from  7.620 22.930 to  7.620 22.955
\plot  7.620 22.955  7.622 22.981 /
\plot  7.622 22.981  7.624 23.006 /
\plot  7.624 23.006  7.628 23.036 /
\plot  7.628 23.036  7.635 23.065 /
\plot  7.635 23.065  7.641 23.093 /
\plot  7.641 23.093  7.652 23.122 /
\plot  7.652 23.122  7.667 23.150 /
\plot  7.667 23.150  7.684 23.175 /
\plot  7.684 23.175  7.703 23.199 /
\plot  7.703 23.199  7.724 23.222 /
\plot  7.724 23.222  7.749 23.241 /
\plot  7.749 23.241  7.779 23.256 /
\plot  7.779 23.256  7.808 23.271 /
\plot  7.808 23.271  7.840 23.281 /
\plot  7.840 23.281  7.876 23.292 /
\plot  7.876 23.292  7.916 23.298 /
\plot  7.916 23.298  7.959 23.307 /
\plot  7.959 23.307  8.003 23.311 /
\plot  8.003 23.311  8.050 23.313 /
\plot  8.050 23.313  8.096 23.315 /
\plot  8.096 23.315  8.143 23.313 /
\plot  8.143 23.313  8.189 23.311 /
\plot  8.189 23.311  8.234 23.307 /
\plot  8.234 23.307  8.276 23.298 /
\plot  8.276 23.298  8.316 23.292 /
\plot  8.316 23.292  8.352 23.281 /
\plot  8.352 23.281  8.384 23.271 /
\plot  8.384 23.271  8.414 23.256 /
\plot  8.414 23.256  8.435 23.245 /
\plot  8.435 23.245  8.454 23.233 /
\plot  8.454 23.233  8.471 23.218 /
\plot  8.471 23.218  8.486 23.201 /
\plot  8.486 23.201  8.501 23.182 /
\plot  8.501 23.182  8.513 23.161 /
\plot  8.513 23.161  8.524 23.135 /
\plot  8.524 23.135  8.534 23.108 /
\plot  8.534 23.108  8.543 23.078 /
\plot  8.543 23.078  8.551 23.044 /
\plot  8.551 23.044  8.556 23.006 /
\plot  8.556 23.006  8.562 22.966 /
\plot  8.562 22.966  8.564 22.924 /
\plot  8.564 22.924  8.568 22.877 /
\plot  8.568 22.877  8.570 22.828 /
\putrule from  8.570 22.828 to  8.570 22.777
\plot  8.570 22.777  8.572 22.722 /
\putrule from  8.572 22.722 to  8.572 22.665
\putrule from  8.572 22.665 to  8.572 22.606
\putrule from  8.572 22.606 to  8.572 22.543
\putrule from  8.572 22.543 to  8.572 22.494
\putrule from  8.572 22.494 to  8.572 22.445
\putrule from  8.572 22.445 to  8.572 22.392
\putrule from  8.572 22.392 to  8.572 22.337
\putrule from  8.572 22.337 to  8.572 22.278
\putrule from  8.572 22.278 to  8.572 22.219
\putrule from  8.572 22.219 to  8.572 22.155
\putrule from  8.572 22.155 to  8.572 22.092
\putrule from  8.572 22.092 to  8.572 22.024
\putrule from  8.572 22.024 to  8.572 21.956
\putrule from  8.572 21.956 to  8.572 21.886
\putrule from  8.572 21.886 to  8.572 21.814
\putrule from  8.572 21.814 to  8.572 21.742
\putrule from  8.572 21.742 to  8.572 21.668
\putrule from  8.572 21.668 to  8.572 21.596
\putrule from  8.572 21.596 to  8.572 21.524
\putrule from  8.572 21.524 to  8.572 21.452
\putrule from  8.572 21.452 to  8.572 21.383
\putrule from  8.572 21.383 to  8.572 21.315
\putrule from  8.572 21.315 to  8.572 21.247
\putrule from  8.572 21.247 to  8.572 21.184
\putrule from  8.572 21.184 to  8.572 21.120
\putrule from  8.572 21.120 to  8.572 21.061
\putrule from  8.572 21.061 to  8.572 21.002
\putrule from  8.572 21.002 to  8.572 20.947
\putrule from  8.572 20.947 to  8.572 20.894
\putrule from  8.572 20.894 to  8.572 20.845
\putrule from  8.572 20.845 to  8.572 20.796
\putrule from  8.572 20.796 to  8.572 20.733
\putrule from  8.572 20.733 to  8.572 20.673
\putrule from  8.572 20.673 to  8.572 20.616
\plot  8.572 20.616  8.570 20.561 /
\putrule from  8.570 20.561 to  8.570 20.511
\plot  8.570 20.511  8.568 20.462 /
\plot  8.568 20.462  8.564 20.415 /
\plot  8.564 20.415  8.562 20.373 /
\plot  8.562 20.373  8.556 20.333 /
\plot  8.556 20.333  8.551 20.295 /
\plot  8.551 20.295  8.543 20.261 /
\plot  8.543 20.261  8.534 20.231 /
\plot  8.534 20.231  8.524 20.204 /
\plot  8.524 20.204  8.513 20.178 /
\plot  8.513 20.178  8.501 20.157 /
\plot  8.501 20.157  8.486 20.138 /
\plot  8.486 20.138  8.471 20.121 /
\plot  8.471 20.121  8.454 20.106 /
\plot  8.454 20.106  8.435 20.094 /
\plot  8.435 20.094  8.414 20.081 /
\plot  8.414 20.081  8.384 20.068 /
\plot  8.384 20.068  8.352 20.058 /
\plot  8.352 20.058  8.316 20.047 /
\plot  8.316 20.047  8.276 20.041 /
\plot  8.276 20.041  8.234 20.032 /
\plot  8.234 20.032  8.189 20.028 /
\plot  8.189 20.028  8.143 20.026 /
\plot  8.143 20.026  8.096 20.024 /
\plot  8.096 20.024  8.050 20.026 /
\plot  8.050 20.026  8.003 20.028 /
\plot  8.003 20.028  7.959 20.032 /
\plot  7.959 20.032  7.916 20.041 /
\plot  7.916 20.041  7.876 20.047 /
\plot  7.876 20.047  7.840 20.058 /
\plot  7.840 20.058  7.808 20.068 /
\plot  7.808 20.068  7.779 20.081 /
\plot  7.779 20.081  7.749 20.098 /
\plot  7.749 20.098  7.724 20.117 /
\plot  7.724 20.117  7.703 20.140 /
\plot  7.703 20.140  7.684 20.163 /
\plot  7.684 20.163  7.667 20.189 /
\plot  7.667 20.189  7.652 20.216 /
\plot  7.652 20.216  7.641 20.246 /
\plot  7.641 20.246  7.635 20.273 /
\plot  7.635 20.273  7.628 20.303 /
\plot  7.628 20.303  7.624 20.333 /
\plot  7.624 20.333  7.622 20.358 /
\plot  7.622 20.358  7.620 20.384 /
\putrule from  7.620 20.384 to  7.620 20.409
\putrule from  7.620 20.409 to  7.620 20.430
\putrule from  7.620 20.430 to  7.620 20.460
\putrule from  7.620 20.460 to  7.620 20.489
\putrule from  7.620 20.489 to  7.620 20.517
\putrule from  7.620 20.517 to  7.620 20.546
\putrule from  7.620 20.546 to  7.620 20.578
\putrule from  7.620 20.578 to  7.620 20.610
\putrule from  7.620 20.610 to  7.620 20.640
\putrule from  7.620 20.640 to  7.620 20.659
\putrule from  7.620 20.659 to  7.620 20.667
\putrule from  7.620 20.667 to  7.620 20.669
}%
%
%
\linethickness= 0.500pt
\setplotsymbol ({\thinlinefont .})
{\color[rgb]{0,0,0}\putrule from 10.469 22.369 to 10.469 22.371
\putrule from 10.469 22.371 to 10.469 22.380
\putrule from 10.469 22.380 to 10.469 22.401
\plot 10.469 22.401 10.471 22.435 /
\putrule from 10.471 22.435 to 10.471 22.479
\putrule from 10.471 22.479 to 10.471 22.530
\plot 10.471 22.530 10.473 22.583 /
\putrule from 10.473 22.583 to 10.473 22.634
\putrule from 10.473 22.634 to 10.473 22.678
\plot 10.473 22.678 10.475 22.720 /
\putrule from 10.475 22.720 to 10.475 22.758
\putrule from 10.475 22.758 to 10.475 22.792
\putrule from 10.475 22.792 to 10.475 22.826
\putrule from 10.475 22.826 to 10.475 22.858
\putrule from 10.475 22.858 to 10.475 22.885
\putrule from 10.475 22.885 to 10.475 22.913
\putrule from 10.475 22.913 to 10.475 22.943
\plot 10.475 22.943 10.473 22.970 /
\plot 10.473 22.970 10.471 23.000 /
\plot 10.471 23.000 10.467 23.029 /
\plot 10.467 23.029 10.463 23.059 /
\plot 10.463 23.059 10.456 23.086 /
\plot 10.456 23.086 10.446 23.114 /
\plot 10.446 23.114 10.435 23.139 /
\plot 10.435 23.139 10.420 23.165 /
\plot 10.420 23.165 10.406 23.188 /
\plot 10.406 23.188 10.389 23.207 /
\plot 10.389 23.207 10.367 23.226 /
\plot 10.367 23.226 10.344 23.243 /
\plot 10.344 23.243 10.319 23.256 /
\plot 10.319 23.256 10.289 23.271 /
\plot 10.289 23.271 10.257 23.281 /
\plot 10.257 23.281 10.221 23.292 /
\plot 10.221 23.292 10.181 23.298 /
\plot 10.181 23.298 10.139 23.307 /
\plot 10.139 23.307 10.094 23.311 /
\plot 10.094 23.311 10.048 23.313 /
\plot 10.048 23.313 10.001 23.315 /
\plot 10.001 23.315  9.955 23.313 /
\plot  9.955 23.313  9.908 23.311 /
\plot  9.908 23.311  9.864 23.307 /
\plot  9.864 23.307  9.821 23.298 /
\plot  9.821 23.298  9.781 23.292 /
\plot  9.781 23.292  9.745 23.281 /
\plot  9.745 23.281  9.713 23.271 /
\plot  9.713 23.271  9.684 23.256 /
\plot  9.684 23.256  9.663 23.245 /
\plot  9.663 23.245  9.644 23.233 /
\plot  9.644 23.233  9.627 23.218 /
\plot  9.627 23.218  9.612 23.201 /
\plot  9.612 23.201  9.597 23.182 /
\plot  9.597 23.182  9.584 23.161 /
\plot  9.584 23.161  9.574 23.135 /
\plot  9.574 23.135  9.563 23.108 /
\plot  9.563 23.108  9.555 23.078 /
\plot  9.555 23.078  9.546 23.044 /
\plot  9.546 23.044  9.542 23.006 /
\plot  9.542 23.006  9.536 22.966 /
\plot  9.536 22.966  9.533 22.924 /
\plot  9.533 22.924  9.529 22.877 /
\plot  9.529 22.877  9.527 22.828 /
\putrule from  9.527 22.828 to  9.527 22.777
\plot  9.527 22.777  9.525 22.722 /
\putrule from  9.525 22.722 to  9.525 22.665
\putrule from  9.525 22.665 to  9.525 22.606
\putrule from  9.525 22.606 to  9.525 22.543
\putrule from  9.525 22.543 to  9.525 22.494
\putrule from  9.525 22.494 to  9.525 22.445
\putrule from  9.525 22.445 to  9.525 22.392
\putrule from  9.525 22.392 to  9.525 22.337
\putrule from  9.525 22.337 to  9.525 22.278
\putrule from  9.525 22.278 to  9.525 22.219
\putrule from  9.525 22.219 to  9.525 22.155
\putrule from  9.525 22.155 to  9.525 22.092
\putrule from  9.525 22.092 to  9.525 22.024
\putrule from  9.525 22.024 to  9.525 21.956
\putrule from  9.525 21.956 to  9.525 21.886
\putrule from  9.525 21.886 to  9.525 21.814
\putrule from  9.525 21.814 to  9.525 21.742
\putrule from  9.525 21.742 to  9.525 21.668
\putrule from  9.525 21.668 to  9.525 21.596
\putrule from  9.525 21.596 to  9.525 21.524
\putrule from  9.525 21.524 to  9.525 21.452
\putrule from  9.525 21.452 to  9.525 21.383
\putrule from  9.525 21.383 to  9.525 21.315
\putrule from  9.525 21.315 to  9.525 21.247
\putrule from  9.525 21.247 to  9.525 21.184
\putrule from  9.525 21.184 to  9.525 21.120
\putrule from  9.525 21.120 to  9.525 21.061
\putrule from  9.525 21.061 to  9.525 21.002
\putrule from  9.525 21.002 to  9.525 20.947
\putrule from  9.525 20.947 to  9.525 20.894
\putrule from  9.525 20.894 to  9.525 20.845
\putrule from  9.525 20.845 to  9.525 20.796
\putrule from  9.525 20.796 to  9.525 20.733
\putrule from  9.525 20.733 to  9.525 20.673
\putrule from  9.525 20.673 to  9.525 20.616
\plot  9.525 20.616  9.527 20.561 /
\putrule from  9.527 20.561 to  9.527 20.511
\plot  9.527 20.511  9.529 20.462 /
\plot  9.529 20.462  9.533 20.415 /
\plot  9.533 20.415  9.536 20.373 /
\plot  9.536 20.373  9.542 20.333 /
\plot  9.542 20.333  9.546 20.295 /
\plot  9.546 20.295  9.555 20.261 /
\plot  9.555 20.261  9.563 20.231 /
\plot  9.563 20.231  9.574 20.204 /
\plot  9.574 20.204  9.584 20.178 /
\plot  9.584 20.178  9.597 20.157 /
\plot  9.597 20.157  9.612 20.138 /
\plot  9.612 20.138  9.627 20.121 /
\plot  9.627 20.121  9.644 20.106 /
\plot  9.644 20.106  9.663 20.094 /
\plot  9.663 20.094  9.684 20.081 /
\plot  9.684 20.081  9.713 20.068 /
\plot  9.713 20.068  9.745 20.058 /
\plot  9.745 20.058  9.781 20.047 /
\plot  9.781 20.047  9.821 20.041 /
\plot  9.821 20.041  9.864 20.032 /
\plot  9.864 20.032  9.908 20.028 /
\plot  9.908 20.028  9.955 20.026 /
\plot  9.955 20.026 10.001 20.024 /
\plot 10.001 20.024 10.048 20.026 /
\plot 10.048 20.026 10.094 20.028 /
\plot 10.094 20.028 10.139 20.032 /
\plot 10.139 20.032 10.181 20.041 /
\plot 10.181 20.041 10.221 20.047 /
\plot 10.221 20.047 10.257 20.058 /
\plot 10.257 20.058 10.289 20.068 /
\plot 10.289 20.068 10.319 20.081 /
\plot 10.319 20.081 10.344 20.096 /
\plot 10.344 20.096 10.367 20.113 /
\plot 10.367 20.113 10.389 20.132 /
\plot 10.389 20.132 10.406 20.151 /
\plot 10.406 20.151 10.422 20.174 /
\plot 10.422 20.174 10.435 20.197 /
\plot 10.435 20.197 10.446 20.225 /
\plot 10.446 20.225 10.456 20.252 /
\plot 10.456 20.252 10.463 20.280 /
\plot 10.463 20.280 10.469 20.309 /
\plot 10.469 20.309 10.473 20.337 /
\plot 10.473 20.337 10.475 20.367 /
\putrule from 10.475 20.367 to 10.475 20.396
\plot 10.475 20.396 10.478 20.424 /
\putrule from 10.478 20.424 to 10.478 20.451
\putrule from 10.478 20.451 to 10.478 20.479
\putrule from 10.478 20.479 to 10.478 20.511
\putrule from 10.478 20.511 to 10.478 20.542
\putrule from 10.478 20.542 to 10.478 20.576
\putrule from 10.478 20.576 to 10.478 20.612
\putrule from 10.478 20.612 to 10.478 20.654
\putrule from 10.478 20.654 to 10.478 20.699
\putrule from 10.478 20.699 to 10.478 20.748
\putrule from 10.478 20.748 to 10.478 20.798
\putrule from 10.478 20.798 to 10.478 20.847
\putrule from 10.478 20.847 to 10.478 20.892
\putrule from 10.478 20.892 to 10.478 20.925
\putrule from 10.478 20.925 to 10.478 20.944
\putrule from 10.478 20.944 to 10.478 20.953
\putrule from 10.478 20.953 to 10.478 20.955
}%
%
%
\linethickness= 0.500pt
\setplotsymbol ({\thinlinefont .})
{\color[rgb]{0,0,0}\putrule from 10.461 22.149 to 10.461 22.147
\putrule from 10.461 22.147 to 10.461 22.136
\putrule from 10.461 22.136 to 10.461 22.111
\plot 10.461 22.111 10.458 22.075 /
\plot 10.458 22.075 10.454 22.035 /
\plot 10.454 22.035 10.448 21.994 /
\plot 10.448 21.994 10.439 21.958 /
\plot 10.439 21.958 10.429 21.922 /
\plot 10.429 21.922 10.414 21.888 /
\plot 10.414 21.888 10.395 21.852 /
\plot 10.395 21.852 10.382 21.829 /
\plot 10.382 21.829 10.365 21.804 /
\plot 10.365 21.804 10.348 21.776 /
\plot 10.348 21.776 10.329 21.747 /
\plot 10.329 21.747 10.308 21.717 /
\plot 10.308 21.717 10.285 21.685 /
\plot 10.285 21.685 10.262 21.651 /
\plot 10.262 21.651 10.238 21.618 /
\plot 10.238 21.618 10.215 21.584 /
\plot 10.215 21.584 10.192 21.550 /
\plot 10.192 21.550 10.171 21.518 /
\plot 10.171 21.518 10.149 21.484 /
\plot 10.149 21.484 10.128 21.455 /
\plot 10.128 21.455 10.111 21.425 /
\plot 10.111 21.425 10.094 21.395 /
\plot 10.094 21.395 10.082 21.368 /
\plot 10.082 21.368 10.065 21.332 /
\plot 10.065 21.332 10.050 21.296 /
\plot 10.050 21.296 10.037 21.262 /
\plot 10.037 21.262 10.029 21.226 /
\plot 10.029 21.226 10.020 21.192 /
\plot 10.020 21.192 10.016 21.158 /
\plot 10.016 21.158 10.014 21.129 /
\putrule from 10.014 21.129 to 10.014 21.099
\plot 10.014 21.099 10.016 21.076 /
\plot 10.016 21.076 10.020 21.052 /
\plot 10.020 21.052 10.027 21.033 /
\plot 10.027 21.033 10.033 21.018 /
\plot 10.033 21.018 10.039 21.006 /
\plot 10.039 21.006 10.048 20.997 /
\plot 10.048 20.997 10.056 20.989 /
\plot 10.056 20.989 10.067 20.980 /
\plot 10.067 20.980 10.077 20.974 /
\plot 10.077 20.974 10.088 20.970 /
\plot 10.088 20.970 10.101 20.968 /
\putrule from 10.101 20.968 to 10.113 20.968
\putrule from 10.113 20.968 to 10.126 20.968
\plot 10.126 20.968 10.141 20.970 /
\plot 10.141 20.970 10.154 20.972 /
\plot 10.154 20.972 10.166 20.976 /
\plot 10.166 20.976 10.179 20.980 /
\plot 10.179 20.980 10.192 20.987 /
\plot 10.192 20.987 10.209 20.997 /
\plot 10.209 20.997 10.228 21.008 /
\plot 10.228 21.008 10.249 21.025 /
\plot 10.249 21.025 10.272 21.044 /
\plot 10.272 21.044 10.300 21.067 /
\plot 10.300 21.067 10.327 21.095 /
\plot 10.327 21.095 10.355 21.118 /
\plot 10.355 21.118 10.372 21.135 /
\plot 10.372 21.135 10.380 21.143 /
\plot 10.380 21.143 10.382 21.146 /
}%
%
%
\linethickness= 0.500pt
\setplotsymbol ({\thinlinefont .})
{\color[rgb]{0,0,0}\putrule from 10.478 20.955 to 10.478 20.957
\putrule from 10.478 20.957 to 10.478 20.970
\putrule from 10.478 20.970 to 10.478 20.997
\putrule from 10.478 20.997 to 10.478 21.035
\putrule from 10.478 21.035 to 10.478 21.078
\putrule from 10.478 21.078 to 10.478 21.116
\putrule from 10.478 21.116 to 10.478 21.150
\putrule from 10.478 21.150 to 10.478 21.177
\putrule from 10.478 21.177 to 10.478 21.203
\putrule from 10.478 21.203 to 10.478 21.224
\putrule from 10.478 21.224 to 10.478 21.247
\putrule from 10.478 21.247 to 10.478 21.268
\putrule from 10.478 21.268 to 10.478 21.289
\plot 10.478 21.289 10.475 21.313 /
\putrule from 10.475 21.313 to 10.475 21.336
\plot 10.475 21.336 10.473 21.357 /
\plot 10.473 21.357 10.471 21.378 /
\plot 10.471 21.378 10.467 21.399 /
\plot 10.467 21.399 10.465 21.421 /
\plot 10.465 21.421 10.458 21.440 /
\plot 10.458 21.440 10.454 21.461 /
\plot 10.454 21.461 10.446 21.482 /
\plot 10.446 21.482 10.437 21.507 /
\plot 10.437 21.507 10.427 21.535 /
\plot 10.427 21.535 10.412 21.569 /
\plot 10.412 21.569 10.397 21.605 /
\plot 10.397 21.605 10.382 21.637 /
\plot 10.382 21.637 10.374 21.660 /
\plot 10.374 21.660 10.367 21.670 /
\putrule from 10.367 21.670 to 10.367 21.673
}%
%
%
\linethickness= 0.500pt
\setplotsymbol ({\thinlinefont .})
{\color[rgb]{0,0,0}\putrule from 10.287 21.812 to 10.287 21.814
\plot 10.287 21.814 10.283 21.823 /
\plot 10.283 21.823 10.272 21.842 /
\plot 10.272 21.842 10.259 21.869 /
\plot 10.259 21.869 10.245 21.901 /
\plot 10.245 21.901 10.232 21.931 /
\plot 10.232 21.931 10.219 21.958 /
\plot 10.219 21.958 10.211 21.984 /
\plot 10.211 21.984 10.200 22.009 /
\plot 10.200 22.009 10.192 22.035 /
\plot 10.192 22.035 10.183 22.062 /
\plot 10.183 22.062 10.175 22.092 /
\plot 10.175 22.092 10.166 22.123 /
\plot 10.166 22.123 10.158 22.155 /
\plot 10.158 22.155 10.152 22.189 /
\plot 10.152 22.189 10.145 22.221 /
\plot 10.145 22.221 10.141 22.250 /
\putrule from 10.141 22.250 to 10.141 22.278
\putrule from 10.141 22.278 to 10.141 22.301
\plot 10.141 22.301 10.145 22.320 /
\plot 10.145 22.320 10.147 22.331 /
\plot 10.147 22.331 10.152 22.339 /
\plot 10.152 22.339 10.158 22.348 /
\plot 10.158 22.348 10.164 22.354 /
\plot 10.164 22.354 10.173 22.360 /
\plot 10.173 22.360 10.181 22.365 /
\plot 10.181 22.365 10.192 22.367 /
\plot 10.192 22.367 10.204 22.369 /
\putrule from 10.204 22.369 to 10.219 22.369
\plot 10.219 22.369 10.232 22.367 /
\plot 10.232 22.367 10.249 22.365 /
\plot 10.249 22.365 10.264 22.360 /
\plot 10.264 22.360 10.281 22.356 /
\plot 10.281 22.356 10.298 22.350 /
\plot 10.298 22.350 10.317 22.344 /
\plot 10.317 22.344 10.336 22.335 /
\plot 10.336 22.335 10.357 22.327 /
\plot 10.357 22.327 10.382 22.314 /
\plot 10.382 22.314 10.408 22.301 /
\plot 10.408 22.301 10.433 22.286 /
\plot 10.433 22.286 10.463 22.269 /
\plot 10.463 22.269 10.492 22.253 /
\plot 10.492 22.253 10.520 22.231 /
\plot 10.520 22.231 10.549 22.212 /
\plot 10.549 22.212 10.577 22.191 /
\plot 10.577 22.191 10.602 22.168 /
\plot 10.602 22.168 10.626 22.147 /
\plot 10.626 22.147 10.649 22.126 /
\plot 10.649 22.126 10.668 22.104 /
\plot 10.668 22.104 10.685 22.081 /
\plot 10.685 22.081 10.700 22.056 /
\plot 10.700 22.056 10.715 22.030 /
\plot 10.715 22.030 10.727 22.001 /
\plot 10.727 22.001 10.738 21.971 /
\plot 10.738 21.971 10.744 21.941 /
\plot 10.744 21.941 10.751 21.910 /
\plot 10.751 21.910 10.753 21.880 /
\plot 10.753 21.880 10.755 21.850 /
\putrule from 10.755 21.850 to 10.755 21.821
\plot 10.755 21.821 10.753 21.795 /
\plot 10.753 21.795 10.751 21.772 /
\plot 10.751 21.772 10.748 21.749 /
\plot 10.748 21.749 10.742 21.723 /
\plot 10.742 21.723 10.736 21.698 /
\plot 10.736 21.698 10.729 21.675 /
\plot 10.729 21.675 10.723 21.649 /
\plot 10.723 21.649 10.712 21.626 /
\plot 10.712 21.626 10.704 21.603 /
\plot 10.704 21.603 10.696 21.582 /
\plot 10.696 21.582 10.685 21.562 /
\plot 10.685 21.562 10.676 21.543 /
\plot 10.676 21.543 10.668 21.527 /
\plot 10.668 21.527 10.660 21.510 /
\plot 10.660 21.510 10.649 21.491 /
\plot 10.649 21.491 10.640 21.469 /
\plot 10.640 21.469 10.628 21.446 /
\plot 10.628 21.446 10.615 21.419 /
\plot 10.615 21.419 10.600 21.391 /
\plot 10.600 21.391 10.588 21.364 /
\plot 10.588 21.364 10.577 21.347 /
\plot 10.577 21.347 10.573 21.338 /
\putrule from 10.573 21.338 to 10.573 21.336
}%
%
%
\put{$-4$
} [lB] at  7.904 19.748
%
%
\put{$+3$
} [lB] at  8.735 21.023
%
%
\linethickness= 0.500pt
\setplotsymbol ({\thinlinefont .})
{\color[rgb]{0,0,0}\circulararc 157.380 degrees from  2.487 21.632 center at  2.869 21.703
}%
%
%
\linethickness= 0.500pt
\setplotsymbol ({\thinlinefont .})
\setdots < 0.0953cm>
{\color[rgb]{0,0,0}\plot  1.429 22.098  1.429 21.241 /
}%
%
%
\linethickness= 0.500pt
\setplotsymbol ({\thinlinefont .})
\setdots < 0.0953cm>
{\color[rgb]{0,0,0}\plot  2.667 21.907  3.048 21.907 /
}%
%
%
\linethickness= 0.500pt
\setplotsymbol ({\thinlinefont .})
\setsolid
{\color[rgb]{0,0,0}\putrule from  4.269 22.149 to  4.269 22.147
\putrule from  4.269 22.147 to  4.269 22.136
\putrule from  4.269 22.136 to  4.269 22.111
\plot  4.269 22.111  4.267 22.075 /
\plot  4.267 22.075  4.263 22.035 /
\plot  4.263 22.035  4.257 21.994 /
\plot  4.257 21.994  4.248 21.958 /
\plot  4.248 21.958  4.238 21.922 /
\plot  4.238 21.922  4.223 21.888 /
\plot  4.223 21.888  4.204 21.852 /
\plot  4.204 21.852  4.191 21.829 /
\plot  4.191 21.829  4.174 21.804 /
\plot  4.174 21.804  4.157 21.776 /
\plot  4.157 21.776  4.138 21.747 /
\plot  4.138 21.747  4.117 21.717 /
\plot  4.117 21.717  4.094 21.685 /
\plot  4.094 21.685  4.070 21.651 /
\plot  4.070 21.651  4.047 21.618 /
\plot  4.047 21.618  4.024 21.584 /
\plot  4.024 21.584  4.000 21.550 /
\plot  4.000 21.550  3.979 21.518 /
\plot  3.979 21.518  3.958 21.484 /
\plot  3.958 21.484  3.937 21.455 /
\plot  3.937 21.455  3.920 21.425 /
\plot  3.920 21.425  3.903 21.395 /
\plot  3.903 21.395  3.890 21.368 /
\plot  3.890 21.368  3.873 21.332 /
\plot  3.873 21.332  3.859 21.296 /
\plot  3.859 21.296  3.846 21.262 /
\plot  3.846 21.262  3.838 21.226 /
\plot  3.838 21.226  3.829 21.192 /
\plot  3.829 21.192  3.825 21.158 /
\plot  3.825 21.158  3.823 21.129 /
\putrule from  3.823 21.129 to  3.823 21.099
\plot  3.823 21.099  3.825 21.076 /
\plot  3.825 21.076  3.829 21.052 /
\plot  3.829 21.052  3.835 21.033 /
\plot  3.835 21.033  3.842 21.018 /
\plot  3.842 21.018  3.848 21.006 /
\plot  3.848 21.006  3.857 20.997 /
\plot  3.857 20.997  3.865 20.989 /
\plot  3.865 20.989  3.876 20.980 /
\plot  3.876 20.980  3.886 20.974 /
\plot  3.886 20.974  3.897 20.970 /
\plot  3.897 20.970  3.909 20.968 /
\putrule from  3.909 20.968 to  3.922 20.968
\putrule from  3.922 20.968 to  3.935 20.968
\plot  3.935 20.968  3.950 20.970 /
\plot  3.950 20.970  3.962 20.972 /
\plot  3.962 20.972  3.975 20.976 /
\plot  3.975 20.976  3.988 20.980 /
\plot  3.988 20.980  4.000 20.987 /
\plot  4.000 20.987  4.017 20.997 /
\plot  4.017 20.997  4.036 21.008 /
\plot  4.036 21.008  4.058 21.025 /
\plot  4.058 21.025  4.081 21.044 /
\plot  4.081 21.044  4.108 21.067 /
\plot  4.108 21.067  4.136 21.095 /
\plot  4.136 21.095  4.163 21.118 /
\plot  4.163 21.118  4.180 21.135 /
\plot  4.180 21.135  4.189 21.143 /
\plot  4.189 21.143  4.191 21.146 /
}%
%
%
\linethickness= 0.500pt
\setplotsymbol ({\thinlinefont .})
{\color[rgb]{0,0,0}\putrule from  4.286 20.955 to  4.286 20.957
\putrule from  4.286 20.957 to  4.286 20.970
\putrule from  4.286 20.970 to  4.286 20.997
\putrule from  4.286 20.997 to  4.286 21.035
\putrule from  4.286 21.035 to  4.286 21.078
\putrule from  4.286 21.078 to  4.286 21.116
\putrule from  4.286 21.116 to  4.286 21.150
\putrule from  4.286 21.150 to  4.286 21.177
\putrule from  4.286 21.177 to  4.286 21.203
\putrule from  4.286 21.203 to  4.286 21.224
\putrule from  4.286 21.224 to  4.286 21.247
\putrule from  4.286 21.247 to  4.286 21.268
\putrule from  4.286 21.268 to  4.286 21.289
\plot  4.286 21.289  4.284 21.313 /
\putrule from  4.284 21.313 to  4.284 21.336
\plot  4.284 21.336  4.282 21.357 /
\plot  4.282 21.357  4.280 21.378 /
\plot  4.280 21.378  4.276 21.399 /
\plot  4.276 21.399  4.274 21.421 /
\plot  4.274 21.421  4.267 21.440 /
\plot  4.267 21.440  4.263 21.461 /
\plot  4.263 21.461  4.255 21.482 /
\plot  4.255 21.482  4.246 21.507 /
\plot  4.246 21.507  4.235 21.535 /
\plot  4.235 21.535  4.221 21.569 /
\plot  4.221 21.569  4.206 21.605 /
\plot  4.206 21.605  4.191 21.637 /
\plot  4.191 21.637  4.183 21.660 /
\plot  4.183 21.660  4.176 21.670 /
\putrule from  4.176 21.670 to  4.176 21.673
}%
%
%
\linethickness= 0.500pt
\setplotsymbol ({\thinlinefont .})
{\color[rgb]{0,0,0}\putrule from  4.096 21.812 to  4.096 21.814
\plot  4.096 21.814  4.092 21.823 /
\plot  4.092 21.823  4.081 21.842 /
\plot  4.081 21.842  4.068 21.869 /
\plot  4.068 21.869  4.053 21.901 /
\plot  4.053 21.901  4.041 21.931 /
\plot  4.041 21.931  4.028 21.958 /
\plot  4.028 21.958  4.020 21.984 /
\plot  4.020 21.984  4.009 22.009 /
\plot  4.009 22.009  4.000 22.035 /
\plot  4.000 22.035  3.992 22.062 /
\plot  3.992 22.062  3.984 22.092 /
\plot  3.984 22.092  3.975 22.123 /
\plot  3.975 22.123  3.967 22.155 /
\plot  3.967 22.155  3.960 22.189 /
\plot  3.960 22.189  3.954 22.221 /
\plot  3.954 22.221  3.950 22.250 /
\putrule from  3.950 22.250 to  3.950 22.278
\putrule from  3.950 22.278 to  3.950 22.301
\plot  3.950 22.301  3.954 22.320 /
\plot  3.954 22.320  3.956 22.331 /
\plot  3.956 22.331  3.960 22.339 /
\plot  3.960 22.339  3.967 22.348 /
\plot  3.967 22.348  3.973 22.354 /
\plot  3.973 22.354  3.981 22.360 /
\plot  3.981 22.360  3.990 22.365 /
\plot  3.990 22.365  4.000 22.367 /
\plot  4.000 22.367  4.013 22.369 /
\putrule from  4.013 22.369 to  4.028 22.369
\plot  4.028 22.369  4.041 22.367 /
\plot  4.041 22.367  4.058 22.365 /
\plot  4.058 22.365  4.072 22.360 /
\plot  4.072 22.360  4.089 22.356 /
\plot  4.089 22.356  4.106 22.350 /
\plot  4.106 22.350  4.125 22.344 /
\plot  4.125 22.344  4.144 22.335 /
\plot  4.144 22.335  4.166 22.327 /
\plot  4.166 22.327  4.191 22.314 /
\plot  4.191 22.314  4.216 22.301 /
\plot  4.216 22.301  4.242 22.286 /
\plot  4.242 22.286  4.271 22.269 /
\plot  4.271 22.269  4.301 22.253 /
\plot  4.301 22.253  4.329 22.231 /
\plot  4.329 22.231  4.358 22.212 /
\plot  4.358 22.212  4.386 22.191 /
\plot  4.386 22.191  4.411 22.168 /
\plot  4.411 22.168  4.434 22.147 /
\plot  4.434 22.147  4.458 22.126 /
\plot  4.458 22.126  4.477 22.104 /
\plot  4.477 22.104  4.494 22.081 /
\plot  4.494 22.081  4.508 22.056 /
\plot  4.508 22.056  4.523 22.030 /
\plot  4.523 22.030  4.536 22.001 /
\plot  4.536 22.001  4.547 21.971 /
\plot  4.547 21.971  4.553 21.941 /
\plot  4.553 21.941  4.559 21.910 /
\plot  4.559 21.910  4.561 21.880 /
\plot  4.561 21.880  4.564 21.850 /
\putrule from  4.564 21.850 to  4.564 21.821
\plot  4.564 21.821  4.561 21.795 /
\plot  4.561 21.795  4.559 21.772 /
\plot  4.559 21.772  4.557 21.749 /
\plot  4.557 21.749  4.551 21.723 /
\plot  4.551 21.723  4.544 21.698 /
\plot  4.544 21.698  4.538 21.675 /
\plot  4.538 21.675  4.532 21.649 /
\plot  4.532 21.649  4.521 21.626 /
\plot  4.521 21.626  4.513 21.603 /
\plot  4.513 21.603  4.504 21.582 /
\plot  4.504 21.582  4.494 21.562 /
\plot  4.494 21.562  4.485 21.543 /
\plot  4.485 21.543  4.477 21.527 /
\plot  4.477 21.527  4.468 21.510 /
\plot  4.468 21.510  4.458 21.491 /
\plot  4.458 21.491  4.449 21.469 /
\plot  4.449 21.469  4.437 21.446 /
\plot  4.437 21.446  4.424 21.419 /
\plot  4.424 21.419  4.409 21.391 /
\plot  4.409 21.391  4.396 21.364 /
\plot  4.396 21.364  4.386 21.347 /
\plot  4.386 21.347  4.381 21.338 /
\putrule from  4.381 21.338 to  4.381 21.336
}%
%
%
\linethickness= 0.500pt
\setplotsymbol ({\thinlinefont .})
{\color[rgb]{0,0,0}\putrule from  1.429 22.669 to  1.429 22.672
\putrule from  1.429 22.672 to  1.429 22.680
\putrule from  1.429 22.680 to  1.429 22.699
\putrule from  1.429 22.699 to  1.429 22.729
\putrule from  1.429 22.729 to  1.429 22.761
\putrule from  1.429 22.761 to  1.429 22.792
\putrule from  1.429 22.792 to  1.429 22.822
\putrule from  1.429 22.822 to  1.429 22.849
\putrule from  1.429 22.849 to  1.429 22.879
\putrule from  1.429 22.879 to  1.429 22.907
\putrule from  1.429 22.907 to  1.429 22.930
\putrule from  1.429 22.930 to  1.429 22.955
\plot  1.429 22.955  1.431 22.981 /
\plot  1.431 22.981  1.433 23.006 /
\plot  1.433 23.006  1.437 23.036 /
\plot  1.437 23.036  1.444 23.065 /
\plot  1.444 23.065  1.450 23.093 /
\plot  1.450 23.093  1.461 23.122 /
\plot  1.461 23.122  1.475 23.150 /
\plot  1.475 23.150  1.492 23.175 /
\plot  1.492 23.175  1.511 23.199 /
\plot  1.511 23.199  1.532 23.222 /
\plot  1.532 23.222  1.558 23.241 /
\plot  1.558 23.241  1.587 23.256 /
\plot  1.587 23.256  1.617 23.271 /
\plot  1.617 23.271  1.649 23.281 /
\plot  1.649 23.281  1.685 23.292 /
\plot  1.685 23.292  1.725 23.298 /
\plot  1.725 23.298  1.767 23.307 /
\plot  1.767 23.307  1.812 23.311 /
\plot  1.812 23.311  1.858 23.313 /
\plot  1.858 23.313  1.905 23.315 /
\plot  1.905 23.315  1.952 23.313 /
\plot  1.952 23.313  1.998 23.311 /
\plot  1.998 23.311  2.043 23.307 /
\plot  2.043 23.307  2.085 23.298 /
\plot  2.085 23.298  2.125 23.292 /
\plot  2.125 23.292  2.161 23.281 /
\plot  2.161 23.281  2.193 23.271 /
\plot  2.193 23.271  2.223 23.256 /
\plot  2.223 23.256  2.248 23.243 /
\plot  2.248 23.243  2.271 23.226 /
\plot  2.271 23.226  2.292 23.207 /
\plot  2.292 23.207  2.309 23.186 /
\plot  2.309 23.186  2.326 23.165 /
\plot  2.326 23.165  2.339 23.139 /
\plot  2.339 23.139  2.349 23.112 /
\plot  2.349 23.112  2.360 23.082 /
\plot  2.360 23.082  2.366 23.053 /
\plot  2.366 23.053  2.373 23.021 /
\plot  2.373 23.021  2.377 22.989 /
\plot  2.377 22.989  2.379 22.957 /
\putrule from  2.379 22.957 to  2.379 22.926
\plot  2.379 22.926  2.381 22.894 /
\putrule from  2.381 22.894 to  2.381 22.860
\putrule from  2.381 22.860 to  2.381 22.828
\putrule from  2.381 22.828 to  2.381 22.794
\putrule from  2.381 22.794 to  2.381 22.761
\putrule from  2.381 22.761 to  2.381 22.725
\putrule from  2.381 22.725 to  2.381 22.684
\putrule from  2.381 22.684 to  2.381 22.640
\putrule from  2.381 22.640 to  2.381 22.591
\putrule from  2.381 22.591 to  2.381 22.536
\putrule from  2.381 22.536 to  2.381 22.479
\putrule from  2.381 22.479 to  2.381 22.418
\putrule from  2.381 22.418 to  2.381 22.358
\putrule from  2.381 22.358 to  2.381 22.301
\putrule from  2.381 22.301 to  2.381 22.257
\putrule from  2.381 22.257 to  2.381 22.223
\putrule from  2.381 22.223 to  2.381 22.202
\putrule from  2.381 22.202 to  2.381 22.195
\putrule from  2.381 22.195 to  2.381 22.193
}%
%
%
\linethickness= 0.500pt
\setplotsymbol ({\thinlinefont .})
{\color[rgb]{0,0,0}\putrule from  2.381 21.622 to  2.381 21.615
\putrule from  2.381 21.615 to  2.381 21.601
\putrule from  2.381 21.601 to  2.381 21.577
\putrule from  2.381 21.577 to  2.381 21.539
\putrule from  2.381 21.539 to  2.381 21.491
\putrule from  2.381 21.491 to  2.381 21.431
\putrule from  2.381 21.431 to  2.381 21.364
\putrule from  2.381 21.364 to  2.381 21.289
\putrule from  2.381 21.289 to  2.381 21.215
\putrule from  2.381 21.215 to  2.381 21.143
\putrule from  2.381 21.143 to  2.381 21.071
\putrule from  2.381 21.071 to  2.381 21.006
\putrule from  2.381 21.006 to  2.381 20.942
\putrule from  2.381 20.942 to  2.381 20.885
\putrule from  2.381 20.885 to  2.381 20.832
\putrule from  2.381 20.832 to  2.381 20.784
\putrule from  2.381 20.784 to  2.381 20.739
\putrule from  2.381 20.739 to  2.381 20.699
\putrule from  2.381 20.699 to  2.381 20.661
\putrule from  2.381 20.661 to  2.381 20.625
\putrule from  2.381 20.625 to  2.381 20.589
\putrule from  2.381 20.589 to  2.381 20.544
\putrule from  2.381 20.544 to  2.381 20.500
\plot  2.381 20.500  2.379 20.458 /
\putrule from  2.379 20.458 to  2.379 20.417
\plot  2.379 20.417  2.377 20.377 /
\plot  2.377 20.377  2.373 20.339 /
\plot  2.373 20.339  2.366 20.301 /
\plot  2.366 20.301  2.360 20.267 /
\plot  2.360 20.267  2.349 20.235 /
\plot  2.349 20.235  2.339 20.206 /
\plot  2.339 20.206  2.326 20.178 /
\plot  2.326 20.178  2.309 20.153 /
\plot  2.309 20.153  2.292 20.132 /
\plot  2.292 20.132  2.271 20.113 /
\plot  2.271 20.113  2.248 20.096 /
\plot  2.248 20.096  2.223 20.081 /
\plot  2.223 20.081  2.193 20.068 /
\plot  2.193 20.068  2.161 20.058 /
\plot  2.161 20.058  2.125 20.047 /
\plot  2.125 20.047  2.085 20.041 /
\plot  2.085 20.041  2.043 20.032 /
\plot  2.043 20.032  1.998 20.028 /
\plot  1.998 20.028  1.952 20.026 /
\plot  1.952 20.026  1.905 20.024 /
\plot  1.905 20.024  1.858 20.026 /
\plot  1.858 20.026  1.812 20.028 /
\plot  1.812 20.028  1.767 20.032 /
\plot  1.767 20.032  1.725 20.041 /
\plot  1.725 20.041  1.685 20.047 /
\plot  1.685 20.047  1.649 20.058 /
\plot  1.649 20.058  1.617 20.068 /
\plot  1.617 20.068  1.587 20.081 /
\plot  1.587 20.081  1.558 20.098 /
\plot  1.558 20.098  1.532 20.117 /
\plot  1.532 20.117  1.511 20.140 /
\plot  1.511 20.140  1.492 20.163 /
\plot  1.492 20.163  1.475 20.189 /
\plot  1.475 20.189  1.461 20.216 /
\plot  1.461 20.216  1.450 20.246 /
\plot  1.450 20.246  1.444 20.273 /
\plot  1.444 20.273  1.437 20.303 /
\plot  1.437 20.303  1.433 20.333 /
\plot  1.433 20.333  1.431 20.358 /
\plot  1.431 20.358  1.429 20.384 /
\putrule from  1.429 20.384 to  1.429 20.409
\putrule from  1.429 20.409 to  1.429 20.430
\putrule from  1.429 20.430 to  1.429 20.460
\putrule from  1.429 20.460 to  1.429 20.489
\putrule from  1.429 20.489 to  1.429 20.517
\putrule from  1.429 20.517 to  1.429 20.546
\putrule from  1.429 20.546 to  1.429 20.578
\putrule from  1.429 20.578 to  1.429 20.610
\putrule from  1.429 20.610 to  1.429 20.640
\putrule from  1.429 20.640 to  1.429 20.659
\putrule from  1.429 20.659 to  1.429 20.667
\putrule from  1.429 20.667 to  1.429 20.669
}%
%
%
\linethickness= 0.500pt
\setplotsymbol ({\thinlinefont .})
{\color[rgb]{0,0,0}\putrule from  4.278 22.369 to  4.278 22.371
\putrule from  4.278 22.371 to  4.278 22.380
\putrule from  4.278 22.380 to  4.278 22.401
\plot  4.278 22.401  4.280 22.435 /
\putrule from  4.280 22.435 to  4.280 22.479
\putrule from  4.280 22.479 to  4.280 22.530
\plot  4.280 22.530  4.282 22.583 /
\putrule from  4.282 22.583 to  4.282 22.634
\putrule from  4.282 22.634 to  4.282 22.678
\plot  4.282 22.678  4.284 22.720 /
\putrule from  4.284 22.720 to  4.284 22.758
\putrule from  4.284 22.758 to  4.284 22.792
\putrule from  4.284 22.792 to  4.284 22.826
\putrule from  4.284 22.826 to  4.284 22.858
\putrule from  4.284 22.858 to  4.284 22.885
\putrule from  4.284 22.885 to  4.284 22.913
\putrule from  4.284 22.913 to  4.284 22.943
\plot  4.284 22.943  4.282 22.970 /
\plot  4.282 22.970  4.280 23.000 /
\plot  4.280 23.000  4.276 23.029 /
\plot  4.276 23.029  4.271 23.059 /
\plot  4.271 23.059  4.265 23.086 /
\plot  4.265 23.086  4.255 23.114 /
\plot  4.255 23.114  4.244 23.139 /
\plot  4.244 23.139  4.229 23.165 /
\plot  4.229 23.165  4.214 23.188 /
\plot  4.214 23.188  4.197 23.207 /
\plot  4.197 23.207  4.176 23.226 /
\plot  4.176 23.226  4.153 23.243 /
\plot  4.153 23.243  4.128 23.256 /
\plot  4.128 23.256  4.098 23.271 /
\plot  4.098 23.271  4.066 23.281 /
\plot  4.066 23.281  4.030 23.292 /
\plot  4.030 23.292  3.990 23.298 /
\plot  3.990 23.298  3.948 23.307 /
\plot  3.948 23.307  3.903 23.311 /
\plot  3.903 23.311  3.857 23.313 /
\plot  3.857 23.313  3.810 23.315 /
\plot  3.810 23.315  3.763 23.313 /
\plot  3.763 23.313  3.717 23.311 /
\plot  3.717 23.311  3.672 23.307 /
\plot  3.672 23.307  3.630 23.298 /
\plot  3.630 23.298  3.590 23.292 /
\plot  3.590 23.292  3.554 23.281 /
\plot  3.554 23.281  3.522 23.271 /
\plot  3.522 23.271  3.493 23.256 /
\plot  3.493 23.256  3.467 23.243 /
\plot  3.467 23.243  3.444 23.226 /
\plot  3.444 23.226  3.423 23.207 /
\plot  3.423 23.207  3.406 23.186 /
\plot  3.406 23.186  3.389 23.165 /
\plot  3.389 23.165  3.376 23.139 /
\plot  3.376 23.139  3.365 23.112 /
\plot  3.365 23.112  3.355 23.082 /
\plot  3.355 23.082  3.349 23.053 /
\plot  3.349 23.053  3.342 23.021 /
\plot  3.342 23.021  3.338 22.989 /
\plot  3.338 22.989  3.336 22.957 /
\putrule from  3.336 22.957 to  3.336 22.926
\plot  3.336 22.926  3.334 22.894 /
\putrule from  3.334 22.894 to  3.334 22.860
\putrule from  3.334 22.860 to  3.334 22.828
\putrule from  3.334 22.828 to  3.334 22.794
\putrule from  3.334 22.794 to  3.334 22.761
\putrule from  3.334 22.761 to  3.334 22.725
\putrule from  3.334 22.725 to  3.334 22.684
\putrule from  3.334 22.684 to  3.334 22.640
\putrule from  3.334 22.640 to  3.334 22.591
\putrule from  3.334 22.591 to  3.334 22.536
\putrule from  3.334 22.536 to  3.334 22.479
\putrule from  3.334 22.479 to  3.334 22.418
\putrule from  3.334 22.418 to  3.334 22.358
\putrule from  3.334 22.358 to  3.334 22.301
\putrule from  3.334 22.301 to  3.334 22.257
\putrule from  3.334 22.257 to  3.334 22.223
\putrule from  3.334 22.223 to  3.334 22.202
\putrule from  3.334 22.202 to  3.334 22.195
\putrule from  3.334 22.195 to  3.334 22.193
}%
%
%
\linethickness= 0.500pt
\setplotsymbol ({\thinlinefont .})
{\color[rgb]{0,0,0}\putrule from  3.334 21.622 to  3.334 21.615
\putrule from  3.334 21.615 to  3.334 21.601
\putrule from  3.334 21.601 to  3.334 21.577
\putrule from  3.334 21.577 to  3.334 21.539
\putrule from  3.334 21.539 to  3.334 21.491
\putrule from  3.334 21.491 to  3.334 21.431
\putrule from  3.334 21.431 to  3.334 21.364
\putrule from  3.334 21.364 to  3.334 21.289
\putrule from  3.334 21.289 to  3.334 21.215
\putrule from  3.334 21.215 to  3.334 21.143
\putrule from  3.334 21.143 to  3.334 21.071
\putrule from  3.334 21.071 to  3.334 21.006
\putrule from  3.334 21.006 to  3.334 20.942
\putrule from  3.334 20.942 to  3.334 20.885
\putrule from  3.334 20.885 to  3.334 20.832
\putrule from  3.334 20.832 to  3.334 20.784
\putrule from  3.334 20.784 to  3.334 20.739
\putrule from  3.334 20.739 to  3.334 20.699
\putrule from  3.334 20.699 to  3.334 20.661
\putrule from  3.334 20.661 to  3.334 20.625
\putrule from  3.334 20.625 to  3.334 20.589
\putrule from  3.334 20.589 to  3.334 20.544
\putrule from  3.334 20.544 to  3.334 20.500
\plot  3.334 20.500  3.336 20.458 /
\putrule from  3.336 20.458 to  3.336 20.417
\plot  3.336 20.417  3.338 20.377 /
\plot  3.338 20.377  3.342 20.339 /
\plot  3.342 20.339  3.349 20.301 /
\plot  3.349 20.301  3.355 20.267 /
\plot  3.355 20.267  3.365 20.235 /
\plot  3.365 20.235  3.376 20.206 /
\plot  3.376 20.206  3.389 20.178 /
\plot  3.389 20.178  3.406 20.153 /
\plot  3.406 20.153  3.423 20.132 /
\plot  3.423 20.132  3.444 20.113 /
\plot  3.444 20.113  3.467 20.096 /
\plot  3.467 20.096  3.493 20.081 /
\plot  3.493 20.081  3.522 20.068 /
\plot  3.522 20.068  3.554 20.058 /
\plot  3.554 20.058  3.590 20.047 /
\plot  3.590 20.047  3.630 20.041 /
\plot  3.630 20.041  3.672 20.032 /
\plot  3.672 20.032  3.717 20.028 /
\plot  3.717 20.028  3.763 20.026 /
\plot  3.763 20.026  3.810 20.024 /
\plot  3.810 20.024  3.857 20.026 /
\plot  3.857 20.026  3.903 20.028 /
\plot  3.903 20.028  3.948 20.032 /
\plot  3.948 20.032  3.990 20.041 /
\plot  3.990 20.041  4.030 20.047 /
\plot  4.030 20.047  4.066 20.058 /
\plot  4.066 20.058  4.098 20.068 /
\plot  4.098 20.068  4.128 20.081 /
\plot  4.128 20.081  4.153 20.096 /
\plot  4.153 20.096  4.176 20.113 /
\plot  4.176 20.113  4.197 20.132 /
\plot  4.197 20.132  4.214 20.151 /
\plot  4.214 20.151  4.231 20.174 /
\plot  4.231 20.174  4.244 20.197 /
\plot  4.244 20.197  4.255 20.225 /
\plot  4.255 20.225  4.265 20.252 /
\plot  4.265 20.252  4.271 20.280 /
\plot  4.271 20.280  4.278 20.309 /
\plot  4.278 20.309  4.282 20.337 /
\plot  4.282 20.337  4.284 20.367 /
\putrule from  4.284 20.367 to  4.284 20.396
\plot  4.284 20.396  4.286 20.424 /
\putrule from  4.286 20.424 to  4.286 20.451
\putrule from  4.286 20.451 to  4.286 20.479
\putrule from  4.286 20.479 to  4.286 20.511
\putrule from  4.286 20.511 to  4.286 20.542
\putrule from  4.286 20.542 to  4.286 20.576
\putrule from  4.286 20.576 to  4.286 20.612
\putrule from  4.286 20.612 to  4.286 20.654
\putrule from  4.286 20.654 to  4.286 20.699
\putrule from  4.286 20.699 to  4.286 20.748
\putrule from  4.286 20.748 to  4.286 20.798
\putrule from  4.286 20.798 to  4.286 20.847
\putrule from  4.286 20.847 to  4.286 20.892
\putrule from  4.286 20.892 to  4.286 20.925
\putrule from  4.286 20.925 to  4.286 20.944
\putrule from  4.286 20.944 to  4.286 20.953
\putrule from  4.286 20.953 to  4.286 20.955
}%
%
%
\put{$-4$
} [lB] at  1.712 19.748
%
%
\put{$+3$
} [lB] at  2.544 21.023
\linethickness=0pt
\putrectangle corners at  1.126 23.340 and 10.780 19.632
\endpicture}

%% file: FindingAnH.tex
\font\thinlinefont=cmr5
\mbox{\beginpicture
\small
\setcoordinatesystem units <0.4cm,0.4cm>
\unitlength=1.04987cm
\linethickness=1pt
\setplotsymbol ({\makebox(0,0)[l]{\tencirc\symbol{'160}}})
\setshadesymbol ({\thinlinefont .})
\setlinear
%
%
\linethickness= 0.500pt
\setplotsymbol ({\thinlinefont .})
{\color[rgb]{0,0,0}\circulararc 106.260 degrees from 12.002 20.003 center at 12.383 20.288
}%
%
%
\linethickness= 0.500pt
\setplotsymbol ({\thinlinefont .})
\setdots < 0.0953cm>
{\color[rgb]{0,0,0}\circulararc 106.260 degrees from 12.764 20.003 center at 12.383 19.717
}%
%
%
\linethickness= 0.500pt
\setplotsymbol ({\thinlinefont .})
\setsolid
{\color[rgb]{0,0,0}\ellipticalarc axes ratio  0.381:0.381  360 degrees 
	from 12.764 20.003 center at 12.383 20.003
}%
%
%
\linethickness= 0.500pt
\setplotsymbol ({\thinlinefont .})
{\color[rgb]{0,0,0}\circulararc 106.260 degrees from 12.002 16.192 center at 12.383 16.478
}%
%
%
\linethickness= 0.500pt
\setplotsymbol ({\thinlinefont .})
\setdots < 0.0953cm>
{\color[rgb]{0,0,0}\circulararc 106.260 degrees from 12.764 16.192 center at 12.383 15.907
}%
%
%
\linethickness= 0.500pt
\setplotsymbol ({\thinlinefont .})
\setsolid
{\color[rgb]{0,0,0}\ellipticalarc axes ratio  0.381:0.381  360 degrees 
	from 12.764 16.192 center at 12.383 16.192
}%
%
%
\linethickness= 0.500pt
\setplotsymbol ({\thinlinefont .})
{\color[rgb]{0,0,0}\putrule from 18.098 20.955 to 18.098 12.383
\plot 18.098 12.383 21.907 15.240 /
\putrule from 21.907 15.240 to 21.907 23.812
\plot 21.907 23.812 18.098 20.955 /
}%
%
%
\linethickness= 0.500pt
\setplotsymbol ({\thinlinefont .})
{\color[rgb]{0,0,0}\putrule from 20.003 20.003 to 18.288 20.003
}%
%
%
\linethickness= 0.500pt
\setplotsymbol ({\thinlinefont .})
{\color[rgb]{0,0,0}\putrule from 20.003 18.098 to 18.479 18.098
}%
%
%
\linethickness= 0.500pt
\setplotsymbol ({\thinlinefont .})
{\color[rgb]{0,0,0}\putrule from 20.003 16.192 to 18.479 16.192
}%
%
%
\linethickness= 0.500pt
\setplotsymbol ({\thinlinefont .})
{\color[rgb]{0,0,0}\putrule from 12.764 20.003 to 17.907 20.003
}%
%
%
\linethickness= 0.500pt
\setplotsymbol ({\thinlinefont .})
{\color[rgb]{0,0,0}\putrule from 12.764 16.192 to 17.907 16.192
}%
%
%
\linethickness= 0.500pt
\setplotsymbol ({\thinlinefont .})
{\color[rgb]{0,0,0}\putrule from 12.002 20.003 to  5.715 20.003
}%
%
%
\linethickness= 0.500pt
\setplotsymbol ({\thinlinefont .})
{\color[rgb]{0,0,0}\putrule from 12.002 16.192 to  5.715 16.192
}%
%
%
\linethickness= 0.500pt
\setplotsymbol ({\thinlinefont .})
{\color[rgb]{0,0,0}\putrule from  3.810 20.955 to  3.810 12.383
\plot  3.810 12.383  7.620 15.240 /
\putrule from  7.620 15.240 to  7.620 16.002
}%
%
%
\linethickness= 0.500pt
\setplotsymbol ({\thinlinefont .})
{\color[rgb]{0,0,0}\putrule from  7.620 20.193 to  7.620 23.812
\plot  7.620 23.812  3.810 20.955 /
}%
%
%
\linethickness= 0.500pt
\setplotsymbol ({\thinlinefont .})
{\color[rgb]{0,0,0}\putrule from  7.620 19.812 to  7.620 18.288
}%
%
%
\linethickness= 0.500pt
\setplotsymbol ({\thinlinefont .})
{\color[rgb]{0,0,0}\putrule from  7.620 17.907 to  7.620 16.383
}%
%
%
\linethickness= 0.500pt
\setplotsymbol ({\thinlinefont .})
{\color[rgb]{0,0,0}\putrule from  5.715 18.098 to  5.719 18.098
\putrule from  5.719 18.098 to  5.726 18.098
\putrule from  5.726 18.098 to  5.743 18.098
\putrule from  5.743 18.098 to  5.766 18.098
\putrule from  5.766 18.098 to  5.800 18.098
\putrule from  5.800 18.098 to  5.844 18.098
\putrule from  5.844 18.098 to  5.903 18.098
\putrule from  5.903 18.098 to  5.973 18.098
\putrule from  5.973 18.098 to  6.056 18.098
\putrule from  6.056 18.098 to  6.151 18.098
\putrule from  6.151 18.098 to  6.259 18.098
\putrule from  6.259 18.098 to  6.378 18.098
\putrule from  6.378 18.098 to  6.505 18.098
\putrule from  6.505 18.098 to  6.640 18.098
\putrule from  6.640 18.098 to  6.782 18.098
\putrule from  6.782 18.098 to  6.928 18.098
\putrule from  6.928 18.098 to  7.078 18.098
\putrule from  7.078 18.098 to  7.228 18.098
\putrule from  7.228 18.098 to  7.381 18.098
\putrule from  7.381 18.098 to  7.533 18.098
\putrule from  7.533 18.098 to  7.681 18.098
\putrule from  7.681 18.098 to  7.827 18.098
\putrule from  7.827 18.098 to  7.969 18.098
\putrule from  7.969 18.098 to  8.107 18.098
\putrule from  8.107 18.098 to  8.240 18.098
\putrule from  8.240 18.098 to  8.367 18.098
\putrule from  8.367 18.098 to  8.490 18.098
\putrule from  8.490 18.098 to  8.608 18.098
\putrule from  8.608 18.098 to  8.721 18.098
\putrule from  8.721 18.098 to  8.826 18.098
\putrule from  8.826 18.098 to  8.928 18.098
\putrule from  8.928 18.098 to  9.025 18.098
\putrule from  9.025 18.098 to  9.119 18.098
\putrule from  9.119 18.098 to  9.205 18.098
\putrule from  9.205 18.098 to  9.290 18.098
\putrule from  9.290 18.098 to  9.368 18.098
\putrule from  9.368 18.098 to  9.445 18.098
\putrule from  9.445 18.098 to  9.519 18.098
\putrule from  9.519 18.098 to  9.588 18.098
\putrule from  9.588 18.098 to  9.654 18.098
\putrule from  9.654 18.098 to  9.720 18.098
\putrule from  9.720 18.098 to  9.781 18.098
\putrule from  9.781 18.098 to  9.842 18.098
\putrule from  9.842 18.098 to  9.940 18.098
\putrule from  9.940 18.098 to 10.035 18.098
\putrule from 10.035 18.098 to 10.126 18.098
\putrule from 10.126 18.098 to 10.213 18.098
\putrule from 10.213 18.098 to 10.295 18.098
\putrule from 10.295 18.098 to 10.376 18.098
\plot 10.376 18.098 10.454 18.095 /
\putrule from 10.454 18.095 to 10.530 18.095
\putrule from 10.530 18.095 to 10.602 18.095
\plot 10.602 18.095 10.672 18.093 /
\plot 10.672 18.093 10.738 18.091 /
\putrule from 10.738 18.091 to 10.801 18.091
\plot 10.801 18.091 10.861 18.089 /
\plot 10.861 18.089 10.918 18.087 /
\plot 10.918 18.087 10.973 18.085 /
\plot 10.973 18.085 11.024 18.081 /
\plot 11.024 18.081 11.072 18.078 /
\plot 11.072 18.078 11.117 18.074 /
\plot 11.117 18.074 11.159 18.070 /
\plot 11.159 18.070 11.199 18.066 /
\plot 11.199 18.066 11.235 18.062 /
\plot 11.235 18.062 11.271 18.057 /
\plot 11.271 18.057 11.305 18.053 /
\plot 11.305 18.053 11.337 18.047 /
\plot 11.337 18.047 11.369 18.040 /
\plot 11.369 18.040 11.398 18.034 /
\plot 11.398 18.034 11.443 18.023 /
\plot 11.443 18.023 11.485 18.013 /
\plot 11.485 18.013 11.529 18.000 /
\plot 11.529 18.000 11.572 17.985 /
\plot 11.572 17.985 11.616 17.971 /
\plot 11.616 17.971 11.663 17.954 /
\plot 11.663 17.954 11.705 17.935 /
\plot 11.705 17.935 11.750 17.915 /
\plot 11.750 17.915 11.794 17.894 /
\plot 11.794 17.894 11.834 17.871 /
\plot 11.834 17.871 11.877 17.850 /
\plot 11.877 17.850 11.915 17.827 /
\plot 11.915 17.827 11.951 17.803 /
\plot 11.951 17.803 11.985 17.780 /
\plot 11.985 17.780 12.016 17.757 /
\plot 12.016 17.757 12.044 17.733 /
\plot 12.044 17.733 12.071 17.708 /
\plot 12.071 17.708 12.097 17.685 /
\plot 12.097 17.685 12.126 17.653 /
\plot 12.126 17.653 12.154 17.619 /
\plot 12.154 17.619 12.181 17.585 /
\plot 12.181 17.585 12.207 17.549 /
\plot 12.207 17.549 12.228 17.511 /
\plot 12.228 17.511 12.251 17.473 /
\plot 12.251 17.473 12.270 17.435 /
\plot 12.270 17.435 12.287 17.397 /
\plot 12.287 17.397 12.302 17.361 /
\plot 12.302 17.361 12.315 17.325 /
\plot 12.315 17.325 12.327 17.293 /
\plot 12.327 17.293 12.336 17.261 /
\plot 12.336 17.261 12.344 17.234 /
\plot 12.344 17.234 12.351 17.209 /
\plot 12.351 17.209 12.359 17.175 /
\plot 12.359 17.175 12.366 17.143 /
\plot 12.366 17.143 12.370 17.113 /
\plot 12.370 17.113 12.374 17.084 /
\plot 12.374 17.084 12.378 17.054 /
\plot 12.378 17.054 12.380 17.026 /
\plot 12.380 17.026 12.383 16.999 /
\putrule from 12.383 16.999 to 12.383 16.974
\putrule from 12.383 16.974 to 12.383 16.948
\putrule from 12.383 16.948 to 12.383 16.923
\putrule from 12.383 16.923 to 12.383 16.899
\putrule from 12.383 16.899 to 12.383 16.876
\putrule from 12.383 16.876 to 12.383 16.847
\putrule from 12.383 16.847 to 12.383 16.815
\putrule from 12.383 16.815 to 12.383 16.777
\putrule from 12.383 16.777 to 12.383 16.734
\putrule from 12.383 16.734 to 12.383 16.688
\putrule from 12.383 16.688 to 12.383 16.645
\putrule from 12.383 16.645 to 12.383 16.609
\putrule from 12.383 16.609 to 12.383 16.586
\putrule from 12.383 16.586 to 12.383 16.576
\putrule from 12.383 16.576 to 12.383 16.573
}%
%
%
\linethickness= 0.500pt
\setplotsymbol ({\thinlinefont .})
{\color[rgb]{0,0,0}\putrule from 17.907 18.098 to 17.903 18.098
\putrule from 17.903 18.098 to 17.894 18.098
\putrule from 17.894 18.098 to 17.875 18.098
\putrule from 17.875 18.098 to 17.850 18.098
\putrule from 17.850 18.098 to 17.810 18.098
\putrule from 17.810 18.098 to 17.759 18.098
\putrule from 17.759 18.098 to 17.695 18.098
\putrule from 17.695 18.098 to 17.617 18.098
\putrule from 17.617 18.098 to 17.528 18.098
\putrule from 17.528 18.098 to 17.427 18.098
\putrule from 17.427 18.098 to 17.316 18.098
\putrule from 17.316 18.098 to 17.196 18.098
\putrule from 17.196 18.098 to 17.071 18.098
\putrule from 17.071 18.098 to 16.940 18.098
\putrule from 16.940 18.098 to 16.806 18.098
\putrule from 16.806 18.098 to 16.673 18.098
\putrule from 16.673 18.098 to 16.538 18.098
\putrule from 16.538 18.098 to 16.406 18.098
\putrule from 16.406 18.098 to 16.277 18.098
\putrule from 16.277 18.098 to 16.152 18.098
\putrule from 16.152 18.098 to 16.032 18.098
\putrule from 16.032 18.098 to 15.915 18.098
\putrule from 15.915 18.098 to 15.805 18.098
\putrule from 15.805 18.098 to 15.699 18.098
\putrule from 15.699 18.098 to 15.600 18.098
\putrule from 15.600 18.098 to 15.505 18.098
\putrule from 15.505 18.098 to 15.416 18.098
\putrule from 15.416 18.098 to 15.331 18.098
\putrule from 15.331 18.098 to 15.251 18.098
\putrule from 15.251 18.098 to 15.174 18.098
\putrule from 15.174 18.098 to 15.102 18.098
\putrule from 15.102 18.098 to 15.035 18.098
\putrule from 15.035 18.098 to 14.969 18.098
\putrule from 14.969 18.098 to 14.906 18.098
\putrule from 14.906 18.098 to 14.846 18.098
\putrule from 14.846 18.098 to 14.787 18.098
\putrule from 14.787 18.098 to 14.732 18.098
\putrule from 14.732 18.098 to 14.647 18.098
\putrule from 14.647 18.098 to 14.565 18.098
\putrule from 14.565 18.098 to 14.486 18.098
\putrule from 14.486 18.098 to 14.408 18.098
\putrule from 14.408 18.098 to 14.334 18.098
\plot 14.334 18.098 14.262 18.100 /
\putrule from 14.262 18.100 to 14.192 18.100
\plot 14.192 18.100 14.125 18.102 /
\plot 14.125 18.102 14.057 18.104 /
\plot 14.057 18.104 13.993 18.106 /
\plot 13.993 18.106 13.932 18.108 /
\plot 13.932 18.108 13.873 18.110 /
\plot 13.873 18.110 13.815 18.114 /
\plot 13.815 18.114 13.763 18.119 /
\plot 13.763 18.119 13.710 18.123 /
\plot 13.710 18.123 13.661 18.129 /
\plot 13.661 18.129 13.614 18.133 /
\plot 13.614 18.133 13.568 18.140 /
\plot 13.568 18.140 13.526 18.148 /
\plot 13.526 18.148 13.485 18.155 /
\plot 13.485 18.155 13.447 18.163 /
\plot 13.447 18.163 13.409 18.172 /
\plot 13.409 18.172 13.371 18.182 /
\plot 13.371 18.182 13.335 18.193 /
\plot 13.335 18.193 13.291 18.205 /
\plot 13.291 18.205 13.248 18.222 /
\plot 13.248 18.222 13.202 18.239 /
\plot 13.202 18.239 13.157 18.256 /
\plot 13.157 18.256 13.111 18.275 /
\plot 13.111 18.275 13.064 18.299 /
\plot 13.064 18.299 13.018 18.320 /
\plot 13.018 18.320 12.973 18.345 /
\plot 12.973 18.345 12.926 18.371 /
\plot 12.926 18.371 12.880 18.396 /
\plot 12.880 18.396 12.838 18.423 /
\plot 12.838 18.423 12.795 18.451 /
\plot 12.795 18.451 12.755 18.479 /
\plot 12.755 18.479 12.717 18.506 /
\plot 12.717 18.506 12.681 18.536 /
\plot 12.681 18.536 12.647 18.563 /
\plot 12.647 18.563 12.617 18.589 /
\plot 12.617 18.589 12.590 18.616 /
\plot 12.590 18.616 12.565 18.644 /
\plot 12.565 18.644 12.541 18.669 /
\plot 12.541 18.669 12.516 18.703 /
\plot 12.516 18.703 12.493 18.735 /
\plot 12.493 18.735 12.471 18.768 /
\plot 12.471 18.768 12.454 18.804 /
\plot 12.454 18.804 12.438 18.840 /
\plot 12.438 18.840 12.425 18.876 /
\plot 12.425 18.876 12.414 18.912 /
\plot 12.414 18.912 12.404 18.951 /
\plot 12.404 18.951 12.397 18.986 /
\plot 12.397 18.986 12.391 19.022 /
\plot 12.391 19.022 12.387 19.056 /
\plot 12.387 19.056 12.385 19.090 /
\putrule from 12.385 19.090 to 12.385 19.122
\plot 12.385 19.122 12.383 19.152 /
\putrule from 12.383 19.152 to 12.383 19.181
\putrule from 12.383 19.181 to 12.383 19.209
\putrule from 12.383 19.209 to 12.383 19.245
\putrule from 12.383 19.245 to 12.383 19.281
\putrule from 12.383 19.281 to 12.383 19.317
\putrule from 12.383 19.317 to 12.383 19.357
\putrule from 12.383 19.357 to 12.383 19.401
\putrule from 12.383 19.401 to 12.383 19.450
\putrule from 12.383 19.450 to 12.383 19.499
\putrule from 12.383 19.499 to 12.383 19.547
\putrule from 12.383 19.547 to 12.383 19.586
\putrule from 12.383 19.586 to 12.383 19.609
\putrule from 12.383 19.609 to 12.383 19.619
\putrule from 12.383 19.619 to 12.383 19.622
}%
%
%
\put{$H$
} [lB] at 15.240 17.216
%
%
\put{$S^2 \times 0$
} [lB] at  3.429 11.502
%
%
\put{$S^2 \times 1 \sim S^2 \times 0$
} [lB] at 17.716 11.502
\linethickness=0pt
\putrectangle corners at  3.397 23.838 and 21.933 11.970
\endpicture}

%% file: FindingAnHSlide.tex
\font\thinlinefont=cmr5
\mbox{\beginpicture
\small
\setcoordinatesystem units <0.4cm,0.4cm>
\unitlength=1.04987cm
\linethickness=1pt
\setplotsymbol ({\makebox(0,0)[l]{\tencirc\symbol{'160}}})
\setshadesymbol ({\thinlinefont .})
\setlinear
%
%
\linethickness= 0.500pt
\setplotsymbol ({\thinlinefont .})
{\color[rgb]{0,0,0}\circulararc 106.260 degrees from 12.002 20.003 center at 12.383 20.288
}%
%
%
\linethickness= 0.500pt
\setplotsymbol ({\thinlinefont .})
\setdots < 0.0953cm>
{\color[rgb]{0,0,0}\circulararc 106.260 degrees from 12.764 20.003 center at 12.383 19.717
}%
%
%
\linethickness= 0.500pt
\setplotsymbol ({\thinlinefont .})
\setsolid
{\color[rgb]{0,0,0}\ellipticalarc axes ratio  0.381:0.381  360 degrees 
	from 12.764 20.003 center at 12.383 20.003
}%
%
%
\linethickness= 0.500pt
\setplotsymbol ({\thinlinefont .})
{\color[rgb]{0,0,0}\circulararc 106.260 degrees from 12.002 16.192 center at 12.383 16.478
}%
%
%
\linethickness= 0.500pt
\setplotsymbol ({\thinlinefont .})
\setdots < 0.0953cm>
{\color[rgb]{0,0,0}\circulararc 106.260 degrees from 12.764 16.192 center at 12.383 15.907
}%
%
%
\linethickness= 0.500pt
\setplotsymbol ({\thinlinefont .})
\setsolid
{\color[rgb]{0,0,0}\ellipticalarc axes ratio  0.381:0.381  360 degrees 
	from 12.764 16.192 center at 12.383 16.192
}%
%
%
\linethickness= 0.500pt
\setplotsymbol ({\thinlinefont .})
{\color[rgb]{0,0,0}\putrule from 18.098 20.955 to 18.098 12.383
\plot 18.098 12.383 21.907 15.240 /
\putrule from 21.907 15.240 to 21.907 23.812
\plot 21.907 23.812 18.098 20.955 /
}%
%
%
\linethickness= 0.500pt
\setplotsymbol ({\thinlinefont .})
{\color[rgb]{0,0,0}\putrule from 20.003 20.003 to 18.288 20.003
}%
%
%
\linethickness= 0.500pt
\setplotsymbol ({\thinlinefont .})
{\color[rgb]{0,0,0}\putrule from 20.003 18.098 to 18.479 18.098
}%
%
%
\linethickness= 0.500pt
\setplotsymbol ({\thinlinefont .})
{\color[rgb]{0,0,0}\putrule from 20.003 16.192 to 18.479 16.192
}%
%
%
\linethickness= 0.500pt
\setplotsymbol ({\thinlinefont .})
{\color[rgb]{0,0,0}\putrule from 12.764 20.003 to 17.907 20.003
}%
%
%
\linethickness= 0.500pt
\setplotsymbol ({\thinlinefont .})
{\color[rgb]{0,0,0}\putrule from 12.764 16.192 to 17.907 16.192
}%
%
%
\linethickness= 0.500pt
\setplotsymbol ({\thinlinefont .})
{\color[rgb]{0,0,0}\putrule from 12.002 20.003 to  5.715 20.003
}%
%
%
\linethickness= 0.500pt
\setplotsymbol ({\thinlinefont .})
{\color[rgb]{0,0,0}\putrule from 12.002 16.192 to  5.715 16.192
}%
%
%
\linethickness= 0.500pt
\setplotsymbol ({\thinlinefont .})
{\color[rgb]{0,0,0}\putrule from  3.810 20.955 to  3.810 12.383
\plot  3.810 12.383  7.620 15.240 /
\putrule from  7.620 15.240 to  7.620 16.002
}%
%
%
\linethickness= 0.500pt
\setplotsymbol ({\thinlinefont .})
{\color[rgb]{0,0,0}\putrule from  7.620 20.193 to  7.620 23.812
\plot  7.620 23.812  3.810 20.955 /
}%
%
%
\linethickness= 0.500pt
\setplotsymbol ({\thinlinefont .})
{\color[rgb]{0,0,0}\putrule from  7.620 19.812 to  7.620 18.288
}%
%
%
\linethickness= 0.500pt
\setplotsymbol ({\thinlinefont .})
{\color[rgb]{0,0,0}\putrule from  7.620 17.907 to  7.620 16.383
}%
%
%
\linethickness= 0.500pt
\setplotsymbol ({\thinlinefont .})
{\color[rgb]{0,0,0}\putrule from  5.715 18.098 to  5.719 18.098
\putrule from  5.719 18.098 to  5.726 18.098
\putrule from  5.726 18.098 to  5.743 18.098
\putrule from  5.743 18.098 to  5.766 18.098
\putrule from  5.766 18.098 to  5.800 18.098
\putrule from  5.800 18.098 to  5.844 18.098
\putrule from  5.844 18.098 to  5.903 18.098
\putrule from  5.903 18.098 to  5.973 18.098
\putrule from  5.973 18.098 to  6.056 18.098
\putrule from  6.056 18.098 to  6.151 18.098
\putrule from  6.151 18.098 to  6.259 18.098
\putrule from  6.259 18.098 to  6.378 18.098
\putrule from  6.378 18.098 to  6.505 18.098
\putrule from  6.505 18.098 to  6.640 18.098
\putrule from  6.640 18.098 to  6.782 18.098
\putrule from  6.782 18.098 to  6.928 18.098
\putrule from  6.928 18.098 to  7.078 18.098
\putrule from  7.078 18.098 to  7.228 18.098
\putrule from  7.228 18.098 to  7.381 18.098
\putrule from  7.381 18.098 to  7.533 18.098
\putrule from  7.533 18.098 to  7.681 18.098
\putrule from  7.681 18.098 to  7.827 18.098
\putrule from  7.827 18.098 to  7.969 18.098
\putrule from  7.969 18.098 to  8.107 18.098
\putrule from  8.107 18.098 to  8.240 18.098
\putrule from  8.240 18.098 to  8.367 18.098
\putrule from  8.367 18.098 to  8.490 18.098
\putrule from  8.490 18.098 to  8.608 18.098
\putrule from  8.608 18.098 to  8.721 18.098
\putrule from  8.721 18.098 to  8.826 18.098
\putrule from  8.826 18.098 to  8.928 18.098
\putrule from  8.928 18.098 to  9.025 18.098
\putrule from  9.025 18.098 to  9.119 18.098
\putrule from  9.119 18.098 to  9.205 18.098
\putrule from  9.205 18.098 to  9.290 18.098
\putrule from  9.290 18.098 to  9.368 18.098
\putrule from  9.368 18.098 to  9.445 18.098
\putrule from  9.445 18.098 to  9.519 18.098
\putrule from  9.519 18.098 to  9.588 18.098
\putrule from  9.588 18.098 to  9.654 18.098
\putrule from  9.654 18.098 to  9.720 18.098
\putrule from  9.720 18.098 to  9.781 18.098
\putrule from  9.781 18.098 to  9.842 18.098
\putrule from  9.842 18.098 to  9.940 18.098
\putrule from  9.940 18.098 to 10.035 18.098
\putrule from 10.035 18.098 to 10.126 18.098
\putrule from 10.126 18.098 to 10.213 18.098
\putrule from 10.213 18.098 to 10.295 18.098
\putrule from 10.295 18.098 to 10.376 18.098
\plot 10.376 18.098 10.454 18.095 /
\putrule from 10.454 18.095 to 10.530 18.095
\putrule from 10.530 18.095 to 10.602 18.095
\plot 10.602 18.095 10.672 18.093 /
\plot 10.672 18.093 10.738 18.091 /
\putrule from 10.738 18.091 to 10.801 18.091
\plot 10.801 18.091 10.861 18.089 /
\plot 10.861 18.089 10.918 18.087 /
\plot 10.918 18.087 10.973 18.085 /
\plot 10.973 18.085 11.024 18.081 /
\plot 11.024 18.081 11.072 18.078 /
\plot 11.072 18.078 11.117 18.074 /
\plot 11.117 18.074 11.159 18.070 /
\plot 11.159 18.070 11.199 18.066 /
\plot 11.199 18.066 11.235 18.062 /
\plot 11.235 18.062 11.271 18.057 /
\plot 11.271 18.057 11.305 18.053 /
\plot 11.305 18.053 11.337 18.047 /
\plot 11.337 18.047 11.369 18.040 /
\plot 11.369 18.040 11.398 18.034 /
\plot 11.398 18.034 11.443 18.023 /
\plot 11.443 18.023 11.485 18.013 /
\plot 11.485 18.013 11.529 18.000 /
\plot 11.529 18.000 11.572 17.985 /
\plot 11.572 17.985 11.616 17.971 /
\plot 11.616 17.971 11.663 17.954 /
\plot 11.663 17.954 11.705 17.935 /
\plot 11.705 17.935 11.750 17.915 /
\plot 11.750 17.915 11.794 17.894 /
\plot 11.794 17.894 11.834 17.871 /
\plot 11.834 17.871 11.877 17.850 /
\plot 11.877 17.850 11.915 17.827 /
\plot 11.915 17.827 11.951 17.803 /
\plot 11.951 17.803 11.985 17.780 /
\plot 11.985 17.780 12.016 17.757 /
\plot 12.016 17.757 12.044 17.733 /
\plot 12.044 17.733 12.071 17.708 /
\plot 12.071 17.708 12.097 17.685 /
\plot 12.097 17.685 12.126 17.653 /
\plot 12.126 17.653 12.154 17.619 /
\plot 12.154 17.619 12.181 17.585 /
\plot 12.181 17.585 12.207 17.549 /
\plot 12.207 17.549 12.228 17.511 /
\plot 12.228 17.511 12.251 17.473 /
\plot 12.251 17.473 12.270 17.435 /
\plot 12.270 17.435 12.287 17.397 /
\plot 12.287 17.397 12.302 17.361 /
\plot 12.302 17.361 12.315 17.325 /
\plot 12.315 17.325 12.327 17.293 /
\plot 12.327 17.293 12.336 17.261 /
\plot 12.336 17.261 12.344 17.234 /
\plot 12.344 17.234 12.351 17.209 /
\plot 12.351 17.209 12.359 17.175 /
\plot 12.359 17.175 12.366 17.143 /
\plot 12.366 17.143 12.370 17.113 /
\plot 12.370 17.113 12.374 17.084 /
\plot 12.374 17.084 12.378 17.054 /
\plot 12.378 17.054 12.380 17.026 /
\plot 12.380 17.026 12.383 16.999 /
\putrule from 12.383 16.999 to 12.383 16.974
\putrule from 12.383 16.974 to 12.383 16.948
\putrule from 12.383 16.948 to 12.383 16.923
\putrule from 12.383 16.923 to 12.383 16.899
\putrule from 12.383 16.899 to 12.383 16.876
\putrule from 12.383 16.876 to 12.383 16.847
\putrule from 12.383 16.847 to 12.383 16.815
\putrule from 12.383 16.815 to 12.383 16.777
\putrule from 12.383 16.777 to 12.383 16.734
\putrule from 12.383 16.734 to 12.383 16.688
\putrule from 12.383 16.688 to 12.383 16.645
\putrule from 12.383 16.645 to 12.383 16.609
\putrule from 12.383 16.609 to 12.383 16.586
\putrule from 12.383 16.586 to 12.383 16.576
\putrule from 12.383 16.576 to 12.383 16.573
}%
%
%
\linethickness= 0.500pt
\setplotsymbol ({\thinlinefont .})
{\color[rgb]{0,0,0}\putrule from 17.907 18.098 to 17.903 18.098
\putrule from 17.903 18.098 to 17.894 18.098
\putrule from 17.894 18.098 to 17.875 18.098
\putrule from 17.875 18.098 to 17.850 18.098
\putrule from 17.850 18.098 to 17.810 18.098
\putrule from 17.810 18.098 to 17.759 18.098
\putrule from 17.759 18.098 to 17.695 18.098
\putrule from 17.695 18.098 to 17.617 18.098
\putrule from 17.617 18.098 to 17.528 18.098
\putrule from 17.528 18.098 to 17.427 18.098
\putrule from 17.427 18.098 to 17.316 18.098
\putrule from 17.316 18.098 to 17.196 18.098
\putrule from 17.196 18.098 to 17.071 18.098
\putrule from 17.071 18.098 to 16.940 18.098
\putrule from 16.940 18.098 to 16.806 18.098
\putrule from 16.806 18.098 to 16.673 18.098
\putrule from 16.673 18.098 to 16.538 18.098
\putrule from 16.538 18.098 to 16.406 18.098
\putrule from 16.406 18.098 to 16.277 18.098
\putrule from 16.277 18.098 to 16.152 18.098
\putrule from 16.152 18.098 to 16.032 18.098
\putrule from 16.032 18.098 to 15.915 18.098
\putrule from 15.915 18.098 to 15.805 18.098
\putrule from 15.805 18.098 to 15.699 18.098
\putrule from 15.699 18.098 to 15.600 18.098
\putrule from 15.600 18.098 to 15.505 18.098
\putrule from 15.505 18.098 to 15.416 18.098
\putrule from 15.416 18.098 to 15.331 18.098
\putrule from 15.331 18.098 to 15.251 18.098
\putrule from 15.251 18.098 to 15.174 18.098
\putrule from 15.174 18.098 to 15.102 18.098
\putrule from 15.102 18.098 to 15.035 18.098
\putrule from 15.035 18.098 to 14.969 18.098
\putrule from 14.969 18.098 to 14.906 18.098
\putrule from 14.906 18.098 to 14.846 18.098
\putrule from 14.846 18.098 to 14.787 18.098
\putrule from 14.787 18.098 to 14.732 18.098
\putrule from 14.732 18.098 to 14.647 18.098
\putrule from 14.647 18.098 to 14.565 18.098
\putrule from 14.565 18.098 to 14.486 18.098
\putrule from 14.486 18.098 to 14.408 18.098
\putrule from 14.408 18.098 to 14.334 18.098
\plot 14.334 18.098 14.262 18.100 /
\putrule from 14.262 18.100 to 14.192 18.100
\plot 14.192 18.100 14.125 18.102 /
\plot 14.125 18.102 14.057 18.104 /
\plot 14.057 18.104 13.993 18.106 /
\plot 13.993 18.106 13.932 18.108 /
\plot 13.932 18.108 13.873 18.110 /
\plot 13.873 18.110 13.815 18.114 /
\plot 13.815 18.114 13.763 18.119 /
\plot 13.763 18.119 13.710 18.123 /
\plot 13.710 18.123 13.661 18.129 /
\plot 13.661 18.129 13.614 18.133 /
\plot 13.614 18.133 13.568 18.140 /
\plot 13.568 18.140 13.526 18.148 /
\plot 13.526 18.148 13.485 18.155 /
\plot 13.485 18.155 13.447 18.163 /
\plot 13.447 18.163 13.409 18.172 /
\plot 13.409 18.172 13.371 18.182 /
\plot 13.371 18.182 13.335 18.193 /
\plot 13.335 18.193 13.291 18.205 /
\plot 13.291 18.205 13.248 18.222 /
\plot 13.248 18.222 13.202 18.239 /
\plot 13.202 18.239 13.157 18.256 /
\plot 13.157 18.256 13.111 18.275 /
\plot 13.111 18.275 13.064 18.299 /
\plot 13.064 18.299 13.018 18.320 /
\plot 13.018 18.320 12.973 18.345 /
\plot 12.973 18.345 12.926 18.371 /
\plot 12.926 18.371 12.880 18.396 /
\plot 12.880 18.396 12.838 18.423 /
\plot 12.838 18.423 12.795 18.451 /
\plot 12.795 18.451 12.755 18.479 /
\plot 12.755 18.479 12.717 18.506 /
\plot 12.717 18.506 12.681 18.536 /
\plot 12.681 18.536 12.647 18.563 /
\plot 12.647 18.563 12.617 18.589 /
\plot 12.617 18.589 12.590 18.616 /
\plot 12.590 18.616 12.565 18.644 /
\plot 12.565 18.644 12.541 18.669 /
\plot 12.541 18.669 12.516 18.703 /
\plot 12.516 18.703 12.493 18.735 /
\plot 12.493 18.735 12.471 18.768 /
\plot 12.471 18.768 12.454 18.804 /
\plot 12.454 18.804 12.438 18.840 /
\plot 12.438 18.840 12.425 18.876 /
\plot 12.425 18.876 12.414 18.912 /
\plot 12.414 18.912 12.404 18.951 /
\plot 12.404 18.951 12.397 18.986 /
\plot 12.397 18.986 12.391 19.022 /
\plot 12.391 19.022 12.387 19.056 /
\plot 12.387 19.056 12.385 19.090 /
\putrule from 12.385 19.090 to 12.385 19.122
\plot 12.385 19.122 12.383 19.152 /
\putrule from 12.383 19.152 to 12.383 19.181
\putrule from 12.383 19.181 to 12.383 19.209
\putrule from 12.383 19.209 to 12.383 19.245
\putrule from 12.383 19.245 to 12.383 19.281
\putrule from 12.383 19.281 to 12.383 19.317
\putrule from 12.383 19.317 to 12.383 19.357
\putrule from 12.383 19.357 to 12.383 19.401
\putrule from 12.383 19.401 to 12.383 19.450
\putrule from 12.383 19.450 to 12.383 19.499
\putrule from 12.383 19.499 to 12.383 19.547
\putrule from 12.383 19.547 to 12.383 19.586
\putrule from 12.383 19.586 to 12.383 19.609
\putrule from 12.383 19.609 to 12.383 19.619
\putrule from 12.383 19.619 to 12.383 19.622
}%
%
%
\linethickness= 0.500pt
\setplotsymbol ({\thinlinefont .})
\setdots < 0.0953cm>
{\color[rgb]{0,0,0}\plot 14.287 20.003 14.287 19.996 /
\plot 14.287 19.996 14.287 19.979 /
\plot 14.287 19.979 14.287 19.954 /
\plot 14.287 19.954 14.287 19.914 /
\plot 14.287 19.914 14.287 19.863 /
\plot 14.287 19.863 14.287 19.804 /
\plot 14.287 19.804 14.287 19.738 /
\plot 14.287 19.738 14.287 19.670 /
\plot 14.287 19.670 14.287 19.605 /
\plot 14.287 19.605 14.287 19.539 /
\plot 14.287 19.539 14.287 19.478 /
\plot 14.287 19.478 14.287 19.423 /
\plot 14.287 19.423 14.287 19.372 /
\plot 14.287 19.372 14.287 19.325 /
\plot 14.287 19.325 14.287 19.283 /
\plot 14.287 19.283 14.287 19.245 /
\plot 14.287 19.245 14.287 19.209 /
\plot 14.287 19.209 14.287 19.177 /
\plot 14.287 19.177 14.287 19.145 /
\plot 14.287 19.145 14.287 19.103 /
\plot 14.287 19.103 14.287 19.063 /
\plot 14.287 19.063 14.287 19.022 /
\plot 14.287 19.022 14.287 18.984 /
\plot 14.287 18.984 14.287 18.946 /
\plot 14.287 18.946 14.287 18.910 /
\plot 14.287 18.910 14.287 18.876 /
\plot 14.287 18.876 14.287 18.845 /
\plot 14.287 18.845 14.287 18.815 /
\plot 14.287 18.815 14.287 18.788 /
\plot 14.287 18.788 14.287 18.764 /
\plot 14.287 18.764 14.287 18.741 /
\plot 14.287 18.741 14.287 18.720 /
\plot 14.287 18.720 14.287 18.701 /
\plot 14.287 18.701 14.287 18.675 /
\plot 14.287 18.675 14.287 18.650 /
\plot 14.287 18.650 14.287 18.627 /
\plot 14.287 18.627 14.287 18.601 /
\plot 14.287 18.601 14.287 18.578 /
\plot 14.287 18.578 14.287 18.555 /
\plot 14.287 18.555 14.287 18.534 /
\plot 14.287 18.534 14.287 18.514 /
\plot 14.287 18.514 14.287 18.495 /
\plot 14.287 18.495 14.287 18.479 /
\plot 14.287 18.479 14.287 18.466 /
\plot 14.287 18.466 14.287 18.453 /
\plot 14.287 18.453 14.285 18.440 /
\plot 14.285 18.440 14.283 18.428 /
\plot 14.283 18.428 14.281 18.413 /
\plot 14.281 18.413 14.277 18.400 /
\plot 14.277 18.400 14.271 18.387 /
\plot 14.271 18.387 14.262 18.375 /
\plot 14.262 18.375 14.252 18.364 /
\plot 14.252 18.364 14.237 18.354 /
\plot 14.237 18.354 14.222 18.343 /
\plot 14.222 18.343 14.205 18.335 /
\plot 14.205 18.335 14.184 18.326 /
\plot 14.184 18.326 14.161 18.320 /
\plot 14.161 18.320 14.137 18.316 /
\plot 14.137 18.316 14.112 18.309 /
\plot 14.112 18.309 14.082 18.307 /
\plot 14.082 18.307 14.048 18.303 /
\plot 14.048 18.303 14.012 18.301 /
\plot 14.012 18.301 13.974 18.301 /
\plot 13.974 18.301 13.932 18.301 /
\plot 13.932 18.301 13.890 18.301 /
\plot 13.890 18.301 13.843 18.303 /
\plot 13.843 18.303 13.796 18.307 /
\plot 13.796 18.307 13.750 18.311 /
\plot 13.750 18.311 13.703 18.318 /
\plot 13.703 18.318 13.659 18.324 /
\plot 13.659 18.324 13.614 18.332 /
\plot 13.614 18.332 13.570 18.341 /
\plot 13.570 18.341 13.526 18.352 /
\plot 13.526 18.352 13.487 18.362 /
\plot 13.487 18.362 13.447 18.373 /
\plot 13.447 18.373 13.407 18.385 /
\plot 13.407 18.385 13.367 18.400 /
\plot 13.367 18.400 13.324 18.415 /
\plot 13.324 18.415 13.284 18.432 /
\plot 13.284 18.432 13.242 18.449 /
\plot 13.242 18.449 13.200 18.468 /
\plot 13.200 18.468 13.157 18.487 /
\plot 13.157 18.487 13.117 18.508 /
\plot 13.117 18.508 13.077 18.527 /
\plot 13.077 18.527 13.041 18.548 /
\plot 13.041 18.548 13.005 18.570 /
\plot 13.005 18.570 12.971 18.591 /
\plot 12.971 18.591 12.939 18.610 /
\plot 12.939 18.610 12.912 18.631 /
\plot 12.912 18.631 12.884 18.650 /
\plot 12.884 18.650 12.859 18.669 /
\plot 12.859 18.669 12.829 18.694 /
\plot 12.829 18.694 12.802 18.720 /
\plot 12.802 18.720 12.774 18.747 /
\plot 12.774 18.747 12.749 18.775 /
\plot 12.749 18.775 12.728 18.802 /
\plot 12.728 18.802 12.704 18.832 /
\plot 12.704 18.832 12.685 18.864 /
\plot 12.685 18.864 12.668 18.895 /
\plot 12.668 18.895 12.653 18.927 /
\plot 12.653 18.927 12.641 18.959 /
\plot 12.641 18.959 12.628 18.989 /
\plot 12.628 18.989 12.620 19.020 /
\plot 12.620 19.020 12.611 19.050 /
\plot 12.611 19.050 12.605 19.082 /
\plot 12.605 19.082 12.598 19.113 /
\plot 12.598 19.113 12.594 19.147 /
\plot 12.594 19.147 12.590 19.183 /
\plot 12.590 19.183 12.586 19.224 /
\plot 12.586 19.224 12.584 19.270 /
\plot 12.584 19.270 12.581 19.321 /
\plot 12.581 19.321 12.579 19.378 /
\plot 12.579 19.378 12.577 19.437 /
\plot 12.577 19.437 12.575 19.494 /
\plot 12.575 19.494 12.575 19.547 /
\plot 12.575 19.547 12.573 19.586 /
\plot 12.573 19.586 12.573 19.609 /
\plot 12.573 19.609 12.573 19.619 /
\plot 12.573 19.619 12.573 19.622 /
}%
%
%
\linethickness= 0.500pt
\setplotsymbol ({\thinlinefont .})
{\color[rgb]{0,0,0}\plot 14.478 20.003 14.478 19.994 /
\plot 14.478 19.994 14.478 19.979 /
\plot 14.478 19.979 14.478 19.950 /
\plot 14.478 19.950 14.478 19.907 /
\plot 14.478 19.907 14.478 19.850 /
\plot 14.478 19.850 14.478 19.780 /
\plot 14.478 19.780 14.478 19.702 /
\plot 14.478 19.702 14.478 19.619 /
\plot 14.478 19.619 14.478 19.533 /
\plot 14.478 19.533 14.478 19.448 /
\plot 14.478 19.448 14.478 19.367 /
\plot 14.478 19.367 14.478 19.291 /
\plot 14.478 19.291 14.478 19.219 /
\plot 14.478 19.219 14.478 19.154 /
\plot 14.478 19.154 14.478 19.094 /
\plot 14.478 19.094 14.478 19.039 /
\plot 14.478 19.039 14.478 18.991 /
\plot 14.478 18.991 14.478 18.944 /
\plot 14.478 18.944 14.478 18.904 /
\plot 14.478 18.904 14.478 18.864 /
\plot 14.478 18.864 14.478 18.828 /
\plot 14.478 18.828 14.478 18.794 /
\plot 14.478 18.794 14.478 18.760 /
\plot 14.478 18.760 14.478 18.726 /
\plot 14.478 18.726 14.480 18.697 /
\plot 14.480 18.697 14.482 18.667 /
\plot 14.482 18.667 14.484 18.637 /
\plot 14.484 18.637 14.489 18.610 /
\plot 14.489 18.610 14.495 18.582 /
\plot 14.495 18.582 14.501 18.557 /
\plot 14.501 18.557 14.512 18.531 /
\plot 14.512 18.531 14.522 18.510 /
\plot 14.522 18.510 14.535 18.487 /
\plot 14.535 18.487 14.550 18.468 /
\plot 14.550 18.468 14.567 18.449 /
\plot 14.567 18.449 14.588 18.432 /
\plot 14.588 18.432 14.609 18.417 /
\plot 14.609 18.417 14.635 18.402 /
\plot 14.635 18.402 14.660 18.390 /
\plot 14.660 18.390 14.690 18.379 /
\plot 14.690 18.379 14.724 18.368 /
\plot 14.724 18.368 14.757 18.360 /
\plot 14.757 18.360 14.796 18.352 /
\plot 14.796 18.352 14.834 18.345 /
\plot 14.834 18.345 14.874 18.339 /
\plot 14.874 18.339 14.918 18.332 /
\plot 14.918 18.332 14.967 18.326 /
\plot 14.967 18.326 15.018 18.322 /
\plot 15.018 18.322 15.071 18.318 /
\plot 15.071 18.318 15.128 18.313 /
\plot 15.128 18.313 15.189 18.309 /
\plot 15.189 18.309 15.253 18.305 /
\plot 15.253 18.305 15.316 18.303 /
\plot 15.316 18.303 15.384 18.299 /
\plot 15.384 18.299 15.452 18.296 /
\plot 15.452 18.296 15.522 18.294 /
\plot 15.522 18.294 15.591 18.292 /
\plot 15.591 18.292 15.661 18.292 /
\plot 15.661 18.292 15.729 18.290 /
\plot 15.729 18.290 15.797 18.290 /
\plot 15.797 18.290 15.864 18.290 /
\plot 15.864 18.290 15.928 18.288 /
\plot 15.928 18.288 15.991 18.288 /
\plot 15.991 18.288 16.053 18.288 /
\plot 16.053 18.288 16.112 18.288 /
\plot 16.112 18.288 16.169 18.288 /
\plot 16.169 18.288 16.224 18.288 /
\plot 16.224 18.288 16.284 18.288 /
\plot 16.284 18.288 16.341 18.288 /
\plot 16.341 18.288 16.398 18.288 /
\plot 16.398 18.288 16.453 18.288 /
\plot 16.453 18.288 16.508 18.288 /
\plot 16.508 18.288 16.563 18.288 /
\plot 16.563 18.288 16.618 18.288 /
\plot 16.618 18.288 16.673 18.288 /
\plot 16.673 18.288 16.726 18.288 /
\plot 16.726 18.288 16.777 18.288 /
\plot 16.777 18.288 16.828 18.288 /
\plot 16.828 18.288 16.878 18.288 /
\plot 16.878 18.288 16.925 18.288 /
\plot 16.925 18.288 16.971 18.288 /
\plot 16.971 18.288 17.016 18.288 /
\plot 17.016 18.288 17.058 18.288 /
\plot 17.058 18.288 17.098 18.288 /
\plot 17.098 18.288 17.137 18.288 /
\plot 17.137 18.288 17.173 18.288 /
\plot 17.173 18.288 17.206 18.288 /
\plot 17.206 18.288 17.240 18.288 /
\plot 17.240 18.288 17.272 18.288 /
\plot 17.272 18.288 17.314 18.288 /
\plot 17.314 18.288 17.357 18.288 /
\plot 17.357 18.288 17.399 18.288 /
\plot 17.399 18.288 17.443 18.288 /
\plot 17.443 18.288 17.490 18.288 /
\plot 17.490 18.288 17.539 18.288 /
\plot 17.539 18.288 17.592 18.288 /
\plot 17.592 18.288 17.647 18.288 /
\plot 17.647 18.288 17.704 18.288 /
\plot 17.704 18.288 17.759 18.288 /
\plot 17.759 18.288 17.810 18.288 /
\plot 17.810 18.288 17.852 18.288 /
\plot 17.852 18.288 17.882 18.288 /
\plot 17.882 18.288 17.899 18.288 /
\plot 17.899 18.288 17.905 18.288 /
\plot 17.905 18.288 17.907 18.288 /
}%
%
%
\linethickness= 0.500pt
\setplotsymbol ({\thinlinefont .})
{\color[rgb]{0,0,0}\plot 18.479 18.288 20.003 18.288 /
}%
%
%
\linethickness= 0.500pt
\setplotsymbol ({\thinlinefont .})
{\color[rgb]{0,0,0}\plot 12.192 19.622 12.192 19.619 /
\plot 12.192 19.619 12.192 19.613 /
\plot 12.192 19.613 12.192 19.596 /
\plot 12.192 19.596 12.192 19.566 /
\plot 12.192 19.566 12.194 19.526 /
\plot 12.194 19.526 12.194 19.480 /
\plot 12.194 19.480 12.196 19.427 /
\plot 12.196 19.427 12.198 19.372 /
\plot 12.198 19.372 12.198 19.321 /
\plot 12.198 19.321 12.203 19.272 /
\plot 12.203 19.272 12.205 19.228 /
\plot 12.205 19.228 12.207 19.188 /
\plot 12.207 19.188 12.211 19.152 /
\plot 12.211 19.152 12.213 19.116 /
\plot 12.213 19.116 12.220 19.084 /
\plot 12.220 19.084 12.224 19.050 /
\plot 12.224 19.050 12.230 19.016 /
\plot 12.230 19.016 12.236 18.984 /
\plot 12.236 18.984 12.245 18.951 /
\plot 12.245 18.951 12.253 18.915 /
\plot 12.253 18.915 12.262 18.879 /
\plot 12.262 18.879 12.275 18.843 /
\plot 12.275 18.843 12.287 18.807 /
\plot 12.287 18.807 12.300 18.768 /
\plot 12.300 18.768 12.315 18.733 /
\plot 12.315 18.733 12.332 18.697 /
\plot 12.332 18.697 12.349 18.663 /
\plot 12.349 18.663 12.366 18.629 /
\plot 12.366 18.629 12.385 18.597 /
\plot 12.385 18.597 12.404 18.567 /
\plot 12.404 18.567 12.425 18.538 /
\plot 12.425 18.538 12.446 18.510 /
\plot 12.446 18.510 12.469 18.483 /
\plot 12.469 18.483 12.493 18.457 /
\plot 12.493 18.457 12.520 18.430 /
\plot 12.520 18.430 12.548 18.402 /
\plot 12.548 18.402 12.577 18.377 /
\plot 12.577 18.377 12.611 18.349 /
\plot 12.611 18.349 12.645 18.322 /
\plot 12.645 18.322 12.681 18.296 /
\plot 12.681 18.296 12.719 18.271 /
\plot 12.719 18.271 12.757 18.248 /
\plot 12.757 18.248 12.795 18.224 /
\plot 12.795 18.224 12.833 18.203 /
\plot 12.833 18.203 12.871 18.182 /
\plot 12.871 18.182 12.910 18.163 /
\plot 12.910 18.163 12.948 18.146 /
\plot 12.948 18.146 12.986 18.129 /
\plot 12.986 18.129 13.020 18.114 /
\plot 13.020 18.114 13.056 18.102 /
\plot 13.056 18.102 13.092 18.087 /
\plot 13.092 18.087 13.130 18.074 /
\plot 13.130 18.074 13.170 18.062 /
\plot 13.170 18.062 13.212 18.047 /
\plot 13.212 18.047 13.255 18.034 /
\plot 13.255 18.034 13.299 18.023 /
\plot 13.299 18.023 13.343 18.011 /
\plot 13.343 18.011 13.390 18.000 /
\plot 13.390 18.000 13.437 17.990 /
\plot 13.437 17.990 13.481 17.979 /
\plot 13.481 17.979 13.528 17.971 /
\plot 13.528 17.971 13.572 17.964 /
\plot 13.572 17.964 13.617 17.956 /
\plot 13.617 17.956 13.661 17.949 /
\plot 13.661 17.949 13.703 17.943 /
\plot 13.703 17.943 13.748 17.939 /
\plot 13.748 17.939 13.788 17.935 /
\plot 13.788 17.935 13.828 17.930 /
\plot 13.828 17.930 13.871 17.928 /
\plot 13.871 17.928 13.913 17.924 /
\plot 13.913 17.924 13.957 17.922 /
\plot 13.957 17.922 14.004 17.920 /
\plot 14.004 17.920 14.053 17.918 /
\plot 14.053 17.918 14.101 17.915 /
\plot 14.101 17.915 14.154 17.913 /
\plot 14.154 17.913 14.205 17.911 /
\plot 14.205 17.911 14.258 17.911 /
\plot 14.258 17.911 14.311 17.909 /
\plot 14.311 17.909 14.366 17.909 /
\plot 14.366 17.909 14.419 17.907 /
\plot 14.419 17.907 14.472 17.907 /
\plot 14.472 17.907 14.525 17.907 /
\plot 14.525 17.907 14.575 17.907 /
\plot 14.575 17.907 14.628 17.907 /
\plot 14.628 17.907 14.679 17.907 /
\plot 14.679 17.907 14.732 17.907 /
\plot 14.732 17.907 14.776 17.907 /
\plot 14.776 17.907 14.821 17.907 /
\plot 14.821 17.907 14.867 17.907 /
\plot 14.867 17.907 14.914 17.907 /
\plot 14.914 17.907 14.965 17.907 /
\plot 14.965 17.907 15.016 17.907 /
\plot 15.016 17.907 15.069 17.907 /
\plot 15.069 17.907 15.124 17.907 /
\plot 15.124 17.907 15.179 17.907 /
\plot 15.179 17.907 15.238 17.907 /
\plot 15.238 17.907 15.297 17.907 /
\plot 15.297 17.907 15.356 17.907 /
\plot 15.356 17.907 15.418 17.907 /
\plot 15.418 17.907 15.481 17.907 /
\plot 15.481 17.907 15.543 17.907 /
\plot 15.543 17.907 15.604 17.907 /
\plot 15.604 17.907 15.668 17.907 /
\plot 15.668 17.907 15.729 17.907 /
\plot 15.729 17.907 15.790 17.907 /
\plot 15.790 17.907 15.852 17.907 /
\plot 15.852 17.907 15.913 17.907 /
\plot 15.913 17.907 15.974 17.907 /
\plot 15.974 17.907 16.036 17.907 /
\plot 16.036 17.907 16.097 17.907 /
\plot 16.097 17.907 16.150 17.907 /
\plot 16.150 17.907 16.205 17.907 /
\plot 16.205 17.907 16.260 17.907 /
\plot 16.260 17.907 16.317 17.907 /
\plot 16.317 17.907 16.379 17.907 /
\plot 16.379 17.907 16.442 17.907 /
\plot 16.442 17.907 16.508 17.907 /
\plot 16.508 17.907 16.578 17.907 /
\plot 16.578 17.907 16.652 17.907 /
\plot 16.652 17.907 16.732 17.907 /
\plot 16.732 17.907 16.815 17.907 /
\plot 16.815 17.907 16.904 17.907 /
\plot 16.904 17.907 16.997 17.907 /
\plot 16.997 17.907 17.094 17.907 /
\plot 17.094 17.907 17.192 17.907 /
\plot 17.192 17.907 17.291 17.907 /
\plot 17.291 17.907 17.391 17.907 /
\plot 17.391 17.907 17.486 17.907 /
\plot 17.486 17.907 17.575 17.907 /
\plot 17.575 17.907 17.657 17.907 /
\plot 17.657 17.907 17.729 17.907 /
\plot 17.729 17.907 17.788 17.907 /
\plot 17.788 17.907 17.835 17.907 /
\plot 17.835 17.907 17.867 17.907 /
\plot 17.867 17.907 17.890 17.907 /
\plot 17.890 17.907 17.901 17.907 /
\plot 17.901 17.907 17.907 17.907 /
}%
%
%
\linethickness= 0.500pt
\setplotsymbol ({\thinlinefont .})
{\color[rgb]{0,0,0}\plot 18.479 17.907 20.003 17.907 /
}%
%
%
\linethickness= 0.500pt
\setplotsymbol ({\thinlinefont .})
{\color[rgb]{0,0,0}\plot 10.478 16.192 10.478 16.199 /
\plot 10.478 16.199 10.478 16.214 /
\plot 10.478 16.214 10.478 16.239 /
\plot 10.478 16.239 10.478 16.275 /
\plot 10.478 16.275 10.478 16.326 /
\plot 10.478 16.326 10.478 16.385 /
\plot 10.478 16.385 10.478 16.453 /
\plot 10.478 16.453 10.478 16.525 /
\plot 10.478 16.525 10.478 16.601 /
\plot 10.478 16.601 10.478 16.675 /
\plot 10.478 16.675 10.478 16.745 /
\plot 10.478 16.745 10.478 16.813 /
\plot 10.478 16.813 10.478 16.876 /
\plot 10.478 16.876 10.478 16.935 /
\plot 10.478 16.935 10.478 16.988 /
\plot 10.478 16.988 10.478 17.039 /
\plot 10.478 17.039 10.478 17.084 /
\plot 10.478 17.084 10.478 17.126 /
\plot 10.478 17.126 10.478 17.166 /
\plot 10.478 17.166 10.478 17.204 /
\plot 10.478 17.204 10.478 17.240 /
\plot 10.478 17.240 10.478 17.289 /
\plot 10.478 17.289 10.478 17.335 /
\plot 10.478 17.335 10.478 17.382 /
\plot 10.478 17.382 10.480 17.427 /
\plot 10.480 17.427 10.480 17.471 /
\plot 10.480 17.471 10.484 17.513 /
\plot 10.484 17.513 10.486 17.556 /
\plot 10.486 17.556 10.490 17.596 /
\plot 10.490 17.596 10.497 17.632 /
\plot 10.497 17.632 10.503 17.666 /
\plot 10.503 17.666 10.511 17.697 /
\plot 10.511 17.697 10.520 17.727 /
\plot 10.520 17.727 10.530 17.752 /
\plot 10.530 17.752 10.543 17.774 /
\plot 10.543 17.774 10.558 17.795 /
\plot 10.558 17.795 10.573 17.812 /
\plot 10.573 17.812 10.592 17.829 /
\plot 10.592 17.829 10.615 17.843 /
\plot 10.615 17.843 10.643 17.858 /
\plot 10.643 17.858 10.672 17.869 /
\plot 10.672 17.869 10.704 17.877 /
\plot 10.704 17.877 10.740 17.884 /
\plot 10.740 17.884 10.776 17.890 /
\plot 10.776 17.890 10.816 17.892 /
\plot 10.816 17.892 10.854 17.892 /
\plot 10.854 17.892 10.894 17.892 /
\plot 10.894 17.892 10.935 17.890 /
\plot 10.935 17.890 10.973 17.886 /
\plot 10.973 17.886 11.011 17.882 /
\plot 11.011 17.882 11.049 17.875 /
\plot 11.049 17.875 11.083 17.869 /
\plot 11.083 17.869 11.115 17.863 /
\plot 11.115 17.863 11.148 17.856 /
\plot 11.148 17.856 11.184 17.848 /
\plot 11.184 17.848 11.220 17.839 /
\plot 11.220 17.839 11.256 17.829 /
\plot 11.256 17.829 11.292 17.818 /
\plot 11.292 17.818 11.331 17.808 /
\plot 11.331 17.808 11.367 17.797 /
\plot 11.367 17.797 11.402 17.784 /
\plot 11.402 17.784 11.436 17.774 /
\plot 11.436 17.774 11.470 17.761 /
\plot 11.470 17.761 11.502 17.750 /
\plot 11.502 17.750 11.532 17.738 /
\plot 11.532 17.738 11.561 17.727 /
\plot 11.561 17.727 11.589 17.716 /
\plot 11.589 17.716 11.620 17.704 /
\plot 11.620 17.704 11.650 17.691 /
\plot 11.650 17.691 11.680 17.678 /
\plot 11.680 17.678 11.712 17.664 /
\plot 11.712 17.664 11.741 17.649 /
\plot 11.741 17.649 11.773 17.634 /
\plot 11.773 17.634 11.803 17.617 /
\plot 11.803 17.617 11.830 17.600 /
\plot 11.830 17.600 11.858 17.583 /
\plot 11.858 17.583 11.883 17.566 /
\plot 11.883 17.566 11.908 17.547 /
\plot 11.908 17.547 11.930 17.530 /
\plot 11.930 17.530 11.951 17.513 /
\plot 11.951 17.513 11.970 17.494 /
\plot 11.970 17.494 11.991 17.471 /
\plot 11.991 17.471 12.012 17.448 /
\plot 12.012 17.448 12.031 17.420 /
\plot 12.031 17.420 12.052 17.393 /
\plot 12.052 17.393 12.069 17.363 /
\plot 12.069 17.363 12.088 17.331 /
\plot 12.088 17.331 12.105 17.300 /
\plot 12.105 17.300 12.120 17.266 /
\plot 12.120 17.266 12.133 17.234 /
\plot 12.133 17.234 12.143 17.204 /
\plot 12.143 17.204 12.152 17.175 /
\plot 12.152 17.175 12.160 17.145 /
\plot 12.160 17.145 12.167 17.115 /
\plot 12.167 17.115 12.173 17.086 /
\plot 12.173 17.086 12.177 17.056 /
\plot 12.177 17.056 12.181 17.024 /
\plot 12.181 17.024 12.186 16.990 /
\plot 12.186 16.990 12.188 16.959 /
\plot 12.188 16.959 12.190 16.927 /
\plot 12.190 16.927 12.190 16.897 /
\plot 12.190 16.897 12.192 16.870 /
\plot 12.192 16.870 12.192 16.842 /
\plot 12.192 16.842 12.192 16.819 /
\plot 12.192 16.819 12.192 16.796 /
\plot 12.192 16.796 12.192 16.770 /
\plot 12.192 16.770 12.192 16.745 /
\plot 12.192 16.745 12.192 16.720 /
\plot 12.192 16.720 12.192 16.692 /
\plot 12.192 16.692 12.192 16.662 /
\plot 12.192 16.662 12.192 16.631 /
\plot 12.192 16.631 12.192 16.603 /
\plot 12.192 16.603 12.192 16.584 /
\plot 12.192 16.584 12.192 16.576 /
\plot 12.192 16.576 12.192 16.573 /
}%
%
%
\linethickness= 0.500pt
\setplotsymbol ({\thinlinefont .})
{\color[rgb]{0,0,0}\plot 10.287 16.192 10.287 16.199 /
\plot 10.287 16.199 10.287 16.216 /
\plot 10.287 16.216 10.287 16.241 /
\plot 10.287 16.241 10.287 16.281 /
\plot 10.287 16.281 10.287 16.332 /
\plot 10.287 16.332 10.287 16.394 /
\plot 10.287 16.394 10.287 16.459 /
\plot 10.287 16.459 10.287 16.529 /
\plot 10.287 16.529 10.287 16.597 /
\plot 10.287 16.597 10.287 16.662 /
\plot 10.287 16.662 10.287 16.726 /
\plot 10.287 16.726 10.287 16.783 /
\plot 10.287 16.783 10.287 16.836 /
\plot 10.287 16.836 10.287 16.885 /
\plot 10.287 16.885 10.287 16.929 /
\plot 10.287 16.929 10.287 16.971 /
\plot 10.287 16.971 10.287 17.010 /
\plot 10.287 17.010 10.287 17.046 /
\plot 10.287 17.046 10.287 17.082 /
\plot 10.287 17.082 10.287 17.124 /
\plot 10.287 17.124 10.287 17.166 /
\plot 10.287 17.166 10.287 17.209 /
\plot 10.287 17.209 10.287 17.251 /
\plot 10.287 17.251 10.287 17.293 /
\plot 10.287 17.293 10.287 17.333 /
\plot 10.287 17.333 10.287 17.374 /
\plot 10.287 17.374 10.287 17.414 /
\plot 10.287 17.414 10.287 17.450 /
\plot 10.287 17.450 10.287 17.486 /
\plot 10.287 17.486 10.287 17.520 /
\plot 10.287 17.520 10.287 17.551 /
\plot 10.287 17.551 10.287 17.579 /
\plot 10.287 17.579 10.287 17.606 /
\plot 10.287 17.606 10.287 17.630 /
\plot 10.287 17.630 10.287 17.653 /
\plot 10.287 17.653 10.287 17.681 /
\plot 10.287 17.681 10.287 17.708 /
\plot 10.287 17.708 10.285 17.733 /
\plot 10.285 17.733 10.285 17.757 /
\plot 10.285 17.757 10.283 17.778 /
\plot 10.283 17.778 10.279 17.799 /
\plot 10.279 17.799 10.272 17.816 /
\plot 10.272 17.816 10.266 17.833 /
\plot 10.266 17.833 10.257 17.846 /
\plot 10.257 17.846 10.249 17.858 /
\plot 10.249 17.858 10.236 17.867 /
\plot 10.236 17.867 10.223 17.875 /
\plot 10.223 17.875 10.211 17.882 /
\plot 10.211 17.882 10.194 17.886 /
\plot 10.194 17.886 10.177 17.890 /
\plot 10.177 17.890 10.154 17.894 /
\plot 10.154 17.894 10.130 17.899 /
\plot 10.130 17.899 10.101 17.901 /
\plot 10.101 17.901 10.071 17.903 /
\plot 10.071 17.903 10.035 17.905 /
\plot 10.035 17.905  9.999 17.905 /
\plot  9.999 17.905  9.959 17.907 /
\plot  9.959 17.907  9.917 17.907 /
\plot  9.917 17.907  9.872 17.907 /
\plot  9.872 17.907  9.828 17.907 /
\plot  9.828 17.907  9.779 17.907 /
\plot  9.779 17.907  9.749 17.907 /
\plot  9.749 17.907  9.718 17.907 /
\plot  9.718 17.907  9.686 17.907 /
\plot  9.686 17.907  9.650 17.907 /
\plot  9.650 17.907  9.612 17.907 /
\plot  9.612 17.907  9.572 17.907 /
\plot  9.572 17.907  9.529 17.907 /
\plot  9.529 17.907  9.483 17.907 /
\plot  9.483 17.907  9.434 17.907 /
\plot  9.434 17.907  9.383 17.907 /
\plot  9.383 17.907  9.330 17.907 /
\plot  9.330 17.907  9.273 17.907 /
\plot  9.273 17.907  9.214 17.907 /
\plot  9.214 17.907  9.152 17.907 /
\plot  9.152 17.907  9.089 17.907 /
\plot  9.089 17.907  9.023 17.907 /
\plot  9.023 17.907  8.958 17.907 /
\plot  8.958 17.907  8.888 17.907 /
\plot  8.888 17.907  8.818 17.907 /
\plot  8.818 17.907  8.748 17.907 /
\plot  8.748 17.907  8.674 17.907 /
\plot  8.674 17.907  8.600 17.907 /
\plot  8.600 17.907  8.524 17.907 /
\plot  8.524 17.907  8.446 17.907 /
\plot  8.446 17.907  8.388 17.907 /
\plot  8.388 17.907  8.329 17.907 /
\plot  8.329 17.907  8.268 17.907 /
\plot  8.268 17.907  8.204 17.907 /
\plot  8.204 17.907  8.139 17.907 /
\plot  8.139 17.907  8.069 17.907 /
\plot  8.069 17.907  7.995 17.907 /
\plot  7.995 17.907  7.916 17.907 /
\plot  7.916 17.907  7.834 17.907 /
\plot  7.834 17.907  7.747 17.907 /
\plot  7.747 17.907  7.654 17.907 /
\plot  7.654 17.907  7.556 17.907 /
\plot  7.556 17.907  7.453 17.907 /
\plot  7.453 17.907  7.345 17.907 /
\plot  7.345 17.907  7.231 17.907 /
\plot  7.231 17.907  7.112 17.907 /
\plot  7.112 17.907  6.989 17.907 /
\plot  6.989 17.907  6.862 17.907 /
\plot  6.862 17.907  6.735 17.907 /
\plot  6.735 17.907  6.608 17.907 /
\plot  6.608 17.907  6.481 17.907 /
\plot  6.481 17.907  6.361 17.907 /
\plot  6.361 17.907  6.244 17.907 /
\plot  6.244 17.907  6.138 17.907 /
\plot  6.138 17.907  6.041 17.907 /
\plot  6.041 17.907  5.956 17.907 /
\plot  5.956 17.907  5.884 17.907 /
\plot  5.884 17.907  5.827 17.907 /
\plot  5.827 17.907  5.783 17.907 /
\plot  5.783 17.907  5.751 17.907 /
\plot  5.751 17.907  5.730 17.907 /
\plot  5.730 17.907  5.719 17.907 /
\plot  5.719 17.907  5.715 17.907 /
}%
%
%
\linethickness= 0.500pt
\setplotsymbol ({\thinlinefont .})
{\color[rgb]{0,0,0}\plot 12.573 16.573 12.573 16.576 /
\plot 12.573 16.576 12.573 16.582 /
\plot 12.573 16.582 12.573 16.599 /
\plot 12.573 16.599 12.573 16.629 /
\plot 12.573 16.629 12.571 16.669 /
\plot 12.571 16.669 12.571 16.715 /
\plot 12.571 16.715 12.569 16.768 /
\plot 12.569 16.768 12.567 16.823 /
\plot 12.567 16.823 12.567 16.874 /
\plot 12.567 16.874 12.562 16.923 /
\plot 12.562 16.923 12.560 16.967 /
\plot 12.560 16.967 12.558 17.007 /
\plot 12.558 17.007 12.554 17.043 /
\plot 12.554 17.043 12.552 17.079 /
\plot 12.552 17.079 12.545 17.111 /
\plot 12.545 17.111 12.541 17.145 /
\plot 12.541 17.145 12.535 17.179 /
\plot 12.535 17.179 12.529 17.211 /
\plot 12.529 17.211 12.522 17.244 /
\plot 12.522 17.244 12.514 17.280 /
\plot 12.514 17.280 12.503 17.316 /
\plot 12.503 17.316 12.493 17.352 /
\plot 12.493 17.352 12.482 17.388 /
\plot 12.482 17.388 12.469 17.427 /
\plot 12.469 17.427 12.457 17.462 /
\plot 12.457 17.462 12.442 17.498 /
\plot 12.442 17.498 12.427 17.532 /
\plot 12.427 17.532 12.412 17.566 /
\plot 12.412 17.566 12.397 17.598 /
\plot 12.397 17.598 12.383 17.628 /
\plot 12.383 17.628 12.368 17.657 /
\plot 12.368 17.657 12.351 17.685 /
\plot 12.351 17.685 12.332 17.716 /
\plot 12.332 17.716 12.311 17.746 /
\plot 12.311 17.746 12.289 17.776 /
\plot 12.289 17.776 12.264 17.808 /
\plot 12.264 17.808 12.239 17.837 /
\plot 12.239 17.837 12.211 17.869 /
\plot 12.211 17.869 12.184 17.899 /
\plot 12.184 17.899 12.154 17.926 /
\plot 12.154 17.926 12.122 17.954 /
\plot 12.122 17.954 12.093 17.979 /
\plot 12.093 17.979 12.061 18.004 /
\plot 12.061 18.004 12.031 18.026 /
\plot 12.031 18.026 12.002 18.047 /
\plot 12.002 18.047 11.970 18.066 /
\plot 11.970 18.066 11.942 18.083 /
\plot 11.942 18.083 11.913 18.098 /
\plot 11.913 18.098 11.883 18.112 /
\plot 11.883 18.112 11.851 18.127 /
\plot 11.851 18.127 11.817 18.142 /
\plot 11.817 18.142 11.781 18.157 /
\plot 11.781 18.157 11.745 18.172 /
\plot 11.745 18.172 11.707 18.184 /
\plot 11.707 18.184 11.669 18.197 /
\plot 11.669 18.197 11.629 18.208 /
\plot 11.629 18.208 11.591 18.218 /
\plot 11.591 18.218 11.551 18.229 /
\plot 11.551 18.229 11.513 18.237 /
\plot 11.513 18.237 11.474 18.244 /
\plot 11.474 18.244 11.436 18.250 /
\plot 11.436 18.250 11.398 18.256 /
\plot 11.398 18.256 11.367 18.260 /
\plot 11.367 18.260 11.335 18.265 /
\plot 11.335 18.265 11.303 18.267 /
\plot 11.303 18.267 11.269 18.271 /
\plot 11.269 18.271 11.231 18.273 /
\plot 11.231 18.273 11.193 18.275 /
\plot 11.193 18.275 11.155 18.277 /
\plot 11.155 18.277 11.113 18.280 /
\plot 11.113 18.280 11.068 18.282 /
\plot 11.068 18.282 11.024 18.284 /
\plot 11.024 18.284 10.975 18.284 /
\plot 10.975 18.284 10.928 18.286 /
\plot 10.928 18.286 10.878 18.286 /
\plot 10.878 18.286 10.827 18.288 /
\plot 10.827 18.288 10.776 18.288 /
\plot 10.776 18.288 10.725 18.288 /
\plot 10.725 18.288 10.672 18.288 /
\plot 10.672 18.288 10.619 18.288 /
\plot 10.619 18.288 10.564 18.288 /
\plot 10.564 18.288 10.509 18.288 /
\plot 10.509 18.288 10.469 18.288 /
\plot 10.469 18.288 10.425 18.288 /
\plot 10.425 18.288 10.382 18.288 /
\plot 10.382 18.288 10.336 18.288 /
\plot 10.336 18.288 10.287 18.288 /
\plot 10.287 18.288 10.236 18.288 /
\plot 10.236 18.288 10.185 18.288 /
\plot 10.185 18.288 10.130 18.288 /
\plot 10.130 18.288 10.073 18.288 /
\plot 10.073 18.288 10.014 18.288 /
\plot 10.014 18.288  9.950 18.288 /
\plot  9.950 18.288  9.887 18.288 /
\plot  9.887 18.288  9.821 18.288 /
\plot  9.821 18.288  9.754 18.288 /
\plot  9.754 18.288  9.684 18.288 /
\plot  9.684 18.288  9.614 18.288 /
\plot  9.614 18.288  9.540 18.288 /
\plot  9.540 18.288  9.468 18.288 /
\plot  9.468 18.288  9.394 18.288 /
\plot  9.394 18.288  9.318 18.288 /
\plot  9.318 18.288  9.243 18.288 /
\plot  9.243 18.288  9.167 18.288 /
\plot  9.167 18.288  9.091 18.288 /
\plot  9.091 18.288  9.015 18.288 /
\plot  9.015 18.288  8.937 18.288 /
\plot  8.937 18.288  8.858 18.288 /
\plot  8.858 18.288  8.780 18.288 /
\plot  8.780 18.288  8.700 18.288 /
\plot  8.700 18.288  8.636 18.288 /
\plot  8.636 18.288  8.570 18.288 /
\plot  8.570 18.288  8.505 18.288 /
\plot  8.505 18.288  8.435 18.288 /
\plot  8.435 18.288  8.365 18.288 /
\plot  8.365 18.288  8.291 18.288 /
\plot  8.291 18.288  8.215 18.288 /
\plot  8.215 18.288  8.134 18.288 /
\plot  8.134 18.288  8.050 18.288 /
\plot  8.050 18.288  7.961 18.288 /
\plot  7.961 18.288  7.870 18.288 /
\plot  7.870 18.288  7.770 18.288 /
\plot  7.770 18.288  7.669 18.288 /
\plot  7.669 18.288  7.561 18.288 /
\plot  7.561 18.288  7.449 18.288 /
\plot  7.449 18.288  7.330 18.288 /
\plot  7.330 18.288  7.207 18.288 /
\plot  7.207 18.288  7.082 18.288 /
\plot  7.082 18.288  6.955 18.288 /
\plot  6.955 18.288  6.824 18.288 /
\plot  6.824 18.288  6.695 18.288 /
\plot  6.695 18.288  6.566 18.288 /
\plot  6.566 18.288  6.441 18.288 /
\plot  6.441 18.288  6.322 18.288 /
\plot  6.322 18.288  6.210 18.288 /
\plot  6.210 18.288  6.109 18.288 /
\plot  6.109 18.288  6.016 18.288 /
\plot  6.016 18.288  5.937 18.288 /
\plot  5.937 18.288  5.870 18.288 /
\plot  5.870 18.288  5.817 18.288 /
\plot  5.817 18.288  5.776 18.288 /
\plot  5.776 18.288  5.747 18.288 /
\plot  5.747 18.288  5.730 18.288 /
\plot  5.730 18.288  5.719 18.288 /
\plot  5.719 18.288  5.715 18.288 /
}%
\linethickness=0pt
\putrectangle corners at  3.785 23.838 and 21.933 12.357
\endpicture}

%% file: OBDAfterH.tex
\font\thinlinefont=cmr5
\mbox{\beginpicture
\small
\setcoordinatesystem units <0.5cm,0.5cm>
\unitlength=1.04987cm
\linethickness=1pt
\setplotsymbol ({\makebox(0,0)[l]{\tencirc\symbol{'160}}})
\setshadesymbol ({\thinlinefont .})
\setlinear
%
%
\linethickness= 0.500pt
\setplotsymbol ({\thinlinefont .})
{\color[rgb]{0,0,0}\circulararc 93.580 degrees from  7.429 19.431 center at  8.547 19.673
}%
%
%
\linethickness= 0.500pt
\setplotsymbol ({\thinlinefont .})
{\color[rgb]{0,0,0}\circulararc 105.600 degrees from  9.716 18.574 center at 10.152 19.581
}%
%
%
\linethickness= 0.500pt
\setplotsymbol ({\thinlinefont .})
\setdashes < 0.1270cm>
{\color[rgb]{0,0,0}\circulararc 126.870 degrees from  8.858 18.574 center at  7.953 18.598
}%
%
%
\linethickness= 0.500pt
\setplotsymbol ({\thinlinefont .})
{\color[rgb]{0,0,0}\circulararc 121.891 degrees from 11.239 19.431 center at 10.716 18.579
}%
%
%
\linethickness= 0.500pt
\setplotsymbol ({\thinlinefont .})
\setsolid
{\color[rgb]{0,0,0}\circulararc 126.870 degrees from  8.954 13.335 center at  9.573 13.645
}%
%
%
\linethickness= 0.500pt
\setplotsymbol ({\thinlinefont .})
{\color[rgb]{0,0,0}\circulararc 115.058 degrees from 10.097 13.144 center at  9.573 12.811
}%
%
%
\linethickness=1pt
\setplotsymbol ({\makebox(0,0)[l]{\tencirc\symbol{'160}}})
\setdashes < 0.1cm>
{\color[rgb]{0,0,0}\circulararc 65.758 degrees from  9.239 18.860 center at  8.793 19.462
}%
%
%
\linethickness= 0.500pt
\setplotsymbol ({\thinlinefont .})
\setsolid
{\color[rgb]{0,0,0}\ellipticalarc axes ratio  1.001:0.334  360 degrees 
	from 10.287 19.956 center at  9.286 19.956
}%
%
%
\linethickness= 0.500pt
\setplotsymbol ({\thinlinefont .})
{\color[rgb]{0,0,0}\ellipticalarc axes ratio  0.425:0.425  360 degrees 
	from  9.665 18.479 center at  9.239 18.479
}%
%
%
\linethickness=1pt
\setplotsymbol ({\makebox(0,0)[l]{\tencirc\symbol{'160}}})
{\color[rgb]{0,0,0}\plot  9.335 16.859  9.335 16.859 /
%
%
\linethickness=1pt
\setplotsymbol ({\makebox(0,0)[l]{\tencirc\symbol{'160}}})
\color[rgb]{0,0,0}\plot  9.335 16.097  9.335 16.097 /
%
%
\linethickness=1pt
\setplotsymbol ({\makebox(0,0)[l]{\tencirc\symbol{'160}}})
\color[rgb]{0,0,0}\plot  9.335 15.335  9.335 15.335 /
%
%
\linethickness=1pt
\setplotsymbol ({\makebox(0,0)[l]{\tencirc\symbol{'160}}})
\color[rgb]{0,0,0}\plot  9.335 14.669  9.335 14.669 /
%
%
\linethickness= 0.500pt
\setplotsymbol ({\thinlinefont .})
\color[rgb]{0,0,0}\putrule from  8.382 20.098 to  8.376 20.098
\plot  8.376 20.098  8.361 20.096 /
\plot  8.361 20.096  8.335 20.094 /
\plot  8.335 20.094  8.302 20.091 /
\plot  8.302 20.091  8.259 20.087 /
\plot  8.259 20.087  8.213 20.081 /
\plot  8.213 20.081  8.166 20.072 /
\plot  8.166 20.072  8.120 20.064 /
\plot  8.120 20.064  8.077 20.055 /
\plot  8.077 20.055  8.037 20.043 /
\plot  8.037 20.043  7.999 20.030 /
\plot  7.999 20.030  7.963 20.015 /
\plot  7.963 20.015  7.927 19.998 /
\plot  7.927 19.998  7.893 19.977 /
\plot  7.893 19.977  7.859 19.954 /
\plot  7.859 19.954  7.830 19.935 /
\plot  7.830 19.935  7.800 19.911 /
\plot  7.800 19.911  7.770 19.886 /
\plot  7.770 19.886  7.739 19.859 /
\plot  7.739 19.859  7.705 19.829 /
\plot  7.705 19.829  7.673 19.797 /
\plot  7.673 19.797  7.639 19.761 /
\plot  7.639 19.761  7.605 19.725 /
\plot  7.605 19.725  7.571 19.685 /
\plot  7.571 19.685  7.537 19.645 /
\plot  7.537 19.645  7.506 19.602 /
\plot  7.506 19.602  7.474 19.560 /
\plot  7.474 19.560  7.444 19.516 /
\plot  7.444 19.516  7.415 19.471 /
\plot  7.415 19.471  7.389 19.427 /
\plot  7.389 19.427  7.364 19.382 /
\plot  7.364 19.382  7.341 19.340 /
\plot  7.341 19.340  7.322 19.296 /
\plot  7.322 19.296  7.303 19.253 /
\plot  7.303 19.253  7.288 19.209 /
\plot  7.288 19.209  7.273 19.169 /
\plot  7.273 19.169  7.262 19.128 /
\plot  7.262 19.128  7.252 19.086 /
\plot  7.252 19.086  7.243 19.044 /
\plot  7.243 19.044  7.237 19.001 /
\plot  7.237 19.001  7.233 18.957 /
\plot  7.233 18.957  7.228 18.910 /
\plot  7.228 18.910  7.226 18.866 /
\plot  7.226 18.866  7.228 18.819 /
\plot  7.228 18.819  7.231 18.773 /
\plot  7.231 18.773  7.235 18.726 /
\plot  7.235 18.726  7.243 18.680 /
\plot  7.243 18.680  7.252 18.633 /
\plot  7.252 18.633  7.262 18.589 /
\plot  7.262 18.589  7.275 18.546 /
\plot  7.275 18.546  7.290 18.504 /
\plot  7.290 18.504  7.307 18.464 /
\plot  7.307 18.464  7.324 18.426 /
\plot  7.324 18.426  7.345 18.390 /
\plot  7.345 18.390  7.366 18.354 /
\plot  7.366 18.354  7.389 18.320 /
\plot  7.389 18.320  7.415 18.288 /
\plot  7.415 18.288  7.442 18.254 /
\plot  7.442 18.254  7.476 18.220 /
\plot  7.476 18.220  7.510 18.186 /
\plot  7.510 18.186  7.548 18.155 /
\plot  7.548 18.155  7.586 18.123 /
\plot  7.586 18.123  7.628 18.089 /
\plot  7.628 18.089  7.673 18.057 /
\plot  7.673 18.057  7.717 18.026 /
\plot  7.717 18.026  7.766 17.994 /
\plot  7.766 17.994  7.813 17.962 /
\plot  7.813 17.962  7.861 17.930 /
\plot  7.861 17.930  7.908 17.899 /
\plot  7.908 17.899  7.954 17.869 /
\plot  7.954 17.869  7.999 17.839 /
\plot  7.999 17.839  8.041 17.810 /
\plot  8.041 17.810  8.084 17.780 /
\plot  8.084 17.780  8.122 17.752 /
\plot  8.122 17.752  8.158 17.725 /
\plot  8.158 17.725  8.191 17.697 /
\plot  8.191 17.697  8.223 17.668 /
\plot  8.223 17.668  8.255 17.636 /
\plot  8.255 17.636  8.287 17.604 /
\plot  8.287 17.604  8.314 17.568 /
\plot  8.314 17.568  8.340 17.530 /
\plot  8.340 17.530  8.365 17.488 /
\plot  8.365 17.488  8.390 17.443 /
\plot  8.390 17.443  8.414 17.393 /
\plot  8.414 17.393  8.439 17.338 /
\plot  8.439 17.338  8.462 17.278 /
\plot  8.462 17.278  8.486 17.217 /
\plot  8.486 17.217  8.507 17.156 /
\plot  8.507 17.156  8.526 17.098 /
\plot  8.526 17.098  8.543 17.048 /
\plot  8.543 17.048  8.556 17.007 /
\plot  8.556 17.007  8.566 16.978 /
\plot  8.566 16.978  8.570 16.963 /
\plot  8.570 16.963  8.572 16.954 /
}%
%
%
\linethickness= 0.500pt
\setplotsymbol ({\thinlinefont .})
{\color[rgb]{0,0,0}\putrule from 10.287 20.098 to 10.293 20.098
\plot 10.293 20.098 10.308 20.096 /
\plot 10.308 20.096 10.334 20.094 /
\plot 10.334 20.094 10.367 20.091 /
\plot 10.367 20.091 10.410 20.087 /
\plot 10.410 20.087 10.456 20.081 /
\plot 10.456 20.081 10.503 20.072 /
\plot 10.503 20.072 10.549 20.064 /
\plot 10.549 20.064 10.592 20.055 /
\plot 10.592 20.055 10.632 20.043 /
\plot 10.632 20.043 10.670 20.030 /
\plot 10.670 20.030 10.706 20.015 /
\plot 10.706 20.015 10.742 19.998 /
\plot 10.742 19.998 10.776 19.977 /
\plot 10.776 19.977 10.812 19.954 /
\plot 10.812 19.954 10.839 19.935 /
\plot 10.839 19.935 10.869 19.911 /
\plot 10.869 19.911 10.899 19.886 /
\plot 10.899 19.886 10.930 19.859 /
\plot 10.930 19.859 10.964 19.829 /
\plot 10.964 19.829 10.996 19.797 /
\plot 10.996 19.797 11.030 19.761 /
\plot 11.030 19.761 11.064 19.725 /
\plot 11.064 19.725 11.098 19.685 /
\plot 11.098 19.685 11.132 19.645 /
\plot 11.132 19.645 11.163 19.602 /
\plot 11.163 19.602 11.195 19.560 /
\plot 11.195 19.560 11.225 19.516 /
\plot 11.225 19.516 11.254 19.471 /
\plot 11.254 19.471 11.280 19.427 /
\plot 11.280 19.427 11.305 19.382 /
\plot 11.305 19.382 11.328 19.340 /
\plot 11.328 19.340 11.347 19.296 /
\plot 11.347 19.296 11.367 19.253 /
\plot 11.367 19.253 11.383 19.209 /
\plot 11.383 19.209 11.396 19.169 /
\plot 11.396 19.169 11.407 19.128 /
\plot 11.407 19.128 11.417 19.086 /
\plot 11.417 19.086 11.426 19.044 /
\plot 11.426 19.044 11.432 19.001 /
\plot 11.432 19.001 11.436 18.957 /
\plot 11.436 18.957 11.441 18.910 /
\plot 11.441 18.910 11.443 18.866 /
\plot 11.443 18.866 11.441 18.819 /
\plot 11.441 18.819 11.438 18.773 /
\plot 11.438 18.773 11.434 18.726 /
\plot 11.434 18.726 11.426 18.680 /
\plot 11.426 18.680 11.417 18.633 /
\plot 11.417 18.633 11.407 18.589 /
\plot 11.407 18.589 11.394 18.546 /
\plot 11.394 18.546 11.379 18.504 /
\plot 11.379 18.504 11.362 18.464 /
\plot 11.362 18.464 11.345 18.426 /
\plot 11.345 18.426 11.324 18.390 /
\plot 11.324 18.390 11.303 18.354 /
\plot 11.303 18.354 11.280 18.320 /
\plot 11.280 18.320 11.256 18.288 /
\plot 11.256 18.288 11.227 18.254 /
\plot 11.227 18.254 11.193 18.220 /
\plot 11.193 18.220 11.159 18.186 /
\plot 11.159 18.186 11.121 18.155 /
\plot 11.121 18.155 11.083 18.123 /
\plot 11.083 18.123 11.041 18.089 /
\plot 11.041 18.089 10.996 18.057 /
\plot 10.996 18.057 10.952 18.026 /
\plot 10.952 18.026 10.903 17.994 /
\plot 10.903 17.994 10.856 17.962 /
\plot 10.856 17.962 10.808 17.930 /
\plot 10.808 17.930 10.761 17.899 /
\plot 10.761 17.899 10.715 17.869 /
\plot 10.715 17.869 10.670 17.839 /
\plot 10.670 17.839 10.628 17.810 /
\plot 10.628 17.810 10.585 17.780 /
\plot 10.585 17.780 10.547 17.752 /
\plot 10.547 17.752 10.511 17.725 /
\plot 10.511 17.725 10.478 17.697 /
\plot 10.478 17.697 10.446 17.668 /
\plot 10.446 17.668 10.414 17.636 /
\plot 10.414 17.636 10.382 17.604 /
\plot 10.382 17.604 10.355 17.568 /
\plot 10.355 17.568 10.329 17.530 /
\plot 10.329 17.530 10.304 17.488 /
\plot 10.304 17.488 10.279 17.443 /
\plot 10.279 17.443 10.255 17.393 /
\plot 10.255 17.393 10.230 17.338 /
\plot 10.230 17.338 10.207 17.278 /
\plot 10.207 17.278 10.183 17.217 /
\plot 10.183 17.217 10.162 17.156 /
\plot 10.162 17.156 10.143 17.098 /
\plot 10.143 17.098 10.126 17.048 /
\plot 10.126 17.048 10.113 17.007 /
\plot 10.113 17.007 10.103 16.978 /
\plot 10.103 16.978 10.099 16.963 /
\plot 10.099 16.963 10.097 16.954 /
}%
%
%
\linethickness= 0.500pt
\setplotsymbol ({\thinlinefont .})
{\color[rgb]{0,0,0}\putrule from  8.763 14.764 to  8.763 14.762
\putrule from  8.763 14.762 to  8.763 14.755
\plot  8.763 14.755  8.761 14.734 /
\plot  8.761 14.734  8.759 14.702 /
\plot  8.759 14.702  8.757 14.656 /
\plot  8.757 14.656  8.752 14.601 /
\plot  8.752 14.601  8.746 14.541 /
\plot  8.746 14.541  8.738 14.482 /
\plot  8.738 14.482  8.729 14.425 /
\plot  8.729 14.425  8.721 14.372 /
\plot  8.721 14.372  8.708 14.323 /
\plot  8.708 14.323  8.695 14.281 /
\plot  8.695 14.281  8.680 14.243 /
\plot  8.680 14.243  8.664 14.207 /
\plot  8.664 14.207  8.642 14.175 /
\plot  8.642 14.175  8.621 14.144 /
\plot  8.621 14.144  8.598 14.118 /
\plot  8.598 14.118  8.572 14.093 /
\plot  8.572 14.093  8.543 14.065 /
\plot  8.543 14.065  8.513 14.040 /
\plot  8.513 14.040  8.481 14.012 /
\plot  8.481 14.012  8.448 13.983 /
\plot  8.448 13.983  8.412 13.953 /
\plot  8.412 13.953  8.374 13.923 /
\plot  8.374 13.923  8.338 13.892 /
\plot  8.338 13.892  8.299 13.858 /
\plot  8.299 13.858  8.263 13.824 /
\plot  8.263 13.824  8.227 13.790 /
\plot  8.227 13.790  8.194 13.754 /
\plot  8.194 13.754  8.164 13.716 /
\plot  8.164 13.716  8.134 13.680 /
\plot  8.134 13.680  8.109 13.640 /
\plot  8.109 13.640  8.086 13.600 /
\plot  8.086 13.600  8.065 13.557 /
\plot  8.065 13.557  8.050 13.523 /
\plot  8.050 13.523  8.039 13.490 /
\plot  8.039 13.490  8.026 13.451 /
\plot  8.026 13.451  8.018 13.413 /
\plot  8.018 13.413  8.009 13.371 /
\plot  8.009 13.371  8.005 13.329 /
\plot  8.005 13.329  8.001 13.284 /
\plot  8.001 13.284  7.999 13.240 /
\putrule from  7.999 13.240 to  7.999 13.191
\plot  7.999 13.191  8.001 13.142 /
\plot  8.001 13.142  8.003 13.094 /
\plot  8.003 13.094  8.012 13.043 /
\plot  8.012 13.043  8.020 12.992 /
\plot  8.020 12.992  8.031 12.941 /
\plot  8.031 12.941  8.043 12.893 /
\plot  8.043 12.893  8.058 12.842 /
\plot  8.058 12.842  8.077 12.793 /
\plot  8.077 12.793  8.096 12.747 /
\plot  8.096 12.747  8.117 12.700 /
\plot  8.117 12.700  8.141 12.658 /
\plot  8.141 12.658  8.166 12.613 /
\plot  8.166 12.613  8.194 12.573 /
\plot  8.194 12.573  8.223 12.533 /
\plot  8.223 12.533  8.255 12.493 /
\plot  8.255 12.493  8.289 12.457 /
\plot  8.289 12.457  8.325 12.418 /
\plot  8.325 12.418  8.365 12.383 /
\plot  8.365 12.383  8.407 12.347 /
\plot  8.407 12.347  8.452 12.313 /
\plot  8.452 12.313  8.501 12.279 /
\plot  8.501 12.279  8.551 12.247 /
\plot  8.551 12.247  8.604 12.213 /
\plot  8.604 12.213  8.661 12.184 /
\plot  8.661 12.184  8.719 12.154 /
\plot  8.719 12.154  8.780 12.126 /
\plot  8.780 12.126  8.841 12.101 /
\plot  8.841 12.101  8.905 12.078 /
\plot  8.905 12.078  8.968 12.057 /
\plot  8.968 12.057  9.034 12.040 /
\plot  9.034 12.040  9.097 12.023 /
\plot  9.097 12.023  9.161 12.010 /
\plot  9.161 12.010  9.227 11.997 /
\plot  9.227 11.997  9.288 11.989 /
\plot  9.288 11.989  9.349 11.985 /
\plot  9.349 11.985  9.411 11.980 /
\putrule from  9.411 11.980 to  9.472 11.980
\plot  9.472 11.980  9.529 11.982 /
\plot  9.529 11.982  9.588 11.985 /
\plot  9.588 11.985  9.646 11.993 /
\plot  9.646 11.993  9.705 12.002 /
\plot  9.705 12.002  9.764 12.014 /
\plot  9.764 12.014  9.821 12.027 /
\plot  9.821 12.027  9.881 12.044 /
\plot  9.881 12.044  9.942 12.065 /
\plot  9.942 12.065 10.001 12.086 /
\plot 10.001 12.086 10.061 12.112 /
\plot 10.061 12.112 10.122 12.139 /
\plot 10.122 12.139 10.181 12.169 /
\plot 10.181 12.169 10.240 12.200 /
\plot 10.240 12.200 10.298 12.234 /
\plot 10.298 12.234 10.353 12.268 /
\plot 10.353 12.268 10.408 12.306 /
\plot 10.408 12.306 10.461 12.344 /
\plot 10.461 12.344 10.509 12.383 /
\plot 10.509 12.383 10.558 12.421 /
\plot 10.558 12.421 10.602 12.461 /
\plot 10.602 12.461 10.643 12.501 /
\plot 10.643 12.501 10.683 12.541 /
\plot 10.683 12.541 10.719 12.581 /
\plot 10.719 12.581 10.753 12.620 /
\plot 10.753 12.620 10.782 12.660 /
\plot 10.782 12.660 10.812 12.700 /
\plot 10.812 12.700 10.839 12.744 /
\plot 10.839 12.744 10.865 12.789 /
\plot 10.865 12.789 10.890 12.833 /
\plot 10.890 12.833 10.911 12.880 /
\plot 10.911 12.880 10.930 12.929 /
\plot 10.930 12.929 10.947 12.977 /
\plot 10.947 12.977 10.962 13.026 /
\plot 10.962 13.026 10.973 13.075 /
\plot 10.973 13.075 10.983 13.125 /
\plot 10.983 13.125 10.990 13.176 /
\plot 10.990 13.176 10.994 13.227 /
\plot 10.994 13.227 10.996 13.276 /
\plot 10.996 13.276 10.994 13.324 /
\plot 10.994 13.324 10.992 13.371 /
\plot 10.992 13.371 10.986 13.418 /
\plot 10.986 13.418 10.977 13.462 /
\plot 10.977 13.462 10.969 13.504 /
\plot 10.969 13.504 10.956 13.545 /
\plot 10.956 13.545 10.943 13.583 /
\plot 10.943 13.583 10.926 13.619 /
\plot 10.926 13.619 10.909 13.652 /
\plot 10.909 13.652 10.890 13.684 /
\plot 10.890 13.684 10.865 13.722 /
\plot 10.865 13.722 10.837 13.758 /
\plot 10.837 13.758 10.806 13.794 /
\plot 10.806 13.794 10.772 13.828 /
\plot 10.772 13.828 10.736 13.862 /
\plot 10.736 13.862 10.698 13.894 /
\plot 10.698 13.894 10.657 13.926 /
\plot 10.657 13.926 10.615 13.957 /
\plot 10.615 13.957 10.573 13.989 /
\plot 10.573 13.989 10.530 14.019 /
\plot 10.530 14.019 10.488 14.046 /
\plot 10.488 14.046 10.448 14.074 /
\plot 10.448 14.074 10.410 14.101 /
\plot 10.410 14.101 10.374 14.127 /
\plot 10.374 14.127 10.340 14.152 /
\plot 10.340 14.152 10.308 14.175 /
\plot 10.308 14.175 10.281 14.201 /
\plot 10.281 14.201 10.255 14.224 /
\plot 10.255 14.224 10.230 14.252 /
\plot 10.230 14.252 10.207 14.281 /
\plot 10.207 14.281 10.188 14.311 /
\plot 10.188 14.311 10.171 14.345 /
\plot 10.171 14.345 10.158 14.383 /
\plot 10.158 14.383 10.145 14.425 /
\plot 10.145 14.425 10.132 14.470 /
\plot 10.132 14.470 10.124 14.520 /
\plot 10.124 14.520 10.116 14.571 /
\plot 10.116 14.571 10.109 14.624 /
\plot 10.109 14.624 10.105 14.671 /
\plot 10.105 14.671 10.101 14.711 /
\plot 10.101 14.711 10.099 14.738 /
\plot 10.099 14.738 10.097 14.755 /
\putrule from 10.097 14.755 to 10.097 14.762
\putrule from 10.097 14.762 to 10.097 14.764
}%
%
%
\linethickness= 0.500pt
\setplotsymbol ({\thinlinefont .})
{\color[rgb]{0,0,0}\plot 11.144 17.240 11.153 17.236 /
\plot 11.153 17.236 11.168 17.225 /
\plot 11.168 17.225 11.193 17.211 /
\plot 11.193 17.211 11.225 17.189 /
\plot 11.225 17.189 11.263 17.164 /
\plot 11.263 17.164 11.299 17.139 /
\plot 11.299 17.139 11.335 17.111 /
\plot 11.335 17.111 11.364 17.084 /
\plot 11.364 17.084 11.392 17.056 /
\plot 11.392 17.056 11.413 17.031 /
\plot 11.413 17.031 11.432 17.001 /
\plot 11.432 17.001 11.449 16.971 /
\plot 11.449 16.971 11.462 16.938 /
\plot 11.462 16.938 11.472 16.908 /
\plot 11.472 16.908 11.481 16.874 /
\plot 11.481 16.874 11.489 16.836 /
\plot 11.489 16.836 11.496 16.798 /
\plot 11.496 16.798 11.502 16.753 /
\plot 11.502 16.753 11.508 16.709 /
\plot 11.508 16.709 11.513 16.662 /
\plot 11.513 16.662 11.517 16.614 /
\plot 11.517 16.614 11.519 16.565 /
\plot 11.519 16.565 11.521 16.516 /
\plot 11.521 16.516 11.523 16.470 /
\putrule from 11.523 16.470 to 11.523 16.423
\plot 11.523 16.423 11.525 16.379 /
\putrule from 11.525 16.379 to 11.525 16.336
\putrule from 11.525 16.336 to 11.525 16.294
\putrule from 11.525 16.294 to 11.525 16.256
\putrule from 11.525 16.256 to 11.525 16.214
\putrule from 11.525 16.214 to 11.525 16.171
\putrule from 11.525 16.171 to 11.525 16.129
\putrule from 11.525 16.129 to 11.525 16.089
\plot 11.525 16.089 11.527 16.049 /
\putrule from 11.527 16.049 to 11.527 16.010
\plot 11.527 16.010 11.529 15.972 /
\plot 11.529 15.972 11.532 15.936 /
\plot 11.532 15.936 11.534 15.903 /
\plot 11.534 15.903 11.538 15.873 /
\plot 11.538 15.873 11.542 15.845 /
\plot 11.542 15.845 11.546 15.822 /
\plot 11.546 15.822 11.551 15.799 /
\plot 11.551 15.799 11.557 15.780 /
\plot 11.557 15.780 11.565 15.759 /
\plot 11.565 15.759 11.574 15.740 /
\plot 11.574 15.740 11.585 15.723 /
\plot 11.585 15.723 11.597 15.706 /
\plot 11.597 15.706 11.612 15.693 /
\plot 11.612 15.693 11.627 15.680 /
\plot 11.627 15.680 11.644 15.670 /
\plot 11.644 15.670 11.661 15.659 /
\plot 11.661 15.659 11.678 15.653 /
\plot 11.678 15.653 11.697 15.646 /
\plot 11.697 15.646 11.714 15.640 /
\plot 11.714 15.640 11.733 15.636 /
\plot 11.733 15.636 11.754 15.634 /
\plot 11.754 15.634 11.779 15.629 /
\plot 11.779 15.629 11.807 15.627 /
\plot 11.807 15.627 11.841 15.625 /
\plot 11.841 15.625 11.879 15.623 /
\putrule from 11.879 15.623 to 11.921 15.623
\plot 11.921 15.623 11.959 15.621 /
\putrule from 11.959 15.621 to 11.987 15.621
\putrule from 11.987 15.621 to 11.999 15.621
\putrule from 11.999 15.621 to 12.002 15.621
}%
%
%
\linethickness= 0.500pt
\setplotsymbol ({\thinlinefont .})
{\color[rgb]{0,0,0}\plot 11.144 14.002 11.153 14.006 /
\plot 11.153 14.006 11.168 14.017 /
\plot 11.168 14.017 11.193 14.031 /
\plot 11.193 14.031 11.225 14.053 /
\plot 11.225 14.053 11.263 14.078 /
\plot 11.263 14.078 11.299 14.103 /
\plot 11.299 14.103 11.335 14.131 /
\plot 11.335 14.131 11.364 14.158 /
\plot 11.364 14.158 11.392 14.186 /
\plot 11.392 14.186 11.413 14.211 /
\plot 11.413 14.211 11.432 14.241 /
\plot 11.432 14.241 11.449 14.271 /
\plot 11.449 14.271 11.462 14.302 /
\plot 11.462 14.302 11.470 14.332 /
\plot 11.470 14.332 11.479 14.362 /
\plot 11.479 14.362 11.487 14.393 /
\plot 11.487 14.393 11.494 14.427 /
\plot 11.494 14.427 11.500 14.463 /
\plot 11.500 14.463 11.504 14.501 /
\plot 11.504 14.501 11.508 14.544 /
\plot 11.508 14.544 11.513 14.584 /
\plot 11.513 14.584 11.517 14.628 /
\plot 11.517 14.628 11.519 14.671 /
\plot 11.519 14.671 11.521 14.715 /
\plot 11.521 14.715 11.523 14.757 /
\putrule from 11.523 14.757 to 11.523 14.798
\plot 11.523 14.798 11.525 14.840 /
\putrule from 11.525 14.840 to 11.525 14.878
\putrule from 11.525 14.878 to 11.525 14.916
\putrule from 11.525 14.916 to 11.525 14.952
\putrule from 11.525 14.952 to 11.525 14.986
\putrule from 11.525 14.986 to 11.525 15.028
\putrule from 11.525 15.028 to 11.525 15.071
\putrule from 11.525 15.071 to 11.525 15.113
\putrule from 11.525 15.113 to 11.525 15.153
\plot 11.525 15.153 11.527 15.193 /
\putrule from 11.527 15.193 to 11.527 15.232
\plot 11.527 15.232 11.529 15.270 /
\plot 11.529 15.270 11.532 15.306 /
\plot 11.532 15.306 11.534 15.339 /
\plot 11.534 15.339 11.538 15.369 /
\plot 11.538 15.369 11.542 15.397 /
\plot 11.542 15.397 11.546 15.420 /
\plot 11.546 15.420 11.551 15.443 /
\plot 11.551 15.443 11.557 15.462 /
\plot 11.557 15.462 11.565 15.483 /
\plot 11.565 15.483 11.574 15.502 /
\plot 11.574 15.502 11.585 15.519 /
\plot 11.585 15.519 11.597 15.536 /
\plot 11.597 15.536 11.612 15.549 /
\plot 11.612 15.549 11.627 15.562 /
\plot 11.627 15.562 11.644 15.572 /
\plot 11.644 15.572 11.661 15.583 /
\plot 11.661 15.583 11.678 15.589 /
\plot 11.678 15.589 11.697 15.596 /
\plot 11.697 15.596 11.714 15.602 /
\plot 11.714 15.602 11.733 15.604 /
\plot 11.733 15.604 11.754 15.608 /
\plot 11.754 15.608 11.779 15.613 /
\plot 11.779 15.613 11.807 15.615 /
\plot 11.807 15.615 11.841 15.617 /
\plot 11.841 15.617 11.879 15.619 /
\putrule from 11.879 15.619 to 11.921 15.619
\plot 11.921 15.619 11.959 15.621 /
\putrule from 11.959 15.621 to 11.987 15.621
\putrule from 11.987 15.621 to 11.999 15.621
\putrule from 11.999 15.621 to 12.002 15.621
}%
%
%
\linethickness=1pt
\setplotsymbol ({\makebox(0,0)[l]{\tencirc\symbol{'160}}})
\setsolid
{\color[rgb]{0,0,0}\putrule from  9.525 19.622 to  9.525 19.624
\putrule from  9.525 19.624 to  9.525 19.636
\putrule from  9.525 19.636 to  9.525 19.664
\putrule from  9.525 19.664 to  9.525 19.702
\putrule from  9.525 19.702 to  9.525 19.744
\putrule from  9.525 19.744 to  9.525 19.782
\putrule from  9.525 19.782 to  9.525 19.816
\putrule from  9.525 19.816 to  9.525 19.844
\putrule from  9.525 19.844 to  9.525 19.869
\putrule from  9.525 19.869 to  9.525 19.890
\putrule from  9.525 19.890 to  9.525 19.918
\putrule from  9.525 19.918 to  9.525 19.945
\putrule from  9.525 19.945 to  9.525 19.973
\putrule from  9.525 19.973 to  9.525 20.000
\putrule from  9.525 20.000 to  9.525 20.028
\putrule from  9.525 20.028 to  9.525 20.051
\putrule from  9.525 20.051 to  9.525 20.077
\putrule from  9.525 20.077 to  9.525 20.098
\putrule from  9.525 20.098 to  9.525 20.119
\putrule from  9.525 20.119 to  9.525 20.144
\putrule from  9.525 20.144 to  9.525 20.172
\putrule from  9.525 20.172 to  9.525 20.206
\putrule from  9.525 20.206 to  9.525 20.242
\putrule from  9.525 20.242 to  9.525 20.269
\putrule from  9.525 20.269 to  9.525 20.286
\putrule from  9.525 20.286 to  9.525 20.288
}%
%
%
\linethickness=1pt
\setplotsymbol ({\makebox(0,0)[l]{\tencirc\symbol{'160}}})
{\color[rgb]{0,0,0}\putrule from  8.954 19.622 to  8.954 19.619
\putrule from  8.954 19.619 to  8.954 19.602
\plot  8.954 19.602  8.956 19.573 /
\putrule from  8.956 19.573 to  8.956 19.535
\plot  8.956 19.535  8.958 19.501 /
\plot  8.958 19.501  8.962 19.469 /
\plot  8.962 19.469  8.964 19.442 /
\plot  8.964 19.442  8.970 19.414 /
\plot  8.970 19.414  8.975 19.393 /
\plot  8.975 19.393  8.981 19.370 /
\plot  8.981 19.370  8.987 19.344 /
\plot  8.987 19.344  8.994 19.317 /
\plot  8.994 19.317  9.004 19.287 /
\plot  9.004 19.287  9.013 19.260 /
\plot  9.013 19.260  9.021 19.232 /
\plot  9.021 19.232  9.032 19.207 /
\plot  9.032 19.207  9.040 19.183 /
\plot  9.040 19.183  9.049 19.160 /
\plot  9.049 19.160  9.057 19.139 /
\plot  9.057 19.139  9.066 19.118 /
\plot  9.066 19.118  9.076 19.097 /
\plot  9.076 19.097  9.087 19.075 /
\plot  9.087 19.075  9.097 19.054 /
\plot  9.097 19.054  9.106 19.035 /
\plot  9.106 19.035  9.116 19.016 /
\plot  9.116 19.016  9.127 18.999 /
\plot  9.127 18.999  9.136 18.984 /
\plot  9.136 18.984  9.144 18.970 /
\plot  9.144 18.970  9.155 18.955 /
\plot  9.155 18.955  9.167 18.940 /
\plot  9.167 18.940  9.182 18.923 /
\plot  9.182 18.923  9.197 18.904 /
\plot  9.197 18.904  9.216 18.885 /
\plot  9.216 18.885  9.231 18.868 /
\plot  9.231 18.868  9.237 18.862 /
\plot  9.237 18.862  9.239 18.860 /
}%
%
%
\put{$g$ tori
} [lB] at 12.097 15.431
%
%
\put{$\alpha_L$
} [lB] at  6.353 19.431
%
%
\put{$\alpha_R$
} [lB] at 11.335 19.431
%
%
\put{$S^1 \times S^1$
} [lB] at 12.106 18
%
\put{$\gamma$
} [lB] at  9.430 20.479

\put{$F_0$} [1B] at 5 16
\linethickness=0pt
\putrectangle corners at  6.921 20.714 and 12.129 11.955
\endpicture}

%% file: Concave1H.1.tex
\font\thinlinefont=cmr5
\mbox{\beginpicture
\small
\setcoordinatesystem units <0.5cm,0.5cm>
\unitlength=1.04987cm
\linethickness=1pt
\setplotsymbol ({\makebox(0,0)[l]{\tencirc\symbol{'160}}})
\setshadesymbol ({\thinlinefont .})
\setlinear
%
%
\linethickness= 0.500pt
\setplotsymbol ({\thinlinefont .})
{\color[rgb]{0,0,0}\circulararc 81.203 degrees from  3.143 21.336 center at  3.580 21.733
}%
%
%
\linethickness= 0.500pt
\setplotsymbol ({\thinlinefont .})
\setdots < 0.0953cm>
{\color[rgb]{0,0,0}\circulararc 75.750 degrees from  4.000 21.241 center at  3.511 20.737
}%
%
%
\linethickness= 0.500pt
\setplotsymbol ({\thinlinefont .})
\setsolid
{\color[rgb]{0,0,0}\ellipticalarc axes ratio  0.425:0.425  360 degrees 
	from  3.950 21.336 center at  3.524 21.336
}%
%
%
\linethickness= 0.500pt
\setplotsymbol ({\thinlinefont .})
{\color[rgb]{0,0,0}\circulararc 81.203 degrees from  5.715 18.955 center at  6.152 19.352
}%
%
%
\linethickness= 0.500pt
\setplotsymbol ({\thinlinefont .})
\setdots < 0.0953cm>
{\color[rgb]{0,0,0}\circulararc 75.750 degrees from  6.572 18.860 center at  6.082 18.356
}%
%
%
\linethickness= 0.500pt
\setplotsymbol ({\thinlinefont .})
\setsolid
{\color[rgb]{0,0,0}\ellipticalarc axes ratio  0.425:0.425  360 degrees 
	from  6.521 18.955 center at  6.096 18.955
}%
%
%
\linethickness= 0.500pt
\setplotsymbol ({\thinlinefont .})
{\color[rgb]{0,0,0}\ellipticalarc axes ratio  0.673:0.673  360 degrees 
	from  2.959 20.574 center at  2.286 20.574
}%
%
%
\linethickness=1pt
\setplotsymbol ({\makebox(0,0)[l]{\tencirc\symbol{'160}}})
{\color[rgb]{0,0,0}\putrule from  1.619 22.193 to  4.286 22.193
}%
%
%
\linethickness=1pt
\setplotsymbol ({\makebox(0,0)[l]{\tencirc\symbol{'160}}})
{\color[rgb]{0,0,0}\putrule from  6.953 20.098 to  6.953 17.431
}%
%
%
\linethickness= 0.500pt
\setplotsymbol ({\thinlinefont .})
\setdots < 0.0953cm>
{\color[rgb]{0,0,0}\plot  3.810 20.955  5.810 19.145 /
}%
%
%
\put{$0$
} [lB] at  2.095 19.3
%
%
\linethickness= 0.500pt
\setplotsymbol ({\thinlinefont .})
\setsolid
{\color[rgb]{0,0,0}\circulararc 81.203 degrees from 13.526 21.336 center at 13.962 21.733
}%
%
%
\linethickness= 0.500pt
\setplotsymbol ({\thinlinefont .})
\setdots < 0.0953cm>
{\color[rgb]{0,0,0}\circulararc 75.750 degrees from 14.383 21.241 center at 13.893 20.737
}%
%
%
\linethickness= 0.500pt
\setplotsymbol ({\thinlinefont .})
\setsolid
{\color[rgb]{0,0,0}\ellipticalarc axes ratio  0.425:0.425  360 degrees 
	from 14.332 21.336 center at 13.906 21.336
}%
%
%
\linethickness= 0.500pt
\setplotsymbol ({\thinlinefont .})
{\color[rgb]{0,0,0}\circulararc 81.203 degrees from 16.097 18.955 center at 16.534 19.352
}%
%
%
\linethickness= 0.500pt
\setplotsymbol ({\thinlinefont .})
\setdots < 0.0953cm>
{\color[rgb]{0,0,0}\circulararc 75.750 degrees from 16.954 18.860 center at 16.465 18.356
}%
%
%
\linethickness= 0.500pt
\setplotsymbol ({\thinlinefont .})
\setsolid
{\color[rgb]{0,0,0}\ellipticalarc axes ratio  0.425:0.425  360 degrees 
	from 16.904 18.955 center at 16.478 18.955
}%
%
%
\linethickness= 0.500pt
\setplotsymbol ({\thinlinefont .})
\setdots < 0.0953cm>
{\color[rgb]{0,0,0}\ellipticalarc axes ratio  1.943:1.943  360 degrees 
	from  4.801 21.241 center at  2.857 21.241
}%
%
%
\linethickness= 0.500pt
\setplotsymbol ({\thinlinefont .})
{\color[rgb]{0,0,0}\ellipticalarc axes ratio  1.943:1.943  360 degrees 
	from 14.802 21.717 center at 12.859 21.717
}%
%
%
\linethickness= 0.500pt
\setplotsymbol ({\thinlinefont .})
{\color[rgb]{0,0,0}\ellipticalarc axes ratio  1.943:1.943  360 degrees 
	from 19.088 18.098 center at 17.145 18.098
}%
%
%
\linethickness= 0.500pt
\setplotsymbol ({\thinlinefont .})
{\color[rgb]{0,0,0}\ellipticalarc axes ratio  1.943:1.943  360 degrees 
	from  8.420 18.669 center at  6.477 18.669
}%
%
%
\linethickness=2pt
\setplotsymbol ({\thinlinefont .})
\setsolid
{\color[rgb]{0,0,0}\plot  8.477 20.479 10.858 20.574 /
%
%
\plot  9.602 20.206 10.859 20.574  9.577 20.840 /
}%
%
%
\linethickness= 0.500pt
\setplotsymbol ({\thinlinefont .})
\setdots < 0.0953cm>
{\color[rgb]{0,0,0}\plot 14.192 20.955 16.192 19.145 /
}%
%
%
\linethickness=1pt
\setplotsymbol ({\makebox(0,0)[l]{\tencirc\symbol{'160}}})
\setsolid
{\color[rgb]{0,0,0}\putrule from 11.525 22.193 to 13.145 22.193
}%
%
%
\linethickness=1pt
\setplotsymbol ({\makebox(0,0)[l]{\tencirc\symbol{'160}}})
{\color[rgb]{0,0,0}\putrule from 13.335 22.193 to 14.573 22.193
}%
%
%
\linethickness=1pt
\setplotsymbol ({\makebox(0,0)[l]{\tencirc\symbol{'160}}})
{\color[rgb]{0,0,0}\putrule from 17.240 18.479 to 17.240 19.717
}%
%
%
\linethickness=1pt
\setplotsymbol ({\makebox(0,0)[l]{\tencirc\symbol{'160}}})
{\color[rgb]{0,0,0}\putrule from 17.240 18.193 to 17.240 16.478
}%
%
%
\linethickness= 0.500pt
\setplotsymbol ({\thinlinefont .})
{\color[rgb]{0,0,0}\putrule from 13.621 21.622 to 13.619 21.622
\plot 13.619 21.622 13.608 21.628 /
\plot 13.608 21.628 13.587 21.641 /
\plot 13.587 21.641 13.551 21.660 /
\plot 13.551 21.660 13.509 21.683 /
\plot 13.509 21.683 13.466 21.709 /
\plot 13.466 21.709 13.426 21.732 /
\plot 13.426 21.732 13.390 21.753 /
\plot 13.390 21.753 13.363 21.774 /
\plot 13.363 21.774 13.339 21.793 /
\plot 13.339 21.793 13.320 21.810 /
\plot 13.320 21.810 13.303 21.827 /
\plot 13.303 21.827 13.291 21.846 /
\plot 13.291 21.846 13.278 21.867 /
\plot 13.278 21.867 13.269 21.888 /
\plot 13.269 21.888 13.261 21.912 /
\plot 13.261 21.912 13.255 21.937 /
\plot 13.255 21.937 13.248 21.967 /
\plot 13.248 21.967 13.244 21.996 /
\plot 13.244 21.996 13.242 22.026 /
\putrule from 13.242 22.026 to 13.242 22.060
\plot 13.242 22.060 13.240 22.092 /
\putrule from 13.240 22.092 to 13.240 22.126
\putrule from 13.240 22.126 to 13.240 22.162
\putrule from 13.240 22.162 to 13.240 22.193
\putrule from 13.240 22.193 to 13.240 22.227
\putrule from 13.240 22.227 to 13.240 22.265
\plot 13.240 22.265 13.238 22.303 /
\plot 13.238 22.303 13.236 22.344 /
\plot 13.236 22.344 13.233 22.384 /
\plot 13.233 22.384 13.229 22.424 /
\plot 13.229 22.424 13.223 22.464 /
\plot 13.223 22.464 13.216 22.504 /
\plot 13.216 22.504 13.208 22.540 /
\plot 13.208 22.540 13.200 22.574 /
\plot 13.200 22.574 13.187 22.604 /
\plot 13.187 22.604 13.174 22.629 /
\plot 13.174 22.629 13.161 22.653 /
\plot 13.161 22.653 13.145 22.674 /
\plot 13.145 22.674 13.125 22.693 /
\plot 13.125 22.693 13.102 22.708 /
\plot 13.102 22.708 13.079 22.720 /
\plot 13.079 22.720 13.053 22.731 /
\plot 13.053 22.731 13.026 22.739 /
\plot 13.026 22.739 12.996 22.746 /
\plot 12.996 22.746 12.969 22.748 /
\plot 12.969 22.748 12.939 22.750 /
\putrule from 12.939 22.750 to 12.910 22.750
\plot 12.910 22.750 12.884 22.748 /
\plot 12.884 22.748 12.859 22.744 /
\plot 12.859 22.744 12.833 22.739 /
\plot 12.833 22.739 12.812 22.733 /
\plot 12.812 22.733 12.783 22.722 /
\plot 12.783 22.722 12.753 22.712 /
\plot 12.753 22.712 12.728 22.697 /
\plot 12.728 22.697 12.702 22.682 /
\plot 12.702 22.682 12.679 22.665 /
\plot 12.679 22.665 12.660 22.646 /
\plot 12.660 22.646 12.641 22.629 /
\plot 12.641 22.629 12.626 22.610 /
\plot 12.626 22.610 12.615 22.591 /
\plot 12.615 22.591 12.605 22.574 /
\plot 12.605 22.574 12.594 22.553 /
\plot 12.594 22.553 12.588 22.528 /
\plot 12.588 22.528 12.584 22.500 /
\plot 12.584 22.500 12.579 22.466 /
\plot 12.579 22.466 12.575 22.430 /
\plot 12.575 22.430 12.573 22.403 /
\putrule from 12.573 22.403 to 12.573 22.386
\putrule from 12.573 22.386 to 12.573 22.384
}%
%
%
\linethickness= 0.500pt
\setplotsymbol ({\thinlinefont .})
{\color[rgb]{0,0,0}\putrule from 12.573 22.098 to 12.573 22.096
\plot 12.573 22.096 12.575 22.090 /
\plot 12.575 22.090 12.577 22.073 /
\plot 12.577 22.073 12.584 22.045 /
\plot 12.584 22.045 12.590 22.005 /
\plot 12.590 22.005 12.598 21.956 /
\plot 12.598 21.956 12.609 21.905 /
\plot 12.609 21.905 12.620 21.852 /
\plot 12.620 21.852 12.630 21.802 /
\plot 12.630 21.802 12.641 21.753 /
\plot 12.641 21.753 12.651 21.711 /
\plot 12.651 21.711 12.660 21.673 /
\plot 12.660 21.673 12.670 21.637 /
\plot 12.670 21.637 12.679 21.605 /
\plot 12.679 21.605 12.689 21.573 /
\plot 12.689 21.573 12.700 21.541 /
\plot 12.700 21.541 12.713 21.507 /
\plot 12.713 21.507 12.728 21.474 /
\plot 12.728 21.474 12.742 21.440 /
\plot 12.742 21.440 12.759 21.406 /
\plot 12.759 21.406 12.778 21.372 /
\plot 12.778 21.372 12.799 21.338 /
\plot 12.799 21.338 12.823 21.306 /
\plot 12.823 21.306 12.846 21.277 /
\plot 12.846 21.277 12.869 21.249 /
\plot 12.869 21.249 12.895 21.224 /
\plot 12.895 21.224 12.920 21.205 /
\plot 12.920 21.205 12.948 21.186 /
\plot 12.948 21.186 12.973 21.173 /
\plot 12.973 21.173 13.003 21.160 /
\plot 13.003 21.160 13.030 21.154 /
\plot 13.030 21.154 13.062 21.150 /
\plot 13.062 21.150 13.096 21.148 /
\plot 13.096 21.148 13.134 21.150 /
\plot 13.134 21.150 13.178 21.156 /
\plot 13.178 21.156 13.229 21.165 /
\plot 13.229 21.165 13.284 21.177 /
\plot 13.284 21.177 13.343 21.190 /
\plot 13.343 21.190 13.401 21.205 /
\plot 13.401 21.205 13.451 21.220 /
\plot 13.451 21.220 13.490 21.230 /
\plot 13.490 21.230 13.513 21.237 /
\plot 13.513 21.237 13.523 21.241 /
\putrule from 13.523 21.241 to 13.526 21.241
}%
%
%
\linethickness= 0.500pt
\setplotsymbol ({\thinlinefont .})
{\color[rgb]{0,0,0}\putrule from 16.764 18.574 to 16.764 18.572
\plot 16.764 18.572 16.770 18.565 /
\plot 16.770 18.565 16.785 18.548 /
\plot 16.785 18.548 16.806 18.523 /
\plot 16.806 18.523 16.834 18.491 /
\plot 16.834 18.491 16.863 18.459 /
\plot 16.863 18.459 16.893 18.432 /
\plot 16.893 18.432 16.925 18.407 /
\plot 16.925 18.407 16.952 18.387 /
\plot 16.952 18.387 16.984 18.373 /
\plot 16.984 18.373 17.016 18.360 /
\plot 17.016 18.360 17.050 18.352 /
\plot 17.050 18.352 17.077 18.347 /
\plot 17.077 18.347 17.109 18.343 /
\plot 17.109 18.343 17.141 18.341 /
\plot 17.141 18.341 17.177 18.339 /
\plot 17.177 18.339 17.215 18.337 /
\plot 17.215 18.337 17.253 18.335 /
\putrule from 17.253 18.335 to 17.295 18.335
\plot 17.295 18.335 17.335 18.332 /
\putrule from 17.335 18.332 to 17.376 18.332
\plot 17.376 18.332 17.418 18.330 /
\plot 17.418 18.330 17.456 18.328 /
\plot 17.456 18.328 17.494 18.326 /
\plot 17.494 18.326 17.530 18.322 /
\plot 17.530 18.322 17.562 18.318 /
\plot 17.562 18.318 17.594 18.311 /
\plot 17.594 18.311 17.621 18.303 /
\plot 17.621 18.303 17.647 18.294 /
\plot 17.647 18.294 17.672 18.284 /
\plot 17.672 18.284 17.695 18.271 /
\plot 17.695 18.271 17.719 18.256 /
\plot 17.719 18.256 17.738 18.241 /
\plot 17.738 18.241 17.757 18.222 /
\plot 17.757 18.222 17.774 18.203 /
\plot 17.774 18.203 17.791 18.184 /
\plot 17.791 18.184 17.803 18.163 /
\plot 17.803 18.163 17.814 18.142 /
\plot 17.814 18.142 17.824 18.121 /
\plot 17.824 18.121 17.831 18.100 /
\plot 17.831 18.100 17.837 18.078 /
\plot 17.837 18.078 17.841 18.057 /
\plot 17.841 18.057 17.843 18.038 /
\putrule from 17.843 18.038 to 17.843 18.017
\plot 17.843 18.017 17.841 17.994 /
\plot 17.841 17.994 17.839 17.968 /
\plot 17.839 17.968 17.833 17.943 /
\plot 17.833 17.943 17.824 17.918 /
\plot 17.824 17.918 17.812 17.892 /
\plot 17.812 17.892 17.799 17.869 /
\plot 17.799 17.869 17.782 17.846 /
\plot 17.782 17.846 17.765 17.824 /
\plot 17.765 17.824 17.746 17.805 /
\plot 17.746 17.805 17.727 17.791 /
\plot 17.727 17.791 17.706 17.776 /
\plot 17.706 17.776 17.685 17.763 /
\plot 17.685 17.763 17.661 17.755 /
\plot 17.661 17.755 17.638 17.746 /
\plot 17.638 17.746 17.609 17.740 /
\plot 17.609 17.740 17.577 17.733 /
\plot 17.577 17.733 17.539 17.729 /
\plot 17.539 17.729 17.496 17.725 /
\plot 17.496 17.725 17.450 17.721 /
\plot 17.450 17.721 17.407 17.719 /
\putrule from 17.407 17.719 to 17.371 17.719
\plot 17.371 17.719 17.348 17.716 /
\putrule from 17.348 17.716 to 17.338 17.716
\putrule from 17.338 17.716 to 17.335 17.716
}%
%
%
\linethickness= 0.500pt
\setplotsymbol ({\thinlinefont .})
{\color[rgb]{0,0,0}\putrule from 17.145 17.716 to 17.143 17.716
\putrule from 17.143 17.716 to 17.132 17.716
\plot 17.132 17.716 17.109 17.719 /
\putrule from 17.109 17.719 to 17.073 17.719
\plot 17.073 17.719 17.031 17.721 /
\plot 17.031 17.721 16.984 17.725 /
\plot 16.984 17.725 16.942 17.729 /
\plot 16.942 17.729 16.904 17.733 /
\plot 16.904 17.733 16.872 17.740 /
\plot 16.872 17.740 16.842 17.746 /
\plot 16.842 17.746 16.819 17.755 /
\plot 16.819 17.755 16.796 17.763 /
\plot 16.796 17.763 16.775 17.776 /
\plot 16.775 17.776 16.753 17.791 /
\plot 16.753 17.791 16.732 17.805 /
\plot 16.732 17.805 16.713 17.824 /
\plot 16.713 17.824 16.692 17.848 /
\plot 16.692 17.848 16.673 17.871 /
\plot 16.673 17.871 16.654 17.896 /
\plot 16.654 17.896 16.637 17.922 /
\plot 16.637 17.922 16.620 17.949 /
\plot 16.620 17.949 16.603 17.977 /
\plot 16.603 17.977 16.588 18.004 /
\plot 16.588 18.004 16.573 18.034 /
\plot 16.573 18.034 16.561 18.059 /
\plot 16.561 18.059 16.548 18.087 /
\plot 16.548 18.087 16.535 18.117 /
\plot 16.535 18.117 16.521 18.150 /
\plot 16.521 18.150 16.504 18.186 /
\plot 16.504 18.186 16.485 18.229 /
\plot 16.485 18.229 16.466 18.277 /
\plot 16.466 18.277 16.447 18.326 /
\plot 16.447 18.326 16.425 18.375 /
\plot 16.425 18.375 16.408 18.417 /
\plot 16.408 18.417 16.396 18.449 /
\plot 16.396 18.449 16.387 18.468 /
\plot 16.387 18.468 16.383 18.476 /
\putrule from 16.383 18.476 to 16.383 18.479
}%
%
%
\put{$0$
} [lB] at 12.668 22.9
%
%
\put{$B_0$
} [lB] at  1.048 22.670
%
%
\put{$B_1$
} [lB] at  7.620 16.574
%
%
\put{$B_0$
} [lB] at 14.002 23.431
%
%
\put{$B_1$
} [lB] at 18.193 19.812
\linethickness=0pt
\putrectangle corners at  0.897 23.749 and 19.105 16.140
\endpicture}

%% file: Concave1H.2.tex
\font\thinlinefont=cmr5
\mbox{\beginpicture
\small
\setcoordinatesystem units <0.5cm,0.5cm>
\unitlength=1.04987cm
\linethickness=1pt
\setplotsymbol ({\makebox(0,0)[l]{\tencirc\symbol{'160}}})
\setshadesymbol ({\thinlinefont .})
\setlinear
%
%
\linethickness= 0.500pt
\setplotsymbol ({\thinlinefont .})
\setdots < 0.0953cm>
{\color[rgb]{0,0,0}\ellipticalarc axes ratio  1.943:1.943  360 degrees 
	from  8.420 18.669 center at  6.477 18.669
}%
%
%
\linethickness= 0.500pt
\setplotsymbol ({\thinlinefont .})
{\color[rgb]{0,0,0}\ellipticalarc axes ratio  1.943:1.943  360 degrees 
	from  4.801 21.241 center at  2.857 21.241
}%
%
%
\linethickness=1pt
\setplotsymbol ({\makebox(0,0)[l]{\tencirc\symbol{'160}}})
\setsolid
{\color[rgb]{0,0,0}\putrule from  1.619 22.193 to  4.286 22.193
}%
%
%
\linethickness=1pt
\setplotsymbol ({\makebox(0,0)[l]{\tencirc\symbol{'160}}})
{\color[rgb]{0,0,0}\putrule from  6.953 20.098 to  6.953 17.431
}%
%
%
\linethickness= 0.500pt
\setplotsymbol ({\thinlinefont .})
{\color[rgb]{0,0,0}\putrule from  2.953 22.003 to  2.953 22.001
\plot  2.953 22.001  2.959 21.994 /
\plot  2.959 21.994  2.972 21.975 /
\plot  2.972 21.975  2.991 21.948 /
\plot  2.991 21.948  3.014 21.916 /
\plot  3.014 21.916  3.040 21.880 /
\plot  3.040 21.880  3.063 21.848 /
\plot  3.063 21.848  3.084 21.819 /
\plot  3.084 21.819  3.105 21.793 /
\plot  3.105 21.793  3.124 21.772 /
\plot  3.124 21.772  3.141 21.751 /
\plot  3.141 21.751  3.160 21.732 /
\plot  3.160 21.732  3.177 21.715 /
\plot  3.177 21.715  3.198 21.698 /
\plot  3.198 21.698  3.219 21.679 /
\plot  3.219 21.679  3.245 21.662 /
\plot  3.245 21.662  3.270 21.645 /
\plot  3.270 21.645  3.300 21.630 /
\plot  3.300 21.630  3.330 21.615 /
\plot  3.330 21.615  3.363 21.605 /
\plot  3.363 21.605  3.397 21.594 /
\plot  3.397 21.594  3.433 21.586 /
\plot  3.433 21.586  3.469 21.579 /
\plot  3.469 21.579  3.509 21.573 /
\plot  3.509 21.573  3.539 21.571 /
\putrule from  3.539 21.571 to  3.573 21.571
\putrule from  3.573 21.571 to  3.607 21.571
\plot  3.607 21.571  3.645 21.573 /
\plot  3.645 21.573  3.685 21.577 /
\plot  3.685 21.577  3.727 21.582 /
\plot  3.727 21.582  3.772 21.590 /
\plot  3.772 21.590  3.816 21.601 /
\plot  3.816 21.601  3.863 21.613 /
\plot  3.863 21.613  3.907 21.628 /
\plot  3.907 21.628  3.954 21.645 /
\plot  3.954 21.645  3.996 21.664 /
\plot  3.996 21.664  4.041 21.685 /
\plot  4.041 21.685  4.081 21.711 /
\plot  4.081 21.711  4.121 21.736 /
\plot  4.121 21.736  4.159 21.764 /
\plot  4.159 21.764  4.189 21.789 /
\plot  4.189 21.789  4.221 21.816 /
\plot  4.221 21.816  4.250 21.846 /
\plot  4.250 21.846  4.282 21.878 /
\plot  4.282 21.878  4.312 21.912 /
\plot  4.312 21.912  4.343 21.948 /
\plot  4.343 21.948  4.375 21.986 /
\plot  4.375 21.986  4.407 22.026 /
\plot  4.407 22.026  4.439 22.066 /
\plot  4.439 22.066  4.470 22.109 /
\plot  4.470 22.109  4.502 22.151 /
\plot  4.502 22.151  4.532 22.195 /
\plot  4.532 22.195  4.561 22.238 /
\plot  4.561 22.238  4.591 22.280 /
\plot  4.591 22.280  4.619 22.322 /
\plot  4.619 22.322  4.646 22.363 /
\plot  4.646 22.363  4.671 22.401 /
\plot  4.671 22.401  4.697 22.439 /
\plot  4.697 22.439  4.722 22.475 /
\plot  4.722 22.475  4.748 22.511 /
\plot  4.748 22.511  4.773 22.549 /
\plot  4.773 22.549  4.801 22.587 /
\plot  4.801 22.587  4.828 22.625 /
\plot  4.828 22.625  4.856 22.661 /
\plot  4.856 22.661  4.885 22.697 /
\plot  4.885 22.697  4.915 22.735 /
\plot  4.915 22.735  4.945 22.769 /
\plot  4.945 22.769  4.976 22.805 /
\plot  4.976 22.805  5.010 22.839 /
\plot  5.010 22.839  5.042 22.871 /
\plot  5.042 22.871  5.076 22.900 /
\plot  5.076 22.900  5.110 22.928 /
\plot  5.110 22.928  5.144 22.953 /
\plot  5.144 22.953  5.177 22.976 /
\plot  5.177 22.976  5.213 23.000 /
\plot  5.213 23.000  5.247 23.019 /
\plot  5.247 23.019  5.283 23.036 /
\plot  5.283 23.036  5.319 23.050 /
\plot  5.319 23.050  5.355 23.065 /
\plot  5.355 23.065  5.393 23.076 /
\plot  5.393 23.076  5.436 23.086 /
\plot  5.436 23.086  5.478 23.097 /
\plot  5.478 23.097  5.522 23.106 /
\plot  5.522 23.106  5.571 23.112 /
\plot  5.571 23.112  5.620 23.116 /
\plot  5.620 23.116  5.671 23.118 /
\plot  5.671 23.118  5.721 23.120 /
\plot  5.721 23.120  5.774 23.118 /
\plot  5.774 23.118  5.825 23.116 /
\plot  5.825 23.116  5.876 23.112 /
\plot  5.876 23.112  5.927 23.106 /
\plot  5.927 23.106  5.977 23.097 /
\plot  5.977 23.097  6.024 23.086 /
\plot  6.024 23.086  6.071 23.076 /
\plot  6.071 23.076  6.115 23.065 /
\plot  6.115 23.065  6.160 23.050 /
\plot  6.160 23.050  6.198 23.036 /
\plot  6.198 23.036  6.238 23.021 /
\plot  6.238 23.021  6.276 23.004 /
\plot  6.276 23.004  6.316 22.987 /
\plot  6.316 22.987  6.356 22.966 /
\plot  6.356 22.966  6.397 22.945 /
\plot  6.397 22.945  6.435 22.921 /
\plot  6.435 22.921  6.475 22.896 /
\plot  6.475 22.896  6.515 22.871 /
\plot  6.515 22.871  6.553 22.843 /
\plot  6.553 22.843  6.589 22.813 /
\plot  6.589 22.813  6.625 22.784 /
\plot  6.625 22.784  6.659 22.754 /
\plot  6.659 22.754  6.691 22.725 /
\plot  6.691 22.725  6.723 22.695 /
\plot  6.723 22.695  6.750 22.663 /
\plot  6.750 22.663  6.775 22.634 /
\plot  6.775 22.634  6.799 22.604 /
\plot  6.799 22.604  6.822 22.572 /
\plot  6.822 22.572  6.843 22.543 /
\plot  6.843 22.543  6.862 22.509 /
\plot  6.862 22.509  6.883 22.473 /
\plot  6.883 22.473  6.900 22.437 /
\plot  6.900 22.437  6.915 22.399 /
\plot  6.915 22.399  6.930 22.358 /
\plot  6.930 22.358  6.943 22.316 /
\plot  6.943 22.316  6.953 22.274 /
\plot  6.953 22.274  6.962 22.227 /
\plot  6.962 22.227  6.968 22.183 /
\plot  6.968 22.183  6.972 22.136 /
\putrule from  6.972 22.136 to  6.972 22.090
\putrule from  6.972 22.090 to  6.972 22.041
\plot  6.972 22.041  6.968 21.994 /
\plot  6.968 21.994  6.964 21.948 /
\plot  6.964 21.948  6.955 21.901 /
\plot  6.955 21.901  6.947 21.857 /
\plot  6.947 21.857  6.934 21.810 /
\plot  6.934 21.810  6.921 21.764 /
\plot  6.921 21.764  6.909 21.728 /
\plot  6.909 21.728  6.894 21.687 /
\plot  6.894 21.687  6.879 21.647 /
\plot  6.879 21.647  6.862 21.607 /
\plot  6.862 21.607  6.843 21.565 /
\plot  6.843 21.565  6.824 21.520 /
\plot  6.824 21.520  6.803 21.476 /
\plot  6.803 21.476  6.780 21.429 /
\plot  6.780 21.429  6.756 21.383 /
\plot  6.756 21.383  6.733 21.334 /
\plot  6.733 21.334  6.708 21.285 /
\plot  6.708 21.285  6.684 21.239 /
\plot  6.684 21.239  6.659 21.190 /
\plot  6.659 21.190  6.634 21.141 /
\plot  6.634 21.141  6.610 21.095 /
\plot  6.610 21.095  6.585 21.050 /
\plot  6.585 21.050  6.562 21.004 /
\plot  6.562 21.004  6.540 20.961 /
\plot  6.540 20.961  6.519 20.919 /
\plot  6.519 20.919  6.498 20.877 /
\plot  6.498 20.877  6.479 20.836 /
\plot  6.479 20.836  6.462 20.796 /
\plot  6.462 20.796  6.441 20.752 /
\plot  6.441 20.752  6.424 20.709 /
\plot  6.424 20.709  6.407 20.665 /
\plot  6.407 20.665  6.390 20.621 /
\plot  6.390 20.621  6.373 20.576 /
\plot  6.373 20.576  6.358 20.530 /
\plot  6.358 20.530  6.344 20.483 /
\plot  6.344 20.483  6.331 20.436 /
\plot  6.331 20.436  6.318 20.390 /
\plot  6.318 20.390  6.306 20.343 /
\plot  6.306 20.343  6.295 20.297 /
\plot  6.295 20.297  6.287 20.250 /
\plot  6.287 20.250  6.278 20.206 /
\plot  6.278 20.206  6.272 20.161 /
\plot  6.272 20.161  6.267 20.119 /
\plot  6.267 20.119  6.261 20.077 /
\plot  6.261 20.077  6.259 20.036 /
\plot  6.259 20.036  6.257 19.998 /
\plot  6.257 19.998  6.255 19.960 /
\putrule from  6.255 19.960 to  6.255 19.922
\putrule from  6.255 19.922 to  6.255 19.882
\plot  6.255 19.882  6.257 19.842 /
\plot  6.257 19.842  6.259 19.799 /
\plot  6.259 19.799  6.263 19.757 /
\plot  6.263 19.757  6.267 19.713 /
\plot  6.267 19.713  6.274 19.670 /
\plot  6.274 19.670  6.280 19.628 /
\plot  6.280 19.628  6.289 19.583 /
\plot  6.289 19.583  6.297 19.541 /
\plot  6.297 19.541  6.308 19.501 /
\plot  6.308 19.501  6.318 19.461 /
\plot  6.318 19.461  6.329 19.423 /
\plot  6.329 19.423  6.339 19.387 /
\plot  6.339 19.387  6.350 19.353 /
\plot  6.350 19.353  6.363 19.321 /
\plot  6.363 19.321  6.373 19.293 /
\plot  6.373 19.293  6.386 19.266 /
\plot  6.386 19.266  6.399 19.241 /
\plot  6.399 19.241  6.416 19.207 /
\plot  6.416 19.207  6.437 19.175 /
\plot  6.437 19.175  6.458 19.143 /
\plot  6.458 19.143  6.483 19.116 /
\plot  6.483 19.116  6.513 19.084 /
\plot  6.513 19.084  6.545 19.054 /
\plot  6.545 19.054  6.581 19.025 /
\plot  6.581 19.025  6.612 18.997 /
\plot  6.612 18.997  6.640 18.976 /
\plot  6.640 18.976  6.659 18.961 /
\plot  6.659 18.961  6.665 18.955 /
\putrule from  6.665 18.955 to  6.668 18.955
}%
%
%
\linethickness= 0.500pt
\setplotsymbol ({\thinlinefont .})
{\color[rgb]{0,0,0}\plot  3.048 20.479  3.054 20.483 /
\plot  3.054 20.483  3.069 20.496 /
\plot  3.069 20.496  3.092 20.515 /
\plot  3.092 20.515  3.126 20.540 /
\plot  3.126 20.540  3.167 20.570 /
\plot  3.167 20.570  3.211 20.604 /
\plot  3.211 20.604  3.255 20.637 /
\plot  3.255 20.637  3.298 20.667 /
\plot  3.298 20.667  3.338 20.697 /
\plot  3.338 20.697  3.376 20.722 /
\plot  3.376 20.722  3.410 20.745 /
\plot  3.410 20.745  3.440 20.764 /
\plot  3.440 20.764  3.469 20.781 /
\plot  3.469 20.781  3.497 20.798 /
\plot  3.497 20.798  3.524 20.811 /
\plot  3.524 20.811  3.552 20.826 /
\plot  3.552 20.826  3.579 20.839 /
\plot  3.579 20.839  3.609 20.849 /
\plot  3.609 20.849  3.639 20.862 /
\plot  3.639 20.862  3.668 20.872 /
\plot  3.668 20.872  3.702 20.883 /
\plot  3.702 20.883  3.734 20.894 /
\plot  3.734 20.894  3.770 20.904 /
\plot  3.770 20.904  3.804 20.915 /
\plot  3.804 20.915  3.840 20.925 /
\plot  3.840 20.925  3.873 20.936 /
\plot  3.873 20.936  3.909 20.947 /
\plot  3.909 20.947  3.943 20.955 /
\plot  3.943 20.955  3.979 20.966 /
\plot  3.979 20.966  4.013 20.976 /
\plot  4.013 20.976  4.049 20.987 /
\plot  4.049 20.987  4.081 20.997 /
\plot  4.081 20.997  4.113 21.008 /
\plot  4.113 21.008  4.149 21.018 /
\plot  4.149 21.018  4.183 21.031 /
\plot  4.183 21.031  4.221 21.044 /
\plot  4.221 21.044  4.259 21.059 /
\plot  4.259 21.059  4.299 21.076 /
\plot  4.299 21.076  4.339 21.090 /
\plot  4.339 21.090  4.379 21.110 /
\plot  4.379 21.110  4.420 21.126 /
\plot  4.420 21.126  4.460 21.146 /
\plot  4.460 21.146  4.498 21.165 /
\plot  4.498 21.165  4.536 21.186 /
\plot  4.536 21.186  4.572 21.205 /
\plot  4.572 21.205  4.606 21.226 /
\plot  4.606 21.226  4.638 21.245 /
\plot  4.638 21.245  4.669 21.266 /
\plot  4.669 21.266  4.699 21.287 /
\plot  4.699 21.287  4.729 21.311 /
\plot  4.729 21.311  4.756 21.334 /
\plot  4.756 21.334  4.784 21.359 /
\plot  4.784 21.359  4.811 21.385 /
\plot  4.811 21.385  4.839 21.412 /
\plot  4.839 21.412  4.866 21.440 /
\plot  4.866 21.440  4.894 21.469 /
\plot  4.894 21.469  4.919 21.499 /
\plot  4.919 21.499  4.947 21.531 /
\plot  4.947 21.531  4.972 21.562 /
\plot  4.972 21.562  4.997 21.594 /
\plot  4.997 21.594  5.023 21.624 /
\plot  5.023 21.624  5.048 21.656 /
\plot  5.048 21.656  5.072 21.685 /
\plot  5.072 21.685  5.093 21.715 /
\plot  5.093 21.715  5.116 21.742 /
\plot  5.116 21.742  5.137 21.770 /
\plot  5.137 21.770  5.160 21.795 /
\plot  5.160 21.795  5.184 21.827 /
\plot  5.184 21.827  5.209 21.857 /
\plot  5.209 21.857  5.237 21.886 /
\plot  5.237 21.886  5.264 21.916 /
\plot  5.264 21.916  5.292 21.948 /
\plot  5.292 21.948  5.321 21.977 /
\plot  5.321 21.977  5.353 22.007 /
\plot  5.353 22.007  5.385 22.037 /
\plot  5.385 22.037  5.414 22.066 /
\plot  5.414 22.066  5.446 22.092 /
\plot  5.446 22.092  5.476 22.117 /
\plot  5.476 22.117  5.508 22.138 /
\plot  5.508 22.138  5.537 22.159 /
\plot  5.537 22.159  5.565 22.178 /
\plot  5.565 22.178  5.592 22.195 /
\plot  5.592 22.195  5.620 22.208 /
\plot  5.620 22.208  5.652 22.223 /
\plot  5.652 22.223  5.683 22.236 /
\plot  5.683 22.236  5.715 22.248 /
\plot  5.715 22.248  5.749 22.257 /
\plot  5.749 22.257  5.783 22.263 /
\plot  5.783 22.263  5.817 22.269 /
\plot  5.817 22.269  5.850 22.274 /
\plot  5.850 22.274  5.884 22.276 /
\putrule from  5.884 22.276 to  5.916 22.276
\plot  5.916 22.276  5.948 22.274 /
\plot  5.948 22.274  5.975 22.272 /
\plot  5.975 22.272  6.001 22.267 /
\plot  6.001 22.267  6.026 22.263 /
\plot  6.026 22.263  6.049 22.257 /
\plot  6.049 22.257  6.066 22.250 /
\plot  6.066 22.250  6.085 22.244 /
\plot  6.085 22.244  6.102 22.236 /
\plot  6.102 22.236  6.119 22.227 /
\plot  6.119 22.227  6.134 22.217 /
\plot  6.134 22.217  6.149 22.204 /
\plot  6.149 22.204  6.162 22.191 /
\plot  6.162 22.191  6.174 22.176 /
\plot  6.174 22.176  6.185 22.159 /
\plot  6.185 22.159  6.193 22.140 /
\plot  6.193 22.140  6.200 22.121 /
\plot  6.200 22.121  6.204 22.100 /
\plot  6.204 22.100  6.208 22.077 /
\plot  6.208 22.077  6.210 22.054 /
\plot  6.210 22.054  6.208 22.028 /
\putrule from  6.208 22.028 to  6.208 22.003
\plot  6.208 22.003  6.204 21.977 /
\plot  6.204 21.977  6.200 21.950 /
\plot  6.200 21.950  6.193 21.922 /
\plot  6.193 21.922  6.187 21.891 /
\plot  6.187 21.891  6.179 21.857 /
\plot  6.179 21.857  6.168 21.821 /
\plot  6.168 21.821  6.155 21.783 /
\plot  6.155 21.783  6.143 21.740 /
\plot  6.143 21.740  6.128 21.700 /
\plot  6.128 21.700  6.113 21.656 /
\plot  6.113 21.656  6.096 21.611 /
\plot  6.096 21.611  6.079 21.567 /
\plot  6.079 21.567  6.062 21.522 /
\plot  6.062 21.522  6.043 21.478 /
\plot  6.043 21.478  6.026 21.435 /
\plot  6.026 21.435  6.007 21.391 /
\plot  6.007 21.391  5.988 21.349 /
\plot  5.988 21.349  5.969 21.304 /
\plot  5.969 21.304  5.952 21.264 /
\plot  5.952 21.264  5.933 21.224 /
\plot  5.933 21.224  5.914 21.184 /
\plot  5.914 21.184  5.895 21.139 /
\plot  5.895 21.139  5.874 21.097 /
\plot  5.874 21.097  5.855 21.050 /
\plot  5.855 21.050  5.834 21.004 /
\plot  5.834 21.004  5.812 20.957 /
\plot  5.812 20.957  5.791 20.908 /
\plot  5.791 20.908  5.770 20.860 /
\plot  5.770 20.860  5.749 20.811 /
\plot  5.749 20.811  5.728 20.762 /
\plot  5.728 20.762  5.709 20.716 /
\plot  5.709 20.716  5.690 20.669 /
\plot  5.690 20.669  5.673 20.623 /
\plot  5.673 20.623  5.656 20.580 /
\plot  5.656 20.580  5.641 20.536 /
\plot  5.641 20.536  5.628 20.496 /
\plot  5.628 20.496  5.616 20.455 /
\plot  5.616 20.455  5.605 20.415 /
\plot  5.605 20.415  5.592 20.371 /
\plot  5.592 20.371  5.582 20.328 /
\plot  5.582 20.328  5.573 20.284 /
\plot  5.573 20.284  5.565 20.240 /
\plot  5.565 20.240  5.556 20.195 /
\plot  5.556 20.195  5.550 20.151 /
\plot  5.550 20.151  5.544 20.104 /
\plot  5.544 20.104  5.537 20.060 /
\plot  5.537 20.060  5.533 20.013 /
\plot  5.533 20.013  5.529 19.969 /
\plot  5.529 19.969  5.527 19.926 /
\plot  5.527 19.926  5.522 19.884 /
\plot  5.522 19.884  5.520 19.844 /
\plot  5.520 19.844  5.518 19.806 /
\plot  5.518 19.806  5.516 19.770 /
\plot  5.516 19.770  5.514 19.734 /
\plot  5.514 19.734  5.512 19.702 /
\plot  5.512 19.702  5.510 19.668 /
\plot  5.510 19.668  5.505 19.630 /
\plot  5.505 19.630  5.501 19.590 /
\plot  5.501 19.590  5.497 19.550 /
\plot  5.497 19.550  5.493 19.511 /
\plot  5.493 19.511  5.486 19.471 /
\plot  5.486 19.471  5.480 19.433 /
\plot  5.480 19.433  5.472 19.395 /
\plot  5.472 19.395  5.465 19.357 /
\plot  5.465 19.357  5.457 19.323 /
\plot  5.457 19.323  5.448 19.289 /
\plot  5.448 19.289  5.440 19.257 /
\plot  5.440 19.257  5.431 19.230 /
\plot  5.431 19.230  5.423 19.202 /
\plot  5.423 19.202  5.414 19.177 /
\plot  5.414 19.177  5.402 19.147 /
\plot  5.402 19.147  5.389 19.120 /
\plot  5.389 19.120  5.376 19.090 /
\plot  5.376 19.090  5.359 19.061 /
\plot  5.359 19.061  5.340 19.027 /
\plot  5.340 19.027  5.319 18.989 /
\plot  5.319 18.989  5.296 18.953 /
\plot  5.296 18.953  5.275 18.917 /
\plot  5.275 18.917  5.256 18.887 /
\plot  5.256 18.887  5.245 18.868 /
\plot  5.245 18.868  5.239 18.862 /
\putrule from  5.239 18.862 to  5.239 18.860
}%
%
%
\linethickness= 0.500pt
\setplotsymbol ({\thinlinefont .})
{\color[rgb]{0,0,0}\putrule from  6.096 22.193 to  6.098 22.193
\putrule from  6.098 22.193 to  6.109 22.193
\putrule from  6.109 22.193 to  6.132 22.193
\plot  6.132 22.193  6.168 22.195 /
\plot  6.168 22.195  6.210 22.197 /
\plot  6.210 22.197  6.257 22.200 /
\plot  6.257 22.200  6.299 22.202 /
\plot  6.299 22.202  6.337 22.204 /
\plot  6.337 22.204  6.369 22.208 /
\plot  6.369 22.208  6.399 22.212 /
\plot  6.399 22.212  6.422 22.219 /
\plot  6.422 22.219  6.445 22.225 /
\plot  6.445 22.225  6.466 22.233 /
\plot  6.466 22.233  6.488 22.242 /
\plot  6.488 22.242  6.507 22.253 /
\plot  6.507 22.253  6.528 22.267 /
\plot  6.528 22.267  6.545 22.280 /
\plot  6.545 22.280  6.564 22.297 /
\plot  6.564 22.297  6.581 22.316 /
\plot  6.581 22.316  6.596 22.333 /
\plot  6.596 22.333  6.608 22.354 /
\plot  6.608 22.354  6.619 22.373 /
\plot  6.619 22.373  6.627 22.394 /
\plot  6.627 22.394  6.636 22.416 /
\plot  6.636 22.416  6.642 22.439 /
\plot  6.642 22.439  6.648 22.462 /
\plot  6.648 22.462  6.653 22.492 /
\plot  6.653 22.492  6.657 22.523 /
\plot  6.657 22.523  6.659 22.562 /
\plot  6.659 22.562  6.661 22.604 /
\plot  6.661 22.604  6.663 22.650 /
\plot  6.663 22.650  6.665 22.693 /
\plot  6.665 22.693  6.668 22.729 /
\putrule from  6.668 22.729 to  6.668 22.752
\putrule from  6.668 22.752 to  6.668 22.763
\putrule from  6.668 22.763 to  6.668 22.765
}%
%
%
\linethickness=1pt
\setplotsymbol ({\makebox(0,0)[l]{\tencirc\symbol{'160}}})
{\color[rgb]{0,0,0}\putrule from  2.191 22.193 to  2.191 22.191
\plot  2.191 22.191  2.189 22.183 /
\plot  2.189 22.183  2.187 22.162 /
\plot  2.187 22.162  2.180 22.130 /
\plot  2.180 22.130  2.176 22.083 /
\plot  2.176 22.083  2.170 22.032 /
\plot  2.170 22.032  2.163 21.979 /
\plot  2.163 21.979  2.159 21.929 /
\plot  2.159 21.929  2.157 21.882 /
\plot  2.157 21.882  2.155 21.838 /
\plot  2.155 21.838  2.157 21.797 /
\plot  2.157 21.797  2.161 21.759 /
\plot  2.161 21.759  2.167 21.723 /
\plot  2.167 21.723  2.176 21.685 /
\plot  2.176 21.685  2.182 21.656 /
\plot  2.182 21.656  2.193 21.626 /
\plot  2.193 21.626  2.203 21.594 /
\plot  2.203 21.594  2.216 21.562 /
\plot  2.216 21.562  2.231 21.531 /
\plot  2.231 21.531  2.250 21.497 /
\plot  2.250 21.497  2.269 21.465 /
\plot  2.269 21.465  2.290 21.431 /
\plot  2.290 21.431  2.316 21.399 /
\plot  2.316 21.399  2.341 21.368 /
\plot  2.341 21.368  2.369 21.338 /
\plot  2.369 21.338  2.398 21.311 /
\plot  2.398 21.311  2.428 21.285 /
\plot  2.428 21.285  2.462 21.262 /
\plot  2.462 21.262  2.496 21.241 /
\plot  2.496 21.241  2.529 21.222 /
\plot  2.529 21.222  2.565 21.207 /
\plot  2.565 21.207  2.603 21.192 /
\plot  2.603 21.192  2.639 21.184 /
\plot  2.639 21.184  2.678 21.175 /
\plot  2.678 21.175  2.718 21.169 /
\plot  2.718 21.169  2.760 21.165 /
\plot  2.760 21.165  2.805 21.160 /
\putrule from  2.805 21.160 to  2.853 21.160
\putrule from  2.853 21.160 to  2.902 21.160
\plot  2.902 21.160  2.953 21.162 /
\plot  2.953 21.162  3.008 21.167 /
\plot  3.008 21.167  3.061 21.171 /
\plot  3.061 21.171  3.116 21.179 /
\plot  3.116 21.179  3.173 21.188 /
\plot  3.173 21.188  3.228 21.196 /
\plot  3.228 21.196  3.283 21.207 /
\plot  3.283 21.207  3.336 21.220 /
\plot  3.336 21.220  3.389 21.232 /
\plot  3.389 21.232  3.440 21.245 /
\plot  3.440 21.245  3.490 21.260 /
\plot  3.490 21.260  3.539 21.273 /
\plot  3.539 21.273  3.588 21.287 /
\plot  3.588 21.287  3.636 21.304 /
\plot  3.636 21.304  3.683 21.321 /
\plot  3.683 21.321  3.732 21.338 /
\plot  3.732 21.338  3.780 21.357 /
\plot  3.780 21.357  3.829 21.376 /
\plot  3.829 21.376  3.880 21.397 /
\plot  3.880 21.397  3.931 21.419 /
\plot  3.931 21.419  3.981 21.442 /
\plot  3.981 21.442  4.030 21.465 /
\plot  4.030 21.465  4.081 21.491 /
\plot  4.081 21.491  4.130 21.516 /
\plot  4.130 21.516  4.176 21.541 /
\plot  4.176 21.541  4.223 21.567 /
\plot  4.223 21.567  4.267 21.592 /
\plot  4.267 21.592  4.310 21.620 /
\plot  4.310 21.620  4.350 21.645 /
\plot  4.350 21.645  4.388 21.670 /
\plot  4.388 21.670  4.424 21.696 /
\plot  4.424 21.696  4.460 21.723 /
\plot  4.460 21.723  4.494 21.749 /
\plot  4.494 21.749  4.528 21.778 /
\plot  4.528 21.778  4.564 21.808 /
\plot  4.564 21.808  4.597 21.840 /
\plot  4.597 21.840  4.633 21.872 /
\plot  4.633 21.872  4.667 21.905 /
\plot  4.667 21.905  4.701 21.941 /
\plot  4.701 21.941  4.737 21.975 /
\plot  4.737 21.975  4.771 22.011 /
\plot  4.771 22.011  4.805 22.047 /
\plot  4.805 22.047  4.841 22.083 /
\plot  4.841 22.083  4.873 22.119 /
\plot  4.873 22.119  4.906 22.155 /
\plot  4.906 22.155  4.940 22.189 /
\plot  4.940 22.189  4.972 22.221 /
\plot  4.972 22.221  5.002 22.250 /
\plot  5.002 22.250  5.033 22.280 /
\plot  5.033 22.280  5.065 22.310 /
\plot  5.065 22.310  5.097 22.335 /
\plot  5.097 22.335  5.129 22.363 /
\plot  5.129 22.363  5.160 22.388 /
\plot  5.160 22.388  5.196 22.413 /
\plot  5.196 22.413  5.232 22.439 /
\plot  5.232 22.439  5.268 22.464 /
\plot  5.268 22.464  5.309 22.487 /
\plot  5.309 22.487  5.349 22.511 /
\plot  5.349 22.511  5.389 22.534 /
\plot  5.389 22.534  5.431 22.555 /
\plot  5.431 22.555  5.474 22.574 /
\plot  5.474 22.574  5.516 22.593 /
\plot  5.516 22.593  5.556 22.608 /
\plot  5.556 22.608  5.596 22.623 /
\plot  5.596 22.623  5.635 22.636 /
\plot  5.635 22.636  5.673 22.648 /
\plot  5.673 22.648  5.709 22.657 /
\plot  5.709 22.657  5.745 22.663 /
\plot  5.745 22.663  5.779 22.669 /
\plot  5.779 22.669  5.812 22.674 /
\plot  5.812 22.674  5.846 22.676 /
\plot  5.846 22.676  5.880 22.678 /
\plot  5.880 22.678  5.914 22.676 /
\plot  5.914 22.676  5.948 22.674 /
\plot  5.948 22.674  5.982 22.669 /
\plot  5.982 22.669  6.016 22.663 /
\plot  6.016 22.663  6.049 22.657 /
\plot  6.049 22.657  6.083 22.648 /
\plot  6.083 22.648  6.115 22.638 /
\plot  6.115 22.638  6.147 22.625 /
\plot  6.147 22.625  6.176 22.612 /
\plot  6.176 22.612  6.204 22.598 /
\plot  6.204 22.598  6.229 22.581 /
\plot  6.229 22.581  6.255 22.566 /
\plot  6.255 22.566  6.278 22.547 /
\plot  6.278 22.547  6.299 22.530 /
\plot  6.299 22.530  6.318 22.511 /
\plot  6.318 22.511  6.339 22.487 /
\plot  6.339 22.487  6.361 22.464 /
\plot  6.361 22.464  6.380 22.439 /
\plot  6.380 22.439  6.397 22.409 /
\plot  6.397 22.409  6.416 22.380 /
\plot  6.416 22.380  6.430 22.348 /
\plot  6.430 22.348  6.445 22.314 /
\plot  6.445 22.314  6.460 22.280 /
\plot  6.460 22.280  6.471 22.244 /
\plot  6.471 22.244  6.481 22.208 /
\plot  6.481 22.208  6.490 22.172 /
\plot  6.490 22.172  6.498 22.138 /
\plot  6.498 22.138  6.502 22.102 /
\plot  6.502 22.102  6.507 22.068 /
\plot  6.507 22.068  6.509 22.037 /
\putrule from  6.509 22.037 to  6.509 22.003
\putrule from  6.509 22.003 to  6.509 21.969
\plot  6.509 21.969  6.507 21.937 /
\plot  6.507 21.937  6.502 21.901 /
\plot  6.502 21.901  6.498 21.865 /
\plot  6.498 21.865  6.492 21.829 /
\plot  6.492 21.829  6.485 21.791 /
\plot  6.485 21.791  6.475 21.751 /
\plot  6.475 21.751  6.466 21.711 /
\plot  6.466 21.711  6.456 21.668 /
\plot  6.456 21.668  6.443 21.628 /
\plot  6.443 21.628  6.430 21.586 /
\plot  6.430 21.586  6.418 21.546 /
\plot  6.418 21.546  6.405 21.505 /
\plot  6.405 21.505  6.392 21.465 /
\plot  6.392 21.465  6.380 21.425 /
\plot  6.380 21.425  6.367 21.383 /
\plot  6.367 21.383  6.354 21.351 /
\plot  6.354 21.351  6.344 21.317 /
\plot  6.344 21.317  6.331 21.281 /
\plot  6.331 21.281  6.320 21.245 /
\plot  6.320 21.245  6.308 21.205 /
\plot  6.308 21.205  6.293 21.165 /
\plot  6.293 21.165  6.280 21.124 /
\plot  6.280 21.124  6.265 21.080 /
\plot  6.265 21.080  6.251 21.038 /
\plot  6.251 21.038  6.236 20.991 /
\plot  6.236 20.991  6.221 20.947 /
\plot  6.221 20.947  6.208 20.900 /
\plot  6.208 20.900  6.193 20.853 /
\plot  6.193 20.853  6.179 20.809 /
\plot  6.179 20.809  6.164 20.762 /
\plot  6.164 20.762  6.149 20.718 /
\plot  6.149 20.718  6.136 20.673 /
\plot  6.136 20.673  6.121 20.629 /
\plot  6.121 20.629  6.109 20.587 /
\plot  6.109 20.587  6.096 20.542 /
\plot  6.096 20.542  6.083 20.502 /
\plot  6.083 20.502  6.073 20.462 /
\plot  6.073 20.462  6.060 20.419 /
\plot  6.060 20.419  6.047 20.377 /
\plot  6.047 20.377  6.035 20.333 /
\plot  6.035 20.333  6.024 20.288 /
\plot  6.024 20.288  6.011 20.242 /
\plot  6.011 20.242  5.999 20.193 /
\plot  5.999 20.193  5.986 20.144 /
\plot  5.986 20.144  5.975 20.096 /
\plot  5.975 20.096  5.965 20.045 /
\plot  5.965 20.045  5.954 19.996 /
\plot  5.954 19.996  5.944 19.945 /
\plot  5.944 19.945  5.933 19.897 /
\plot  5.933 19.897  5.925 19.846 /
\plot  5.925 19.846  5.918 19.799 /
\plot  5.918 19.799  5.910 19.751 /
\plot  5.910 19.751  5.905 19.704 /
\plot  5.905 19.704  5.899 19.660 /
\plot  5.899 19.660  5.895 19.613 /
\plot  5.895 19.613  5.893 19.571 /
\plot  5.893 19.571  5.891 19.526 /
\plot  5.891 19.526  5.889 19.478 /
\plot  5.889 19.478  5.886 19.429 /
\plot  5.886 19.429  5.889 19.380 /
\putrule from  5.889 19.380 to  5.889 19.329
\plot  5.889 19.329  5.893 19.279 /
\plot  5.893 19.279  5.895 19.228 /
\plot  5.895 19.228  5.901 19.175 /
\plot  5.901 19.175  5.905 19.124 /
\plot  5.905 19.124  5.914 19.071 /
\plot  5.914 19.071  5.920 19.020 /
\plot  5.920 19.020  5.929 18.970 /
\plot  5.929 18.970  5.939 18.921 /
\plot  5.939 18.921  5.950 18.872 /
\plot  5.950 18.872  5.961 18.828 /
\plot  5.961 18.828  5.971 18.785 /
\plot  5.971 18.785  5.984 18.745 /
\plot  5.984 18.745  5.994 18.707 /
\plot  5.994 18.707  6.007 18.671 /
\plot  6.007 18.671  6.020 18.637 /
\plot  6.020 18.637  6.032 18.605 /
\plot  6.032 18.605  6.049 18.570 /
\plot  6.049 18.570  6.066 18.534 /
\plot  6.066 18.534  6.085 18.500 /
\plot  6.085 18.500  6.107 18.468 /
\plot  6.107 18.468  6.128 18.438 /
\plot  6.128 18.438  6.151 18.411 /
\plot  6.151 18.411  6.174 18.383 /
\plot  6.174 18.383  6.200 18.360 /
\plot  6.200 18.360  6.225 18.337 /
\plot  6.225 18.337  6.253 18.318 /
\plot  6.253 18.318  6.280 18.301 /
\plot  6.280 18.301  6.306 18.284 /
\plot  6.306 18.284  6.333 18.271 /
\plot  6.333 18.271  6.361 18.258 /
\plot  6.361 18.258  6.386 18.250 /
\plot  6.386 18.250  6.413 18.239 /
\plot  6.413 18.239  6.445 18.231 /
\plot  6.445 18.231  6.479 18.224 /
\plot  6.479 18.224  6.515 18.218 /
\plot  6.515 18.218  6.555 18.212 /
\plot  6.555 18.212  6.602 18.208 /
\plot  6.602 18.208  6.653 18.203 /
\plot  6.653 18.203  6.710 18.201 /
\plot  6.710 18.201  6.769 18.199 /
\plot  6.769 18.199  6.826 18.197 /
\plot  6.826 18.197  6.879 18.195 /
\plot  6.879 18.195  6.917 18.193 /
\putrule from  6.917 18.193 to  6.941 18.193
\putrule from  6.941 18.193 to  6.951 18.193
\putrule from  6.951 18.193 to  6.953 18.193
}%
%
%
\put{$B_0$
} [lB] at  1.048 22.670
%
%
\put{$B_1$
} [lB] at  7.620 16.574
%
%
\put{$A$
} [lB] at  1.5 21.336
\linethickness=0pt
\putrectangle corners at  0.897 23.199 and  8.437 16.542
\endpicture}

%% file: 2HandleToRound1Handle.tex
\font\thinlinefont=cmr5
\mbox{\beginpicture
\small
\setcoordinatesystem units <0.7cm,0.7cm>
\unitlength=1.04987cm
\linethickness=1pt
\setplotsymbol ({\makebox(0,0)[l]{\tencirc\symbol{'160}}})
\setshadesymbol ({\thinlinefont .})
\setlinear
%
%
\linethickness= 0.500pt
\setplotsymbol ({\thinlinefont .})
\setdots < 0.0953cm>
{\color[rgb]{0,0,0}\circulararc 180.000 degrees from  5.118 20.240 center at  7.413 20.240
}%
%
%
\linethickness= 0.500pt
\setplotsymbol ({\thinlinefont .})
{\color[rgb]{0,0,0}\circulararc 180.106 degrees from 10.914 17.316 center at 13.208 17.314
}%
%
%
\linethickness= 0.500pt
\setplotsymbol ({\thinlinefont .})
{\color[rgb]{0,0,0}\circulararc 180.000 degrees from 10.914 20.240 center at 13.208 20.240
}%
%
%
\linethickness= 0.500pt
\setplotsymbol ({\thinlinefont .})
{\color[rgb]{0,0,0}\circulararc 179.947 degrees from  8.043 14.446 center at 10.339 14.447
}%
%
%
\linethickness= 0.500pt
\setplotsymbol ({\thinlinefont .})
{\color[rgb]{0,0,0}\circulararc 180.106 degrees from  5.118 17.316 center at  7.413 17.314
}%
%
%
\linethickness=1pt
\setplotsymbol ({\makebox(0,0)[l]{\tencirc\symbol{'160}}})
\setsolid
{\color[rgb]{0,0,0}\putrule from  5.118 20.240 to  9.707 20.240
}%
%
%
\linethickness=1pt
\setplotsymbol ({\makebox(0,0)[l]{\tencirc\symbol{'160}}})
{\color[rgb]{0,0,0}\putrule from 10.914 17.316 to 15.502 17.316
}%
%
%
\linethickness= 0.500pt
\setplotsymbol ({\thinlinefont .})
{\color[rgb]{0,0,0}\putrule from  5.692 18.519 to  9.133 18.519
}%
%
%
\linethickness= 0.500pt
\setplotsymbol ({\thinlinefont .})
{\color[rgb]{0,0,0}\putrule from 11.430 18.519 to 14.872 18.519
}%
%
%
\linethickness=1pt
\setplotsymbol ({\makebox(0,0)[l]{\tencirc\symbol{'160}}})
{\color[rgb]{0,0,0}\putrule from 10.914 20.240 to 15.502 20.240
}%
%
%
\linethickness=1pt
\setplotsymbol ({\makebox(0,0)[l]{\tencirc\symbol{'160}}})
{\color[rgb]{0,0,0}\putrule from  5.118 17.316 to  9.707 17.316
}%
%
%
\linethickness=1pt
\setplotsymbol ({\makebox(0,0)[l]{\tencirc\symbol{'160}}})
{\color[rgb]{0,0,0}\putrule from  8.043 14.446 to 12.634 14.446
}%
%
%
\linethickness= 0.500pt
\setplotsymbol ({\thinlinefont .})
\setdots < 0.1429cm>
{\color[rgb]{0,0,0}\circulararc 71.712 degrees from  7.127 16.510 center at  6.839 16.112
}%
%
%
\linethickness= 0.500pt
\setplotsymbol ({\thinlinefont .})
\setsolid
{\color[rgb]{0,0,0}\circulararc 53.169 degrees from 13.551 19.611 center at 13.150 18.809
}%
%
%
\linethickness= 0.500pt
\setplotsymbol ({\thinlinefont .})
\setdots < 0.1429cm>
{\color[rgb]{0,0,0}\circulararc 70.427 degrees from 14.010 19.380 center at 13.723 18.974
}%
%
%
\linethickness= 0.500pt
\setplotsymbol ({\thinlinefont .})
\setsolid
{\color[rgb]{0,0,0}\circulararc 70.427 degrees from 13.437 19.380 center at 13.723 19.787
}%
%
%
\linethickness= 0.500pt
\setplotsymbol ({\thinlinefont .})
\setdots < 0.1429cm>
{\color[rgb]{0,0,0}\circulararc 70.427 degrees from 12.863 19.380 center at 12.576 18.974
}%
%
%
\linethickness= 0.500pt
\setplotsymbol ({\thinlinefont .})
\setsolid
{\color[rgb]{0,0,0}\circulararc 70.427 degrees from 12.289 19.380 center at 12.576 19.787
}%
%
%
\linethickness= 0.500pt
\setplotsymbol ({\thinlinefont .})
{\color[rgb]{0,0,0}\circulararc 70.427 degrees from  7.698 16.510 center at  7.985 16.916
}%
%
%
\linethickness= 0.500pt
\setplotsymbol ({\thinlinefont .})
\setdots < 0.1429cm>
{\color[rgb]{0,0,0}\circulararc 71.960 degrees from  8.272 16.510 center at  7.985 16.115
}%
%
%
\linethickness= 0.500pt
\setplotsymbol ({\thinlinefont .})
\setsolid
{\color[rgb]{0,0,0}\circulararc 70.183 degrees from  6.551 16.510 center at  6.839 16.920
}%
%
%
\linethickness= 0.500pt
\setplotsymbol ({\thinlinefont .})
{\color[rgb]{0,0,0}\circulararc 53.100 degrees from  7.813 16.741 center at  7.412 15.938
}%
%
%
\linethickness= 0.500pt
\setplotsymbol ({\thinlinefont .})
{\color[rgb]{0,0,0}\ellipticalarc axes ratio  0.286:0.286  360 degrees 
	from 14.010 19.380 center at 13.724 19.380
}%
%
%
\linethickness= 0.500pt
\setplotsymbol ({\thinlinefont .})
{\color[rgb]{0,0,0}\ellipticalarc axes ratio  0.286:0.286  360 degrees 
	from  8.272 16.510 center at  7.986 16.510
}%
%
%
\linethickness= 0.500pt
\setplotsymbol ({\thinlinefont .})
{\color[rgb]{0,0,0}\ellipticalarc axes ratio  0.286:0.286  360 degrees 
	from  7.125 16.510 center at  6.839 16.510
}%
%
%
\linethickness= 0.500pt
\setplotsymbol ({\thinlinefont .})
{\color[rgb]{0,0,0}\ellipticalarc axes ratio  0.286:0.286  360 degrees 
	from 12.863 19.380 center at 12.577 19.380
}%
%
%
\linethickness=1pt
\setplotsymbol ({\makebox(0,0)[l]{\tencirc\symbol{'160}}})
{\color[rgb]{0,0,0}\putrule from  9.995 19.380 to 10.569 19.380
%
%
\plot 10.264 19.304 10.569 19.380 10.264 19.456 /
}%
%
%
\linethickness=1pt
\setplotsymbol ({\makebox(0,0)[l]{\tencirc\symbol{'160}}})
{\color[rgb]{0,0,0}\putrule from  9.995 16.510 to 10.569 16.510
%
%
\plot 10.264 16.434 10.569 16.510 10.264 16.586 /
}%
%
%
\linethickness= 0.500pt
\setplotsymbol ({\thinlinefont .})
{\color[rgb]{0,0,0}\plot  9.100 13.153  9.110 13.593 /
%
%
\plot  9.179 13.287  9.110 13.593  9.027 13.290 /
}%
%
%
\linethickness=1pt
\setplotsymbol ({\makebox(0,0)[l]{\tencirc\symbol{'160}}})
{\color[rgb]{0,0,0}\putrule from 15.697 19.380 to 16.271 19.380
%
%
\plot 15.966 19.304 16.271 19.380 15.966 19.456 /
}%
%
%
\linethickness=1pt
\setplotsymbol ({\makebox(0,0)[l]{\tencirc\symbol{'160}}})
{\color[rgb]{0,0,0}\putrule from 15.697 16.552 to 16.271 16.552
%
%
\plot 15.966 16.476 16.271 16.552 15.966 16.629 /
}%
%
%
\linethickness= 0.500pt
\setplotsymbol ({\thinlinefont .})
\setdots < 0.1429cm>
{\color[rgb]{0,0,0}\plot 12.666 16.106 13.841 16.165 /
}%
%
%
\linethickness= 0.500pt
\setplotsymbol ({\thinlinefont .})
\setsolid
{\color[rgb]{0,0,0}\plot  9.711 13.894 10.858 13.881 /
}%
%
%
\linethickness= 0.500pt
\setplotsymbol ({\thinlinefont .})
{\color[rgb]{0,0,0}\plot 14.400 16.222 14.362 16.739 /
%
%
\plot 14.460 16.440 14.362 16.739 14.308 16.429 /
}%
%
%
\linethickness= 0.500pt
\setplotsymbol ({\thinlinefont .})
\setdots < 0.1429cm>
{\color[rgb]{0,0,0}\plot  7.127 16.510  7.698 16.510 /
}%
%
%
\linethickness= 0.500pt
\setplotsymbol ({\thinlinefont .})
{\color[rgb]{0,0,0}\plot 12.863 19.380 13.437 19.380 /
}%
%
%
\linethickness= 0.500pt
\setplotsymbol ({\thinlinefont .})
\setsolid
{\color[rgb]{0,0,0}\plot 12.124 16.040 12.171 16.440 /
%
%
\plot 12.211 16.129 12.171 16.440 12.060 16.146 /
}%
%
%
\linethickness= 0.500pt
\setplotsymbol ({\thinlinefont .})
\setdots < 0.1429cm>
{\color[rgb]{0,0,0}\plot  9.743 13.384 10.924 13.422 /
}%
%
%
\linethickness= 0.500pt
\setplotsymbol ({\thinlinefont .})
\setsolid
{\color[rgb]{0,0,0}\plot 11.481 13.441 11.451 13.938 /
}%
%
%
\linethickness= 0.500pt
\setplotsymbol ({\thinlinefont .})
{\color[rgb]{0,0,0}\plot 11.481 13.506 11.434 13.983 /
%
%
\plot 11.540 13.687 11.434 13.983 11.388 13.672 /
}%
%
%
\linethickness= 0.500pt
\setplotsymbol ({\thinlinefont .})
{\color[rgb]{0,0,0}\putrule from 12.601 16.633 to 13.807 16.633
}%
%
%
\linethickness= 0.500pt
\setplotsymbol ({\thinlinefont .})
{\color[rgb]{0,0,0}\putrule from  6.608 16.341 to  6.604 16.341
\plot  6.604 16.341  6.587 16.339 /
\plot  6.587 16.339  6.562 16.336 /
\plot  6.562 16.336  6.538 16.334 /
\plot  6.538 16.334  6.524 16.328 /
\plot  6.524 16.328  6.513 16.322 /
\plot  6.513 16.322  6.505 16.311 /
\plot  6.505 16.311  6.500 16.296 /
\plot  6.500 16.296  6.496 16.279 /
\plot  6.496 16.279  6.494 16.258 /
\putrule from  6.494 16.258 to  6.494 16.237
\putrule from  6.494 16.237 to  6.494 16.216
\putrule from  6.494 16.216 to  6.494 16.199
\putrule from  6.494 16.199 to  6.494 16.180
\putrule from  6.494 16.180 to  6.494 16.159
\putrule from  6.494 16.159 to  6.494 16.135
\putrule from  6.494 16.135 to  6.494 16.112
\putrule from  6.494 16.112 to  6.494 16.089
\putrule from  6.494 16.089 to  6.494 16.066
\putrule from  6.494 16.066 to  6.494 16.042
\putrule from  6.494 16.042 to  6.494 16.019
\putrule from  6.494 16.019 to  6.494 15.994
\plot  6.494 15.994  6.496 15.966 /
\putrule from  6.496 15.966 to  6.496 15.939
\plot  6.496 15.939  6.498 15.909 /
\plot  6.498 15.909  6.502 15.883 /
\plot  6.502 15.883  6.507 15.858 /
\plot  6.507 15.858  6.513 15.833 /
\plot  6.513 15.833  6.521 15.811 /
\plot  6.521 15.811  6.530 15.788 /
\plot  6.530 15.788  6.543 15.765 /
\plot  6.543 15.765  6.557 15.742 /
\plot  6.557 15.742  6.572 15.723 /
\plot  6.572 15.723  6.589 15.706 /
\plot  6.589 15.706  6.608 15.691 /
\plot  6.608 15.691  6.627 15.680 /
\plot  6.627 15.680  6.644 15.672 /
\plot  6.644 15.672  6.661 15.665 /
\plot  6.661 15.665  6.682 15.661 /
\plot  6.682 15.661  6.703 15.657 /
\plot  6.703 15.657  6.727 15.655 /
\plot  6.727 15.655  6.752 15.653 /
\plot  6.752 15.653  6.778 15.651 /
\putrule from  6.778 15.651 to  6.801 15.651
\putrule from  6.801 15.651 to  6.824 15.651
\putrule from  6.824 15.651 to  6.850 15.651
\putrule from  6.850 15.651 to  6.864 15.651
\putrule from  6.864 15.651 to  6.879 15.651
\putrule from  6.879 15.651 to  6.896 15.651
\putrule from  6.896 15.651 to  6.917 15.651
\putrule from  6.917 15.651 to  6.938 15.651
\putrule from  6.938 15.651 to  6.964 15.651
\putrule from  6.964 15.651 to  6.991 15.651
\putrule from  6.991 15.651 to  7.023 15.651
\putrule from  7.023 15.651 to  7.059 15.651
\putrule from  7.059 15.651 to  7.099 15.651
\putrule from  7.099 15.651 to  7.142 15.651
\putrule from  7.142 15.651 to  7.188 15.651
\putrule from  7.188 15.651 to  7.239 15.651
\putrule from  7.239 15.651 to  7.292 15.651
\putrule from  7.292 15.651 to  7.351 15.651
\putrule from  7.351 15.651 to  7.413 15.651
\putrule from  7.413 15.651 to  7.453 15.651
\putrule from  7.453 15.651 to  7.493 15.651
\putrule from  7.493 15.651 to  7.537 15.651
\putrule from  7.537 15.651 to  7.584 15.651
\putrule from  7.584 15.651 to  7.633 15.651
\putrule from  7.633 15.651 to  7.688 15.651
\putrule from  7.688 15.651 to  7.747 15.651
\putrule from  7.747 15.651 to  7.811 15.651
\putrule from  7.811 15.651 to  7.880 15.651
\putrule from  7.880 15.651 to  7.957 15.651
\putrule from  7.957 15.651 to  8.037 15.651
\putrule from  8.037 15.651 to  8.122 15.651
\putrule from  8.122 15.651 to  8.213 15.651
\putrule from  8.213 15.651 to  8.308 15.651
\putrule from  8.308 15.651 to  8.407 15.651
\putrule from  8.407 15.651 to  8.507 15.651
\putrule from  8.507 15.651 to  8.606 15.651
\putrule from  8.606 15.651 to  8.704 15.651
\putrule from  8.704 15.651 to  8.795 15.651
\putrule from  8.795 15.651 to  8.877 15.651
\putrule from  8.877 15.651 to  8.949 15.651
\putrule from  8.949 15.651 to  9.011 15.651
\putrule from  9.011 15.651 to  9.057 15.651
\putrule from  9.057 15.651 to  9.093 15.651
\putrule from  9.093 15.651 to  9.116 15.651
\putrule from  9.116 15.651 to  9.127 15.651
\putrule from  9.127 15.651 to  9.133 15.651
}%
%
%
\linethickness= 0.500pt
\setplotsymbol ({\thinlinefont .})
{\color[rgb]{0,0,0}\putrule from  5.692 15.651 to  5.694 15.651
\putrule from  5.694 15.651 to  5.709 15.651
\putrule from  5.709 15.651 to  5.736 15.651
\putrule from  5.736 15.651 to  5.770 15.651
\putrule from  5.770 15.651 to  5.802 15.651
\putrule from  5.802 15.651 to  5.831 15.651
\putrule from  5.831 15.651 to  5.857 15.651
\putrule from  5.857 15.651 to  5.882 15.651
\putrule from  5.882 15.651 to  5.903 15.651
\putrule from  5.903 15.651 to  5.927 15.651
\putrule from  5.927 15.651 to  5.952 15.651
\putrule from  5.952 15.651 to  5.977 15.651
\putrule from  5.977 15.651 to  6.005 15.651
\putrule from  6.005 15.651 to  6.030 15.651
\putrule from  6.030 15.651 to  6.054 15.651
\putrule from  6.054 15.651 to  6.077 15.651
\putrule from  6.077 15.651 to  6.096 15.651
\putrule from  6.096 15.651 to  6.113 15.651
\putrule from  6.113 15.651 to  6.126 15.651
\putrule from  6.126 15.651 to  6.136 15.651
\putrule from  6.136 15.651 to  6.147 15.651
\plot  6.147 15.651  6.157 15.653 /
\plot  6.157 15.653  6.166 15.655 /
\plot  6.166 15.655  6.172 15.659 /
\plot  6.172 15.659  6.179 15.663 /
\plot  6.179 15.663  6.185 15.670 /
\plot  6.185 15.670  6.189 15.676 /
\plot  6.189 15.676  6.193 15.684 /
\plot  6.193 15.684  6.195 15.695 /
\plot  6.195 15.695  6.200 15.708 /
\putrule from  6.200 15.708 to  6.200 15.723
\plot  6.200 15.723  6.202 15.740 /
\plot  6.202 15.740  6.204 15.759 /
\plot  6.204 15.759  6.206 15.782 /
\putrule from  6.206 15.782 to  6.206 15.809
\putrule from  6.206 15.809 to  6.206 15.839
\plot  6.206 15.839  6.208 15.871 /
\putrule from  6.208 15.871 to  6.208 15.905
\putrule from  6.208 15.905 to  6.208 15.939
\putrule from  6.208 15.939 to  6.208 15.977
\putrule from  6.208 15.977 to  6.208 16.013
\putrule from  6.208 16.013 to  6.208 16.051
\putrule from  6.208 16.051 to  6.208 16.087
\putrule from  6.208 16.087 to  6.208 16.123
\putrule from  6.208 16.123 to  6.208 16.161
\putrule from  6.208 16.161 to  6.208 16.203
\plot  6.208 16.203  6.210 16.245 /
\plot  6.210 16.245  6.212 16.288 /
\putrule from  6.212 16.288 to  6.212 16.330
\plot  6.212 16.330  6.217 16.372 /
\plot  6.217 16.372  6.219 16.413 /
\plot  6.219 16.413  6.223 16.451 /
\plot  6.223 16.451  6.227 16.485 /
\plot  6.227 16.485  6.234 16.516 /
\plot  6.234 16.516  6.240 16.544 /
\plot  6.240 16.544  6.246 16.569 /
\plot  6.246 16.569  6.255 16.593 /
\plot  6.255 16.593  6.267 16.614 /
\plot  6.267 16.614  6.280 16.631 /
\plot  6.280 16.631  6.293 16.645 /
\plot  6.293 16.645  6.310 16.658 /
\plot  6.310 16.658  6.325 16.667 /
\plot  6.325 16.667  6.342 16.675 /
\plot  6.342 16.675  6.358 16.679 /
\plot  6.358 16.679  6.375 16.681 /
\plot  6.375 16.681  6.390 16.684 /
\putrule from  6.390 16.684 to  6.405 16.684
\putrule from  6.405 16.684 to  6.418 16.684
\putrule from  6.418 16.684 to  6.443 16.684
\putrule from  6.443 16.684 to  6.466 16.684
\putrule from  6.466 16.684 to  6.490 16.684
\putrule from  6.490 16.684 to  6.511 16.684
\putrule from  6.511 16.684 to  6.528 16.684
\putrule from  6.528 16.684 to  6.543 16.684
\putrule from  6.543 16.684 to  6.557 16.684
\putrule from  6.557 16.684 to  6.572 16.684
\putrule from  6.572 16.684 to  6.591 16.684
\putrule from  6.591 16.684 to  6.606 16.684
\putrule from  6.606 16.684 to  6.608 16.684
}%
%
%
\linethickness= 0.500pt
\setplotsymbol ({\thinlinefont .})
{\color[rgb]{0,0,0}\putrule from 10.844 14.372 to 10.844 14.370
\putrule from 10.844 14.370 to 10.844 14.353
\plot 10.844 14.353 10.846 14.323 /
\putrule from 10.846 14.323 to 10.846 14.285
\plot 10.846 14.285 10.848 14.249 /
\plot 10.848 14.249 10.850 14.218 /
\putrule from 10.850 14.218 to 10.850 14.190
\plot 10.850 14.190 10.852 14.163 /
\putrule from 10.852 14.163 to 10.852 14.141
\putrule from 10.852 14.141 to 10.852 14.116
\plot 10.852 14.116 10.854 14.091 /
\plot 10.854 14.091 10.856 14.063 /
\putrule from 10.856 14.063 to 10.856 14.034
\plot 10.856 14.034 10.858 14.002 /
\plot 10.858 14.002 10.861 13.972 /
\putrule from 10.861 13.972 to 10.861 13.942
\plot 10.861 13.942 10.863 13.915 /
\plot 10.863 13.915 10.865 13.887 /
\putrule from 10.865 13.887 to 10.865 13.858
\plot 10.865 13.858 10.867 13.830 /
\plot 10.867 13.830 10.869 13.799 /
\plot 10.869 13.799 10.871 13.767 /
\plot 10.871 13.767 10.873 13.735 /
\plot 10.873 13.735 10.875 13.701 /
\plot 10.875 13.701 10.878 13.672 /
\plot 10.878 13.672 10.882 13.642 /
\plot 10.882 13.642 10.884 13.614 /
\plot 10.884 13.614 10.886 13.587 /
\plot 10.886 13.587 10.888 13.564 /
\plot 10.888 13.564 10.890 13.538 /
\plot 10.890 13.538 10.894 13.513 /
\plot 10.894 13.513 10.897 13.485 /
\plot 10.897 13.485 10.901 13.460 /
\plot 10.901 13.460 10.905 13.437 /
\plot 10.905 13.437 10.909 13.413 /
\plot 10.909 13.413 10.914 13.392 /
\plot 10.914 13.392 10.918 13.373 /
\plot 10.918 13.373 10.922 13.356 /
\plot 10.922 13.356 10.928 13.335 /
\plot 10.928 13.335 10.937 13.314 /
\plot 10.937 13.314 10.943 13.295 /
\plot 10.943 13.295 10.954 13.276 /
\plot 10.954 13.276 10.962 13.257 /
\plot 10.962 13.257 10.973 13.242 /
\plot 10.973 13.242 10.983 13.227 /
\plot 10.983 13.227 10.992 13.216 /
\plot 10.992 13.216 11.007 13.202 /
\plot 11.007 13.202 11.021 13.187 /
\plot 11.021 13.187 11.038 13.174 /
\plot 11.038 13.174 11.055 13.164 /
\plot 11.055 13.164 11.072 13.155 /
\plot 11.072 13.155 11.089 13.149 /
\plot 11.089 13.149 11.104 13.142 /
\plot 11.104 13.142 11.121 13.138 /
\plot 11.121 13.138 11.142 13.134 /
\plot 11.142 13.134 11.165 13.130 /
\plot 11.165 13.130 11.182 13.125 /
\putrule from 11.182 13.125 to 11.184 13.125
}%
%
%
\linethickness= 0.500pt
\setplotsymbol ({\thinlinefont .})
{\color[rgb]{0,0,0}\putrule from  8.215 16.684 to  8.219 16.684
\putrule from  8.219 16.684 to  8.236 16.684
\putrule from  8.236 16.684 to  8.266 16.684
\putrule from  8.266 16.684 to  8.291 16.684
\putrule from  8.291 16.684 to  8.312 16.684
\putrule from  8.312 16.684 to  8.329 16.684
\putrule from  8.329 16.684 to  8.348 16.684
\putrule from  8.348 16.684 to  8.367 16.684
\putrule from  8.367 16.684 to  8.388 16.684
\putrule from  8.388 16.684 to  8.412 16.684
\putrule from  8.412 16.684 to  8.433 16.684
\putrule from  8.433 16.684 to  8.454 16.684
\putrule from  8.454 16.684 to  8.471 16.684
\putrule from  8.471 16.684 to  8.490 16.684
\plot  8.490 16.684  8.511 16.681 /
\plot  8.511 16.681  8.530 16.679 /
\plot  8.530 16.679  8.551 16.675 /
\plot  8.551 16.675  8.572 16.671 /
\plot  8.572 16.671  8.592 16.665 /
\plot  8.592 16.665  8.608 16.654 /
\plot  8.608 16.654  8.621 16.645 /
\plot  8.621 16.645  8.634 16.635 /
\plot  8.634 16.635  8.644 16.622 /
\plot  8.644 16.622  8.657 16.609 /
\plot  8.657 16.609  8.668 16.593 /
\plot  8.668 16.593  8.676 16.576 /
\plot  8.676 16.576  8.685 16.559 /
\plot  8.685 16.559  8.689 16.542 /
\plot  8.689 16.542  8.693 16.527 /
\putrule from  8.693 16.527 to  8.693 16.510
\putrule from  8.693 16.510 to  8.693 16.495
\plot  8.693 16.495  8.689 16.478 /
\plot  8.689 16.478  8.685 16.461 /
\plot  8.685 16.461  8.676 16.447 /
\plot  8.676 16.447  8.668 16.430 /
\plot  8.668 16.430  8.657 16.415 /
\plot  8.657 16.415  8.644 16.400 /
\plot  8.644 16.400  8.634 16.387 /
\plot  8.634 16.387  8.621 16.377 /
\plot  8.621 16.377  8.608 16.368 /
\plot  8.608 16.368  8.592 16.360 /
\plot  8.592 16.360  8.572 16.353 /
\plot  8.572 16.353  8.551 16.349 /
\plot  8.551 16.349  8.530 16.345 /
\plot  8.530 16.345  8.511 16.343 /
\plot  8.511 16.343  8.490 16.341 /
\putrule from  8.490 16.341 to  8.471 16.341
\putrule from  8.471 16.341 to  8.454 16.341
\putrule from  8.454 16.341 to  8.433 16.341
\putrule from  8.433 16.341 to  8.412 16.341
\putrule from  8.412 16.341 to  8.388 16.341
\putrule from  8.388 16.341 to  8.367 16.341
\putrule from  8.367 16.341 to  8.348 16.341
\putrule from  8.348 16.341 to  8.329 16.341
\putrule from  8.329 16.341 to  8.312 16.341
\putrule from  8.312 16.341 to  8.291 16.341
\putrule from  8.291 16.341 to  8.266 16.341
\putrule from  8.266 16.341 to  8.236 16.341
\putrule from  8.236 16.341 to  8.219 16.341
\putrule from  8.219 16.341 to  8.215 16.341
}%
%
%
\linethickness= 0.500pt
\setplotsymbol ({\thinlinefont .})
{\color[rgb]{0,0,0}\putrule from 11.373 15.651 to 11.375 15.651
\putrule from 11.375 15.651 to 11.390 15.651
\putrule from 11.390 15.651 to 11.417 15.651
\putrule from 11.417 15.651 to 11.451 15.651
\putrule from 11.451 15.651 to 11.485 15.651
\putrule from 11.485 15.651 to 11.515 15.651
\putrule from 11.515 15.651 to 11.544 15.651
\putrule from 11.544 15.651 to 11.572 15.651
\putrule from 11.572 15.651 to 11.597 15.651
\putrule from 11.597 15.651 to 11.623 15.651
\putrule from 11.623 15.651 to 11.652 15.651
\putrule from 11.652 15.651 to 11.684 15.651
\putrule from 11.684 15.651 to 11.716 15.651
\putrule from 11.716 15.651 to 11.748 15.651
\putrule from 11.748 15.651 to 11.779 15.651
\putrule from 11.779 15.651 to 11.809 15.651
\putrule from 11.809 15.651 to 11.834 15.651
\putrule from 11.834 15.651 to 11.860 15.651
\putrule from 11.860 15.651 to 11.883 15.651
\putrule from 11.883 15.651 to 11.904 15.651
\plot 11.904 15.651 11.927 15.653 /
\plot 11.927 15.653 11.946 15.655 /
\plot 11.946 15.655 11.968 15.659 /
\plot 11.968 15.659 11.987 15.663 /
\plot 11.987 15.663 12.002 15.672 /
\plot 12.002 15.672 12.016 15.682 /
\plot 12.016 15.682 12.029 15.693 /
\plot 12.029 15.693 12.040 15.708 /
\plot 12.040 15.708 12.048 15.723 /
\plot 12.048 15.723 12.057 15.740 /
\plot 12.057 15.740 12.065 15.759 /
\plot 12.065 15.759 12.071 15.780 /
\plot 12.071 15.780 12.080 15.805 /
\plot 12.080 15.805 12.086 15.833 /
\plot 12.086 15.833 12.093 15.862 /
\plot 12.093 15.862 12.099 15.892 /
\plot 12.099 15.892 12.103 15.922 /
\plot 12.103 15.922 12.109 15.951 /
\plot 12.109 15.951 12.114 15.983 /
\plot 12.114 15.983 12.118 16.013 /
\plot 12.118 16.013 12.122 16.040 /
\plot 12.122 16.040 12.124 16.068 /
\plot 12.124 16.068 12.129 16.095 /
\plot 12.129 16.095 12.133 16.127 /
\plot 12.133 16.127 12.137 16.159 /
\plot 12.137 16.159 12.141 16.190 /
\plot 12.141 16.190 12.145 16.224 /
\plot 12.145 16.224 12.152 16.256 /
\plot 12.152 16.256 12.156 16.288 /
\plot 12.156 16.288 12.160 16.320 /
\plot 12.160 16.320 12.164 16.351 /
\plot 12.164 16.351 12.167 16.379 /
\plot 12.167 16.379 12.171 16.406 /
\plot 12.171 16.406 12.175 16.434 /
\plot 12.175 16.434 12.179 16.463 /
\plot 12.179 16.463 12.184 16.495 /
\plot 12.184 16.495 12.190 16.525 /
\plot 12.190 16.525 12.194 16.557 /
\plot 12.194 16.557 12.200 16.586 /
\plot 12.200 16.586 12.207 16.616 /
\plot 12.207 16.616 12.213 16.643 /
\plot 12.213 16.643 12.220 16.671 /
\plot 12.220 16.671 12.228 16.694 /
\plot 12.228 16.694 12.234 16.715 /
\plot 12.234 16.715 12.243 16.734 /
\plot 12.243 16.734 12.251 16.749 /
\plot 12.251 16.749 12.260 16.764 /
\plot 12.260 16.764 12.270 16.779 /
\plot 12.270 16.779 12.281 16.789 /
\plot 12.281 16.789 12.294 16.800 /
\plot 12.294 16.800 12.308 16.808 /
\plot 12.308 16.808 12.321 16.815 /
\plot 12.321 16.815 12.336 16.821 /
\plot 12.336 16.821 12.353 16.823 /
\plot 12.353 16.823 12.368 16.825 /
\plot 12.368 16.825 12.383 16.823 /
\plot 12.383 16.823 12.397 16.821 /
\plot 12.397 16.821 12.414 16.819 /
\plot 12.414 16.819 12.433 16.813 /
\plot 12.433 16.813 12.452 16.804 /
\plot 12.452 16.804 12.474 16.792 /
\plot 12.474 16.792 12.495 16.777 /
\plot 12.495 16.777 12.516 16.760 /
\plot 12.516 16.760 12.535 16.741 /
\plot 12.535 16.741 12.554 16.717 /
\plot 12.554 16.717 12.571 16.696 /
\plot 12.571 16.696 12.584 16.671 /
\plot 12.584 16.671 12.596 16.645 /
\plot 12.596 16.645 12.605 16.622 /
\plot 12.605 16.622 12.611 16.597 /
\plot 12.611 16.597 12.617 16.571 /
\plot 12.617 16.571 12.622 16.542 /
\plot 12.622 16.542 12.626 16.510 /
\plot 12.626 16.510 12.630 16.476 /
\plot 12.630 16.476 12.634 16.442 /
\plot 12.634 16.442 12.636 16.406 /
\plot 12.636 16.406 12.639 16.372 /
\plot 12.639 16.372 12.641 16.339 /
\plot 12.641 16.339 12.643 16.305 /
\plot 12.643 16.305 12.645 16.273 /
\plot 12.645 16.273 12.647 16.239 /
\plot 12.647 16.239 12.649 16.205 /
\plot 12.649 16.205 12.653 16.169 /
\plot 12.653 16.169 12.656 16.133 /
\plot 12.656 16.133 12.660 16.095 /
\plot 12.660 16.095 12.666 16.057 /
\plot 12.666 16.057 12.670 16.021 /
\plot 12.670 16.021 12.677 15.987 /
\plot 12.677 15.987 12.683 15.953 /
\plot 12.683 15.953 12.689 15.924 /
\plot 12.689 15.924 12.694 15.896 /
\plot 12.694 15.896 12.702 15.871 /
\plot 12.702 15.871 12.708 15.847 /
\plot 12.708 15.847 12.715 15.826 /
\plot 12.715 15.826 12.723 15.805 /
\plot 12.723 15.805 12.734 15.784 /
\plot 12.734 15.784 12.747 15.765 /
\plot 12.747 15.765 12.761 15.748 /
\plot 12.761 15.748 12.776 15.731 /
\plot 12.776 15.731 12.793 15.718 /
\plot 12.793 15.718 12.814 15.706 /
\plot 12.814 15.706 12.835 15.695 /
\plot 12.835 15.695 12.857 15.687 /
\plot 12.857 15.687 12.882 15.680 /
\plot 12.882 15.680 12.903 15.674 /
\plot 12.903 15.674 12.926 15.670 /
\plot 12.926 15.670 12.952 15.668 /
\plot 12.952 15.668 12.979 15.663 /
\plot 12.979 15.663 13.011 15.661 /
\plot 13.011 15.661 13.047 15.659 /
\plot 13.047 15.659 13.085 15.657 /
\plot 13.085 15.657 13.128 15.655 /
\plot 13.128 15.655 13.172 15.653 /
\putrule from 13.172 15.653 to 13.221 15.653
\plot 13.221 15.653 13.269 15.651 /
\putrule from 13.269 15.651 to 13.322 15.651
\putrule from 13.322 15.651 to 13.377 15.651
\putrule from 13.377 15.651 to 13.432 15.651
\putrule from 13.432 15.651 to 13.492 15.651
\putrule from 13.492 15.651 to 13.551 15.651
\putrule from 13.551 15.651 to 13.591 15.651
\putrule from 13.591 15.651 to 13.631 15.651
\putrule from 13.631 15.651 to 13.676 15.651
\putrule from 13.676 15.651 to 13.720 15.651
\putrule from 13.720 15.651 to 13.769 15.651
\putrule from 13.769 15.651 to 13.822 15.651
\putrule from 13.822 15.651 to 13.879 15.651
\putrule from 13.879 15.651 to 13.938 15.651
\putrule from 13.938 15.651 to 14.004 15.651
\putrule from 14.004 15.651 to 14.076 15.651
\putrule from 14.076 15.651 to 14.152 15.651
\putrule from 14.152 15.651 to 14.232 15.651
\putrule from 14.232 15.651 to 14.315 15.651
\putrule from 14.315 15.651 to 14.402 15.651
\putrule from 14.402 15.651 to 14.491 15.651
\putrule from 14.491 15.651 to 14.577 15.651
\putrule from 14.577 15.651 to 14.660 15.651
\putrule from 14.660 15.651 to 14.738 15.651
\putrule from 14.738 15.651 to 14.806 15.651
\putrule from 14.806 15.651 to 14.865 15.651
\putrule from 14.865 15.651 to 14.912 15.651
\putrule from 14.912 15.651 to 14.946 15.651
\putrule from 14.946 15.651 to 14.969 15.651
\putrule from 14.969 15.651 to 14.980 15.651
\putrule from 14.980 15.651 to 14.986 15.651
}%
%
%
\linethickness= 0.500pt
\setplotsymbol ({\thinlinefont .})
{\color[rgb]{0,0,0}\plot 14.300 16.825 14.315 16.811 /
\plot 14.315 16.811 14.328 16.792 /
\plot 14.328 16.792 14.340 16.770 /
\plot 14.340 16.770 14.351 16.749 /
\plot 14.351 16.749 14.362 16.724 /
\plot 14.362 16.724 14.368 16.696 /
\plot 14.368 16.696 14.374 16.669 /
\plot 14.374 16.669 14.381 16.641 /
\plot 14.381 16.641 14.385 16.614 /
\plot 14.385 16.614 14.389 16.588 /
\plot 14.389 16.588 14.393 16.561 /
\plot 14.393 16.561 14.398 16.533 /
\plot 14.398 16.533 14.402 16.504 /
\plot 14.402 16.504 14.404 16.472 /
\plot 14.404 16.472 14.406 16.440 /
\plot 14.406 16.440 14.408 16.408 /
\plot 14.408 16.408 14.410 16.377 /
\plot 14.410 16.377 14.408 16.347 /
\putrule from 14.408 16.347 to 14.408 16.320
\plot 14.408 16.320 14.406 16.290 /
\plot 14.406 16.290 14.402 16.264 /
\plot 14.402 16.264 14.395 16.235 /
\plot 14.395 16.235 14.389 16.207 /
\plot 14.389 16.207 14.381 16.178 /
\plot 14.381 16.178 14.370 16.150 /
\plot 14.370 16.150 14.357 16.123 /
\plot 14.357 16.123 14.343 16.097 /
\plot 14.343 16.097 14.328 16.074 /
\plot 14.328 16.074 14.311 16.055 /
\plot 14.311 16.055 14.296 16.036 /
\plot 14.296 16.036 14.277 16.021 /
\plot 14.277 16.021 14.256 16.008 /
\plot 14.256 16.008 14.232 15.996 /
\plot 14.232 15.996 14.207 15.985 /
\plot 14.207 15.985 14.182 15.977 /
\plot 14.182 15.977 14.154 15.972 /
\plot 14.154 15.972 14.127 15.968 /
\putrule from 14.127 15.968 to 14.101 15.968
\putrule from 14.101 15.968 to 14.076 15.968
\plot 14.076 15.968 14.050 15.972 /
\plot 14.050 15.972 14.031 15.977 /
\plot 14.031 15.977 14.010 15.983 /
\plot 14.010 15.983 13.991 15.991 /
\plot 13.991 15.991 13.970 16.002 /
\plot 13.970 16.002 13.949 16.013 /
\plot 13.949 16.013 13.928 16.027 /
\plot 13.928 16.027 13.909 16.044 /
\plot 13.909 16.044 13.892 16.061 /
\plot 13.892 16.061 13.875 16.080 /
\plot 13.875 16.080 13.860 16.101 /
\plot 13.860 16.101 13.847 16.123 /
\plot 13.847 16.123 13.837 16.144 /
\plot 13.837 16.144 13.826 16.169 /
\plot 13.826 16.169 13.820 16.197 /
\plot 13.820 16.197 13.811 16.224 /
\plot 13.811 16.224 13.805 16.256 /
\plot 13.805 16.256 13.801 16.290 /
\plot 13.801 16.290 13.796 16.324 /
\plot 13.796 16.324 13.794 16.360 /
\plot 13.794 16.360 13.792 16.394 /
\putrule from 13.792 16.394 to 13.792 16.427
\putrule from 13.792 16.427 to 13.792 16.459
\plot 13.792 16.459 13.794 16.491 /
\plot 13.794 16.491 13.796 16.518 /
\plot 13.796 16.518 13.799 16.552 /
\plot 13.799 16.552 13.803 16.586 /
\plot 13.803 16.586 13.807 16.618 /
\plot 13.807 16.618 13.815 16.650 /
\plot 13.815 16.650 13.822 16.681 /
\plot 13.822 16.681 13.830 16.711 /
\plot 13.830 16.711 13.841 16.736 /
\plot 13.841 16.736 13.851 16.760 /
\plot 13.851 16.760 13.862 16.779 /
\plot 13.862 16.779 13.873 16.796 /
\plot 13.873 16.796 13.885 16.813 /
\plot 13.885 16.813 13.898 16.825 /
\plot 13.898 16.825 13.915 16.838 /
\plot 13.915 16.838 13.932 16.849 /
\plot 13.932 16.849 13.951 16.859 /
\plot 13.951 16.859 13.972 16.868 /
\plot 13.972 16.868 13.991 16.874 /
\plot 13.991 16.874 14.014 16.878 /
\plot 14.014 16.878 14.036 16.880 /
\plot 14.036 16.880 14.057 16.883 /
\putrule from 14.057 16.883 to 14.080 16.883
\putrule from 14.080 16.883 to 14.103 16.883
\putrule from 14.103 16.883 to 14.131 16.883
\plot 14.131 16.883 14.158 16.878 /
\plot 14.158 16.878 14.186 16.874 /
\plot 14.186 16.874 14.213 16.868 /
\plot 14.213 16.868 14.239 16.859 /
\plot 14.239 16.859 14.262 16.849 /
\plot 14.262 16.849 14.281 16.838 /
\plot 14.281 16.838 14.300 16.825 /
}%
%
%
\linethickness= 0.500pt
\setplotsymbol ({\thinlinefont .})
{\color[rgb]{0,0,0}\putrule from  9.686 14.368 to  9.686 14.366
\putrule from  9.686 14.366 to  9.686 14.355
\plot  9.686 14.355  9.688 14.334 /
\plot  9.688 14.334  9.690 14.307 /
\putrule from  9.690 14.307 to  9.690 14.281
\plot  9.690 14.281  9.692 14.258 /
\plot  9.692 14.258  9.694 14.235 /
\putrule from  9.694 14.235 to  9.694 14.211
\plot  9.694 14.211  9.696 14.190 /
\plot  9.696 14.190  9.699 14.169 /
\plot  9.699 14.169  9.701 14.144 /
\plot  9.701 14.144  9.703 14.118 /
\plot  9.703 14.118  9.705 14.091 /
\putrule from  9.705 14.091 to  9.705 14.063
\plot  9.705 14.063  9.707 14.036 /
\plot  9.707 14.036  9.709 14.010 /
\plot  9.709 14.010  9.711 13.985 /
\plot  9.711 13.985  9.713 13.964 /
\putrule from  9.713 13.964 to  9.713 13.936
\plot  9.713 13.936  9.716 13.911 /
\plot  9.716 13.911  9.718 13.883 /
\putrule from  9.718 13.883 to  9.718 13.856
\plot  9.718 13.856  9.720 13.830 /
\putrule from  9.720 13.830 to  9.720 13.805
\plot  9.720 13.805  9.722 13.784 /
\plot  9.722 13.784  9.724 13.760 /
\putrule from  9.724 13.760 to  9.724 13.739
\plot  9.724 13.739  9.726 13.716 /
\putrule from  9.726 13.716 to  9.726 13.693
\plot  9.726 13.693  9.728 13.665 /
\plot  9.728 13.665  9.730 13.638 /
\plot  9.730 13.638  9.732 13.610 /
\plot  9.732 13.610  9.735 13.583 /
\plot  9.735 13.583  9.737 13.555 /
\putrule from  9.737 13.555 to  9.737 13.532
\plot  9.737 13.532  9.739 13.509 /
\plot  9.739 13.509  9.741 13.481 /
\plot  9.741 13.481  9.743 13.451 /
\plot  9.743 13.451  9.745 13.422 /
\plot  9.745 13.422  9.749 13.390 /
\plot  9.749 13.390  9.751 13.360 /
\plot  9.751 13.360  9.754 13.331 /
\plot  9.754 13.331  9.758 13.301 /
\plot  9.758 13.301  9.760 13.274 /
\plot  9.760 13.274  9.764 13.246 /
\plot  9.764 13.246  9.766 13.219 /
\plot  9.766 13.219  9.771 13.189 /
\plot  9.771 13.189  9.777 13.159 /
\plot  9.777 13.159  9.781 13.130 /
\plot  9.781 13.130  9.787 13.100 /
\plot  9.787 13.100  9.794 13.073 /
\plot  9.794 13.073  9.802 13.047 /
\plot  9.802 13.047  9.809 13.024 /
\plot  9.809 13.024  9.817 13.003 /
\plot  9.817 13.003  9.823 12.986 /
\plot  9.823 12.986  9.832 12.969 /
\plot  9.832 12.969  9.840 12.954 /
\plot  9.840 12.954  9.853 12.939 /
\plot  9.853 12.939  9.864 12.924 /
\plot  9.864 12.924  9.878 12.910 /
\plot  9.878 12.910  9.893 12.897 /
\plot  9.893 12.897  9.910 12.884 /
\plot  9.910 12.884  9.929 12.874 /
\plot  9.929 12.874  9.948 12.863 /
\plot  9.948 12.863  9.970 12.855 /
\plot  9.970 12.855  9.991 12.848 /
\plot  9.991 12.848 10.012 12.842 /
\plot 10.012 12.842 10.033 12.838 /
\plot 10.033 12.838 10.058 12.831 /
\plot 10.058 12.831 10.086 12.827 /
\plot 10.086 12.827 10.118 12.823 /
\plot 10.118 12.823 10.152 12.819 /
\plot 10.152 12.819 10.190 12.816 /
\plot 10.190 12.816 10.230 12.812 /
\plot 10.230 12.812 10.270 12.810 /
\plot 10.270 12.810 10.315 12.808 /
\plot 10.315 12.808 10.359 12.806 /
\plot 10.359 12.806 10.406 12.804 /
\plot 10.406 12.804 10.452 12.802 /
\plot 10.452 12.802 10.503 12.799 /
\putrule from 10.503 12.799 to 10.539 12.799
\putrule from 10.539 12.799 to 10.577 12.799
\plot 10.577 12.799 10.615 12.797 /
\putrule from 10.615 12.797 to 10.657 12.797
\putrule from 10.657 12.797 to 10.702 12.797
\plot 10.702 12.797 10.746 12.795 /
\putrule from 10.746 12.795 to 10.795 12.795
\putrule from 10.795 12.795 to 10.844 12.795
\putrule from 10.844 12.795 to 10.892 12.795
\plot 10.892 12.795 10.943 12.793 /
\putrule from 10.943 12.793 to 10.994 12.793
\putrule from 10.994 12.793 to 11.045 12.793
\putrule from 11.045 12.793 to 11.093 12.793
\putrule from 11.093 12.793 to 11.142 12.793
\putrule from 11.142 12.793 to 11.189 12.793
\putrule from 11.189 12.793 to 11.235 12.793
\putrule from 11.235 12.793 to 11.278 12.793
\plot 11.278 12.793 11.320 12.791 /
\putrule from 11.320 12.791 to 11.360 12.791
\putrule from 11.360 12.791 to 11.400 12.791
\putrule from 11.400 12.791 to 11.447 12.791
\putrule from 11.447 12.791 to 11.494 12.791
\putrule from 11.494 12.791 to 11.540 12.791
\putrule from 11.540 12.791 to 11.587 12.791
\putrule from 11.587 12.791 to 11.637 12.791
\putrule from 11.637 12.791 to 11.688 12.791
\putrule from 11.688 12.791 to 11.745 12.791
\putrule from 11.745 12.791 to 11.805 12.791
\putrule from 11.805 12.791 to 11.864 12.791
\putrule from 11.864 12.791 to 11.921 12.791
\putrule from 11.921 12.791 to 11.976 12.791
\putrule from 11.976 12.791 to 12.018 12.791
\putrule from 12.018 12.791 to 12.050 12.791
\putrule from 12.050 12.791 to 12.069 12.791
\putrule from 12.069 12.791 to 12.076 12.791
\putrule from 12.076 12.791 to 12.078 12.791
}%
%
%
\put{$K$
} [lB] at  9.248 18.462
%
%
\put{$B$
} [lB] at  7.127 20.299
%
%
\put{$-1$
} [lB] at  7.239 16.912
%
%
\put{$-1$
} [lB] at 10.080 13.983
%
%
\put{$\Sigma$
} [lB] at  6.945 19.306
%
%
\put{$-1$
} [lB] at 13.034 16.722
%
%
\put{$-1$
} [lB] at 12.8 19.8
%
%
\linethickness=1pt
\setplotsymbol ({\makebox(0,0)[l]{\tencirc\symbol{'160}}})
{\color[rgb]{1,1,1}\put{\makebox(0,0)[l]{\circle*{ 0.195}}} at  9.167 14.446
}%
%
%
\linethickness=1pt
\setplotsymbol ({\makebox(0,0)[l]{\tencirc\symbol{'160}}})
{\color[rgb]{1,1,1}\put{\makebox(0,0)[l]{\circle*{ 0.144}}} at 11.343 14.448
}%
%
%
\linethickness=3pt
\setplotsymbol ({\makebox(0,0)[l]{\tencirc\symbol{'162}}})
\color[rgb]{0,0,0}\plot 12.666 16.106 12.666 16.106 /
%
%
\linethickness=3pt
\setplotsymbol ({\makebox(0,0)[l]{\tencirc\symbol{'162}}})
\color[rgb]{0,0,0}\plot  9.743 13.373  9.743 13.373 /
%
%
\linethickness=3pt
\setplotsymbol ({\makebox(0,0)[l]{\tencirc\symbol{'162}}})
\color[rgb]{0,0,0}\plot 10.918 13.418 10.918 13.418 /
%
%
\linethickness=3pt
\setplotsymbol ({\makebox(0,0)[l]{\tencirc\symbol{'162}}})
\color[rgb]{0,0,0}\plot 13.841 16.165 13.841 16.165 /
%
%
\put{$U$
} [lB] at 10.846 14.772
%
%
\put{$K'$
} [lB] at  9.174 14.802
%
%
\put{$K'$
} [lB] at 12.137 16.876
%
%
\put{$U$
} [lB] at 13.915 16.944
%
%
\linethickness= 0.500pt
\setplotsymbol ({\thinlinefont .})
{\color[rgb]{0,0,0}\putrule from  8.484 12.791 to  8.488 12.791
\putrule from  8.488 12.791 to  8.507 12.791
\plot  8.507 12.791  8.539 12.793 /
\putrule from  8.539 12.793 to  8.566 12.793
\putrule from  8.566 12.793 to  8.592 12.793
\plot  8.592 12.793  8.613 12.795 /
\putrule from  8.613 12.795 to  8.630 12.795
\putrule from  8.630 12.795 to  8.647 12.795
\putrule from  8.647 12.795 to  8.666 12.795
\putrule from  8.666 12.795 to  8.687 12.795
\putrule from  8.687 12.795 to  8.708 12.795
\putrule from  8.708 12.795 to  8.729 12.795
\putrule from  8.729 12.795 to  8.748 12.795
\putrule from  8.748 12.795 to  8.767 12.795
\putrule from  8.767 12.795 to  8.786 12.795
\plot  8.786 12.795  8.807 12.793 /
\putrule from  8.807 12.793 to  8.829 12.793
\putrule from  8.829 12.793 to  8.850 12.793
\putrule from  8.850 12.793 to  8.871 12.793
\putrule from  8.871 12.793 to  8.892 12.793
\putrule from  8.892 12.793 to  8.911 12.793
\putrule from  8.911 12.793 to  8.928 12.793
\putrule from  8.928 12.793 to  8.945 12.793
\plot  8.945 12.793  8.960 12.795 /
\plot  8.960 12.795  8.977 12.797 /
\plot  8.977 12.797  8.994 12.802 /
\plot  8.994 12.802  9.009 12.810 /
\plot  9.009 12.810  9.021 12.821 /
\plot  9.021 12.821  9.034 12.833 /
\plot  9.034 12.833  9.045 12.848 /
\plot  9.045 12.848  9.051 12.865 /
\plot  9.051 12.865  9.059 12.884 /
\plot  9.059 12.884  9.066 12.905 /
\plot  9.066 12.905  9.072 12.933 /
\plot  9.072 12.933  9.078 12.965 /
\plot  9.078 12.965  9.083 12.996 /
\plot  9.083 12.996  9.089 13.032 /
\plot  9.089 13.032  9.093 13.068 /
\plot  9.093 13.068  9.095 13.106 /
\plot  9.095 13.106  9.100 13.144 /
\putrule from  9.100 13.144 to  9.100 13.174
\plot  9.100 13.174  9.102 13.206 /
\plot  9.102 13.206  9.104 13.238 /
\putrule from  9.104 13.238 to  9.104 13.271
\plot  9.104 13.271  9.106 13.310 /
\plot  9.106 13.310  9.108 13.348 /
\putrule from  9.108 13.348 to  9.108 13.388
\plot  9.108 13.388  9.110 13.428 /
\putrule from  9.110 13.428 to  9.110 13.468
\plot  9.110 13.468  9.112 13.509 /
\putrule from  9.112 13.509 to  9.112 13.547
\plot  9.112 13.547  9.114 13.585 /
\putrule from  9.114 13.585 to  9.114 13.621
\plot  9.114 13.621  9.116 13.657 /
\putrule from  9.116 13.657 to  9.116 13.691
\plot  9.116 13.691  9.119 13.727 /
\putrule from  9.119 13.727 to  9.119 13.765
\plot  9.119 13.765  9.121 13.803 /
\plot  9.121 13.803  9.123 13.841 /
\plot  9.123 13.841  9.125 13.879 /
\putrule from  9.125 13.879 to  9.125 13.919
\plot  9.125 13.919  9.127 13.957 /
\plot  9.127 13.957  9.129 13.995 /
\plot  9.129 13.995  9.131 14.034 /
\plot  9.131 14.034  9.133 14.067 /
\plot  9.133 14.067  9.136 14.101 /
\plot  9.136 14.101  9.138 14.133 /
\plot  9.138 14.133  9.142 14.165 /
\plot  9.142 14.165  9.144 14.201 /
\plot  9.144 14.201  9.146 14.237 /
\plot  9.146 14.237  9.150 14.273 /
\plot  9.150 14.273  9.155 14.309 /
\plot  9.155 14.309  9.159 14.345 /
\plot  9.159 14.345  9.165 14.381 /
\plot  9.165 14.381  9.169 14.415 /
\plot  9.169 14.415  9.178 14.448 /
\plot  9.178 14.448  9.184 14.476 /
\plot  9.184 14.476  9.193 14.503 /
\plot  9.193 14.503  9.201 14.527 /
\plot  9.201 14.527  9.210 14.548 /
\plot  9.210 14.548  9.220 14.569 /
\plot  9.220 14.569  9.233 14.586 /
\plot  9.233 14.586  9.246 14.603 /
\plot  9.246 14.603  9.260 14.620 /
\plot  9.260 14.620  9.275 14.635 /
\plot  9.275 14.635  9.294 14.647 /
\plot  9.294 14.647  9.311 14.658 /
\plot  9.311 14.658  9.328 14.669 /
\plot  9.328 14.669  9.347 14.675 /
\plot  9.347 14.675  9.364 14.681 /
\plot  9.364 14.681  9.381 14.685 /
\plot  9.381 14.685  9.398 14.690 /
\plot  9.398 14.690  9.426 14.692 /
\plot  9.426 14.692  9.451 14.694 /
\plot  9.451 14.694  9.478 14.692 /
\plot  9.478 14.692  9.508 14.690 /
\plot  9.508 14.690  9.533 14.685 /
\plot  9.533 14.685  9.557 14.679 /
\plot  9.557 14.679  9.578 14.673 /
\plot  9.578 14.673  9.595 14.666 /
\plot  9.595 14.666  9.610 14.660 /
\plot  9.610 14.660  9.622 14.654 /
\plot  9.622 14.654  9.633 14.645 /
\plot  9.633 14.645  9.644 14.635 /
\plot  9.644 14.635  9.650 14.626 /
\plot  9.650 14.626  9.656 14.618 /
\plot  9.656 14.618  9.660 14.607 /
\plot  9.660 14.607  9.663 14.599 /
\plot  9.663 14.599  9.667 14.586 /
\plot  9.667 14.586  9.669 14.569 /
\plot  9.669 14.569  9.673 14.548 /
\plot  9.673 14.548  9.675 14.527 /
\plot  9.675 14.527  9.677 14.510 /
\putrule from  9.677 14.510 to  9.677 14.508
}%
%
%
\linethickness= 0.500pt
\setplotsymbol ({\thinlinefont .})
{\color[rgb]{0,0,0}\putrule from 11.032 13.168 to 11.034 13.168
\plot 11.034 13.168 11.047 13.161 /
\plot 11.047 13.161 11.072 13.153 /
\plot 11.072 13.153 11.104 13.144 /
\plot 11.104 13.144 11.134 13.136 /
\plot 11.134 13.136 11.157 13.130 /
\plot 11.157 13.130 11.180 13.125 /
\plot 11.180 13.125 11.199 13.123 /
\putrule from 11.199 13.123 to 11.220 13.123
\plot 11.220 13.123 11.242 13.125 /
\plot 11.242 13.125 11.265 13.128 /
\plot 11.265 13.128 11.288 13.134 /
\plot 11.288 13.134 11.311 13.142 /
\plot 11.311 13.142 11.335 13.153 /
\plot 11.335 13.153 11.354 13.166 /
\plot 11.354 13.166 11.373 13.178 /
\plot 11.373 13.178 11.386 13.193 /
\plot 11.386 13.193 11.398 13.208 /
\plot 11.398 13.208 11.411 13.227 /
\plot 11.411 13.227 11.424 13.246 /
\plot 11.424 13.246 11.436 13.269 /
\plot 11.436 13.269 11.447 13.295 /
\plot 11.447 13.295 11.455 13.320 /
\plot 11.455 13.320 11.462 13.348 /
\plot 11.462 13.348 11.468 13.373 /
\plot 11.468 13.373 11.474 13.401 /
\plot 11.474 13.401 11.477 13.424 /
\plot 11.477 13.424 11.479 13.449 /
\plot 11.479 13.449 11.481 13.477 /
\putrule from 11.481 13.477 to 11.481 13.506
\putrule from 11.481 13.506 to 11.481 13.538
\putrule from 11.481 13.538 to 11.481 13.572
\plot 11.481 13.572 11.479 13.606 /
\plot 11.479 13.606 11.477 13.640 /
\plot 11.477 13.640 11.474 13.674 /
\plot 11.474 13.674 11.470 13.708 /
\plot 11.470 13.708 11.468 13.741 /
\plot 11.468 13.741 11.464 13.775 /
\plot 11.464 13.775 11.462 13.805 /
\plot 11.462 13.805 11.458 13.835 /
\plot 11.458 13.835 11.453 13.866 /
\plot 11.453 13.866 11.449 13.898 /
\plot 11.449 13.898 11.445 13.932 /
\plot 11.445 13.932 11.441 13.968 /
\plot 11.441 13.968 11.434 14.004 /
\plot 11.434 14.004 11.430 14.038 /
\plot 11.430 14.038 11.424 14.074 /
\plot 11.424 14.074 11.417 14.108 /
\plot 11.417 14.108 11.411 14.139 /
\plot 11.411 14.139 11.407 14.171 /
\plot 11.407 14.171 11.400 14.201 /
\plot 11.400 14.201 11.394 14.230 /
\plot 11.394 14.230 11.388 14.264 /
\plot 11.388 14.264 11.381 14.298 /
\plot 11.381 14.298 11.373 14.330 /
\plot 11.373 14.330 11.362 14.364 /
\plot 11.362 14.364 11.354 14.398 /
\plot 11.354 14.398 11.343 14.431 /
\plot 11.343 14.431 11.333 14.461 /
\plot 11.333 14.461 11.320 14.491 /
\plot 11.320 14.491 11.307 14.516 /
\plot 11.307 14.516 11.295 14.541 /
\plot 11.295 14.541 11.282 14.563 /
\plot 11.282 14.563 11.267 14.580 /
\plot 11.267 14.580 11.250 14.601 /
\plot 11.250 14.601 11.229 14.618 /
\plot 11.229 14.618 11.208 14.635 /
\plot 11.208 14.635 11.184 14.649 /
\plot 11.184 14.649 11.159 14.660 /
\plot 11.159 14.660 11.134 14.671 /
\plot 11.134 14.671 11.108 14.677 /
\plot 11.108 14.677 11.083 14.681 /
\plot 11.083 14.681 11.062 14.683 /
\putrule from 11.062 14.683 to 11.041 14.683
\putrule from 11.041 14.683 to 11.015 14.683
\plot 11.015 14.683 10.990 14.679 /
\plot 10.990 14.679 10.964 14.671 /
\plot 10.964 14.671 10.941 14.662 /
\plot 10.941 14.662 10.920 14.652 /
\plot 10.920 14.652 10.901 14.641 /
\plot 10.901 14.641 10.886 14.630 /
\plot 10.886 14.630 10.873 14.618 /
\plot 10.873 14.618 10.863 14.601 /
\plot 10.863 14.601 10.854 14.582 /
\plot 10.854 14.582 10.848 14.556 /
\plot 10.848 14.556 10.846 14.529 /
\plot 10.846 14.529 10.844 14.512 /
\putrule from 10.844 14.512 to 10.844 14.508
}%
\linethickness=0pt
\putrectangle corners at  5.072 20.451 and 16.317 12.133
\endpicture}

%% file: 1R2Lhandles.tex
\font\thinlinefont=cmr5
\mbox{\beginpicture
\small
\setcoordinatesystem units <1.04987cm,1.04987cm>
\unitlength=1.04987cm
\linethickness=1pt
\setplotsymbol ({\makebox(0,0)[l]{\tencirc\symbol{'160}}})
\setshadesymbol ({\thinlinefont .})
\setlinear
%
%
\linethickness= 0.500pt
\setplotsymbol ({\thinlinefont .})
{\color[rgb]{0,0,0}\circulararc 106.260 degrees from  5.524 23.812 center at  5.906 24.098
}%
%
%
\linethickness= 0.500pt
\setplotsymbol ({\thinlinefont .})
\setdots < 0.0953cm>
{\color[rgb]{0,0,0}\circulararc 106.260 degrees from  6.287 23.812 center at  5.906 23.527
}%
%
%
\linethickness= 0.500pt
\setplotsymbol ({\thinlinefont .})
\setsolid
{\color[rgb]{0,0,0}\ellipticalarc axes ratio  0.381:0.381  360 degrees 
	from  6.287 23.812 center at  5.905 23.812
}%
%
%
\linethickness= 0.500pt
\setplotsymbol ({\thinlinefont .})
{\color[rgb]{0,0,0}\circulararc 106.260 degrees from  5.524 21.717 center at  5.906 22.003
}%
%
%
\linethickness= 0.500pt
\setplotsymbol ({\thinlinefont .})
\setdots < 0.0953cm>
{\color[rgb]{0,0,0}\circulararc 106.260 degrees from  6.287 21.717 center at  5.906 21.431
}%
%
%
\linethickness= 0.500pt
\setplotsymbol ({\thinlinefont .})
\setsolid
{\color[rgb]{0,0,0}\ellipticalarc axes ratio  0.381:0.381  360 degrees 
	from  6.287 21.717 center at  5.905 21.717
}%
%
%
\linethickness= 0.500pt
\setplotsymbol ({\thinlinefont .})
{\color[rgb]{0,0,0}\circulararc 106.260 degrees from 11.430 21.717 center at 11.811 22.003
}%
%
%
\linethickness= 0.500pt
\setplotsymbol ({\thinlinefont .})
\setdots < 0.0953cm>
{\color[rgb]{0,0,0}\circulararc 106.260 degrees from 12.192 21.717 center at 11.811 21.431
}%
%
%
\linethickness= 0.500pt
\setplotsymbol ({\thinlinefont .})
\setsolid
{\color[rgb]{0,0,0}\ellipticalarc axes ratio  0.381:0.381  360 degrees 
	from 12.192 21.717 center at 11.811 21.717
}%
%
%
\linethickness= 0.500pt
\setplotsymbol ({\thinlinefont .})
{\color[rgb]{0,0,0}\circulararc 106.260 degrees from 11.430 23.812 center at 11.811 24.098
}%
%
%
\linethickness= 0.500pt
\setplotsymbol ({\thinlinefont .})
\setdots < 0.0953cm>
{\color[rgb]{0,0,0}\circulararc 106.260 degrees from 12.192 23.812 center at 11.811 23.527
}%
%
%
\linethickness= 0.500pt
\setplotsymbol ({\thinlinefont .})
\setsolid
{\color[rgb]{0,0,0}\ellipticalarc axes ratio  0.381:0.381  360 degrees 
	from 12.192 23.812 center at 11.811 23.812
}%
%
%
\linethickness= 0.500pt
\setplotsymbol ({\thinlinefont .})
{\color[rgb]{0,0,0}\circulararc 106.260 degrees from  6.858 23.241 center at  7.239 23.527
}%
%
%
\linethickness= 0.500pt
\setplotsymbol ({\thinlinefont .})
\setdots < 0.0953cm>
{\color[rgb]{0,0,0}\circulararc 106.260 degrees from  7.620 23.241 center at  7.239 22.955
}%
%
%
\linethickness= 0.500pt
\setplotsymbol ({\thinlinefont .})
\setsolid
{\color[rgb]{0,0,0}\ellipticalarc axes ratio  0.381:0.381  360 degrees 
	from  7.620 23.241 center at  7.239 23.241
}%
%
%
\linethickness= 0.500pt
\setplotsymbol ({\thinlinefont .})
{\color[rgb]{0,0,0}\circulararc 106.260 degrees from  8.001 23.241 center at  8.382 23.527
}%
%
%
\linethickness= 0.500pt
\setplotsymbol ({\thinlinefont .})
\setdots < 0.0953cm>
{\color[rgb]{0,0,0}\circulararc 106.260 degrees from  8.763 23.241 center at  8.382 22.955
}%
%
%
\linethickness= 0.500pt
\setplotsymbol ({\thinlinefont .})
\setsolid
{\color[rgb]{0,0,0}\ellipticalarc axes ratio  0.381:0.381  360 degrees 
	from  8.763 23.241 center at  8.382 23.241
}%
%
%
\linethickness=5pt
\setplotsymbol ({\makebox(0,0)[l]{\tencirc\symbol{'164}}})
{\color[rgb]{0,0,0}\plot 12.954 23.241 12.954 23.241 /
%
%
\linethickness=5pt
\setplotsymbol ({\makebox(0,0)[l]{\tencirc\symbol{'164}}})
\color[rgb]{0,0,0}\plot 14.287 23.241 14.287 23.241 /
%
%
\linethickness= 0.500pt
\setplotsymbol ({\thinlinefont .})
\setdashes < 0.1270cm>
\color[rgb]{0,0,0}\plot 12.954 23.241 14.287 23.241 /
}%
%
%
\linethickness= 0.500pt
\setplotsymbol ({\thinlinefont .})
\setsolid
{\color[rgb]{0,0,0}\putrule from 12.954 22.860 to 14.287 22.860
}%
%
%
\linethickness= 1pt
\setplotsymbol ({\thinlinefont .})
{\color[rgb]{0,0,0}\putrule from  6.287 23.812 to  6.608 23.812
}%
%
%
\linethickness= 0.500pt
\setplotsymbol ({\thinlinefont .})
\setdots < 0.0953cm>
{\color[rgb]{0,0,0}\plot  9.525 23.812  9.525 21.717 /
}%
%
%
\linethickness= 0.500pt
\setplotsymbol ({\thinlinefont .})
{\color[rgb]{0,0,0}\plot 11.049 23.812 11.049 21.717 /
}%
%
%
\linethickness= 0.500pt
\setplotsymbol ({\thinlinefont .})
{\color[rgb]{0,0,0}\plot 15.240 23.812 15.240 21.717 /
}%
%
%
\linethickness= 0.500pt
\setplotsymbol ({\thinlinefont .})
\setsolid
{\color[rgb]{0,0,0}\putrule from 12.954 23.050 to 12.954 23.812
%
%
\plot 13.018 23.559 12.954 23.812 12.890 23.559 /
}%
%
%
\linethickness= 0.500pt
\setplotsymbol ({\thinlinefont .})
{\color[rgb]{0,0,0}\putrule from 14.287 23.050 to 14.287 23.812
%
%
\plot 14.351 23.559 14.287 23.812 14.224 23.559 /
}%
%
%
\linethickness= 0.500pt
\setplotsymbol ({\thinlinefont .})
\setdots < 0.0953cm>
{\color[rgb]{0,0,0}\plot  5.144 23.812  5.144 21.717 /
}%
%
%
\linethickness=1pt
\setplotsymbol ({\makebox(0,0)[l]{\tencirc\symbol{'160}}})
\setsolid
{\color[rgb]{0,0,0}\putrule from  6.731 23.812 to  8.907 23.812
}%
%
%
\linethickness=1pt
\setplotsymbol ({\makebox(0,0)[l]{\tencirc\symbol{'160}}})
{\color[rgb]{0,0,0}\putrule from  9.017 23.812 to  9.525 23.812
}%
%
%
\linethickness=1pt
\setplotsymbol ({\makebox(0,0)[l]{\tencirc\symbol{'160}}})
{\color[rgb]{0,0,0}\putrule from  5.524 23.812 to  5.144 23.812
}%
%
%
\linethickness=1pt
\setplotsymbol ({\makebox(0,0)[l]{\tencirc\symbol{'160}}})
{\color[rgb]{0,0,0}\putrule from  5.524 21.717 to  5.144 21.717
}%
%
%
\linethickness=1pt
\setplotsymbol ({\makebox(0,0)[l]{\tencirc\symbol{'160}}})
{\color[rgb]{0,0,0}\putrule from  6.287 21.717 to  9.525 21.717
}%
%
%
\linethickness=1pt
\setplotsymbol ({\makebox(0,0)[l]{\tencirc\symbol{'160}}})
{\color[rgb]{0,0,0}\putrule from 11.049 23.812 to 11.430 23.812
}%
%
%
\linethickness=1pt
\setplotsymbol ({\makebox(0,0)[l]{\tencirc\symbol{'160}}})
{\color[rgb]{0,0,0}\putrule from 12.192 23.812 to 12.895 23.812
}%
%
%
\linethickness=1pt
\setplotsymbol ({\makebox(0,0)[l]{\tencirc\symbol{'160}}})
{\color[rgb]{0,0,0}\putrule from 13.030 23.812 to 14.228 23.812
}%
%
%
\linethickness=1pt
\setplotsymbol ({\makebox(0,0)[l]{\tencirc\symbol{'160}}})
{\color[rgb]{0,0,0}\putrule from 14.340 23.812 to 15.240 23.812
}%
%
%
\linethickness=1pt
\setplotsymbol ({\makebox(0,0)[l]{\tencirc\symbol{'160}}})
{\color[rgb]{0,0,0}\putrule from 15.240 21.717 to 12.192 21.717
}%
%
%
\linethickness=1pt
\setplotsymbol ({\makebox(0,0)[l]{\tencirc\symbol{'160}}})
{\color[rgb]{0,0,0}\putrule from 11.049 21.717 to 11.430 21.717
}%
%
%
\linethickness= 0.500pt
\setplotsymbol ({\thinlinefont .})
{\color[rgb]{0,0,0}\putrule from  5.905 23.431 to  5.905 23.429
\putrule from  5.905 23.429 to  5.905 23.419
\plot  5.905 23.419  5.908 23.398 /
\plot  5.908 23.398  5.912 23.362 /
\plot  5.912 23.362  5.916 23.319 /
\plot  5.916 23.319  5.922 23.277 /
\plot  5.922 23.277  5.929 23.237 /
\plot  5.929 23.237  5.939 23.201 /
\plot  5.939 23.201  5.950 23.173 /
\plot  5.950 23.173  5.965 23.150 /
\plot  5.965 23.150  5.982 23.131 /
\plot  5.982 23.131  6.001 23.114 /
\plot  6.001 23.114  6.020 23.103 /
\plot  6.020 23.103  6.043 23.093 /
\plot  6.043 23.093  6.068 23.084 /
\plot  6.068 23.084  6.096 23.078 /
\plot  6.096 23.078  6.128 23.072 /
\plot  6.128 23.072  6.160 23.070 /
\plot  6.160 23.070  6.191 23.067 /
\plot  6.191 23.067  6.223 23.070 /
\plot  6.223 23.070  6.255 23.072 /
\plot  6.255 23.072  6.287 23.078 /
\plot  6.287 23.078  6.314 23.084 /
\plot  6.314 23.084  6.339 23.093 /
\plot  6.339 23.093  6.363 23.103 /
\plot  6.363 23.103  6.382 23.114 /
\plot  6.382 23.114  6.399 23.127 /
\plot  6.399 23.127  6.413 23.144 /
\plot  6.413 23.144  6.428 23.161 /
\plot  6.428 23.161  6.439 23.180 /
\plot  6.439 23.180  6.449 23.201 /
\plot  6.449 23.201  6.458 23.224 /
\plot  6.458 23.224  6.464 23.249 /
\plot  6.464 23.249  6.469 23.275 /
\plot  6.469 23.275  6.473 23.300 /
\plot  6.473 23.300  6.475 23.326 /
\plot  6.475 23.326  6.477 23.351 /
\putrule from  6.477 23.351 to  6.477 23.374
\putrule from  6.477 23.374 to  6.477 23.400
\putrule from  6.477 23.400 to  6.477 23.423
\putrule from  6.477 23.423 to  6.477 23.448
\putrule from  6.477 23.448 to  6.477 23.476
\putrule from  6.477 23.476 to  6.477 23.506
\putrule from  6.477 23.506 to  6.477 23.537
\putrule from  6.477 23.537 to  6.477 23.573
\putrule from  6.477 23.573 to  6.477 23.614
\putrule from  6.477 23.614 to  6.477 23.656
\putrule from  6.477 23.656 to  6.477 23.694
\putrule from  6.477 23.694 to  6.477 23.726
\putrule from  6.477 23.726 to  6.477 23.747
\putrule from  6.477 23.747 to  6.477 23.755
\putrule from  6.477 23.755 to  6.477 23.757
}%
%
%
\linethickness= 0.500pt
\setplotsymbol ({\thinlinefont .})
{\color[rgb]{0,0,0}\putrule from  6.477 23.868 to  6.477 23.870
\putrule from  6.477 23.870 to  6.477 23.884
\plot  6.477 23.884  6.479 23.914 /
\putrule from  6.479 23.914 to  6.479 23.956
\plot  6.479 23.956  6.483 24.001 /
\plot  6.483 24.001  6.485 24.041 /
\plot  6.485 24.041  6.490 24.075 /
\plot  6.490 24.075  6.494 24.102 /
\plot  6.494 24.102  6.500 24.122 /
\plot  6.500 24.122  6.509 24.138 /
\plot  6.509 24.138  6.513 24.147 /
\plot  6.513 24.147  6.519 24.153 /
\plot  6.519 24.153  6.524 24.158 /
\plot  6.524 24.158  6.530 24.164 /
\plot  6.530 24.164  6.536 24.168 /
\plot  6.536 24.168  6.543 24.170 /
\plot  6.543 24.170  6.549 24.172 /
\plot  6.549 24.172  6.557 24.174 /
\putrule from  6.557 24.174 to  6.564 24.174
\plot  6.564 24.174  6.572 24.172 /
\plot  6.572 24.172  6.581 24.170 /
\plot  6.581 24.170  6.587 24.166 /
\plot  6.587 24.166  6.596 24.162 /
\plot  6.596 24.162  6.602 24.158 /
\plot  6.602 24.158  6.608 24.149 /
\plot  6.608 24.149  6.615 24.141 /
\plot  6.615 24.141  6.621 24.132 /
\plot  6.621 24.132  6.625 24.122 /
\plot  6.625 24.122  6.632 24.111 /
\plot  6.632 24.111  6.636 24.098 /
\plot  6.636 24.098  6.642 24.075 /
\plot  6.642 24.075  6.648 24.047 /
\plot  6.648 24.047  6.653 24.016 /
\plot  6.653 24.016  6.657 23.982 /
\plot  6.657 23.982  6.661 23.944 /
\plot  6.661 23.944  6.663 23.904 /
\plot  6.663 23.904  6.665 23.863 /
\putrule from  6.665 23.863 to  6.665 23.823
\plot  6.665 23.823  6.668 23.785 /
\putrule from  6.668 23.785 to  6.668 23.749
\putrule from  6.668 23.749 to  6.668 23.715
\putrule from  6.668 23.715 to  6.668 23.686
\putrule from  6.668 23.686 to  6.668 23.652
\putrule from  6.668 23.652 to  6.668 23.620
\plot  6.668 23.620  6.670 23.592 /
\putrule from  6.670 23.592 to  6.670 23.565
\plot  6.670 23.565  6.674 23.539 /
\plot  6.674 23.539  6.676 23.518 /
\plot  6.676 23.518  6.682 23.499 /
\plot  6.682 23.499  6.689 23.484 /
\plot  6.689 23.484  6.697 23.474 /
\plot  6.697 23.474  6.708 23.463 /
\plot  6.708 23.463  6.723 23.453 /
\plot  6.723 23.453  6.744 23.446 /
\plot  6.744 23.446  6.771 23.442 /
\plot  6.771 23.442  6.807 23.438 /
\plot  6.807 23.438  6.850 23.434 /
\plot  6.850 23.434  6.883 23.431 /
\putrule from  6.883 23.431 to  6.900 23.431
\putrule from  6.900 23.431 to  6.905 23.431
}%
%
%
\linethickness= 0.500pt
\setplotsymbol ({\thinlinefont .})
{\color[rgb]{0,0,0}\putrule from  6.909 23.050 to  6.905 23.050
\putrule from  6.905 23.050 to  6.888 23.050
\plot  6.888 23.050  6.852 23.048 /
\plot  6.852 23.048  6.811 23.044 /
\plot  6.811 23.044  6.773 23.040 /
\plot  6.773 23.040  6.746 23.036 /
\plot  6.746 23.036  6.723 23.029 /
\plot  6.723 23.029  6.708 23.019 /
\plot  6.708 23.019  6.695 23.006 /
\plot  6.695 23.006  6.687 22.991 /
\plot  6.687 22.991  6.678 22.970 /
\plot  6.678 22.970  6.674 22.947 /
\plot  6.674 22.947  6.670 22.919 /
\plot  6.670 22.919  6.668 22.890 /
\putrule from  6.668 22.890 to  6.668 22.860
\putrule from  6.668 22.860 to  6.668 22.828
\putrule from  6.668 22.828 to  6.668 22.809
\putrule from  6.668 22.809 to  6.668 22.788
\plot  6.668 22.788  6.665 22.767 /
\plot  6.665 22.767  6.663 22.746 /
\plot  6.663 22.746  6.661 22.722 /
\plot  6.661 22.722  6.657 22.697 /
\plot  6.657 22.697  6.651 22.674 /
\plot  6.651 22.674  6.642 22.650 /
\plot  6.642 22.650  6.632 22.629 /
\plot  6.632 22.629  6.617 22.608 /
\plot  6.617 22.608  6.602 22.589 /
\plot  6.602 22.589  6.585 22.572 /
\plot  6.585 22.572  6.564 22.555 /
\plot  6.564 22.555  6.540 22.543 /
\plot  6.540 22.543  6.513 22.530 /
\plot  6.513 22.530  6.483 22.521 /
\plot  6.483 22.521  6.449 22.511 /
\plot  6.449 22.511  6.411 22.502 /
\plot  6.411 22.502  6.371 22.496 /
\plot  6.371 22.496  6.329 22.490 /
\plot  6.329 22.490  6.287 22.483 /
\plot  6.287 22.483  6.244 22.479 /
\plot  6.244 22.479  6.202 22.473 /
\plot  6.202 22.473  6.162 22.468 /
\plot  6.162 22.468  6.124 22.464 /
\plot  6.124 22.464  6.090 22.458 /
\plot  6.090 22.458  6.060 22.454 /
\plot  6.060 22.454  6.032 22.447 /
\plot  6.032 22.447  6.009 22.441 /
\plot  6.009 22.441  5.988 22.432 /
\plot  5.988 22.432  5.971 22.424 /
\plot  5.971 22.424  5.956 22.413 /
\plot  5.956 22.413  5.941 22.403 /
\plot  5.941 22.403  5.931 22.392 /
\plot  5.931 22.392  5.922 22.380 /
\plot  5.922 22.380  5.916 22.367 /
\plot  5.916 22.367  5.912 22.354 /
\plot  5.912 22.354  5.910 22.339 /
\plot  5.910 22.339  5.908 22.327 /
\plot  5.908 22.327  5.905 22.314 /
\putrule from  5.905 22.314 to  5.905 22.301
\putrule from  5.905 22.301 to  5.905 22.288
\putrule from  5.905 22.288 to  5.905 22.267
\putrule from  5.905 22.267 to  5.905 22.242
\putrule from  5.905 22.242 to  5.905 22.214
\putrule from  5.905 22.214 to  5.905 22.181
\putrule from  5.905 22.181 to  5.905 22.145
\putrule from  5.905 22.145 to  5.905 22.117
\putrule from  5.905 22.117 to  5.905 22.100
\putrule from  5.905 22.100 to  5.905 22.098
}%
%
%
\linethickness= 0.500pt
\setplotsymbol ({\thinlinefont .})
{\color[rgb]{0,0,0}\putrule from  8.727 23.431 to  8.729 23.431
\putrule from  8.729 23.431 to  8.748 23.431
\plot  8.748 23.431  8.780 23.434 /
\plot  8.780 23.434  8.820 23.438 /
\plot  8.820 23.438  8.854 23.442 /
\plot  8.854 23.442  8.882 23.446 /
\plot  8.882 23.446  8.901 23.453 /
\plot  8.901 23.453  8.915 23.463 /
\plot  8.915 23.463  8.926 23.474 /
\plot  8.926 23.474  8.932 23.484 /
\plot  8.932 23.484  8.939 23.499 /
\plot  8.939 23.499  8.945 23.518 /
\plot  8.945 23.518  8.949 23.539 /
\plot  8.949 23.539  8.951 23.565 /
\putrule from  8.951 23.565 to  8.951 23.592
\plot  8.951 23.592  8.954 23.620 /
\putrule from  8.954 23.620 to  8.954 23.652
\putrule from  8.954 23.652 to  8.954 23.686
\putrule from  8.954 23.686 to  8.954 23.715
\putrule from  8.954 23.715 to  8.954 23.749
\putrule from  8.954 23.749 to  8.954 23.785
\plot  8.954 23.785  8.956 23.823 /
\putrule from  8.956 23.823 to  8.956 23.863
\plot  8.956 23.863  8.958 23.904 /
\plot  8.958 23.904  8.960 23.944 /
\plot  8.960 23.944  8.964 23.982 /
\plot  8.964 23.982  8.968 24.016 /
\plot  8.968 24.016  8.973 24.047 /
\plot  8.973 24.047  8.979 24.075 /
\plot  8.979 24.075  8.985 24.098 /
\plot  8.985 24.098  8.989 24.111 /
\plot  8.989 24.111  8.996 24.122 /
\plot  8.996 24.122  9.000 24.132 /
\plot  9.000 24.132  9.006 24.141 /
\plot  9.006 24.141  9.013 24.149 /
\plot  9.013 24.149  9.019 24.158 /
\plot  9.019 24.158  9.025 24.162 /
\plot  9.025 24.162  9.034 24.168 /
\plot  9.034 24.168  9.040 24.170 /
\plot  9.040 24.170  9.049 24.172 /
\plot  9.049 24.172  9.057 24.174 /
\putrule from  9.057 24.174 to  9.064 24.174
\plot  9.064 24.174  9.072 24.172 /
\putrule from  9.072 24.172 to  9.078 24.172
\plot  9.078 24.172  9.085 24.168 /
\plot  9.085 24.168  9.091 24.164 /
\plot  9.091 24.164  9.097 24.160 /
\plot  9.097 24.160  9.102 24.153 /
\plot  9.102 24.153  9.108 24.149 /
\plot  9.108 24.149  9.112 24.141 /
\plot  9.112 24.141  9.121 24.126 /
\plot  9.121 24.126  9.127 24.105 /
\plot  9.127 24.105  9.131 24.079 /
\plot  9.131 24.079  9.136 24.047 /
\plot  9.136 24.047  9.138 24.007 /
\plot  9.138 24.007  9.142 23.965 /
\putrule from  9.142 23.965 to  9.142 23.925
\plot  9.142 23.925  9.144 23.895 /
\putrule from  9.144 23.895 to  9.144 23.882
\putrule from  9.144 23.882 to  9.144 23.880
}%
%
%
\linethickness= 0.500pt
\setplotsymbol ({\thinlinefont .})
{\color[rgb]{0,0,0}\putrule from  9.144 23.751 to  9.144 23.749
\putrule from  9.144 23.749 to  9.144 23.741
\putrule from  9.144 23.741 to  9.144 23.719
\putrule from  9.144 23.719 to  9.144 23.688
\putrule from  9.144 23.688 to  9.144 23.650
\putrule from  9.144 23.650 to  9.144 23.609
\putrule from  9.144 23.609 to  9.144 23.569
\putrule from  9.144 23.569 to  9.144 23.533
\putrule from  9.144 23.533 to  9.144 23.499
\putrule from  9.144 23.499 to  9.144 23.470
\putrule from  9.144 23.470 to  9.144 23.440
\putrule from  9.144 23.440 to  9.144 23.412
\putrule from  9.144 23.412 to  9.144 23.383
\putrule from  9.144 23.383 to  9.144 23.351
\plot  9.144 23.351  9.142 23.319 /
\putrule from  9.142 23.319 to  9.142 23.288
\plot  9.142 23.288  9.138 23.254 /
\plot  9.138 23.254  9.133 23.222 /
\plot  9.133 23.222  9.129 23.190 /
\plot  9.129 23.190  9.121 23.161 /
\plot  9.121 23.161  9.112 23.135 /
\plot  9.112 23.135  9.102 23.112 /
\plot  9.102 23.112  9.089 23.093 /
\plot  9.089 23.093  9.074 23.076 /
\plot  9.074 23.076  9.055 23.063 /
\plot  9.055 23.063  9.034 23.053 /
\plot  9.034 23.053  9.009 23.046 /
\plot  9.009 23.046  8.975 23.042 /
\plot  8.975 23.042  8.937 23.040 /
\putrule from  8.937 23.040 to  8.892 23.040
\plot  8.892 23.040  8.843 23.042 /
\plot  8.843 23.042  8.797 23.044 /
\plot  8.797 23.044  8.757 23.046 /
\plot  8.757 23.046  8.731 23.048 /
\plot  8.731 23.048  8.721 23.050 /
\putrule from  8.721 23.050 to  8.719 23.050
}%
%
%
\linethickness= 0.500pt
\setplotsymbol ({\thinlinefont .})
{\color[rgb]{0,0,0}\putrule from  7.239 22.860 to  7.239 22.858
\plot  7.239 22.858  7.241 22.849 /
\plot  7.241 22.849  7.243 22.830 /
\plot  7.243 22.830  7.247 22.803 /
\plot  7.247 22.803  7.256 22.771 /
\plot  7.256 22.771  7.264 22.741 /
\plot  7.264 22.741  7.277 22.714 /
\plot  7.277 22.714  7.292 22.689 /
\plot  7.292 22.689  7.311 22.663 /
\plot  7.311 22.663  7.334 22.638 /
\plot  7.334 22.638  7.353 22.619 /
\plot  7.353 22.619  7.377 22.598 /
\plot  7.377 22.598  7.404 22.578 /
\plot  7.404 22.578  7.434 22.555 /
\plot  7.434 22.555  7.465 22.534 /
\plot  7.465 22.534  7.501 22.513 /
\plot  7.501 22.513  7.537 22.492 /
\plot  7.537 22.492  7.578 22.473 /
\plot  7.578 22.473  7.616 22.456 /
\plot  7.616 22.456  7.656 22.443 /
\plot  7.656 22.443  7.696 22.430 /
\plot  7.696 22.430  7.734 22.422 /
\plot  7.734 22.422  7.772 22.418 /
\plot  7.772 22.418  7.811 22.416 /
\plot  7.811 22.416  7.849 22.418 /
\plot  7.849 22.418  7.887 22.422 /
\plot  7.887 22.422  7.925 22.430 /
\plot  7.925 22.430  7.965 22.443 /
\plot  7.965 22.443  8.005 22.456 /
\plot  8.005 22.456  8.043 22.473 /
\plot  8.043 22.473  8.084 22.492 /
\plot  8.084 22.492  8.120 22.513 /
\plot  8.120 22.513  8.156 22.534 /
\plot  8.156 22.534  8.187 22.555 /
\plot  8.187 22.555  8.217 22.578 /
\plot  8.217 22.578  8.244 22.598 /
\plot  8.244 22.598  8.268 22.619 /
\plot  8.268 22.619  8.287 22.638 /
\plot  8.287 22.638  8.310 22.663 /
\plot  8.310 22.663  8.329 22.689 /
\plot  8.329 22.689  8.344 22.714 /
\plot  8.344 22.714  8.357 22.741 /
\plot  8.357 22.741  8.365 22.771 /
\plot  8.365 22.771  8.374 22.803 /
\plot  8.374 22.803  8.378 22.830 /
\plot  8.378 22.830  8.380 22.849 /
\plot  8.380 22.849  8.382 22.858 /
\putrule from  8.382 22.858 to  8.382 22.860
}%
%
%
\linethickness= 0.500pt
\setplotsymbol ({\thinlinefont .})
{\color[rgb]{0,0,0}\putrule from 11.811 23.431 to 11.811 23.429
\putrule from 11.811 23.429 to 11.811 23.419
\plot 11.811 23.419 11.813 23.398 /
\plot 11.813 23.398 11.817 23.362 /
\plot 11.817 23.362 11.824 23.319 /
\plot 11.824 23.319 11.832 23.277 /
\plot 11.832 23.277 11.843 23.237 /
\plot 11.843 23.237 11.855 23.201 /
\plot 11.855 23.201 11.870 23.173 /
\plot 11.870 23.173 11.889 23.150 /
\plot 11.889 23.150 11.910 23.131 /
\plot 11.910 23.131 11.938 23.114 /
\plot 11.938 23.114 11.966 23.103 /
\plot 11.966 23.103 11.995 23.093 /
\plot 11.995 23.093 12.029 23.084 /
\plot 12.029 23.084 12.067 23.078 /
\plot 12.067 23.078 12.107 23.072 /
\plot 12.107 23.072 12.150 23.070 /
\plot 12.150 23.070 12.192 23.067 /
\plot 12.192 23.067 12.234 23.070 /
\plot 12.234 23.070 12.277 23.072 /
\plot 12.277 23.072 12.317 23.078 /
\plot 12.317 23.078 12.355 23.084 /
\plot 12.355 23.084 12.389 23.093 /
\plot 12.389 23.093 12.418 23.103 /
\plot 12.418 23.103 12.446 23.114 /
\plot 12.446 23.114 12.469 23.127 /
\plot 12.469 23.127 12.490 23.144 /
\plot 12.490 23.144 12.507 23.161 /
\plot 12.507 23.161 12.522 23.180 /
\plot 12.522 23.180 12.537 23.201 /
\plot 12.537 23.201 12.548 23.224 /
\plot 12.548 23.224 12.556 23.249 /
\plot 12.556 23.249 12.562 23.273 /
\plot 12.562 23.273 12.567 23.300 /
\plot 12.567 23.300 12.569 23.326 /
\plot 12.569 23.326 12.571 23.351 /
\plot 12.571 23.351 12.573 23.374 /
\putrule from 12.573 23.374 to 12.573 23.398
\putrule from 12.573 23.398 to 12.573 23.421
\putrule from 12.573 23.421 to 12.573 23.448
\putrule from 12.573 23.448 to 12.573 23.476
\putrule from 12.573 23.476 to 12.573 23.506
\putrule from 12.573 23.506 to 12.573 23.537
\putrule from 12.573 23.537 to 12.573 23.571
\putrule from 12.573 23.571 to 12.573 23.611
\putrule from 12.573 23.611 to 12.573 23.652
\putrule from 12.573 23.652 to 12.573 23.692
\putrule from 12.573 23.692 to 12.573 23.721
\putrule from 12.573 23.721 to 12.573 23.743
\putrule from 12.573 23.743 to 12.573 23.751
\putrule from 12.573 23.751 to 12.573 23.753
}%
%
%
\linethickness= 0.500pt
\setplotsymbol ({\thinlinefont .})
{\color[rgb]{0,0,0}\putrule from 12.573 23.870 to 12.573 23.872
\putrule from 12.573 23.872 to 12.573 23.887
\plot 12.573 23.887 12.575 23.916 /
\plot 12.575 23.916 12.579 23.956 /
\plot 12.579 23.956 12.584 24.001 /
\plot 12.584 24.001 12.590 24.041 /
\plot 12.590 24.041 12.598 24.075 /
\plot 12.598 24.075 12.609 24.102 /
\plot 12.609 24.102 12.622 24.124 /
\plot 12.622 24.124 12.636 24.138 /
\plot 12.636 24.138 12.649 24.149 /
\plot 12.649 24.149 12.664 24.158 /
\plot 12.664 24.158 12.683 24.164 /
\plot 12.683 24.164 12.700 24.170 /
\plot 12.700 24.170 12.721 24.172 /
\plot 12.721 24.172 12.742 24.174 /
\plot 12.742 24.174 12.764 24.172 /
\plot 12.764 24.172 12.785 24.168 /
\plot 12.785 24.168 12.806 24.164 /
\plot 12.806 24.164 12.827 24.155 /
\plot 12.827 24.155 12.844 24.145 /
\plot 12.844 24.145 12.863 24.130 /
\plot 12.863 24.130 12.878 24.115 /
\plot 12.878 24.115 12.890 24.098 /
\plot 12.890 24.098 12.903 24.079 /
\plot 12.903 24.079 12.912 24.056 /
\plot 12.912 24.056 12.922 24.028 /
\plot 12.922 24.028 12.929 23.999 /
\plot 12.929 23.999 12.935 23.965 /
\plot 12.935 23.965 12.941 23.929 /
\plot 12.941 23.929 12.946 23.891 /
\plot 12.946 23.891 12.948 23.848 /
\plot 12.948 23.848 12.952 23.806 /
\putrule from 12.952 23.806 to 12.952 23.764
\plot 12.952 23.764 12.954 23.719 /
\putrule from 12.954 23.719 to 12.954 23.677
\putrule from 12.954 23.677 to 12.954 23.633
\putrule from 12.954 23.633 to 12.954 23.590
\putrule from 12.954 23.590 to 12.954 23.552
\putrule from 12.954 23.552 to 12.954 23.512
\putrule from 12.954 23.512 to 12.954 23.470
\putrule from 12.954 23.470 to 12.954 23.427
\putrule from 12.954 23.427 to 12.954 23.383
\putrule from 12.954 23.383 to 12.954 23.338
\putrule from 12.954 23.338 to 12.954 23.292
\putrule from 12.954 23.292 to 12.954 23.245
\putrule from 12.954 23.245 to 12.954 23.199
\putrule from 12.954 23.199 to 12.954 23.154
\putrule from 12.954 23.154 to 12.954 23.110
\putrule from 12.954 23.110 to 12.954 23.070
\putrule from 12.954 23.070 to 12.954 23.029
\putrule from 12.954 23.029 to 12.954 22.991
\putrule from 12.954 22.991 to 12.954 22.957
\putrule from 12.954 22.957 to 12.954 22.924
\putrule from 12.954 22.924 to 12.954 22.894
\putrule from 12.954 22.894 to 12.954 22.866
\putrule from 12.954 22.866 to 12.954 22.839
\plot 12.954 22.839 12.952 22.813 /
\plot 12.952 22.813 12.950 22.786 /
\plot 12.950 22.786 12.946 22.761 /
\plot 12.946 22.761 12.941 22.735 /
\plot 12.941 22.735 12.935 22.710 /
\plot 12.935 22.710 12.929 22.686 /
\plot 12.929 22.686 12.918 22.665 /
\plot 12.918 22.665 12.905 22.644 /
\plot 12.905 22.644 12.893 22.625 /
\plot 12.893 22.625 12.876 22.608 /
\plot 12.876 22.608 12.859 22.591 /
\plot 12.859 22.591 12.838 22.576 /
\plot 12.838 22.576 12.816 22.564 /
\plot 12.816 22.564 12.791 22.553 /
\plot 12.791 22.553 12.764 22.543 /
\plot 12.764 22.543 12.732 22.534 /
\plot 12.732 22.534 12.698 22.526 /
\plot 12.698 22.526 12.662 22.517 /
\plot 12.662 22.517 12.622 22.511 /
\plot 12.622 22.511 12.577 22.504 /
\plot 12.577 22.504 12.531 22.498 /
\plot 12.531 22.498 12.482 22.494 /
\plot 12.482 22.494 12.433 22.487 /
\plot 12.433 22.487 12.383 22.483 /
\plot 12.383 22.483 12.332 22.479 /
\plot 12.332 22.479 12.283 22.475 /
\plot 12.283 22.475 12.234 22.471 /
\plot 12.234 22.471 12.188 22.468 /
\plot 12.188 22.468 12.143 22.464 /
\plot 12.143 22.464 12.103 22.460 /
\plot 12.103 22.460 12.067 22.456 /
\plot 12.067 22.456 12.033 22.451 /
\plot 12.033 22.451 12.002 22.447 /
\plot 12.002 22.447 11.970 22.441 /
\plot 11.970 22.441 11.942 22.435 /
\plot 11.942 22.435 11.919 22.428 /
\plot 11.919 22.428 11.898 22.420 /
\plot 11.898 22.420 11.879 22.411 /
\plot 11.879 22.411 11.862 22.401 /
\plot 11.862 22.401 11.849 22.390 /
\plot 11.849 22.390 11.836 22.380 /
\plot 11.836 22.380 11.828 22.369 /
\plot 11.828 22.369 11.822 22.356 /
\plot 11.822 22.356 11.817 22.346 /
\plot 11.817 22.346 11.815 22.333 /
\plot 11.815 22.333 11.813 22.322 /
\plot 11.813 22.322 11.811 22.310 /
\putrule from 11.811 22.310 to 11.811 22.299
\putrule from 11.811 22.299 to 11.811 22.288
\putrule from 11.811 22.288 to 11.811 22.267
\putrule from 11.811 22.267 to 11.811 22.242
\putrule from 11.811 22.242 to 11.811 22.214
\putrule from 11.811 22.214 to 11.811 22.181
\putrule from 11.811 22.181 to 11.811 22.145
\putrule from 11.811 22.145 to 11.811 22.117
\putrule from 11.811 22.117 to 11.811 22.100
\putrule from 11.811 22.100 to 11.811 22.098
}%
%
%
\linethickness= 0.500pt
\setplotsymbol ({\thinlinefont .})
{\color[rgb]{0,0,0}\putrule from 14.859 23.753 to 14.859 23.751
\putrule from 14.859 23.751 to 14.859 23.743
\putrule from 14.859 23.743 to 14.859 23.721
\putrule from 14.859 23.721 to 14.859 23.692
\putrule from 14.859 23.692 to 14.859 23.652
\putrule from 14.859 23.652 to 14.859 23.611
\putrule from 14.859 23.611 to 14.859 23.571
\putrule from 14.859 23.571 to 14.859 23.537
\putrule from 14.859 23.537 to 14.859 23.506
\putrule from 14.859 23.506 to 14.859 23.476
\putrule from 14.859 23.476 to 14.859 23.448
\putrule from 14.859 23.448 to 14.859 23.421
\putrule from 14.859 23.421 to 14.859 23.393
\putrule from 14.859 23.393 to 14.859 23.366
\putrule from 14.859 23.366 to 14.859 23.336
\putrule from 14.859 23.336 to 14.859 23.305
\putrule from 14.859 23.305 to 14.859 23.271
\putrule from 14.859 23.271 to 14.859 23.235
\putrule from 14.859 23.235 to 14.859 23.199
\putrule from 14.859 23.199 to 14.859 23.163
\putrule from 14.859 23.163 to 14.859 23.127
\putrule from 14.859 23.127 to 14.859 23.091
\putrule from 14.859 23.091 to 14.859 23.055
\putrule from 14.859 23.055 to 14.859 23.019
\putrule from 14.859 23.019 to 14.859 22.987
\putrule from 14.859 22.987 to 14.859 22.953
\putrule from 14.859 22.953 to 14.859 22.919
\plot 14.859 22.919 14.857 22.883 /
\plot 14.857 22.883 14.855 22.847 /
\plot 14.855 22.847 14.851 22.811 /
\plot 14.851 22.811 14.846 22.773 /
\plot 14.846 22.773 14.840 22.737 /
\plot 14.840 22.737 14.831 22.703 /
\plot 14.831 22.703 14.821 22.672 /
\plot 14.821 22.672 14.810 22.644 /
\plot 14.810 22.644 14.796 22.617 /
\plot 14.796 22.617 14.781 22.593 /
\plot 14.781 22.593 14.764 22.574 /
\plot 14.764 22.574 14.745 22.557 /
\plot 14.745 22.557 14.721 22.543 /
\plot 14.721 22.543 14.696 22.528 /
\plot 14.696 22.528 14.668 22.517 /
\plot 14.668 22.517 14.637 22.511 /
\plot 14.637 22.511 14.605 22.504 /
\plot 14.605 22.504 14.573 22.500 /
\putrule from 14.573 22.500 to 14.542 22.500
\plot 14.542 22.500 14.510 22.502 /
\plot 14.510 22.502 14.478 22.507 /
\plot 14.478 22.507 14.450 22.513 /
\plot 14.450 22.513 14.425 22.521 /
\plot 14.425 22.521 14.402 22.532 /
\plot 14.402 22.532 14.383 22.543 /
\plot 14.383 22.543 14.366 22.555 /
\plot 14.366 22.555 14.351 22.572 /
\plot 14.351 22.572 14.336 22.589 /
\plot 14.336 22.589 14.326 22.608 /
\plot 14.326 22.608 14.315 22.629 /
\plot 14.315 22.629 14.307 22.653 /
\plot 14.307 22.653 14.300 22.678 /
\plot 14.300 22.678 14.296 22.703 /
\plot 14.296 22.703 14.292 22.731 /
\plot 14.292 22.731 14.290 22.758 /
\plot 14.290 22.758 14.287 22.784 /
\putrule from 14.287 22.784 to 14.287 22.809
\putrule from 14.287 22.809 to 14.287 22.835
\putrule from 14.287 22.835 to 14.287 22.860
\putrule from 14.287 22.860 to 14.287 22.890
\putrule from 14.287 22.890 to 14.287 22.919
\putrule from 14.287 22.919 to 14.287 22.951
\putrule from 14.287 22.951 to 14.287 22.985
\putrule from 14.287 22.985 to 14.287 23.019
\putrule from 14.287 23.019 to 14.287 23.055
\putrule from 14.287 23.055 to 14.287 23.091
\putrule from 14.287 23.091 to 14.287 23.129
\putrule from 14.287 23.129 to 14.287 23.165
\putrule from 14.287 23.165 to 14.287 23.201
\putrule from 14.287 23.201 to 14.287 23.237
\putrule from 14.287 23.237 to 14.287 23.273
\putrule from 14.287 23.273 to 14.287 23.305
\putrule from 14.287 23.305 to 14.287 23.338
\putrule from 14.287 23.338 to 14.287 23.372
\putrule from 14.287 23.372 to 14.287 23.408
\putrule from 14.287 23.408 to 14.287 23.446
\putrule from 14.287 23.446 to 14.287 23.487
\putrule from 14.287 23.487 to 14.287 23.527
\putrule from 14.287 23.527 to 14.287 23.567
\putrule from 14.287 23.567 to 14.287 23.607
\putrule from 14.287 23.607 to 14.287 23.645
\putrule from 14.287 23.645 to 14.287 23.681
\putrule from 14.287 23.681 to 14.287 23.715
\putrule from 14.287 23.715 to 14.287 23.749
\putrule from 14.287 23.749 to 14.287 23.781
\putrule from 14.287 23.781 to 14.287 23.812
\putrule from 14.287 23.812 to 14.287 23.842
\putrule from 14.287 23.842 to 14.287 23.872
\plot 14.287 23.872 14.290 23.904 /
\plot 14.290 23.904 14.292 23.933 /
\plot 14.292 23.933 14.296 23.961 /
\plot 14.296 23.961 14.300 23.990 /
\plot 14.300 23.990 14.307 24.016 /
\plot 14.307 24.016 14.315 24.041 /
\plot 14.315 24.041 14.326 24.064 /
\plot 14.326 24.064 14.336 24.083 /
\plot 14.336 24.083 14.351 24.100 /
\plot 14.351 24.100 14.366 24.117 /
\plot 14.366 24.117 14.383 24.130 /
\plot 14.383 24.130 14.406 24.143 /
\plot 14.406 24.143 14.434 24.155 /
\plot 14.434 24.155 14.465 24.164 /
\plot 14.465 24.164 14.499 24.170 /
\plot 14.499 24.170 14.535 24.174 /
\plot 14.535 24.174 14.573 24.177 /
\putrule from 14.573 24.177 to 14.611 24.177
\plot 14.611 24.177 14.647 24.174 /
\plot 14.647 24.174 14.681 24.168 /
\plot 14.681 24.168 14.713 24.162 /
\plot 14.713 24.162 14.740 24.153 /
\plot 14.740 24.153 14.764 24.141 /
\plot 14.764 24.141 14.787 24.126 /
\plot 14.787 24.126 14.806 24.107 /
\plot 14.806 24.107 14.821 24.079 /
\plot 14.821 24.079 14.834 24.047 /
\plot 14.834 24.047 14.842 24.009 /
\plot 14.842 24.009 14.851 23.967 /
\plot 14.851 23.967 14.855 23.927 /
\plot 14.855 23.927 14.857 23.897 /
\plot 14.857 23.897 14.859 23.884 /
\putrule from 14.859 23.884 to 14.859 23.882
}%
%
%
\put{$-1$
} [lB] at 13.335 22.479
%
%
\put{$-1$
} [lB] at  7.620 22.098
%
%
\put{$+1$
} [lB] at  6.287 24.289
%
%
\put{$+1$
} [lB] at 12.478 24.289
\linethickness=0pt
\putrectangle corners at  5.097 24.543 and 15.287 21.321
\endpicture}